\documentstyle[10pt]{amsart}
\RequirePackage{amssymb}
\RequirePackage[T1]{fontenc}

\newtheorem{theorem}{Theorem}[section]
\newtheorem{lemma}[theorem]{Lemma}
\newtheorem{proposition}[theorem]{Proposition}
\newtheorem{corollary}[theorem]{Corollary} 
\newtheorem{definition}[theorem]{Definition}

\theoremstyle{remark}
\newtheorem{remark}[theorem]{Remark}
\theoremstyle{example}
\newtheorem{example}[theorem]{Example}

\newcommand{\x}{\operatorname{x}}

\newcommand{\ad}{{\operatorname{ad}}}

\newcommand{\Tr}{\text{Tr}}

\newcommand{\diag}{\text{diag}}

\newcommand{\tr}{\text{tr}\,} 
\newcommand{\nc}{\newcommand}
\nc{\Symm}{{\on{Sym}}}

\newcommand{\solu}[1]{\begin{sol}{\bf (\ref{#1})}}

\newcommand{\on}{\operatorname}

\nc{\cE}{{\cal E}}

\newcommand{\kk}{{\bf k}}
\newcommand{\ZZ}{{\mathbb Z}}

\newcommand{\QQ}{{\mathbb Q}}
\newcommand{\CC}{{\mathbb C}}
\newcommand{\RR}{{\mathbb{R}}}
\newcommand{\NN}{{\mathbb N}}

\renewcommand{\sl}{{\mathfrak sl}}
\nc{\SL}{{\mathfrak sl}}
\nc{\gt}{{\mathfrak gt}}
\nc{\grt}{{\mathfrak grt}}
\nc{\gtm}{{\mathfrak gtm}}
\nc{\grtm}{{\mathfrak grtm}}
\nc{\gtmd}{{\mathfrak gtmd}}
\nc{\grtmd}{{\mathfrak grtmd}}
\nc{\pb}{{\mathfrak{pb}}}
\renewcommand{\t}{{\mathfrak{t}}}
\nc{\vt}{{\ov{\t}}}

\renewcommand{\d}{{\mathfrak d}}
\nc{\g}{{\mathfrak{g}}}

\nc{\G}{{\mathfrak{G}}}
\nc{\HH}{{\mathfrak H}}

\newcommand{\h}{{\mathfrak{h}}}

\newcommand{\n}{{\mathfrak{n}}}
\newcommand{\m}{{\mathfrak{m}}}

\newcommand{\SG}{{\mathfrak{S}}}
\newcommand{\f}{{\mathfrak{f}}}

\nc{\wh}{\widehat}\nc{\wt}{\widetilde} 
\renewcommand{\i}{\on{i}}

\newcommand{\ov}{\overline}

\newcommand{\zz}{{\bf z}}
\newcommand{\nn}{{\bf n}}
\newcommand{\mm}{{\bf m}}

\newcommand{\ben}{\begin{enumerate}}
\newcommand{\een}{\end{enumerate}}

\newcommand{\cO}{{\mathcal O}}

\newcommand{\PP}{{\mathbb{P}}}

\newcommand{\cC}{{\mathcal C}}
\newcommand{\cA}{{\mathcal A}}

\newcommand{\cB}{{\mathcal B}}

\begin{document}

\title[Universal KZB equations]{Universal KZB equations I: the elliptic
case}

\author{Damien Calaque}
\address{D.C.: Universit\'e Lyon 1, Institut Camille Jordan, CNRS UMR5201, 
F-69622 Villeurbanne, France}
\email{calaque@@math.univ-lyon1.fr}

\author{Benjamin Enriquez}
\address{B.E.: ULP Strasbourg 1 et IRMA (CNRS), 7 rue Ren\'e Descartes, 
F-67084 Strasbourg, France}
\email{enriquez@@math.u-strasbg.fr}

\author{Pavel Etingof}
\address{P.E.: Department of Mathematics, MIT, Cambridge, MA 02139, USA} 
\email{etingof@@math.mit.edu}

\dedicatory{To Yuri Ivanovich Manin on his 70th birthday}

\begin{abstract} 
We define a universal version of the Knizhnik-Zamolodchikov-Bernard (KZB) 
connection in genus $1$. This is a flat connection over a principal bundle 
on the moduli space of elliptic curves with marked points. It restricts to a 
flat connection on configuration spaces of points on elliptic curves, which 
can be used for proving the formality of the pure braid groups on genus 1 
surfaces. We study the monodromy of this connection and show that it gives 
rise to a relation between the KZ associator and a generating series for 
iterated integrals of Eisenstein forms. We show that the universal
KZB connection realizes as the usual KZB connection for simple Lie 
algebras, and that in the $\sl_n$ case this realization factors through 
the Cherednik algebras. This leads us to define a functor from the category 
of equivariant $D$-modules on $\sl_n$ to that of modules over the Cherednik 
algebra, and to compute the character of irreducible equivariant $D$-modules 
over $\sl_n$ which are supported on the nilpotent cone. 
\end{abstract}

\maketitle

\tableofcontents

\section*{Introduction}
The KZ system was introduced in \cite{KZ} as a system of equations
satisfied by correlation functions in conformal field theory. It was then
realized that this system has a universal version (\cite{Dr:Gal}).
The monodromy of this system leads to representations of the braid
groups, which can be used for proving the that the pure braid groups,
which are the fundamental groups of the configuration spaces of $\CC$,
are formal (i.e., their Lie algebras are isomorphic with their associated
graded Lie algebras, which is a holonomy Lie algebra and thus has an
explicit presentation).
This fact was first proved in the framework of minimal model theory
(\cite{Su,Ko}). These results gave rise to Drinfeld's theory of associators
and quasi-Hopf algebras (\cite{Dr:QH,Dr:Gal}); one of the purposes of
this work was to give an algebraic construction of the formality isomorphisms,
and indeed one of its by-products is the fact that these
isomorphisms can be defined over $\QQ$.

In the case of configuration spaces over surfaces of genus $\geq 1$, similar 
Lie algebra isomorphisms were constructed by Bezrukavnikov (\cite{Bez}), 
using results of Kriz 
(\cite{Kr}). In this series of papers, we will show that this result can be 
reproved using a suitable flat connection over configuration spaces. This 
connection is a universal version of the KZB connection (\cite{Be1,Be2}), 
which is the higher genus analogue of the KZ connection.

In this paper, we focus on the case of genus $1$. We define the 
universal KZB connection (Section \ref{sect:1}), and rederive from there
the formality result (Section \ref{sect:2}).  As in the integrable case of the 
KZB connection, the universal KZB connection extends from the configuration 
spaces $\bar C(E_{\tau},n)/S_{n}$ to the moduli space ${\cal M}_{1,[n]}$
of elliptic curves with $n$ unordered marked points (Section \ref{sect:4}). 
This means that: (a) the connection can be extended to the directions of 
variation of moduli, and  (b) it is modular invariant.

This connection then gives rise to a monodromy morphism $\gamma_{n}:\Gamma_{1,[n]}
\to {\bold G}_{n} \rtimes S_{n}$, which we analyze in Section \ref{sect:5}. 
The images of most generators can be expressed using the KZ associator, 
but the image $\tilde\Theta$ of the $S$-transformation expresses using iterated 
integrals of Eisenstein series. The relations between generators give rise to 
relations between $\tilde\Theta$ and the KZ associator, identities 
(\ref{gamma12}). This identity may be viewed as an elliptic analogue
of the pentagon identity, as it is a ``de Rham'' analogue of the relation 
6AS in \cite{HLS} (in \cite{Ma}, the question was asked of the existence 
of this kind of identity).

In Section \ref{sect:5bis}, we investigate how to algebraically 
construct a morphism $\Gamma_{1,[n]}\to {\bold G}_{n}\rtimes S_{n}$. 
We show that a morphism $\overline{\on{B}}_{1,n}
\to \on{exp}(\hat{\bar\t}_{1,n})\rtimes S_{n}$ can be constructed
using an associator only (here $\overline{\on{B}}_{1,n}$ is the reduced 
braid group of $n$ points on the torus). \cite{Dr:Gal} then implies that 
the formality isomorphism can be defined over $\QQ$. In the last part of 
Section \ref{sect:5bis}, we develop the analogue of the theory of 
quasitriangular quasibialgebras (QTQBA's), namely elliptic structures over 
QTQBA's. These structures give rise to representations of 
$\overline{\on{B}}_{1,n}$, and they can be modified by twist. We hope that 
in the case of a simple Lie algebra, and using suitable twists, the elliptic 
structure given in Section \ref{5:6} will give rise to elliptic structures 
over the quantum group $U_{q}(\g)$ (where $q\in \CC^{\times}$) or over the 
Lusztig quantum group (when $q$ is a root of unity), yielding back the 
representations of $\overline{\on{B}}_{1,n}$ from conformal field
theory.

In Section \ref{sect:6}, we show that the universal KZB connection indeed 
specializes to the ordinary KZB connection.

Sections 7-9 are dedicated applications of the ideas of the preceding 
sections (in particular, Section 6) to representation theory of 
Cherednik algebras.

More precisely, In Section 7, we construct a 
homomorphism from the Lie algebra $\bar {\mathfrak {t}}_{1,n}\rtimes 
{\mathfrak {d}}$ to the rational Cherednik algebra $H_n(k)$ of type 
$A_{n-1}$. 
This allows us to consider the elliptic KZB connection with values in 
representations of the rational Cherednik algebra. The monodromy of this 
connection then gives representations of the true Cherednik algebra
(i.e. the double affine Hecke algebra). In particular, this gives a simple 
way of constructing an isomorphism between the rational Cherednik 
algebra and the double affine Hecke algebra, with formal deformation 
parameters.

In Section 8, we consider the special representation $V_N$ of the
rational Cherednik algebra $H_n(k)$, $k=N/n$, for which the
elliptic KZB connection is the KZB connection for (holomorphic)
$n$-point correlation functions of the WZW model for ${\rm SL}_N(\CC)$ on
the elliptic curve, when the marked points are labeled by the
vector representation ${\Bbb C}^N$. This representation is realized
in the space of equivariant polynomial functions on
${\mathfrak{sl}}_N$ with values in $({\Bbb C}^N)^{\otimes n}$, and
we show that it is irreducible, and calculate its character.

In Section 9, we generalize the construction of Section 8, by
replacing, in the construction of $V_N$, the space of polynomial
functions on ${\mathfrak{sl}}_N$ with an arbitrary $D$-module on
${\mathfrak{sl}}_N$.  This gives rise to an exact functor from
the category of (equivariant) $D$-modules on ${\mathfrak{sl}}_N$ to
the category of representations of $H_n(N/n)$. We study this
functor in detail. In particular, we show that this functor maps
$D$-modules concentrated on the nilpotent cone to modules from
category ${\mathcal O}_-$ of highest weight modules over the
Cherednik algebra, and is closely related to the Gan-Ginzburg
functor, \cite{GG1}. Using these facts, we show that it maps
irreducible $D$-modules on the nilpotent cone to irreducible
representations of the Cherednik algebra, and determine their
highest weights. As an application, we compute the decomposition of 
cuspidal $D$-modules into irreducible representations of $\on{SL}_N(\CC)$. 
Finally, we describe the generalization of the above result 
to the trigonometric case (which involves $D$-modules on the group 
and trigonometric Cherednik algebras), and point out several directions 
for generalization.

\section{Bundles with flat connections on (reduced) configuration spaces} 
\label{sect:1}

\subsection{The Lie algebras $\t_{1,n}$ and $\bar\t_{1,n}$}

Let $n\geq 1$ be an integer and $\kk$ be a field of characteristic zero. 
We define $\t_{1,n}^\kk$ as the Lie algebra with generators 
$x_{i},y_{i}$ ($i=1,...,n$) and $t_{ij}$ ($i\neq j\in \{1,...,n\}$) 
and relations 
\begin{equation} \label{inf:pure:braid}
t_{ij} = t_{ji}, \quad [t_{ij},t_{ik} + t_{jk}] = 0, \quad 
[t_{ij},t_{kl}] = 0, 
\end{equation}
$$
[x_i,y_j] = t_{ij}, \quad [x_i,x_j] = [y_i,y_j] = 0, \quad
[x_i,y_i] = -\sum_{j|j\neq i} t_{ij}, 
$$
$$
[x_i,t_{jk}] = [y_i,t_{jk}] = 0 , \quad [x_i+x_j,t_{ij}] = [y_i+y_j,t_{ij}] = 0. 
$$
($i,j,k,l$ are distinct). In this Lie algebra, $\sum_i x_i$ 
and $\sum_i y_i$ are central; we then define 
$\bar\t_{1,n}^\kk := \t_{1,n}^\kk / (\sum_{i} x_{i}, \sum_{i} y_{i})$. 
Both $\t_{1,n}^\kk$ and $\bar\t_{1,n}^\kk$ are positively graded, where 
$\on{deg}(x_{i}) = \on{deg}(y_{i})=1$.

The symmetric group $S_{n}$ acts by automorphisms of $\t_{1,n}^\kk$
by $\sigma(x_{i}):= x_{\sigma(i)}$, $\sigma(y_{i}):= y_{\sigma(i)}$, 
$\sigma(t_{ij}):= t_{\sigma(i)\sigma(j)}$; this induces an action of
$S_{n}$ by automorphisms of $\bar\t_{1,n}^\kk$.

We will set $\t_{1,n}:= \t_{1,n}^\CC$, 
$\bar\t_{1,n}:= \bar\t_{1,n}^\CC$ in Sections \ref{sect:1} to \ref{sect:5}.

\subsection{Bundles with flat connections over $C(E,n)$ and $\bar C(E,n)$}

Let $E$ be an elliptic curve, $C(E,n)$ be the configuration space 
$E^{n} - \{$diagonals$\}$ ($n\geq 1$) and $\bar C(E,n) := C(E,n)/E$ 
be the reduced configuration space. We will define 
a\footnote{We will denote by $\hat\g$ or $\g^{\wedge}$ the degree completion 
of a positively graded Lie algebra $\g$.} 
$\on{exp}(\hat{\bar\t}_{1,n})$-principal bundle with a flat 
(holomorphic) connection $(\bar P_{E,n},\bar\nabla_{E,n}) \to \bar C(E,n)$. 
For this, we define a $\on{exp}(\hat\t_{1,n})$-principal bundle with a 
flat connection $(P_{E,n},\nabla_{E,n}) \to C(E,n)$. Its image under 
the natural morphism $\on{exp}(\hat\t_{n}) \to \on{exp}(\hat{\bar\t}_{n})$ is 
a $\exp(\hat{\bar\t}_{1,n})$-bundle with connection 
$(\tilde P_{E,n},\tilde\nabla_{E,n})\to C(E,n)$, and we then prove that 
$(\tilde P_{E,n},\tilde\nabla_{E,n})$ is the pull-back of a pair 
$(\bar P_{E,n},\bar\nabla_{E,n})$ under the canonical projection 
$C(E,n)\to\bar C(E,n)$.

For this, we fix a uniformization $E \simeq E_{\tau}$, where for 
$\tau\in{\mathfrak H}$, ${\mathfrak H} := \{\tau\in\CC | \Im(\tau)>0\}$, 
$E_{\tau} := \CC/\Lambda_{\tau}$ and $\Lambda_\tau:= \ZZ+\ZZ\tau$.

We then have $C(E_\tau,n)=(\CC^{n}-\on{Diag}_{n,\tau})/\Lambda_{\tau}^{n}$, 
where $\on{Diag}_{n,\tau} := \{\zz = (z_{1},...,z_{n})\in\CC^{n} | 
z_{ij}:= z_{i} - z_{j}\in \Lambda_{\tau}$ for some $i\neq j\}$. We define 
$P_{\tau,n}$ as the restriction to $C(E_\tau,n)$ of the bundle over 
$\CC^{n}/\Lambda_{\tau}^{n}$ for which a section on 
$U \subset \CC^{n}/\Lambda_{\tau}^{n}$ is a regular map 
$f : \pi^{-1}(U)\to \on{exp}(\hat\t_{1,n})$, such 
that\footnote{We set $\i := \sqrt{-1}$, leaving $i$ for indices.}
$f(\zz + \delta_{i}) = f(\zz)$, $f(\zz + \tau\delta_{i}) = 
e^{-2\pi\i x_{i}}f(\zz)$ (here $\pi : \CC^{n} \to \CC^{n}/\Lambda_{\tau}^{n}$ 
is the canonical projection and $\delta_{i}$ is the $i$th vector of 
the canonical basis of $\CC^{n}$).

The bundle $\tilde P_{\tau,n}\to C(E_\tau,n)$ derived from $P_{\tau,n}$
is the pull-back of a bundle $\bar P_{\tau,n}\to \bar C(E_\tau,n)$ 
since the $e^{-2\pi\i\bar x_{i}}\in \on{exp}(\hat{\bar\t}_{1,n})$
commute pairwise and their product is $1$. Here $x\mapsto \bar x$ is the 
map $\hat\t_{1,n}\to\hat{\bar\t}_{1,n}$.

A flat connection $\nabla_{\tau,n}$ on $P_{\tau,n}$ is then the same 
as an equivariant flat connection over the trivial bundle over 
$\CC^{n} - \on{Diag}_{n,\tau}$, i.e., a connection of the form 
$$
\nabla_{\tau,n} := \on{d} - \sum_{i=1}^{n} K_{i}(\zz|\tau)\on{d}z_{i}, 
$$ 
where $K_{i}(-|\tau) : \CC^{n} \to \hat\t_{1,n}$ is holomorphic on 
$\CC^n - \on{Diag}_{n,\tau}$, such that:

(a) $K_{i}(\zz + \delta_{j}|\tau) = K_{i}(\zz|\tau)$, 
$K_{i}(\zz + \tau\delta_{j}|\tau) 
= e^{-2\pi\i\on{ad}(x_{j})}(K_{i}(\zz|\tau))$,

(b) $[\partial/\partial z_{i} - K_{i}(\zz|\tau),\partial /\partial z_{j} 
- K_{j}(\zz|\tau)] = 0$ for any $i,j$.

$\nabla_{\tau,n}$ then induces a flat connection $\tilde\nabla_{\tau,n}$ 
on $\tilde P_{\tau,n}$. Then $\tilde\nabla_{\tau,n}$ is the pull-back of 
a (necessarily flat) connection on $\bar P_{\tau,n}$ iff:

(c) $\bar K_{i}(\zz|\tau) = \bar K_{i}(\zz + u(\sum_{i}\delta_{i})|\tau)$
and $\sum_{i}\bar K_{i}(\zz|\tau)=0$ for $\zz\in\CC^{n} - \on{Diag}_{n,\tau}$, 
$u\in\CC$.

In order to define the $K_{i}(\zz|\tau)$, we first recall some facts on 
theta-functions. 
There is a unique holomorphic function $\CC\times\HH \to \CC$, $(z,\tau)\mapsto \theta(z|\tau)$, such that 
$\{z | \theta(z|\tau) = 0\} = \Lambda_{\tau}$, $\theta(z+1|\tau)
= -\theta(z|\tau) = \theta(-z|\tau)$ and $\theta(z+\tau|\tau) = - e^{-\pi\i\tau} e^{-2\pi\i z}\theta(z|\tau)$, 
and $\theta_{z}(0|\tau)=1$. We have $\theta(z|\tau+1) = \theta(z|\tau)$, while $\theta(-z/\tau|-1/\tau)
= - (1/\tau) e^{(\pi\i/\tau) z^2} \theta(z|\tau)$. 
If $\eta(\tau) = q^{1/24}\prod_{n\geq 1} (1-q^n)$ where $q = e^{2\pi\i\tau}$, 
and if we set $\vartheta(z|\tau) := \eta(\tau)^3 \theta(z|\tau)$, then 
$\partial_\tau\vartheta = (1/4\pi\i) \partial_z^2\vartheta$.

Let us set 
$$
k(z,x|\tau) := {{\theta(z+x|\tau)}\over{\theta(z|\tau)\theta(x|\tau)}} 
- {1\over x}. 
$$
When $\tau$ is fixed, $k(z,x|\tau)$
belongs to $\on{Hol}(\CC - \Lambda_\tau)[[x]]$. Substituting $x=\ad x_i$, 
we get a linear map $\t_{1,n}\to (\t_{1,n}\otimes\on{Hol}(\CC 
- \Lambda_\tau))^{\wedge}$, and taking the image of $t_{ij}$, we define 
$$
K_{ij}(z|\tau) := k(z,\on{ad}x_i|\tau)(t_{ij})
=\big( {{\theta(z + \on{ad}(x_i)|\tau)}\over{\theta(z|\tau)}}
{{\on{ad}(x_i)}\over{\theta(\on{ad}(x_i)|\tau)}} - 1\big)(y_j); 
$$
it is a holomorphic function on $\CC - \Lambda_\tau$ with values in 
$\hat\t_{1,n}$.

Now set $\zz:=(z_1,\ldots,z_n)$, $z_{ij}:=z_i-z_j$ and define 
$$
K_i(\zz|\tau):= -y_i + \sum_{j|j\neq i} K_{ij}(z_{ij}|\tau). 
$$

Let us check that the $K_{i}(\zz|\tau)$ satisfy condition (c). 
We have clearly  $K_{i}(\zz + u(\sum_{i} \delta_{i})) = K_{i}(\zz)$. 
We have $k(z,x|\tau) + k(-z,-x|\tau) = 0$, so $K_{ij}(z|\tau) + K_{ji}(-z|\tau) = 0$, so that
  $\sum_{i} K_{i}(\zz|\tau) = -\sum_{i} y_{i}$,  which implies $\sum_{i} \bar K_{i}(\zz|\tau)=0$.

\begin{lemma}
$K_{i}(\zz+\delta_{j}|\tau) = K_{i}(\zz|\tau)$ and 
$K_{i}(\zz+\tau\delta_{j}|\tau) = e^{-2\pi\i \on{ad}x_{j}}(K_{i}(\zz|\tau))$, i.e., the 
$K_{i}(\zz|\tau)$ satisfy condition (a). 
\end{lemma}

{\em Proof.} We have $k(z\pm 1,x|\tau) = k(z,x|\tau)$ so for any $j$, 
$K_{i}(\zz + \delta_{j}|\tau) = K_{i}(\zz|\tau)$. We have 
$k(z\pm \tau,x|\tau) = e^{\mp 2\pi\i x}k(z,x|\tau) + (e^{\mp 2\pi\i x}-1)/x$, 
so if $j\neq i$, 
$K_{i}(\zz + \tau\delta_{j}|\tau) = \sum_{j'\neq i,j} K_{ij'}(z_{ij'}|\tau) + 
e^{2\pi\i\on{ad}x_{i}} K_{ij}(z_{ij}|\tau) 
+ {{e^{2\pi\i\on{ad}x_{i}}-1}\over {\on{ad}x_{i}}}(t_{ij}) - y_{i}$. 
Then 
$$
{{e^{2\pi\i\on{ad}x_{i}}-1}\over{\on{ad}x_{i}}}(t_{ij})
= {{1-e^{-2\pi\i\on{ad}x_{j}}}\over{\on{ad}x_{j}}}(t_{ij})
= (1-e^{-2\pi\i\on{ad}x_{j}})(y_{i}),
$$ 
$e^{2\pi\i\on{ad}x_{i}} (K_{ij}(z_{ij}|\tau)) 
= e^{-2\pi\i\on{ad}x_{j}} (K_{ij}(z_{ij}|\tau))$
and for $j'\neq i,j$, $K_{ij'}(z_{ij'}|\tau) = e^{-2\pi\i\on{ad}x_{j}}(K_{ij'}(z_{ij'}|\tau))$, 
so $K_{i}(\zz + \tau\delta_{j}|\tau) 
= e^{-2\pi\i\on{ad}x_{j}}(K_{i}(\zz|\tau))$.

Now $K_{i}(\zz + \tau\delta_{i}|\tau) 
= -\sum_{i}y_{i} - \sum_{j|j\neq i}K_{j}(\zz + \tau\delta_{i}|\tau)
= -\sum_{i}y_{i} - e^{-2\pi\i\on{ad}x_{i}}(\sum_{j|j\neq i} K_{j}(\zz|\tau))
= e^{-2\pi\i\on{ad}x_{i}}(-\sum_{i}y_{i} - \sum_{j|j\neq i}K_{j}(\zz|\tau)) = 
e^{-2\pi\i\on{ad}x_{i}}K_{i}(\zz|\tau)$ (the first and last equality follow 
from the proof of (c), the second equality has just been proved, the 
third equality follows from the centrality of $\sum_{i} y_{i}$). 
\hfill \qed \medskip

\begin{proposition}
$[\partial/\partial z_i-K_i(\zz|\tau),
\partial/\partial z_j-K_j(\zz|\tau)]=0$, 
i.e., the $K_{i}(\zz|\tau)$ satisfy condition (b). 
\end{proposition}

{\em Proof.} For $i\neq j$, let us set $K_{ij} := K_{ij}(z_{ij}|\tau)$. 
Recall that $K_{ij} + K_{ji}=0$, therefore if 
$\partial_i := \partial/\partial z_i$
$$
\partial_i K_{ij} - \partial_j K_{ji} = 0, \quad 
[y_i - K_{ij} , y_j - K_{ji}] = -[K_{ij},y_i + y_j].
$$

Moreover, if $i,j,k,l$ are distinct, then $[K_{ik}, K_{jl}] = 0$.
It follows that if $i\neq j$, 
\begin{align*}
& [\partial_i - K_i(\zz|\tau), \partial_j - K_j(\zz|\tau)] 
\\ & = [y_i + y_j,K_{ij}]
+ \sum_{k|k\neq i,j} \big( [K_{ik},K_{jk}] + [K_{ij},K_{jk}]
+ [K_{ij},K_{ik}] + [y_j,K_{ik}] - [y_i,K_{jk}] \big) . 
\end{align*}

Let us assume for a while that if $k\notin \{i,j\}$, then 
\begin{equation} \label{eq:3}
-[y_i,K_{jk}] - [y_j,K_{ki}] - [y_k,K_{ij}]
+ [K_{ji},K_{ki}] + [K_{kj},K_{ij}]
+ [K_{ik},K_{jk}] = 0 
\end{equation} 
(this is the universal version of the classical dynamical Yang-Baxter equation).

Then (\ref{eq:3}) implies that 
$$
[\partial_i - K_i(\zz|\tau), \partial_j - K_j(\zz|\tau)] = [y_i + y_j,K_{ij}]
+ \sum_{k|k\neq i,j} [y_k,K_{ij}] = [\sum_k y_k,K_{ij}] = 0 
$$
(as $\sum_k y_k$ is central), which proves the proposition.

Let us now prove (\ref{eq:3}). 
If $f(x)\in \CC[[x]]$, then 
$$
[y_k,f(\ad x_i)(t_{ij})] = {{f(\ad x_i) - f(-\ad x_j)}\over{\ad x_i + 
\ad x_j}}[-t_{ki},t_{ij}], 
$$
$$
[y_i,f(\ad x_j)(t_{jk})] = {{f(\ad x_j) - f(-\ad x_k)}\over{\ad x_j + 
\ad x_k}}[-t_{ij},t_{jk}] = 
{{f(\ad x_j) - f(\ad x_i + \ad x_j)}\over{-\ad x_i}}[-t_{ij},t_{jk}] , 
$$
$$
[y_j,f(\ad x_k)(t_{ki})] = {{f(\ad x_k) - f(-\ad x_i)}\over{\ad x_k + 
\ad x_i}}[-t_{jk},t_{ki}] =
{{f(-\ad x_i - \ad x_j) - f(-\ad x_i)}\over{- \ad x_j}}
[-t_{jk},t_{ki}] . 
$$
The first identity is proved as follows: 
\begin{align*}
& [y_k,(\ad x_i)^n(t_{ij})] = -\sum_{s=0}^{n-1}
(\ad x_i)^s (\ad  t_{ki}) (\ad x_i)^{n-1-s}(t_{ij})
= -\sum_{s=0}^{n-1}
(\ad x_i)^s (\ad  t_{ki}) (-\ad x_j)^{n-1-s}(t_{ij})
\\ & 
= -\sum_{s=0}^{n-1}
(\ad x_i)^s (-\ad x_j)^{n-1-s} (\ad  t_{ki})(t_{ij})
= f(\ad x_i,-\ad x_j)([-t_{ki},t_{ij}]), 
\end{align*}
where $f(u,v) = (u^n-v^n)/(u-v)$. The two next identities 
follow from this one and from the fact that $x_i+x_j+x_k$
commutes with $t_{ij},t_{ik},t_{jk}$.

Then, if we write $k(z,x)$ instead of $k(z,x|\tau)$, the l.h.s. of (\ref{eq:3}) is equal to 
\begin{align*}
& \Big( k(z_{ij},-\ad x_j)k(z_{ik},\ad x_i + \ad x_j) 
- k(z_{ij},\ad x_i)k(z_{jk},\ad x_i + \ad x_j) 
+ k(z_{ik},\ad x_i)k(z_{jk},\ad x_j) 
\\ & 
+ {{k(z_{jk},\ad x_j) - k(z_{jk},\ad x_i + \ad x_j)}\over{\ad x_i}} 
+ {{k(z_{ik},\ad x_i) - k(z_{ij},\ad x_i + \ad x_j)}\over{\ad x_j}} 
\\ & - {{k(z_{ij},\ad x_i) - k(z_{ij},-\ad x_j)}\over{\ad x_i + \ad x_j}}
\Big) [t_{ij},t_{ik}]. 
\end{align*}

So (\ref{eq:3}) follows from the identity 
\begin{align*}
& k(z,-v)k(z',u+v) - k(z,u) k(z'-z,u+v) + k(z',u) k(z'-z,v) 
\\ & + {{k(z'-z,v) - k(z'-z,u+v)}\over{u}}
+ {{k(z',u) - k(z',u+v)}\over{v}} 
- {{k(z,u) - k(z,-v)}\over{u+v}} = 0, 
\end{align*}
where $u,v$ are formal variables, which is a consequence of the 
theta-functions identity 
\begin{align} \label{theta:id}
& \nonumber \big( k(z,-v) - {1\over v}\big) \big( k(z',u+v) + {1\over {u+v}}\big) 
- \big( k(z,u) + {1\over u}\big) \big( k(z'-z,u+v) + {1\over{u+v}}\big) 
\\ & + \big( k(z',u) + {1\over u}\big) \big( k(z'-z,v) + {1\over v}\big) = 0. 
\end{align}
\hfill \qed \medskip

We have therefore proved:

\begin{theorem}
$(P_{\tau,n},\nabla_{\tau,n})$ is a flat connection on $C(E_\tau,n)$, 
and the induced flat connection 
$(\tilde P_{\tau,n},\tilde\nabla_{\tau,n})$ is the pull-back of a 
unique flat connection 
$(\bar P_{\tau,n},\bar\nabla_{\tau,n})$ on $\bar C(E_\tau,n)$. 
\end{theorem}

\subsection{Bundles with flat connections on $C(E,n)/S_{n}$ and $\bar C(E,n)/S_{n}$}

The group $S_{n}$ acts freely by automorphisms of $C(E,n)$ by $\sigma(z_{1},...,z_{n}):= 
(z_{\sigma^{-1}(1)},...,z_{\sigma^{-1}(n)})$. This descends to a free action of 
$S_{n}$ on $\bar C(E,n)$. We set $C(E,[n]):= C(E,n)/S_{n}$, $\bar C(E,[n]):= 
\bar C(E,n)/S_{n}$.

We will show that $(P_{\tau,n},\nabla_{\tau,n})$
induces a bundle with flat connection $(P_{\tau,[n]},\nabla_{\tau,[n]})$
on $C(E_{\tau},[n])$ with group $\on{exp}(\hat\t_{1,n})\rtimes S_{n}$, 
and similarly $(\bar P_{\tau,n},\bar\nabla_{\tau,n})$ induces 
$(\bar P_{\tau,[n]},\bar\nabla_{\tau,[n]})$ on 
$\bar C(E_{\tau},[n])$ with group $\on{exp}(\hat{\bar\t}_{1,n})\rtimes S_{n}$.

We define $P_{\tau,[n]} \to C(E_{\tau},[n])$ by the condition that 
a section of $U\subset C(E_{\tau},[n])$ is a regular map $\pi^{-1}(U)\to 
\on{exp}(\hat\t_{1,n})\rtimes S_{n}$, satisfying again $f(\zz+\delta_{i}) = f(\zz)$, 
$f(\zz + \tau\delta_{i}) = e^{-2\pi\i x_{i}}f(\zz)$ and the additional requirement 
$f(\sigma\zz) = \sigma f(\zz)$ 
(where $\tilde\pi : \CC^{n} - \on{Diag}_{\tau,n} \to C(E_{\tau},[n])$ is the canonical 
projection). It is clear that $\nabla_{\tau,n}$ is $S_{n}$-invariant, which implies 
that it defines a flat connection  $\nabla_{\tau,[n]}$ on  $C(E_{\tau},[n])$.

The bundle $\bar P(E_{\tau},[n])\to \bar C(E_{\tau},[n])$ is defined 
by the additional requirement $f(\zz + u(\sum_{i}\delta_{i})) = f(\zz)$ 
and $\bar\nabla_{\tau,n}$ then induces a flat connection 
$\bar \nabla_{\tau,[n]}$ on $\bar C(E_{\tau},[n])$.

\section{Formality of pure braid groups on the torus} \label{sect:2}

\subsection{Reminders on Malcev Lie algebras}

Let $\kk$ be a field of characteristic $0$ and let $\g$ be a pronilpotent
$\kk$-Lie algebra. Set $\g^1 = \g$, $\g^{k+1} = [\g,\g^k]$; then 
$\g = \g^1 \supset \g^2...$ is a decreasing filtration of $\g$. The 
associated graded Lie algebra is 
$\on{gr}(\g) := \oplus_{k\geq 1}\g^k/\g^{k+1}$; we also consider its 
completion $\hat{\on{gr}}(\g) := \hat{\oplus}_{k\geq 1}\g^k/\g^{k+1}$ 
(here $\hat\oplus$ is the direct product).  We say that $\g$ is formal 
iff there exists an isomorphism of filtered Lie algebras $\g\simeq
\hat{\on{gr}}(\g)$, whose associated graded morphism is the identity. 
We will use the following fact: if $\g$ is a pronilpotent Lie algebra, 
$\t$ is a positively graded Lie algebra, and there exists an isomorphism 
$\g\simeq \hat\t$ of filtered Lie algebras, then $\g$ is formal, 
and the associated graded morphism $\on{gr}(\g)\to \t$ is an 
isomorphism of graded Lie algebras.

If $\Gamma$ is a finitely generated group, there exists a unique 
pair $(\Gamma(\kk),i_{\Gamma})$ of a prounipotent algebraic group 
$\Gamma(\kk)$ and a group morphism $i_{\Gamma} : \Gamma\to \Gamma(\kk)$, 
which is initial in the category of all pairs $(U,j)$, where $U$ is 
prounipotent $\kk$-algebraic group and $j : \Gamma\to U$ is a group morphism.

We denote by $\on{Lie}(\Gamma)_\kk$ the Lie algebra of $\Gamma(\kk)$. 
Then we have $\Gamma(\kk) = \on{exp}(\on{Lie}(\Gamma)_\kk)$; 
$\on{Lie}(\Gamma)_\kk$ is a pronilpotent Lie algebra. We have 
$\on{Lie}(\Gamma)_\kk = \on{Lie}(\Gamma)_\QQ\otimes \kk$. We say that 
$\Gamma$ is formal iff $\on{Lie}(\Gamma)_\CC$ is formal (one can 
show that this implies that $\on{Lie}(\Gamma)_\QQ$ is formal).

When $\Gamma$ is presented by generators $g_{1},...,g_{n}$ and relations 
$R_{i}(g_{1},...,g_{n})$ ($i=1,...,p$),  $\on{Lie}(\Gamma)_\QQ$ is the 
quotient of the topologically free Lie algebra $\hat\f_{n}$ generated by 
$\gamma_{1},...,\gamma_{n}$ by the topological ideal generated by 
$\on{log}(R_{i}(e^{\gamma_{1}},...,e^{\gamma_{n}}))$ ($i=1,...,p$).

The decreasing filtration of $\hat\f_{n}$ is $\hat\f_{n} = 
(\hat\f_{n})^{1} \supset (\hat\f_{n})^{2} \supset...$, where 
$(\hat\f_{n})^{k}$ is the part of $\hat\f_{n}$ of degree $\geq k$ in the 
generators $\gamma_{1},...,\gamma_{n}$. The image of this 
filtration by the projection is map is the decreasing filtration 
$\on{Lie}(\Gamma)_\QQ = \on{Lie}(\Gamma)_\QQ^{1}\supset 
\on{Lie}(\Gamma)_\QQ^{2} \supset...$ of $\on{Lie}(\Gamma)_\QQ$.

\subsection{Presentation of $\on{PB}_{1,n}$}
For $\tau\in\HH$, let $U_\tau \subset \CC^n - \on{Diag}_{n,\tau}$ be the 
open subset of 
all $\zz = (z_1,...,z_n)$, of the form $z_i = a_i + \tau b_i$, where 
$0<a_1<...<a_n<1$ and $0<b_1<...<b_n<1$. 
If $\zz_0=(z_{1}^{0},...,z_{n}^{0})\in U_\tau$, its image $\overline\zz_0$ in $E_\tau^n$
actually belongs to the configuration space $C(E_\tau,n)$.

The pure braid group of $n$ points on the torus $\on{PB}_{1,n}$
may be viewed as $\on{PB}_{1,n} = \pi_1(C(E_{\tau},n),\overline\zz_0)$. Denote
by $X_i,Y_i\in \on{PB}_{1,n}$ the classes of the projection of the 
paths $[0,1]\ni t \mapsto \zz_{0} - t \delta_i$ and 
$[0,1]\ni t \mapsto \zz_{0} - t \tau \delta_i$.

Set $A_i := X_i...X_n$, $B_i := Y_i...Y_n$ for $i=1,...,n$. 
According to \cite{Bi1}, $A_i,B_i$ ($i=1,...,n$) generate $\on{PB}_{1,n}$
and a presentation of $\on{PB}_{1,n}$ is, in terms of these generators: 
$$
(A_i,A_j) = (B_i,B_j)=1\; (\on{any\ }i,j),\quad   (A_{1},B_{j}) = (B_{1},A_{j})=1\; 
(\on{any\ } j),
$$
$$
(B_{k},A_{k}A_{j}^{-1}) = 
(B_{k}B_{j}^{-1},A_{k}) = C_{jk}\; (j\leq k),\quad (A_{i},C_{jk}) = (B_{i},C_{jk})=1\; 
(i\leq j\leq k), 
$$
where $(g,h) = ghg^{-1}h^{-1}$.

\subsection{Alternative presentations of $\t_{1,n}$}

We now give two variants of the defining presentation of $\t_{1,n}$. Presentation 
(A) below is the original presentation in \cite{Bez}, and presentation (B) will be suited
to the comparison with the above presentation of $\on{PB}_{1,n}$.

\begin{lemma} \label{lemma:pres}
$\t_{1,n}$ admits the following presentations:

(A) generators are $x_i$, $y_i$
($i = 1,...,n$), relations are $[x_i,y_j] = [x_j,y_i]$ ($i\neq j$), $[x_i,x_j] 
= [y_i,y_j]=0$ (any $i,j$), $[\sum_j x_j,y_i] = [\sum_j y_j,x_i]=0$ (any $i$), $[x_i,[x_j,y_k]]
= [y_i,[y_j,x_k]]=0$ ($i,j,k$ are distinct);

(B) generators are $a_{i}$, $b_{i}$ ($i=1,...,n$), relations are 
$[a_{i},a_{j}] = [b_{i},b_{j}]=0$ (any $i,j$), 
$[a_{1},b_{j}]=[b_{1},a_{j}]=0$ (any $j$), $[a_{j},b_{k}]=[a_{k},b_{j}]$ (any $i,j$), 
$[a_{i},c_{jk}] = [b_{i},c_{jk}]=0$ ($i\leq j\leq k$), where $c_{jk} = [b_{k},a_{k}-a_{j}]$.

The isomorphism of presentations (A) and (B) is $a_{i} = \sum_{j=i}^{n} x_{j}$, 
$b_{i} = \sum_{j=i}^{n} y_{j}$. 
\end{lemma}

{\em Proof.} Let us prove that the initial relations for $x_{i},y_{i},t_{ij}$ imply the
relations (A) for $x_{i},y_{i}$. Let us assume the initial relations. 
If $i\neq j$, since $[x_i,y_j] = t_{ij}$ and $t_{ij} = t_{ji}$, 
we get $[x_i,y_j] = [x_j,y_i]$. The relations $[x_i,x_j] 
= [y_i,y_j]=0$ (any $i,j$) are contained in the initial relations. For any $i$, since 
$[x_i,y_i] = -\sum_{j|j\neq i}t_{ij}$ and $[x_j,y_i] = t_{ji} = t_{ij}$ ($j\neq i$), 
we get $[\sum_j x_j,y_i] = 0$. Similarly, $[\sum_j y_j,x_i]=0$ (for any $i$). 
If $i,j,k$ are distinct, since $[x_j,y_k]=t_{jk}$ and $[x_i,t_{jk}]=0$, we get $[x_i,[x_j,y_k]]=0$
and similarly we prove $[x_i,[y_j,x_k]]=0$.

Let us now prove that the relations (A) for $x_{i},y_{i}$ imply the initial relations
for $x_i,y_i$ and $t_{ij}:= [x_i,y_j]$ ($i\neq j$). Assume the relations (A). 
If $i\neq j$, since $[x_i,y_j]=[x_j,y_i]$, we have $t_{ij} = t_{ji}$. The 
relation $t_{ij} = [x_i,y_j]$ ($i\neq j$) is clear and $[x_i,x_j]=[y_i,y_j]=0$
(any $i,j$) are already in relations (A). Since for any $i$, $[\sum_j x_j,y_i]=0$, 
we get $[x_i,y_i] = -\sum_{j|j\neq i} [x_j,y_i] = -\sum_{j|j\neq i} t_{ji}
= -\sum_{j|j\neq i} t_{ij}$. If $i,j,k$ are distinct, the relations $[x_i,[x_j,y_k]] = 
[y_i,[y_j,x_k]]=0$ imply $[x_i,t_{jk}]=[y_i,t_{jk}]=0$. If $i\neq j$, since 
$[\sum_k x_k,x_i] = [\sum_k x_k,y_j]=0$, we get $[\sum_k x_k,t_{ij}]=0$
and $[x_k,t_{ij}]=0$ for $k\notin\{i,j\}$ then implies $[x_i+x_j,t_{ij}]=0$. 
One proves similarly $[y_i+y_j,t_{ij}]=0$. We have already shown that 
$[x_i,t_{kl}] = [y_j,t_{kl}]=0$ for $i,j,k,l$ distinct, which implies 
$[[x_i,y_j],t_{kl}]=0$, i.e., $[t_{ij},t_{kl}]=0$. 
If $i,j,k$ are distinct, we have shown that $[t_{ij},y_k]=0$ and $[t_{ij},x_i+x_j]=0$, 
which implies $[t_{ij},[x_i+x_j,y_k]]=0$, i.e., 
$[t_{ij},t_{ik}+t_{jk}]=0$.

Let us prove that the relations (A) for $x_{i},y_{i}$ imply relations 
(B) for $a_{i}:= \sum_{j=i}^{n}x_{j}$, $b_{i}:= \sum_{j=i}^{n}y_{j}$. 
Summing up the relations $[x_{i'},x_{j'}]=[y_{i'},y_{j'}]=0$ and $[x_{i'},y_{j'}]=
[x_{j'},y_{i'}]$ for $i'=i,...,n$ and $j'=j,...,n$, we get 
$[a_{i},a_{j}]=[b_{i},b_{j}]=0$ and $[a_{i},b_{j}]=[a_{j},b_{i}]$ (for any $i,j$). 
Summing up $[\sum_{j}x_{j},y_{i'}]=[\sum_{j}y_{j},x_{i'}]=0$ for 
$i'=i,...,n$, we get $[a_{1},b_{i}]=[a_{i},b_{1}]=0$ (for any $i$). Finally, 
$c_{jk} = \sum_{\alpha=j}^{k-1}\sum_{\beta=k}^{n} t_{\alpha\beta}$
(in terms of the initial presentation) so the relations $[x_{i'},t_{\alpha\beta}]=0$
for $i'\neq\alpha,\beta$ and $[x_{\alpha}+x_{\beta},t_{\alpha\beta}]=0$
imply $[a_{i},c_{jk}]=0$ for $i\leq j\leq k$. Similarly, one shows
$[b_{i},c_{jk}]=0$ for $i\leq j\leq k$.

Let us prove that the relations (B) for $a_{i},b_{i}$ imply relations (A)
for $x_{i}:= a_{i}-a_{i+1}$, $y_{i}:= b_{i}-b_{i+1}$ (with the convention 
$a_{n+1}=b_{n+1}=0$). As before, $[a_{i},a_{j}]=[b_{i},b_{j}]=0$, $[a_{i},b_{j}]
= [a_{j},b_{i}]$ imply $[x_{i},x_{j}]=[y_{i},y_{j}]=0$, $[x_{i},y_{j}]
= [x_{j},y_{i}]$ (for any $i,j$). We set $t_{ij}:= [x_{i},y_{j}]$ 
for $i\neq j$, then we have $t_{ij}=t_{ji}$. 
We have for $j<k$, $t_{jk} = c_{jk}-c_{j,k+1}-c_{j+1,k}+c_{j+1,k+1}$ (we set $c_{i,n+1}:=0$), 
so $[a_{i},c_{jk}]=0$ implies $[\sum_{i'=i}^{n} x_{i'},t_{jk}]=0$ for $i\leq j< k$. 
When $i<j< k$, the difference between this relation and its analogue of 
$(i+1,j,k)$ gives $[x_{i},t_{jk}] = 0$ for $i<j<k$. This can be rewritten 
$[x_{i},[x_{j},y_{k}]]=0$ and since $[x_{i},x_{j}]=0$, we get $[x_{j},[x_{i},y_{k}]]=0$, 
so $[x_{j},t_{ik}]=0$ and by changing indices, $[x_{i},t_{jk}]=0$ for 
$j<i<k$. Rewriting again  $[x_{i},t_{jk}] = 0$ for $i<j< k$ as $[x_{i},[y_{j},x_{k}]]=0$
and using $[x_{i},x_{k}]=0$, we get $[x_{k},[x_{i},y_{j}]]=0$. i.e., $[x_{k},t_{ij}]=0$, 
which we rewrite $[x_{i},t_{jk}]=0$ for $j<k<i$. Finally, $[x_{i},t_{jk}]=0$ for 
$j<k$ and $i\notin\{j,k\}$, which implies $[x_{i},t_{jk}]=0$ for $i,j,k$ different. 
One proves similarly $[y_{i},t_{jk}]=0$ for $i,j,k$ different. 
\hfill \qed \medskip

\subsection{The formality of $\on{PB}_{1,n}$}

The flat connection $\on{d} - \sum_{i=1}^n K_i(\zz|\tau)\on{d}z_i$
gives rise to a monodromy representation $\mu_{\zz_{0},\tau} : \on{PB}_{1,n}
=\pi_1(C,\overline\zz_0)\to\exp(\hat\t_{1,n})$, 
which factors through a morphism $\mu_{\zz_{0},\tau}(\CC):\on{PB}_{1,n}(\CC)
\to\exp(\hat\t_{1,n})$. Let $\on{Lie}(\mu_{\zz_{0},\tau}):
\on{Lie}(\on{PB}_{1,n})_\CC \to\hat\t_{1,n}$ be  the corresponding morphism 
between pronilpotent Lie algebras.

\begin{proposition}
$\on{Lie}(\mu_{\zz_{0},\tau})$ is an isomorphism of filtered Lie algebras,
so that $\on{PB}_{1,n}$ is formal. 
\end{proposition}

{\em Proof.} As we have seen, $\on{Lie}(\on{PB}_{1,n})_\CC$ (denoted
$\on{Lie}(\on{PB}_{1,n})$ in this proof) is the 
quotient of the 
topologically free Lie algebra generated by $\alpha_{i},\beta_{i}$
($i=1,...,n$) by the topological ideal generated by $[\alpha_{i},\alpha_{j}]$, 
$[\beta_{i},\beta_{j}]$, $[\alpha_{1},\beta_{j}]$, $[\beta_{1},\alpha_{j}]$, 
$\on{log}(e^{\beta_{k}},e^{\alpha_{k}-\alpha_{j}}) 
- \on{log}(e^{\beta_{k}-\beta_{j}},e^{\alpha_{k}})$, 
$[\alpha_{i},\gamma_{jk}]$, $[\beta_{i},\gamma_{jk}]$
where $\gamma_{jk} = \on{log}(e^{\beta_{k}},e^{\alpha_{k}-\alpha_{j}}) $.

This presentation and the above presentation (B) of $\t_{1,n}$ imply that there is a 
morphism of graded Lie algebras $p_n:\t_{1,n}\to
\on{gr}\on{Lie}(\on{PB}_{1,n})$ 
defined by 
$a_i\mapsto [\alpha_{i}]$, $b_{i}\mapsto [\beta_{i}]$, where 
$\alpha\mapsto [\alpha]$ is the projection map $\on{Lie}(\on{PB}_{1,n})\to 
\on{gr}_{1}\on{Lie}(\on{PB}_{1,n})$.

$p_{n}$ is surjective because $\on{gr}\on{Lie}\Gamma$ is generated in degree 
$1$ (as the associated graded of any quotient of a topologically free Lie 
algebra).

There is a unique derivation $\tilde\Delta_{0}\in \on{Der}(\t_{1,n})$, such 
that $\tilde\Delta_{0}(x_{i}) = y_{i}$ and $\tilde\Delta_{0}(y_{i})=0$. 
This derivation gives rise to a one-parameter group of automorphisms
of $\on{Der}(\t_{1,n})$, defined by $\on{exp}(s\tilde\Delta_{0})(x_{i}):= 
x_{i} + s y_{i}$, $\on{exp}(s\tilde\Delta_{0})(y_{i}) = y_{i}$.

$\on{Lie}(\mu_{\zz_{0},\tau})$ induces a morphism 
$\on{gr} \on{Lie}(\mu_{\zz_{0},\tau}) : 
\on{gr}\on{Lie}(\on{PB}_{1,n})\to \t_{1,n}$. We will now 
prove that 
\begin{equation} \label{formality}
\textrm{gr}\on{Lie}(\mu_{\zz_{0},\tau})\circ p_n
=\on{exp}(- {{\tau}\over{2\pi\i}}\tilde\Delta_{0}) \circ w, 
\end{equation}
where $w$ is the automorphism of 
$\t_{1,n}$ defined by $w(a_{i}) = -b_{i}$, $w(b_{i}) = 2\pi\i a_{i}$.

$\mu_{\zz_{0},\tau}$ is defined as follows. Let $F_{\zz_{0}}(\zz)$ be the solution of 
$(\partial/\partial z_{i})F_{\zz_{0}}(\zz) = K_{i}(\zz|\tau)F_{\zz_{0}}(\zz)$, 
$F_{\zz_{0}}(\zz_{0})=1$ on $U_{\tau}$; let $H_{\tau} := \{\zz = (z_{1},...,z_{n}) | 
z_{i} = a_{i} + \tau b_{i}, 0<a_{1}<...<a_{n}<1\}$ and 
$V_{\tau} := \{\zz = (z_{1},...,z_{n}) | 
z_{i} = a_{i} + \tau b_{i}, 0<b_{1}<...<b_{n}<1\}$; let 
$F_{\zz_{0}}^{H}$ and $F_{\zz_{0}}^{V}$ be the analytic prolongations of 
$F_{\zz_{0}}$ to $H_{\tau}$ and $V_{\tau}$; then 
$$
F_{\zz_{0}}^{H}(\zz + \delta_{i}) = F_{\zz_{0}}^{H}(\zz)\mu_{\zz_{0},\tau}(X_{i}), 
\quad 
e^{2\pi\i x_{i}}F_{\zz_{0}}^{V}(\zz + \tau\delta_{i}) = 
F_{\zz_{0}}^{V}(\zz)\mu_{\zz_{0},\tau}(Y_{i}). 
$$ 
We have $\on{log}F_{\zz_{0}}(\zz) = -\sum_{i} (z_{i} - z_{i}^{0})y_{i}$ + terms of 
degree $\geq 2$, where $\t_{1,n}$ is graded by $\on{deg}(x_{i}) = \on{deg}(y_{i})=1$, 
which implies that $\on{log}\mu_{\zz_{0},\tau}(X_{i}) = -y_{i}$ + terms of degree $\geq 2$, 
$\on{log}\mu_{\zz_{0},\tau}(Y_{i}) = 2\pi\i x_{i} - \tau y_{i}$ + terms of degree $\geq 2$.
  Therefore $\on{Lie}(\mu_{\zz_{0},\tau})(\alpha_{i}) = 
\on{log}\mu_{\zz_{0},\tau}(A_{i}) = - b_{i}$ + terms of degree $\geq 2$, 
$\on{Lie}(\mu_{\zz_{0},\tau})(\beta_{i}) = 
\on{log}\mu_{\zz_{0},\tau}(B_{i}) = 2\pi\i a_{i} - \tau b_{i}$ + terms of degree $\geq 2$.
So $\on{gr}\on{Lie}(\mu_{\zz_{0},\tau})([\alpha_{i}]) = -b_{i}$, 
$\on{gr}\on{Lie}(\mu_{\zz_{0},\tau})([\beta_{i}]) = 
2\pi\i a_{i} - \tau b_{i}$.

It follows that $\on{gr}\on{Lie}(\mu_{\zz_{0},\tau}) \circ p_{n}$ is the endomorphism 
$a_{i}\mapsto -b_{i}$, $b_{i}\mapsto 2\pi\i a_{i} - \tau b_{i}$ of $\t_{1,n}$, 
which is the automorphism $\on{exp}( - {{\tau}\over{2\pi\i}}\tilde\Delta_{0}) 
\circ w$; this proves (\ref{formality}).

Since we already proved that $p_{n}$ is surjective,  it follows that 
$\on{gr}\on{Lie}(\mu_{\zz_{0},\tau})$ and $p_{n}$ are both isomorphisms. 
As $\on{Lie}(\on{PB}_{1,n})$ and $\hat\t_{1,n}$ are both complete and separated, 
$\on{Lie}(\mu_{\zz_{0},\tau})$ is bijective, and since it is a morphism, it is an 
isomorphism of filtered Lie algebras. 
\hfill \qed \medskip

\subsection{The formality of $\overline{\on{PB}}_{1,n}$}

Let $\zz_{0}\in U_{\tau}$ and $[\zz_{0}]\in \bar C(E_{\tau},n)$
be its image. We set $\overline{\on{PB}}_{1,n} := 
\pi_{1}(\bar C(E_{\tau},n),[\zz_{0}])$. Then $\overline{\on{PB}}_{1,n}$ 
is the quotient of $\on{PB}_{1,n}$ by its central 
subgroup (isomorphic to $\ZZ^{2}$) generated by $A_{1}$ and $B_{1}$. 
We have $\mu_{\zz_{0},\tau}(A_{1}) = e^{-\sum_{i}y_{i}}$
and $\mu_{\zz_{0},\tau}(B_{1}) = e^{2\pi\i \sum_{i}x_{i} 
- \tau\sum_{i}y_{i}}$, so $\on{Lie}(\mu_{\zz_{0},\tau})(\alpha_{1}) = -a_{1}$, 
$\on{Lie}(\mu_{\zz_{0},\tau})(\beta_{1}) = 2\pi\i a_{1} - \tau b_{1}$, which 
implies that $\on{Lie}(\mu_{\zz_{0},\tau})$ induces an isomorphism 
between $\on{Lie}(\overline{\on{PB}}_{1,n})_\CC$ and $\bar\t_{1,n}$. 
In particular, $\overline{\on{PB}}_{1,n}$ is formal.

\begin{remark}
Let $\on{Diag}_{n}:=\{(\zz,\tau)\in\CC^{n}\times\HH | \zz\in\on{Diag}_{n,\tau}\}$ and
let $U \subset (\CC^n\times\HH) - \on{Diag}_n$ be the set of all 
$(\zz,\tau)$ such that $\zz\in U_\tau$. 
Each element of $U$ gives rise to a Lie algebra isomorphism 
$\mu_{\zz,\tau} : \on{Lie}(\on{PB}_{1,n}) \simeq \hat\t_{1,n}$. 
For an infinitesimal  $(\on{d}\zz,\on{d}\tau)$, the composition 
$\mu_{\zz + \on{d}\zz,\tau + \on{d}\tau} \circ \mu_{\zz,\tau}^{-1}$ 
is then an infinitesimal automorphism of $\hat\t_{1,n}$. 
This defines a flat connection over $U$ with values in the trivial Lie 
algebra bundle 
with Lie algebra $\on{Der}(\hat\t_{1,n})$. When $\on{d}\tau=0$, the infinitesimal 
automorphism has the form  $\on{exp}(\sum_{i} K_{i}(\zz|\tau)\on{d}z_{i})$, 
so the connection has the form $\on{d} - \sum_{i} \on{ad}(K_{i}(\zz|\tau))
\on{d}z_{i}  - \tilde\Delta(\zz|\tau)\on{d}\tau$, where $\tilde\Delta : 
U \to \on{Der}(\hat\t_{1,n})$ 
is a meromorphic map with poles at $\on{Diag}_{n}$. In the next section, 
we determine
a map $\Delta : (\CC^{n}\times\HH) - \on{Diag}_{n} \to \on{Der}(\hat\t_{1,n})$ 
with the same flatness properties as $\tilde\Delta(\zz|\tau)$. 
\end{remark}

\subsection{The isomorphisms $\on{B}_{1,n}(\CC)\simeq \on{exp}(\hat\t_{1,n})
\rtimes S_{n}$, $\overline{\on{B}}_{1,n}(\CC) \simeq 
\on{exp}(\hat{\bar\t}_{1,n})\rtimes S_{n}$}

Let $\zz_{0}$ be as above; we define $\on{B}_{1,n}:= \pi_{1}(C(E_{\tau},[n]),[\zz_{0}])$
and $\overline{\on{B}}_{1,n}:= \pi_{1}(\bar C(E_{\tau},[n]),[\overline\zz_{0}])$, 
where $x\mapsto [x]$ is the canonical projection $C(E_{\tau},n)\to C(E_{\tau},[n])$
or $\bar C(E_{\tau},n)\to \bar C(E_{\tau},[n])$.

We have an exact sequence $1\to \on{PB}_{1,n} \to \on{B}_{1,n} \to S_{n}\to 1$, 
We then define groups $\on{B}_{1,n}(\CC)$ fitting in an exact sequence
$1\to \on{PB}_{1,n}(\CC) \to \on{B}_{1,n}(\CC) \to S_{n}\to 1$
as follows: the morphism $\on{B}_{1,n}\to \on{Aut}(\on{PB}_{1,n})$
extends to $\on{B}_{1,n} \to \on{Aut}(\on{PB}_{1,n}(\CC))$; we then 
construct the semidirect product $\on{PB}_{1,n}(\CC)\rtimes \on{B}_{1,n}$; 
then $\on{PB}_{1,n}$ embeds diagonally as a normal subgroup of this 
semidirect product, and $\on{B}_{1,n}(\CC)$ is defined as the quotient 
$(\on{PB}_{1,n}(\CC)\rtimes \on{B}_{1,n}) / \on{PB}_{1,n}$.

The monodromy of $\nabla_{\tau,[n]}$ then gives rise to a group morphism 
$\on{B}_{1,n}\to \on{exp}(\hat\t_{1,n})\rtimes S_{n}$, which factors 
through $\on{B}_{1,n}(\CC)\to \on{exp}(\hat\t_{1,n})\rtimes S_{n}$. 
Since this map commutes with the natural morphisms to $S_{n}$ and using the 
isomorphism $\on{PB}_{1,n}(\CC)\simeq \on{exp}(\hat\t_{1,n})$, we obtain that 
$\on{B}_{1,n}(\CC)\to \on{exp}(\hat\t_{1,n})\rtimes S_{n}$ is an isomorphism.

Similarly, starting from the exact sequence 
$1\to \overline{\on{PB}}_{1,n} \to \overline{\on{B}}_{1,n}\to S_{n}\to 1$
one defines a group $\overline{\on{B}}_{1,n}(\CC)$ fitting in an exact sequence 
$1\to \overline{\on{PB}}_{1,n} \to \overline{\on{B}}_{1,n}(\CC)\to S_{n}\to 1$
together with an isomorphism 
$\overline{\on{B}}_{1,n}(\CC)\to \on{exp}(\hat{\bar\t}_{1,n})\rtimes S_{n}$.

\section{Bundles with flat connection on ${\cal M}_{1,n}$ and 
${\cal M}_{1,[n]}$} \label{sect:4}

We first define Lie algebras of derivations of $\bar\t_{1,n}$
and a related group ${\bold G}_{n}$. We then define a principal 
${\bold G}_{n}$-bundle with flat connection of 
${\cal M}_{1,n}$ and a principal ${\bold G}_{n}\rtimes S_{n}$-bundle 
with flat connection on the moduli space ${\cal M}_{1,[n]}$ of elliptic 
curves with $n$ unordered marked points.

\subsection{Derivations of the Lie algebras $\t_{1,n}$ and $\bar\t_{1,n}$ 
and associated groups}

Let $\d$ be the Lie algebra with generators $\Delta_0,d,X$ and $\delta_{2m}$
($m\geq 1$), and relations: 
$$
[d,X] = 2X, \quad [d,\Delta_0] = -2\Delta_0, \quad [X,\Delta_0]=d, 
$$
$$
[\delta_{2m},X] = 0, \quad [d,\delta_{2m}] = 2m\delta_{2m}, \quad 
\on{ad}(\Delta_0)^{2m+1}(\delta_{2m}) = 0. 
$$

\begin{proposition}
We have a Lie algebra morphism $\d\to \on{Der}(\t_{1,n})$,
denoted by $\xi\mapsto \tilde\xi$, 
such that
$$
\tilde d(x_i) = x_i, \tilde d(y_i)=-y_i, \tilde d(t_{ij}) = 0, \quad
\tilde X(x_i) = 0, \tilde X(y_i)=x_i, \tilde X(t_{ij}) = 0, 
$$ 
$$
\tilde \Delta_0(x_i)=y_i, \tilde \Delta_0(y_i)=0, \tilde \Delta_0(t_{ij})=0, 
$$
$$
\tilde\delta_{2m}(x_i)=0, \tilde \delta_{2m}(t_{ij}) = [t_{ij},
(\on{ad}x_i)^{2m}(t_{ij})], \tilde \delta_{2m}(y_i) = \sum_{j|j\neq i}
{1\over 2} \sum_{p+q=2m-1} [(\on{ad}x_i)^p(t_{ij}),(-\on{ad}x_i)^q(t_{ij})].
$$
This induces a Lie algebra morphism $\d\to \on{Der}(\bar\t_{1,n})$. 
\end{proposition}

{\em Proof.}
The fact that $\tilde\Delta_0,\tilde d,\tilde X$ are derivations and commute
according to the Lie bracket of $\sl_2$ is clear.

Let us prove that  $\tilde\delta_{2m}$ is a derivation. 
We have $\tilde\delta_{2m}(t_{ij}) = [t_{ij},\sum_{i<j} 
(\ad x_i)^{2m}(t_{ij})]$, which implies that $\tilde\delta_{2m}$ preserves 
the infinitesimal pure braid identities. It clearly preserves the relations
$[x_i,x_j]=0$, $[x_i,y_j]=t_{ij}$, $[x_k,t_{ij}] = 0$, 
$[x_i+x_j,t_{ij}]=0$.

Let us prove that $\tilde\delta_{2m}$ preserves the relation $[y_k,t_{ij}] = 0$.
\begin{align*}
& [\tilde\delta_{2m}(y_k),t_{ij}] 
= {1\over 2}\sum_{p+q = 2m-1}(-1)^q 
[[(\ad x_k)^p(t_{ki}), (\ad x_k)^q(t_{ki})] 
+ [(\ad x_k)^p(t_{kj}), (\ad x_k)^q(t_{kj})],t_{ij}] 
\\ & = 
{1\over 2}\sum_{p+q = 2m-1}(-1)^{q+1} 
[[(\ad x_k)^p(t_{ki}), (\ad x_k)^q(t_{kj})] 
+ [(\ad x_k)^p(t_{kj}), (\ad x_k)^q(t_{ki})],t_{ij}] 
\\ & = 
\sum_{p+q = 2m-1}(-1)^{q+1} 
[[(\ad x_k)^p(t_{ki}), (\ad x_k)^q(t_{kj})],t_{ij}]
= 
[t_{ij}, 
\sum_{p+q=2m-1} (-1)^p (\ad x_i)^p(\ad x_j)^q([t_{ki},t_{kj}])].
\end{align*}

On the other hand, $[y_k,\tilde\delta_{2m}(t_{ij})] = 
[y_k,[t_{ij},(\ad x_i)^{2m}(t_{ij})]]
= [t_{ij},[y_k,(\ad x_i)^{2m}(t_{ij})]]$. 
Now 
\begin{align*}
& [y_k,(\ad x_i)^{2m}(t_{ij})] = 
-\sum_{\alpha + \beta = 2m-1} (\ad x_i)^\alpha 
\big( [t_{ki} , (\ad x_i)^\beta(t_{ij}) ] \big) 
\\ & = -\sum_{\alpha + \beta = 2m-1} (\ad x_i)^\alpha
[t_{ki},(-\ad x_j)^\beta(t_{ij})]
= -\sum_{\alpha + \beta = 2m-1}
(\ad x_i)^\alpha (-\ad x_j)^\beta([t_{ki},t_{ij}])
\\ & = \sum_{p+q = 2m-1}
(-1)^{p+1} (\ad x_i)^p(\ad x_j)^q ([t_{ki},t_{kj}]). 
\end{align*}
Hence $[\tilde\delta_{2m}(y_k),t_{ij}] + [y_k,\tilde\delta_{2m}(t_{ij})] = 0$.

Let us prove that $\tilde\delta_{2m}$ preserves the relation $[y_i,y_j]=0$, 
i.e., that $[\tilde\delta_{2m}(y_i),y_j] + [y_i,\tilde\delta_{2m}(y_j)] = 0$.

We have 
\begin{align*}
& [y_i,\tilde\delta_{2m}(y_j)] = {1\over 2}[y_i,\sum_{p+q = 2m-1} (-1)^q 
[(\ad x_j)^p(t_{ji}),(\ad x_j)^q(t_{ji})]]
\\ & + {1\over 2}\sum_{k\neq i,j}
[y_i,\sum_{p+q = 2m-1} (-1)^q 
[(\ad x_j)^p(t_{jk}),(\ad x_j)^q(t_{jk})]]. 
\end{align*}
Now 
\begin{align} \label{part}
& {1\over 2}[y_i,\sum_{p+q = 2m-1} (-1)^q 
[(\ad x_j)^p(t_{ji}),(\ad x_j)^q(t_{ji})]] - (i\leftrightarrow j)
\\ & \nonumber = - {1\over 2}[y_i+y_j,\sum_{p+q = 2m-1} (-1)^q 
[(\ad x_i)^p(t_{ij}),(\ad x_i)^q(t_{ij})]]
\\ &  \nonumber 
=  \sum_{p+q = 2m-1} (-1)^{q+1}
  [[y_i + y_j, (\ad x_i)^p(t_{ij})],(\ad x_i)^q(t_{ij})].
\end{align}

A computation similar to the above computation of 
$[y_k,(\ad x_i)^{2m}(t_{ij})]$ yields
$$
[y_i+y_j,(\ad x_i)^p(t_{ij})] = (-1)^p \sum_{\alpha+\beta=p-1} 
[(\ad x_k)^\alpha(t_{ik}),(\ad x_j)^\beta(t_{jk})], 
$$
so 
$$
(\ref{part}) = \sum_{\alpha+\beta+\gamma=2m-2}
[(\ad x_i)^\alpha(t_{ij}),[(\ad x_k)^\beta(t_{ik}),
(\ad x_j)^\gamma(t_{jk})]]. 
$$

If now $k\neq i,j$, then 
$$
[y_i,{1\over 2}\sum_{p+q = 2m-1} (-1)^q 
[(\ad x_j)^p(t_{jk}),(\ad x_j)^q(t_{jk})]] = 
\sum_{p+q = 2m-1} (-1)^q [[y_i,(\ad x_j)^p(t_{jk})],(\ad x_j)^q(t_{jk})]. 
$$ 
As we have seen, 
\begin{align*}
& [y_j,(\ad x_i)^p(t_{ik})]  = (-1)^p \sum_{\alpha + \beta = p-1}
(-\ad x_i)^\alpha (\ad x_k)^\beta[t_{ij},t_{ik}]
\\ & 
= (-1)^{p+1} \sum_{\alpha + \beta = p-1} [(-\ad x_i)^\alpha(t_{ij}),
(\ad x_k)^\beta(t_{jk})]
\end{align*}
So we get 
\begin{align*}
& [y_i,{1\over 2}\sum_{p+q = 2m-1} (-1)^q 
[(\ad x_j)^p(t_{jk}),(\ad x_j)^q(t_{jk})]] 
\\ & 
= \sum_{\alpha + \beta+\gamma = 2m-2}
[[(\ad x_i)^{\alpha}(t_{ij}),(\ad x_k)^{\beta}(t_{ik})],
(\ad x_j)^\gamma(t_{jk})]
\end{align*}
therefore
\begin{align*}
& [y_i,{1\over 2} \sum_{p+q=2m-1}(-1)^q
[(\ad x_j)^p(t_{jk}),(\ad x_j)^q(t_{jk})]] - (i\leftrightarrow j) 
\\ & = \sum_{\alpha+\beta+\gamma = 2m-2} [(\ad x_i)^\alpha(t_{ij}), 
[(\ad x_k)^{\beta}(t_{ik}),(\ad x_j)^{\gamma}(t_{jk})]]. 
\end{align*}
Therefore $[y_i,\tilde\delta_{2m}(y_j)] + [\tilde\delta_{2m}(y_i),y_j] = 0$.

Since $\tilde\delta_{2m}(\sum_i x_i) = \tilde\delta_{2m}(\sum_i y_i) = 0$ 
and $\sum_i x_i$
and $\sum_i y_i$ are central, $\tilde\delta_{2m}$ preserves the relations
$[\sum_i x_i,y_j] = 0$ and $[\sum_k x_k,t_{ij}] = [\sum_k y_k,t_{ij}] = 0$.
It follows that $\tilde\delta_{2m}$ preserves the relations $[x_i+x_j,t_{ij}]
= [y_i + y_j,t_{ij}] = 0$ and $[x_i,y_i] = -\sum_{j|j\neq i} t_{ij}$.  All this 
proves that $\tilde\delta_{2m}$ is a derivation.

Let us show that 
$\on{ad}(\tilde\Delta_0)^{2m+1}(\tilde\delta_{2m}) = 0$ for $m\geq 1$. 
We have 
\begin{align*}
& \on{ad}(\tilde\Delta_0)^{2m+1}(\tilde\delta_{2m})(x_i) = - (2m+1) 
\tilde\Delta_0^{2m} \circ \tilde\delta_{2m} \circ \tilde\Delta_0(x_i) = 
- (2m+1) \tilde\Delta_0^{2m} \circ \tilde\delta_{2m} (y_i) 
\\
  & = - (2m+1) \tilde\Delta_0^{2m} 
(\sum_{j|j\neq i} {1\over 2} \sum_{p+q=2m-1}[(\on{ad}x_i)^p(t_{ij}),
(-\on{ad}x_i)^q(t_{ij})])= 0;
\end{align*} 
the last part of this computation 
implies that $\on{ad}(\tilde\Delta_0)^{2m+1}(\tilde\delta_{2m})(y_i) = 0$, 
therefore $\on{ad}(\tilde\Delta_0)^{2m+1}(\tilde\delta_{2m}) = 0$.

We have clearly $[\tilde X,\tilde\delta_{2m}]=0$ and 
$[\tilde d,\tilde\delta_{2m}]= 2m\tilde\delta_{2m}$. It follows that we 
have a Lie algebra morphism $\d\to\on{Der}(\t_{1,n})$. 
Since $\tilde d,\tilde\Delta_{0},\tilde X$ and $\tilde\delta_{2m}$ all map 
$\CC(\sum_{i}x_{i})\oplus \CC(\sum_{i}y_{i})$ to itself, this induces a 
Lie algebra morphism $\d\to\on{Der}(\bar\t_{1,n})$. 
\hfill \qed \medskip

Let $e,f,h$ be the standard basis of $\sl_2$. Then we have a 
Lie algebra morphism $\d\to \sl_2$, defined by $\delta_{2n}\mapsto 0$, 
$d\mapsto h$, $X\mapsto e$, $\Delta_0\mapsto f$. We denote by 
$\d_+\subset \d$ its kernel.

Since the morphism $\d\to \sl_2$ has a section (given by 
$e,f,h\mapsto X,\Delta_0,d$), we have a semidirect product 
decomposition $\d = \d_+ \rtimes \sl_2$.

We then have 
$$
\bar\t_{1,n}\rtimes \d = (\bar\t_{1,n}\rtimes \d_+) \rtimes \sl_2. 
$$

\begin{lemma} 
$\bar\t_{1,n} \rtimes \d_+$ is positively graded. 
\end{lemma}

{\em Proof.} We define compatible $\ZZ^2$-gradings of $\d$ and 
$\bar\t_{1,n}$ by 
$\on{deg}(\Delta_0) = (-1,1)$, $\on{deg}(d)=(0,0)$, $\on{deg}(X) = (1,-1)$, 
$\on{deg}(\delta_{2m}) = (2m+1,1)$, $\on{deg}(x_i) = (1,0)$, 
$\on{deg}(y_i)=(0,1)$, 
$\on{deg}(t_{ij}) = (1,1)$.

We define the support of $\d$ (resp., $\bar\t_{1,n}$) as the subset of 
$\ZZ^2$ of indices for which the corresponding component 
of $\d$ (resp., $\bar\t_{1,n}$) is nonzero.

Since the $\bar x_i$ on one hand, the $\bar y_i$ on the other hand generate 
abelian Lie subalgebras of $\bar\t_{1,n}$, the support of $\bar\t_{1,n}$
is contained in $\NN_{>0}^2 \cup \{(1,0),(0,1)\}$.

On the other hand, $\d_+$
is generated by the $\on{ad}(\Delta_0)^p(\delta_{2m})$, which 
all have degrees in $\NN_{>0}^{2}$. 
It follows that the support of $\d_+$ is contained in 
$\NN_{>0}^2$.

Therefore the support of $\bar\t_{1,n} \rtimes \d_+$ is contained in 
$\NN_{>0}^2\cup \{(1,0),(0,1)\}$, so this Lie algebra
is positively graded. \hfill \qed \medskip

\begin{lemma}
$\bar\t_{1,n}\rtimes\d_+$ is a sum of finite dimensional 
$\sl_2$-modules; $\d_+$ is a sum of irreducible 
odd dimensional $\sl_2$-modules. 
\end{lemma}

{\em Proof.} A generating space for $\bar\t_{1,n}$
is $\sum_{i} (\CC \bar x_i\oplus \CC\bar y_i)$, which is a sum of
finite dimensional $\sl_{2}$-modules, so $\bar\t_{1,n}$ is a sum of 
finite dimensional $\sl_2$-modules.

A generating space for $\d_+$ is the sum over $m\geq 1$ of its 
$\sl_2$-submodules generated by the $\delta_{2m}$, which are zero or 
irreducible odd dimensional, therefore $\d_+$ is a sum of odd dimensional 
$\sl_2$-modules. (In fact, the $\sl_2$-submodule generated 
by $\delta_{2m}$ is nonzero, as it follows from the construction of the 
above morphism $\d_+\to \on{Der}(\bar\t_{1,n})$ that $\delta_{2m}\neq 0$.)
\hfill \qed \medskip

It follows that $\bar\t_{1,n}$, $\bar\d_+$ and $\bar\t_{1,n}\rtimes \d_+$
integrate to $\on{SL}_2(\CC)$-modules (while $\bar\d_+$ even integrates
to a $\on{PSL}_2(\CC)$-module).

We can form in particular the semidirect products 
$$
{\bold G}_n := 
\on{exp}((\bar\t_{1,n} \rtimes \d_+)^\wedge) \rtimes \on{SL}_2(\CC) 
$$
and $\on{exp}(\hat\d_+)\rtimes \on{PSL}_2(\CC)$; we have morphisms 
${\bold G}_n \to\on{exp}(\hat\d_+)\rtimes \on{PSL}_2(\CC)$ (this is 
a 2-covering if $n=1$ since $\bar\t_{1,1}=0$).

Observe that the action of $S_{n}$ by automorphisms of $\bar\t_{1,n}$
extends to an action on $\bar\t_{1,n}\rtimes\d$, where the action on $\d$
is trivial. This gives rise to an action of $S_{n}$ by automorphisms of 
${\bold G}_{n}$.

\subsection{Bundle with flat connection on ${\cal M}_{1,n}$}

The semidirect product $((\ZZ^n)^2\times\CC)\rtimes \on{SL}_2(\ZZ)$
acts on $(\CC^{n}\times\HH) - \on{Diag}_{n}$ by 
$(\nn,\mm,u) * (\zz,\tau) := (\nn+ \tau\mm + u(\sum_i \delta_i),\tau)$
for $(\nn,\mm,u)\in (\ZZ^n)^2 \times\CC$ and 
$\bigl( \begin{smallmatrix} \alpha & \beta \\ \gamma & \delta
\end{smallmatrix}\bigl) * 
(\zz,\tau) := ({\zz\over{\gamma\tau+\delta}},{{\alpha\tau+\beta}\over
{\gamma\tau+\delta}})$ for 
$\bigl( \begin{smallmatrix} \alpha & \beta \\ \gamma & \delta
\end{smallmatrix}\bigl) \in
\on{SL}_2(\ZZ)$
(here $\on{Diag}_{n} := \{(\zz,\tau)\in \CC^{n}\times\HH |$ 
for some $i\neq j, z_{ij}\in \Lambda_\tau\}$). 
The quotient then identifies with the moduli space ${\cal M}_{1,n}$
of elliptic curves with $n$ marked points.

Set ${\bold G}_n := \on{exp}((\bar\t_{1,n}\rtimes\d_+)^\wedge)\rtimes
\on{SL}_2(\CC)$. We will define a principal ${\bold G}_n$-bundle 
with flat connection 
$({\cal P}_n,\nabla_{{\cal P}_n})$ over ${\cal M}_{1,n}$.

For $u\in\CC^{\times}$, $u^{d} := \bigl( \begin{smallmatrix} u & 0 \\
  0 & u^{-1}\end{smallmatrix}\bigl)\in \on{SL}_{2}(\CC)\subset {\bold G}_{n}$
  and for $v\in \CC$, $e^{vX}:= \bigl( \begin{smallmatrix} 1 & v \\
   0 & 1 \end{smallmatrix}\bigl) \in \on{SL}_{2}(\CC)\subset {\bold G}_{n}$. 
Since $[X,\bar x_{i}]=0$, we consistently set $\on{exp}
(aX + \sum_{i}b_{i}\bar x_{i}) := 
\on{exp}(aX)\on{exp}(\sum_{i} b_{i}\bar x_{i})$.

\begin{proposition} \label{prop:bundle}
There exists a unique principal ${\bold G}_{n}$-bundle ${\cal P}_{n}$
over ${\cal M}_{1,n}$, such that a section of 
$U\subset {\cal M}_{1,n}$ is a function $f : \pi^{-1}(U) 
\to {\bold G}_n$ (where $\pi : (\CC^{n}\times\HH) - \on{Diag}_{n} \to 
{\cal M}_{1,n}$ is the canonical projection), such that 
$f(\zz+\delta_i|\tau) = f(\zz + u(\sum_i\delta_i)|\tau) = f(\zz|\tau)$, 
$f(\zz+\tau\delta_i|\tau) = e^{-2\pi\i\bar x_i} f(\zz|\tau)$, 
$f(\zz|\tau+1) = f(\zz|\tau)$ and 
$f({\zz\over\tau}|-{1\over\tau}) = \tau^{d}\on{exp}
({{2\pi\i}\over\tau}(\sum_{i} z_{i}\bar x_{i} + X)) f(\zz|\tau)$. 
\end{proposition}

{\em Proof.} Let $c_{\tilde g} : \CC^{n}\times\HH \to {\bold G}_{n}$ be a 
family of holomorphic functions (where $\tilde g\in 
((\ZZ^{n})^{2}\times \CC)\rtimes \on{SL}_{2}(\ZZ)$)
satisfying the cocycle condition $c_{\tilde g\tilde g'}(\zz|\tau) 
= c_{\tilde g}(\tilde g'*(\zz|\tau))c_{\tilde g'}(\zz|\tau)$. 
Then there exists a unique principal ${\bold G}_{n}$-bundle over 
${\cal M}_{1,n}$
such that a section of $U\subset {\cal M}_{1,n}$ is a function 
$f : \pi^{-1}(U)\to {\bold G}_{n}$ such that 
$f(\tilde g * (\zz|\tau)) = c_{\tilde g}(\zz|\tau)f(\zz|\tau)$.

We will now prove that there is a unique cocycle such that 
$c_{(u,0,0)} = c_{(0,\delta_{i},0)} = 1$, 
$c_{(0,0,\delta_{i})} = e^{-2\pi\i\bar x_{i}}$, $c_{S} = 1$ and 
$c_{T}(\zz|\tau) = \tau^{d} \on{exp}({{2\pi\i}\over\tau}
(\sum_{i}z_{i}\bar x_{i} + X))$, where 
$S = \bigl(\begin{smallmatrix} 1 & 1 \\ 0 & 1 \end{smallmatrix}\bigl)$, 
$T = \bigl(\begin{smallmatrix} 0 & -1 \\  1 & 0 \end{smallmatrix}\bigl)$.

Such a cocycle is the same as a family of functions 
$c_{g} : \CC^{n}\times\HH\to {\bold G}_{n}$
(where $g\in\on{SL}_{2}(\ZZ)$), satisfying the cocycle conditions 
$c_{gg'}(\zz|\tau) = c_{g}(g' * (\zz|\tau))
c_{g'}(\zz|\tau)$ for $g,g'\in\on{SL}_{2}(\ZZ)$, and 
$c_g(\zz + \delta_i|\tau) = e^{2\pi\i \gamma\bar x_i}c_g(\zz|\tau)$, 
$c_g(\zz + \tau\delta_i|\tau) = e^{-2\pi\i \delta\bar x_i} 
c_g(\zz|\tau) e^{2\pi\i\bar x_i}$
and $c_g(\zz + u(\sum_i\delta_i)|\tau)=c_g(\zz|\tau)$ for 
$g = \big(\begin{smallmatrix}\alpha & \beta \\
  \gamma & \delta\end{smallmatrix}\big)\in\on{SL}_{2}(\ZZ)$.

\begin{lemma}
There exists a unique family of functions 
$c_{g} : \CC^{n}\times \HH\to {\bold G}_{n}$ such that 
$c_{gg'}(\zz|\tau) = c_{g}(g'*(\zz|\tau))
c_{g'}(\zz|\tau)$ for $g,g'\in\on{SL}_{2}(\ZZ)$, with 
$$
c_{S}(\zz|\tau) = 1, \quad 
c_{T}(\zz|\tau) = \tau^d e^{(2\pi\i/\tau)(\sum_j z_j \bar x_j + X)}. 
$$
\end{lemma}

{\em Proof.} $\on{SL}_{2}(\ZZ)$ is the group generated by $\tilde S$, 
$\tilde T$ and relations $\tilde T^{4} = 1$, $(\tilde S\tilde T)^{3}
=\tilde T^{2}$, $\tilde S\tilde T^2 = \tilde T^{2}\tilde S$. 
Let $\langle \tilde S,\tilde T\rangle$ be the free group with generators 
$\tilde S$, $\tilde T$; then there is a unique family of maps 
$c_{\tilde g} : \CC^{n}\times\HH \to {\bold G}_{n}$, 
$\tilde g\in \langle \tilde S,\tilde T\rangle$ satisfying the cocycle 
conditions (w.r.t. the action of  $\langle \tilde S,\tilde T\rangle$ 
on $\CC^{n}\times\HH$ through its quotient $\on{SL}_{2}(\ZZ)$)
and $c_{\tilde S} = c_{S}$, $c_{\tilde T} = c_{T}$. 
It remains to show that $c_{\tilde T^{4}} = 1$, 
$c_{(\tilde S\tilde T)^{3}} = c_{\tilde T^{2}}$ and 
$c_{\tilde S\tilde T^{2}} = c_{\tilde T^{2}\tilde S}$.

For this, we show that $c_{\tilde T^{2}}(\zz|\tau) = (-1)^{d}$. We have 
$c_{\tilde T^{2}}(\zz|\tau) = c_{T}(\zz/\tau|-1/\tau) c_{T}(\zz|\tau) = 
(-\tau)^{-d} \on{exp}(-2\pi\i\tau(\sum_{j} (z_{j}/\tau)\bar x_{j} + X)) 
\tau^{d}\on{exp}({{2\pi\i}\over\tau}(\sum_{j} z_{j}\bar x_{j} + X)) = (-1)^{d}$
since $\tau^{d}X\tau^{-d} = \tau^{2}X$, $\tau^{d}\bar x_{i}\tau^{-d} = \tau\bar x_{i}$.

Since $((-1)^{d})^{2} = 1^{d}=1$, we get $c_{\tilde T^{4}}=1$. Since $c_{\tilde S}$ and 
$c_{\tilde T^{2}}$ are both constant and commute, we also get $c_{\tilde S\tilde T^{2}} = c_{\tilde T^{2}\tilde S}$.

We finally have $c_{\tilde S\tilde T}(\zz|\tau) = c_{T}(\zz|\tau)$ while $\tilde S\tilde T = \bigl( \begin{smallmatrix} 
1  & -1\\  1 & 0\end{smallmatrix}\bigl)$, $(\tilde S\tilde T)^{2} = \bigl( \begin{smallmatrix} 
0  & -1\\  1 & -1\end{smallmatrix}\bigl)$ so 
\begin{align*}
& c_{(\tilde S\tilde T)^{3}}(\zz|\tau) 
= c_{T}({\zz\over{\tau-1}}|{1\over{1-\tau}}) c_{T}({\zz\over\tau}|{{\tau-1}\over\tau}) c_{T}(\zz|\tau)
= ({1\over{1-\tau}})^{d} \on{exp}(-2\pi\i\sum z_{j}\bar x_{j} + 2\pi\i (1-\tau)X) 
\\ & 
({{\tau-1}\over\tau})^{d} \on{exp}({{2\pi\i}\over{\tau-1}} \sum_{j}z_{j}\bar x_{j}
+ 2\pi\i {{\tau}\over{\tau-1}}X) \tau^{d} 
\on{exp}({{2\pi\i}\over\tau}(\sum_{j} z_{j}\bar x_{j} + X))
\\
  & = (-1)^{d} \on{exp}({{2\pi\i}\over {1-\tau}}(\sum_{j}z_{j}\bar x_{j} + X))
  \on{exp}({{2\pi\i}\over{\tau(\tau-1)}}(\sum_{j} z_{j}\bar x_{j} + X))
\on{exp}({{2\pi\i}\over\tau}(\sum_{j} z_{j}\bar x_{j} + X)) =(-1)^{d}, 
\end{align*}
so $c_{(\tilde S\tilde T)^{3}} = c_{\tilde T^{2}}$. \hfill \qed\medskip

{\em End of proof of Proposition \ref{prop:bundle}.} 
We now check that the maps $c_{g}$ satisfy the remaining conditions, i.e., 
$c(\zz + u(\sum_{i}\delta_{i})|\tau) = c_{g}(\zz|\tau)$, 
$c_{g}(\zz + \delta_{i}|\tau) = e^{2\pi\i\gamma\bar x_{i}}c_{g}(\zz|\tau)$, 
$c_{g}(\zz+\tau\delta_{i}|\tau) = e^{-2\pi\i \delta\bar x_{i}}c_{g}(\zz|\tau)e^{2\pi\i \bar x_{i}}$. 
The cocycle identity $c_{gg'}(\zz|\tau) = c_{g}(g'*(\zz|\tau))c_{g'}(\zz|\tau)$ implies that 
it suffices to prove these identities for $g = S$ and $g=T$. They are trivially satisfied if $g=S$. 
When $g=T$, the first identity follows from $\sum_{i}\bar x_{i}=0$, the third identity 
follows from the fact that $(X,\bar x_{1},...,\bar x_{n})$ is a commutative family, the second 
identity follows from the same fact together with $\tau^{d}\bar x_{i}\tau^{-d} = \tau\bar x_{i}$. 
\hfill \qed\medskip

Set 
$$
g(z,x|\tau) := {{\theta(z+x|\tau)}\over{\theta(z|\tau)\theta(x|\tau)}}
\big( {\theta'\over\theta}(z+x|\tau) - {\theta'\over\theta}(x|\tau)\big)
+{1\over x^{2}} = k_x(z,x|\tau),
$$
(we set $f'(z|\tau):= (\partial/\partial z)f(z|\tau)$).

We have $g(z,x|\tau)\in \on{Hol}((\CC\times \HH) - \on{Diag}_{1})[[x]]$, therefore
$g(z,\on{ad}\bar x_{i}|\tau)$ is a linear map $\bar\t_{1,n} \to 
(\on{Hol}((\CC\times\HH) - \on{Diag}_{1})\otimes \bar\t_{1,n})^{\wedge}$, 
so $g(z,\on{ad}\bar x_{i}|\tau)(\bar t_{ij}) \in
(\on{Hol}((\CC\times\HH) - \on{Diag}_{1})\otimes \bar\t_{1,n})^{\wedge}$. 
Therefore $$g(\zz|\tau) := \sum_{i<j} g(z_{ij},\on{ad}\bar x_i|\tau)(\bar t_{ij})$$
is a meromorphic function $\CC^{n}\times\HH\to \hat{\bar\t}_{1,n}$ with only poles at $\on{Diag}_{n}$.

We set 
$$
\bar\Delta(\zz|\tau) := - {1\over{2\pi\i}}\Delta_0  -{1\over{2\pi\i}} \sum_{n\geq 1} a_{2n} E_{2n+2}(\tau)
\delta_{2n} +{1\over{2\pi\i}} g(\zz|\tau) , 
$$
where $a_{2n} = - (2n+1) B_{2n+2} (2\i\pi)^{2n+2}/(2n+2)!$ and $B_n$
are the Bernoulli numbers given by $x/(e^x-1) = \sum_{r\geq 0} 
(B_r/r!)x^r$. This is a meromorphic function $\CC^{n}\times\HH\to 
(\bar\t_{1,n}\rtimes\d_{+})^{\wedge}\rtimes\n_{+}
\subset \on{Lie}({\bold G}_{1,n})$ (where $\n_{+} = \CC\Delta_{0}\subset \sl_{2}$)
with only poles at $\on{Diag}_{n}$.

For $\psi(x) = \sum_{n\geq 1}b_{2n}x^{2n}$, we set $\delta_{\psi}:= \sum_{n\geq 1}
b_{2n}\delta_{2n}$, $\Delta_{\psi}:= \Delta_{0} + 
\sum_{n\geq 1} b_{2n}\delta_{2n}$. If we set 
$$
\varphi(x|\tau) = - x^{-2} - (\theta'/\theta)'(x|\tau) 
+ (x^{-2} + (\theta'/\theta)'(x|\tau))_{|x=0} 
= g(0,0|\tau)-g(0,x|\tau), 
$$
then $\varphi(x|\tau) = \sum_{n\geq 1} a_{2n} E_{2n+2}(\tau)x^{2n}$, 
so that 
$$
\bar\Delta(\zz|\tau) 
= -{1\over{2\pi\i}}\Delta_{\varphi(*|\tau)} 
+ {1\over{2\pi\i}} g(\zz|\tau).
$$

\begin{theorem} \label{thm:nabla}
There is a unique flat connection $\nabla_{{\cal P}_{n}}$ on ${\cal P}_{n}$, whose pull-back 
to $(\CC^{n} \times\HH) - \on{Diag}_{n}$ is the connection 
$$
\on{d} - \bar\Delta(\zz|\tau)\on{d}\tau - \sum_i \bar K_i(\zz|\tau)\on{d}z_i 
$$
on the trivial ${\bold G}_{n}$-bundle. 
\end{theorem}

{\em Proof.} We should check that the connection 
$\on{d} - \bar\Delta(\zz|\tau)\on{d}\tau - \sum_i \bar K_i(\zz|\tau)\on{d}z_i$
is equivariant and flat, which is expressed as follows (taking into account that we already checked
the equivariance and flatness of $\on{d} - \sum_{i} \bar K_{i}(\zz|\tau)\on{d}z_{i}$ for any $\tau$):

(equivariance) for $g = \bigl( \begin{smallmatrix}\alpha & \beta \\
  \gamma & \delta\end{smallmatrix}\bigl)\in \on{SL}_{2}(\ZZ)$
\begin{equation} \label{equiv:K}
{1\over{\gamma\tau+\delta}} \bar K_i({{\zz}\over{\gamma\tau+\delta}} | {{\alpha\tau+\beta}
\over{\gamma\tau+\delta}}) = \on{Ad}(c_g(\zz|\tau))(\bar K_i(\zz|\tau))
+ [(\partial/\partial z_i)c_g(\zz|\tau)]c_g(\zz|\tau)^{-1},
\end{equation} 
\begin{equation} \label{shift:Delta}
\bar\Delta(\zz + \delta_i|\tau) = 
\bar\Delta(\zz + u(\sum_i\delta_i)|\tau) = \bar\Delta(\zz|\tau), \quad 
\bar\Delta(\zz + \tau\delta_i|\tau) =e^{-2\pi\i\on{ad}x_i}
(\bar\Delta(\zz|\tau) - \bar K_i(\zz|\tau)),
\end{equation} 
\begin{align} \label{equiv:Delta}
{1\over{(\gamma\tau+\delta)^2}} \bar\Delta({{\zz}\over{\gamma\tau+\delta}}|
{{\alpha\tau+\beta}\over{\gamma\tau+\delta}}) & = 
\on{Ad}(c_g(\zz|\tau))(\bar\Delta(\zz|\tau))
+ {\gamma\over{\gamma z + \delta}} \sum_{i=1}^n z_i \on{Ad}(c_g(\zz|\tau))(\bar K_i(\zz|\tau)) \nonumber
\\ & + [({\partial\over{\partial\tau}} + {\gamma\over{\gamma\tau + \delta}} \sum_{i=1}^n z_i{{\partial}\over
{\partial z_i}})c_g(\zz|\tau)]c_g(\zz|\tau)^{-1}, 
\end{align}

(flatness) $[\partial/\partial\tau - \bar\Delta(\zz|\tau),
\partial/\partial z_i - \bar K_i(\zz|\tau)]=0$.

Let us now check the equivariance identity (\ref{equiv:K}) for 
$\bar K_{i}(\zz|\tau)$.  The cocycle identity $c_{gg'}(\zz|\tau) 
= c_{g}(g'*(\zz|\tau))c_{g'}(\zz|\tau)$ implies that it suffices 
to check it when $g = S$ and $g = T$. When $g=S$, this is the 
identity $\bar K_{i}(\zz|\tau+1) = \bar K_{i}(\zz|\tau)$, which follows 
from the identity $\theta(z|\tau+1) = \theta(z|\tau)$. When $g=T$, 
we have to check the identity 
\begin{equation} \label{mod:K}
{1\over \tau}\bar K_{i}({\zz\over\tau}|-{1\over\tau}) 
= \on{Ad}(\tau^{d} e^{{{2\pi\i}\over\tau}(\sum_{i}z_{i}\bar x_{i} + X)})
(\bar K_{i}(\zz|\tau)) + 2\pi\i\bar x_{i}. 
\end{equation}
We have 
\begin{align*}
& 2\pi\i\bar x_{i} - \on{Ad}(e^{2\pi\i(\sum_{i} z_{i}\bar x_{i} + X)})
(\bar y_{i}/\tau)
\\ & 
= - \on{Ad}(e^{2\pi\i(\sum_{i}z_{i}\bar x_{i})})(\bar y_{i}/\tau)
\quad \on{(as} \; \on{Ad}(e^{2\pi\i\tau X})(\bar y_{i}/\tau) 
= \bar y_{i}/\tau + 2\pi\i \bar x_{i})
\\ & 
= - {{\bar y_{i}}\over\tau} - {{e^{2\pi\i\on{ad}(\sum_{k}z_{k}\bar x_{k})}-1}
\over{\on{ad}(\sum_{k}z_{k}\bar x_{k})}}
([\sum_{j}z_{j}\bar x_{j},{{\bar y_{i}}\over\tau}]) 
= - {{\bar y_{i}}\over\tau} - {{e^{2\pi\i\on{ad}(\sum_{k}z_{k}\bar x_{k})}-1}
\over{\on{ad}(\sum_{k}z_{k}\bar x_{k})}}(\sum_{j|j\neq i} 
{{z_{ji}}\over\tau}\bar t_{ij})
\\ & 
= - {{\bar y_{i}}\over\tau} - \sum_{j|j\neq i} 
{{e^{2\pi\i\on{ad}(\sum_{k}z_{k}\bar x_{k})}-1}
\over{\on{ad}(\sum_{k}z_{k}\bar x_{k})}}({{z_{ji}}\over\tau}\bar t_{ij})
= - {{\bar y_{i}}\over\tau} - \sum_{j|j\neq i} 
{{e^{2\pi\i\on{ad}(z_{ij}\bar x_{i})}-1}
\over{\on{ad}(z_{ij}\bar x_{i})}}({{z_{ji}}\over\tau}\bar t_{ij})
\\
  & = - {{\bar y_{i}}\over\tau} +\sum_{j|j\neq i}
  {{e^{2\pi\i\on{ad}(z_{ij}\bar x_{i})}-1}
\over{\on{ad}(\bar x_{i})}}({{\bar t_{ij}}\over\tau}), 
\end{align*}
therefore 
\begin{equation} \label{15}
{1\over\tau}(\sum_{j} {{e^{2\pi\i z_{ij}\on{ad}\bar x_{i}}-1}
\over{\on{ad}\bar x_{i}}}(\bar t_{ij}) - \bar y_{i}) 
= - \on{Ad}(\tau^{d}e^{{{2\pi\i}\over\tau}
(\sum_{i}z_{i}\bar x_{i}+X)})(\bar y_{i}) + 2\pi\i\bar x_{i}. 
\end{equation}
We have $\theta(z/\tau|-1/\tau) = (1/\tau)e^{(\pi\i/\tau)z^{2}}\theta(z|\tau)$, therefore
\begin{equation} \label{16}
{1\over\tau}k({z\over\tau},x|-{1\over\tau}) = e^{2\pi\i zx}k(z,\tau x|\tau) 
+ {{e^{2\pi\i zx}-1}\over{x\tau}}. 
\end{equation}
Substituting $(z,x)=(z_{ij},\on{ad}\bar x_{i})$ ($j\neq i$), 
applying to $\bar t_{ij}$, summing over $j$ and adding up 
identity (\ref{15}), we get 
\begin{align*}
& {1\over\tau}(\sum_{j|j\neq i} k({{z_{ij}}\over\tau},\on{ad}\bar x_{i}|
-{1\over\tau})(\bar t_{ij}) - \bar y_{i})
\\ & 
= \sum_{j|j\neq i} e^{2\pi\i z_{ij}\on{ad}\bar x_{i}}k(z_{ij},
\tau\on{ad}\bar x_{i}|\tau)(\bar t_{ij}) 
- \on{Ad}(\tau^{d}e^{{{2\pi\i}\over\tau}(\sum_{i}z_{i}\bar x_{i}+X)})
(\bar y_{i}) + 2\pi\i\bar x_{i}. 
\end{align*}
Since $e^{2\pi\i z_{ij}\on{ad}\bar x_{i}}k(z_{ij},
\tau\on{ad}\bar x_{i}|\tau)(\bar t_{ij}) 
= \on{Ad}(\tau^{d}e^{(2\pi\i/\tau)(\sum_{i}z_{i}\bar x_{i}+X)})
\big(k(z_{ij},\on{ad}\bar x_{i})(\bar t_{ij})\big)$, 
this implies (\ref{mod:K}). This ends the proof of (\ref{equiv:K}).

Let us now check the shift identities (\ref{shift:Delta}) in 
$\bar\Delta(\zz|\tau)$. The first part is 
immediate; let us check the last identity. We have 
$k(z+\tau,x|\tau) = e^{-2\pi\i x}g(z,x|\tau) 
+ (e^{-2\pi\i x}-1)/x$, 
therefore $g(z+\tau,x|\tau) = e^{-2\pi\i x}g(z,x|\tau)
-2\pi\i e^{-2\pi\i x}k(z,x|\tau) 
+ {1\over x}({{1-e^{-2\pi\i x}}\over x} - 2\pi\i e^{-2\pi\i x})$. 
Substituting $(z,x) = (z_{ij},\on{ad}\bar x_{i})$
($j\neq i$), applying to $\bar t_{ij}$, summing up and adding up 
$\sum_{k,l|k,l\neq j}g(z_{kl},\on{ad}\bar x_{k}|\tau)(\bar t_{kl})$, we get 
\begin{align*}
& g(\zz+\tau\delta_{i}|\tau)
\\ & 
= e^{-2\pi\i\on{ad}\bar x_{i}}(g(\zz|\tau)) 
- 2\pi\i e^{-2\pi\i\on{ad}\bar x_{i}}
(\bar K_{i}(\zz|\tau)+\bar y_{i}) + \sum_{j|j\neq i} 
{1\over{\on{ad}\bar x_{i}}}
({{1-e^{-2\pi\i\on{ad}\bar x_{i}}}\over{\on{ad}\bar x_{i}}} 
- 2\pi\i e^{-2\pi\i\on{ad}\bar x_{i}})(\bar t_{ij})
\\ & 
= e^{-2\pi\i\on{ad}\bar x_{i}}(g(\zz|\tau)) 
- 2\pi\i e^{-2\pi\i\on{ad}\bar x_{i}}
(\bar K_{i}(\zz|\tau)+\bar y_{i}) 
- ({{1-e^{-2\pi\i\on{ad}\bar x_{i}}}\over{\on{ad}\bar x_{i}}} 
- 2\pi\i e^{-2\pi\i\on{ad}\bar x_{i}})(\bar y_{i})
\\ & 
= e^{-2\pi\i\on{ad}\bar x_{i}}(g(\zz|\tau)) 
- 2\pi\i e^{-2\pi\i\on{ad}\bar x_{i}}
(\bar K_{i}(\zz|\tau)) 
- {{1-e^{-2\pi\i\on{ad}\bar x_{i}}}\over{\on{ad}\bar x_{i}}}(\bar y_{i}); 
\end{align*}
on the other hand, we have $e^{-2\pi\i\on{ad}\bar x_{i}}(\Delta_{0}) 
= \Delta_{0} + {{1-e^{-2\pi\i\on{ad}\bar x_{i}}}
\over{\on{ad}\bar x_{i}}}(\bar y_{i})$ (as $[\Delta_{0},\bar x_{i}] 
= \bar y_{i}$), therefore 
$g(\zz + \delta_{i}|\tau) - \Delta_{0} = e^{-2\pi\i\on{ad}\bar x_{i}}
(g(\zz|\tau) - \Delta_{0} - 2\pi\i\bar K_{i}(\zz|\tau))$. Since the 
$\delta_{2n}$ commute with 
$\bar x_{i}$, we get $\bar\Delta(\zz + \tau\delta_{i}|\tau) 
= e^{-2\pi\i\on{ad}\bar x_{i}}
(\bar\Delta(\zz|\tau) - \bar K_{i}(\zz|\tau))$, as wanted.

Let us now check the equivariance identities (\ref{equiv:Delta}) for 
$\bar\Delta(\zz|\tau)$. As above,
the cocycle identities imply that it suffices to check (\ref{equiv:Delta}) 
for $g=S,T$. When $g=S$, 
this identity follows from $\sum_{i}\bar K_{i}(\zz|\tau)=0$. 
When $g=T$, it is written 
\begin{equation} \label{mod:Delta}
{1\over\tau^{2}}\bar\Delta({\zz\over\tau}|-{1\over\tau}) = 
\on{Ad}(c_{T}(\zz|\tau))
\Big(\bar\Delta(\zz|\tau) + {1\over\tau}\sum_{i} z_{i}\bar K_{i}(\zz|\tau)\Big) + {d\over\tau} 
-2\pi\i X. 
\end{equation}
The modularity identity (\ref{16}) for $k(z,x|\tau)$ implies that 
$$
{1\over\tau^{2}}g({z\over\tau},x|-{1\over\tau}) = e^{2\pi\i zx}g(z,\tau x|\tau) 
+ {{2\pi\i z}\over\tau} e^{2\pi\i zx}k(z,\tau x|\tau) + {{1-e^{2\pi\i zx}}
\over{\tau^{2}x^{2}}}
+ {{2\pi\i z}\over\tau^{2}}{{e^{2\pi\i zx}}\over x}. 
$$
This implies 
\begin{align*}
& {1\over\tau^{2}}\sum_{i<j} g({{z_{ij}}\over\tau},\on{ad}\bar x_{i}|
-{1\over\tau})(\bar t_{ij})
= \sum_{i<j} e^{2\pi\i z_{ij}\on{ad}\bar x_{i}} g(z_{ij},
\tau\on{ad}\bar x_{i}|\tau)(\bar t_{ij})
\\
  & + \sum_{i<j} {{2\pi\i}\over\tau} z_{ij} e^{2\pi\i z_{ij}\on{ad}\bar x_{i}} 
k(z_{ij},\tau \on{ad}\bar x_{i}|\tau)(\bar t_{ij})
+ \sum_{i<j} ({{1 - e^{2\pi\i z_{ij}\on{ad}\bar x_{i}}}
\over{\tau^{2}(\on{ad}\bar x_{i})^{2}}} 
+ {{2\pi\i z_{ij}}\over{\tau^{2}}}{{e^{2\pi\i z_{ij}\on{ad}\bar x_{i}}}
\over{\on{ad}\bar x_{i}}})
(\bar t_{ij}). 
\end{align*}
We compute as above
$$
\sum_{i<j}e^{2\pi\i z_{ij}\on{ad}\bar x_{i}}g(z_{ij},\tau\on{ad}\bar x_{i}|\tau)(\bar t_{ij}) = 
\on{Ad}(\tau^{d}e^{{{2\pi i}\over\tau}(\sum_{i}z_{i}\bar x_{i}+X)})(g(\zz|\tau)), 
$$
$$
\sum_{i<j} {{2\pi\i}\over\tau}z_{ij}e^{2\pi\i z_{ij}\on{ad}\bar x_{i}} 
k(z_{ij},\tau\on{ad}\bar x_{i}|\tau)(\bar t_{ij})
= \sum_{i}{{2\pi\i}\over\tau}z_{i}(\sum_{j|j\neq i} e^{2\pi\i z_{ij}\on{ad}\bar x_{i}}
k(z_{ij},\tau\on{ad}\bar x_{i}|\tau)(\bar t_{ij}))
$$
(using $k(z,x|\tau)+k(-z,-x|\tau)=0$) and 
$$
\sum_{i<j}e^{2\pi\i z_{ij}\on{ad}\bar x_{i}}k(z_{ij},\tau\on{ad}\bar x_{i}|\tau)(\bar t_{ij})
= \on{Ad}(\tau^{d}e^{{{2\pi i}\over\tau}(\sum_{i}z_{i}\bar x_{i}+X)})(\bar K_{i}(\zz|\tau)+\bar y_{i}). 
$$
Therefore 
\begin{align*}
{1\over\tau^{2}} g({\zz\over\tau}|-{1\over\tau})
& = \on{Ad}(c_{T}(\zz|\tau))\Big( g(\zz|\tau) + {{2\pi\i}\over\tau}\sum_{i}z_{i}\bar K_{i}(\zz|\tau)
+ {{2\pi\i}\over\tau}\sum_{i}z_{i}\bar y_{i}\Big)
\\ & 
+ \sum_{i<j} ({{1-e^{2\pi\i z_{ij}\on{ad}\bar x_{i}}}
\over{\tau^{2}(\on{ad}\bar x_{i})^{2}}} + {{2\pi\i z_{ij}}\over \tau^{2}}
{{e^{2\pi\i z_{ij}\on{ad}\bar x_{i}}}\over{\on{ad}\bar x_{i}}})(\bar t_{ij}), 
\end{align*}
which implies 
\begin{align*}
& 
{1\over\tau^{2}}\bar\Delta({\zz\over\tau}|-{1\over\tau}) 
= \on{Ad}(c_{T}(\zz|\tau))
\big( \bar\Delta(\zz|\tau) + {1\over\tau}\sum_{i}\bar K_{i}(\zz|\tau)\big)
\\ & 
+ \on{Ad}(c_{T}(\zz|\tau))({1\over\tau}\sum_{i}z_{i}\bar y_{i}) 
+ {1\over{2\pi\i}}
\sum_{i<j} ({{1-e^{2\pi\i z_{ij}\on{ad}\bar x_{i}}}
\over{\tau^{2}(\on{ad}\bar x_{i})^{2}}} + {{2\pi\i z_{ij}}\over \tau^{2}}
{{e^{2\pi\i z_{ij}\on{ad}\bar x_{i}}}\over{\on{ad}\bar x_{i}}})(\bar t_{ij})
\\
  & + {1\over{2\pi\i}}\big(\on{Ad}(c_{T}(\zz|\tau))(\Delta_{\varphi(*|\tau)}) 
- {1\over\tau^{2}}\Delta_{\varphi(*|-1/\tau)}\big) . 
\end{align*}
To prove (\ref{mod:Delta}), it then suffices to prove 
\begin{align} \label{remainingbis}
& \on{Ad}(c_{T}(\zz|\tau))({1\over\tau}\sum_{i}z_{i}\bar y_{i}) 
+ {1\over{2\pi\i}}
\sum_{i<j} ({{1-e^{2\pi\i z_{ij}\on{ad}\bar x_{i}}}
\over{\tau^{2}(\on{ad}\bar x_{i})^{2}}} + {{2\pi\i z_{ij}}\over \tau^{2}}
{{e^{2\pi\i z_{ij}\on{ad}\bar x_{i}}}\over{\on{ad}\bar x_{i}}})(\bar t_{ij})
\nonumber \\
  & + {1\over{2\pi\i}}\big(\on{Ad}(c_{T}(\zz|\tau))(\Delta_{\varphi(*|\tau)}) 
- {1\over\tau^{2}}\Delta_{\varphi(*|-1/\tau)}\big)  = {d\over\tau} - 2\pi\i X. 
\end{align}
We compute 
$$
\on{Ad}(c_{T}(\zz|\tau))({1\over\tau}\sum_{i}z_{i}\bar y_{i}) = 
{1\over\tau^{2}} \sum_{i}z_{i}\bar y_{i} 
+ {{2\pi\i}\over\tau}\sum_{i}z_{i}\bar x_{i}
+\sum_{i<j} (-{1\over\tau^{2}}) z_{ij} 
{{e^{2\pi\i z_{ij}\on{ad}\bar x_{i}}-1}\over{\on{ad}\bar x_{i}}}(\bar t_{ij}). 
$$
We also have $\on{Ad}(c_{T}(\zz|\tau))(E_{2n+2}(\tau)\delta_{2n}) 
= {1\over\tau^{2}}E_{2n+2}(-{1\over\tau})\delta_{2n}$
since $[\delta_{2n},\bar x_{i}]=[\delta_{2n},X]=0$ and 
$[d,\delta_{2n}]=2n\delta_{2n}$, and since $E_{2n+2}(-1/\tau)=
\tau^{2n+2}E_{2n+2}(\tau)$. This implies 
$$
\on{Ad}(c_{T}(\zz|\tau))(\delta_{\varphi(*|\tau)}) = \delta_{\varphi(*|-1/\tau)}. 
$$
We now compute $\on{Ad}(c_{T}(\zz|\tau))(\Delta_{0}) - (1/\tau^{2})\Delta_{0}$. 
We have 
$\on{Ad}(c_{T}(\zz|\tau))(\Delta_{0}) = 
\on{Ad}(e^{2\pi\i \sum_{i} z_{i}\bar x_{i}}) \circ $ $
\on{Ad}(\tau^{d}e^{(2\pi\i/\tau)X})(\Delta_{0})$, and 
$\on{Ad}(\tau^{d}e^{(2\pi\i/\tau)X})(\Delta_{0})
= (1/\tau^{2})\Delta_{0} + (2\pi\i/\tau)d - (2\pi\i)^{2}X$. Now 
$\on{Ad}(e^{2\pi\i\sum_{i}z_{i}\bar x_{i}})(X)=X$, 
$\on{Ad}(e^{2\pi\i\sum_{i}z_{i}\bar x_{i}})(d) 
= d - 2\pi\i \sum_{i}z_{i}\bar x_{i}$. We now compute 
\begin{align*}
& \on{Ad}(e^{2\pi\i\sum_{i}z_{i}\bar x_{i}})(\Delta_{0}) = \Delta_{0} 
+ {{e^{2\pi\i\sum_{i}z_{i}\on{ad}\bar x_{i}}-1}\over
{2\pi\i\on{ad}(\sum_{i}z_{i}\bar x_{i})}}
([2\pi\i\sum_{i}z_{i}\bar x_{i},\Delta_{0}])
\\
& = \Delta_{0} 
- {{e^{2\pi\i\sum_{i}z_{i}\on{ad}\bar x_{i}}-1}\over
{\on{ad}(\sum_{i}z_{i}\bar x_{i})}}
( \sum_{i} z_{i}\bar y_{i})
= \Delta_{0} 
- \sum_{i} {{e^{2\pi\i\sum_{j|j\neq i}z_{ji}\on{ad}\bar x_{j}}-1}\over
{\on{ad}(\sum_{j|j\neq i}z_{ji}\bar x_{j})}}
( z_{i}\bar y_{i})
  \\ & 
= \Delta_{0} - \sum_{i} \Big(2\pi\i z_{i}\bar y_{i} + {1\over
{\on{ad}(\sum_{j|j\neq i} z_{ji}\bar x_{j})}} 
({{e^{2\pi\i\sum_{j|j\neq i} z_{ji}\on{ad}\bar x_{j}}-1}\over
{\on{ad}(\sum_{j|j\neq i}z_{ji}\bar x_{j})}} 
- 2\pi\i)([\sum_{j|j\neq i}z_{ji}\bar x_{j},z_{i}\bar y_{i}]) \Big)
\\ & 
= \Delta_{0} - \sum_{i} 2\pi\i z_{i}\bar y_{i} - \sum_{i\neq j} 
\Big( {1\over{\on{ad}(\bar x_{j})}} 
({{e^{2\pi\i z_{ji}\on{ad}\bar x_{j}}-1}\over{\on{ad}(z_{ji}\bar x_{j})}} 
- 2\pi\i)(z_{i}\bar t_{ij}) \Big); 
\end{align*}
the last sum decomposes as 
\begin{align*}
& \sum_{i<j} {1\over{\on{ad}(\bar x_{j})}} 
({{e^{2\pi\i z_{ji}\on{ad}\bar x_{j}}-1}\over{\on{ad}(z_{ji}\bar x_{j})}} 
- 2\pi\i)(z_{i}\bar t_{ij})  + \sum_{i>j}  {1\over{\on{ad}(\bar x_{j})}} 
({{e^{2\pi\i z_{ji}\on{ad}\bar x_{j}}-1}\over{\on{ad}(z_{ji}\bar x_{j})}} 
- 2\pi\i)(z_{i}\bar t_{ij}) 
\\ & = 
\sum_{i<j}  {1\over{\on{ad}(\bar x_{j})}} 
({{e^{2\pi\i z_{ji}\on{ad}\bar x_{j}}-1}\over{\on{ad}(z_{ji}\bar x_{j})}} 
- 2\pi\i)(z_{i}\bar t_{ij})  +   {1\over{\on{ad}(\bar x_{i})}} 
({{e^{2\pi\i z_{ij}\on{ad}\bar x_{i}}-1}\over{\on{ad}(z_{ij}\bar x_{i})}} 
- 2\pi\i)(z_{j}\bar t_{ij})
\\ & 
= \sum_{i<j}   {1\over{\on{ad}(\bar x_{i})}} 
({{e^{2\pi\i z_{ij}\on{ad}\bar x_{i}}-1}\over{\on{ad}(z_{ij}\bar x_{i})}} 
- 2\pi\i)(z_{ji}\bar t_{ij}), 
\end{align*}
so 
$$
\on{Ad}(e^{2\pi\i\sum_{i}z_{i}\bar x_{i}})(\Delta_{0}) = \Delta_{0} 
-2\pi \i\sum_{i}z_{i}\bar y_{i} 
- \sum_{i<j}   {1\over{\on{ad}(\bar x_{i})}} 
({{e^{2\pi\i z_{ij}\on{ad}\bar x_{i}}-1}\over{\on{ad}(z_{ij}\bar x_{i})}} 
- 2\pi\i)(z_{ji}\bar t_{ij}), 
$$
and finally 
\begin{align*}
& \on{Ad}(c_{T}(\zz|\tau))(\Delta_{\varphi(*|\tau)}) 
- {1\over\tau^{2}}\Delta_{\varphi(*|-1/\tau)} 
\\ & = 
-{{2\pi \i}\over{\tau^{2}}}\sum_{i}z_{i}\bar y_{i} 
- {1\over\tau^{2}}\sum_{i<j}   {1\over{\on{ad}(\bar x_{i})}} 
({{e^{2\pi\i z_{ij}\on{ad}\bar x_{i}}-1}\over{\on{ad}(z_{ij}\bar x_{i})}} 
- 2\pi\i)(z_{ji}\bar t_{ij}) + {{2\pi\i}\over\tau}(d 
- 2\pi\i\sum_{i}z_{i}\bar x_{i}) - (2\pi\i)^{2}X, 
\end{align*}
which implies (\ref{remainingbis}). This proves (\ref{mod:Delta}) 
and therefore (\ref{equiv:Delta}).

We then prove that flatness identity $[\partial/\partial\tau 
- \bar\Delta(\zz|\tau),
\partial/\partial z_{i} - \bar K_{i}(\zz|\tau)]=0$.  For this, we will prove 
that $(\partial/\partial\tau)\bar K_{i}(\zz|\tau) 
= (\partial/\partial\tau)\bar\Delta(\zz|\tau)$, and that 
$[\bar\Delta(\zz|\tau),\bar K_{i}(\zz|\tau)]=0$.

Let us first prove 
\begin{equation} \label{crossed}
(\partial/\partial\tau)\bar K_{i}(\zz|\tau) 
= (\partial/\partial z_{i})\bar\Delta(\zz|\tau). 
\end{equation}
We have $(\partial/\partial\tau)\bar K_{i}(\zz|\tau) 
=\sum_{j|j\neq i} (\partial_{\tau}k)(z_{ij},\on{ad}\bar x_{i}|\tau)(\bar t_{ij})$
and $ (\partial/\partial z_{i})\bar\Delta(\zz|\tau) = (2\pi\i)^{-1}\sum_{j|j\neq i}$ $
(\partial_{z}g)(z_{ij},\on{ad}\bar x_{i})(\bar t_{ij})$ (where 
$\partial_{\tau}:= \partial/\partial\tau$, $\partial_{z} = \partial/\partial z$)
so it suffices to prove the identity 
$(\partial_{\tau}k)(z,x|\tau) = (2\pi\i)^{-1}(\partial_{z}g)(z,x|\tau)$, i.e., 
$(\partial_{\tau}k)(z,x|\tau) = (2\pi\i)^{-1}(\partial_{z}\partial_{x}k)(z,x|\tau)$. 
In this identity, $k(z,x|\tau)$ may be replaced by $\tilde k(z,x|\tau):= k(z,x|\tau) + 1/x
= \theta(z+x|\tau)/(\theta(z|\tau)\theta(x|\tau))$. Dividing by $\tilde k(z,x|\tau)$, 
the wanted identity is rewritten as 
$$
2\pi\i \big( {{\partial_{\tau}\theta}\over\theta}(z+x|\tau) 
- {{\partial_{\tau}\theta}\over\theta}(z|\tau) 
- {{\partial_{\tau}\theta}\over\theta}(x|\tau)\big) = 
({{\theta'}\over\theta})'(z+x|\tau) + \big({\theta'\over\theta}(z+x|\tau) 
- {\theta'\over\theta}(z|\tau)\big)
\big({\theta'\over\theta}(z+x|\tau) - {\theta'\over\theta}(x|\tau)\big) 
$$
(recall that $f'(z|\tau) = \partial_{z}f(z|\tau)$), 
or taking into account the heat equation 
$4\pi\i (\partial_{\tau}\theta/\theta)(z|\tau) = 
(\theta''/\theta)(z|\tau) - 12\pi\i (\partial_{\tau}\eta/\eta)(\tau)$, 
as follows
\begin{align} \label{id:theta}
& 2 \Big( {\theta'\over\theta}(z|\tau){\theta'\over\theta}(x|\tau)
- {\theta'\over\theta}(x|\tau){\theta'\over\theta}(z+x|\tau)
- {\theta'\over\theta}(z|\tau){\theta'\over\theta}(z+x|\tau) \Big) 
\\
& + {\theta''\over\theta}(z|\tau)
  \nonumber 
+ {\theta''\over\theta}(x|\tau) + {\theta''\over\theta}(z+x|\tau) 
- 12 \pi \i {{\partial_\tau \eta}\over \eta}(\tau)=0
\end{align}
Let us prove (\ref{id:theta}). Denote its l.h.s. by $F(z,x|\tau)$. Since $\theta(z|\tau)$
is odd w.r.t. $z$, $F(z,x|\tau)$ is invariant under the permutation of $z,x,-z-x$. 
The identities $(\theta'/\theta)(z+\tau|\tau) = (\theta'/\theta)(z|\tau) - 2\pi\i$ and
$(\theta''/\theta)(z+\tau|\tau) = (\theta''/\theta)(z|\tau) 
- 4\pi\i(\theta'/\theta)(z|\tau) +(2\pi\i)^{2} $ imply that 
$F(z,x|\tau)$ is elliptic in $z,x$ (w.r.t. the lattice $\Lambda_{\tau}$). The possible 
poles of $F(z,x|\tau)$  as a function of $z$ are simple at $z = 0$ and $z = -x$ 
(mod $\Lambda_{\tau}$), but one checks that $F(z,x|\tau)$ is regular 
at these points, so it is constant in $z$. By the $\SG_{3}$-symmetry, it is also constant in 
$x$, hence it is a function of $\tau$ only: $F(z,x|\tau) = F(\tau)$.

To compute this function,  we compute 
$F(z,0|\tau) = [-2(\theta'/\theta)' -2(\theta'/\theta)^{2}
+2\theta''/\theta](z|\tau)
+ (\theta''/\theta)(0|\tau) - 12\pi\i (\partial_{\tau}\eta/\theta)(\tau)$, hence 
$F(\tau) = (\theta''/\theta)(0|\tau) - 12\pi\i(\partial_{\tau}\eta/\eta)(\tau)$. 
The above heat equation then implies that 
$F(\tau) = 4\pi\i (\partial_{\tau}\theta/\theta)(0|\tau)$. 
Now $\theta'(0|\tau)=1$ implies that
$\theta(z|\tau)$ has the expansion $\theta(z|\tau) 
= z + \sum_{n\geq 2}a_{n}(\tau)z^{n}$ as $z\to 0$, 
which implies $(\partial_{\tau}\theta/\theta)(0|\tau)=0$. So $F(\tau)=0$, 
which implies (\ref{id:theta}) and therefore (\ref{crossed}).

We now prove
\begin{equation} \label{comm}
[\bar\Delta(\zz|\tau),\bar K_{i}(\zz|\tau)]=0. 
\end{equation}
Since $\tau$ is constant in what follows, we will write $k(z,x)$, $g(z,x)$, $\varphi$ instead of 
$k(z,x|\tau)$, $g(z,x|\tau)$, $\varphi(*|\tau)$. 
For $i\neq j$, let us set $g_{ij}:= g(z_{ij},\on{ad}\bar x_{i})(\bar t_{ij})$. 
Since $g(z,x|\tau) = g(-z,-x|\tau)$, we have $g_{ij} = g_{ji}$. 
Recall that $\bar K_{ij} = k(z_{ij},\on{ad}\bar x_{i})(t_{ij})$.

We have 
\begin{align} \label{start}
& 2\pi\i[\bar\Delta(\zz|\tau),\bar K_{i}(\zz|\tau)] = 
[-\Delta_{\varphi} + \sum_{i,j|i<j}g_{ij}, -\bar y_{i}+\sum_{j|j\neq i} \bar K_{ij}]
\\ & \nonumber
= [\Delta_{\varphi},\bar y_{i}] + 
\sum_{j|j\neq i} \Big( - [\Delta_{\varphi},\bar K_{ij}] + [\bar y_{i},g_{ij}]
+ [g_{ij},\bar K_{ij}]\Big)
\\ & \nonumber 
+ \sum_{j,k|j\neq i,k\neq i, j<k} \big( [\bar y_{i},g_{jk}] + [g_{ik}+g_{jk},\bar K_{ij}]
  + [g_{ij}+g_{jk},\bar K_{ik}]\big). 
\end{align}

One computes 
\begin{equation} \label{partial:1}
[\Delta_{\varphi},\bar y_{i}] = \sum_{\alpha} 
[f_{\alpha}(\on{ad}\bar x_{i})(\bar t_{ij}),
g_{\alpha}(-\on{ad}\bar x_{i})(\bar t_{ij})], \quad \on{where} \; 
\sum_{\alpha}f_{\alpha}(u)g_{\alpha}(v) 
= {1\over 2}{{\varphi(u) - \varphi(v)}\over{u-v}}. 
\end{equation}

If $f(x)\in\CC[[x]]$, then 
\begin{align*}
& [\Delta_{0},f(\on{ad}\bar x_{i})(\bar t_{ij})] 
- [\bar y_{i},f'(\on{ad}\bar x_{i})(\bar t_{ij})]
= \sum_{\alpha} [h_{\alpha}(\on{ad}\bar x_{i})(\bar t_{ij}),
k_{\alpha}(\on{ad}\bar x_{i})(\bar t_{ij})]
\\
  & + \sum_{k|k\neq i,j} {{f(\on{ad}\bar x_{i})-f(-\on{ad}\bar x_{j})
  - f'(-\on{ad}\bar x_{j})(\on{ad}\bar x_{i}+\on{ad}\bar x_{j})}
\over{(\on{ad}\bar x_{i} + \on{ad}\bar x_{j})^{2}}}([\bar t_{ij},\bar t_{ik}]), 
\end{align*}
where 
$$
\sum_{\alpha}h_{\alpha}(u)k_{\alpha}(v) = {1\over 2} \Big(
{1\over v^{2}}\big(f(u+v) - f(u) - vf'(u)\big)
- {1\over u^{2}}\big(f(u+v) - f(v) - uf'(v)\big) \Big). 
$$
Since $g(z,x) = k_{x}(z,x)$, we get 
\begin{align} \label{partial:2}
& -[\Delta_{0},\bar K_{ij}] + [\bar y_{i},g_{ij}] = -\sum_{\alpha} 
[f_{\alpha}^{ij}(\on{ad}\bar x_{i})(\bar t_{ij}), g_{\alpha}^{ij}
(\on{ad}\bar x_{i})(\bar t_{ij})]
\\ & \nonumber 
+ \sum_{k|k\neq i,j} {{k(z_{ij},\on{ad}\bar x_{i}) - k(z_{ij},-\on{ad}\bar x_{j}) 
- (\on{ad}\bar x_{i}+ \on{ad}\bar x_{j})k_{x}(z_{ij},-\on{ad}\bar x_{j})}
\over{(\on{ad}\bar x_{i}+\on{ad}\bar x_{j})^{2}}}([\bar t_{ij},\bar t_{jk}]), 
\end{align}
where 
$$
\sum_{\alpha}f^{ij}_{\alpha}(u)g^{ij}_{\alpha}(v) = {1\over 2} \Big(
{1\over v^{2}}\big(k(z_{ij},u+v) - k(z_{ij},u) - vk_{x}(z_{ij},u)\big)
- {1\over u^{2}}\big(k(z_{ij},u+v) - k(z_{ij},v) - u k_{x}(z_{ij},v)\big) \Big). 
$$

For $f(x)\in\CC[[x]]$, we have 
$$
[\delta_{\varphi},f(\on{ad}\bar x_{i})(\bar t_{ij})] = \sum_{\alpha}
[l_{\alpha}(\on{ad}\bar x_{i})(\bar t_{ij}),m_{\alpha}(\on{ad}\bar x_{i})
(\bar t_{ij})], \quad \on{where}
\sum_{\alpha} l_{\alpha}(u)m_{\alpha}(v) = f(u+v)\varphi(v), 
$$
therefore 
\begin{equation} \label{partial:3}
-[\delta_{\varphi},\bar K_{ij}] = 
-\sum_{\alpha}[l_{\alpha}^{ij}(\on{ad}\bar x_{i})(\bar t_{ij}),
m_{\alpha}^{ij}(\on{ad}\bar x_{i})(\bar t_{ij})], \quad \on{where}\; 
\sum_{\alpha} l_{\alpha}^{ij}(u)m_{\alpha}^{ij}(v) = k(z_{ij},u+v)\varphi(v). 
\end{equation}

For $j,k\neq i$ and $j<k$, we have 
$$
[\bar y_{i},g_{jk}] + [g_{ik}+g_{jk},\bar K_{ij}] + [g_{ij}+g_{jk},\bar K_{ik}]
= [\bar y_{i},g_{jk}] - [g_{ki},\bar K_{ji}] - [g_{ji},\bar K_{ki}]
+ [g_{jk},\bar K_{ij}] +[g_{jk},\bar K_{ik}], 
$$
and since for any $f(x)\in\CC[[x]]$, 
$[\bar y_{i},f(\on{ad}\bar x_{i})(\bar t_{jk})] = - {{f(\on{ad}\bar x_{j}) 
- f(-\on{ad}\bar x_{k})}
\over{\on{ad}\bar x_{j} + \on{ad}\bar x_{k}}}([\bar t_{ij},\bar t_{jk}])$, we get 
\begin{align} \label{partial:4}
& [\bar y_{i},g_{jk}] + [g_{ik}+g_{jk},\bar K_{ij}] + [g_{ij}+g_{jk},\bar K_{ik}]
\\ & \nonumber
= \Big(-{{g(z_{jk},\on{ad}\bar x_{j}) - g(z_{jk},-\on{ad}\bar x_{k})}\over
{\on{ad}\bar x_{j} + \on{ad}\bar x_{k}}} - g(z_{ki},\on{ad}\bar x_{k})
k(z_{ji},\on{ad}\bar x_{j})
+ g(z_{ji},\on{ad}\bar x_{j})k(z_{ki},\on{ad}\bar x_{k})
\\ & \nonumber
- g(z_{kj},\on{ad}\bar x_{k})k(z_{ij},\on{ad}\bar x_{i})
  + g(z_{jk},\on{ad}\bar x_{j})k(z_{ik},\on{ad}\bar x_{i})\Big)([\bar t_{ij},\bar t_{jk}]). 
\end{align}
Summing up (\ref{partial:1}), (\ref{partial:2}), (\ref{partial:3}) and (\ref{partial:4}), 
(\ref{start}) gives
\begin{align*}
& 2\pi\i[\bar\Delta(\zz|\tau),\bar K_{i}(\zz|\tau)] 
\\ & = \sum_{j|j\neq i} \sum_{\alpha}
[F^{ij}_{\alpha}(\on{ad}\bar x_{i})(\bar t_{ij}),G^{ij}_{\alpha}
(\on{ad}\bar x_{i})(\bar t_{ij})]
+ \sum_{j,k|j\neq i,k\neq i} H(z_{ij},z_{ik},-\on{ad}\bar x_{j},-\on{ad}\bar x_{k})
([t_{ij},t_{jk}]), 
\end{align*}
where 
$\sum_{\alpha} F_{\alpha}^{ij}(u)G_{\alpha}^{ij}(v) = L(z_{ij},u,v)$, 
\begin{align*}
& L(z,u,v) = {1\over 2}{{\varphi(u)-\varphi(v)}\over{u+v}} 
+ {1\over 2}k(z,u+v) (\varphi(u) - \varphi(v)) + {1\over 2}(g(z,u)k(z,v) - k(z,u)g(z,v))
\\
  & - {1\over 2} \Big({1\over v^{2}} \big(k(z,u+v)-k(z,u) - vk_{x}(z,u)\big)
  - {1\over u^{2}} \big(k(z,u+v)-k(z,v) - uk_{x}(z,v)\big)  \Big)\end{align*}
and 
\begin{align*}
& H(z,z',u,v) = {1\over v^{2}} \big( k(z,u+v)-k(z,u)-vk_{x}(z,u) \big)
  - {1\over u^{2}} \big( k(z',u+v)-k(z',v)-uk_{x}(z',v) \big) 
\\
  & + {1\over{u+v}} \big( g(z'-z,-u)-g(z'-z,v) \big)  - g(-z',-v)k(-z,-u) + g(-z,-u)k(-z',-v) 
\\
  & - g(z-z',-v)k(z,u+v) + g(z'-z,-u)k(z',u+v). 
\end{align*}
Explicit computation shows that $H(z,z',u,v)=0$, which implies that $L(z,u,v)=0$ since 
$L(z,u,v) = -{1\over 2} H(z,z,u,v)$. This proves (\ref{comm}). 
\hfill \qed \medskip

\begin{remark}
Define $\Delta(\zz|\tau)$ by the same formula as $\bar\Delta(\zz|\tau)$, 
replacing $\bar x_{i},\bar y_{i}$ by $x_{i},y_{i}$. Then 
$\on{d} - \Delta(\zz|\tau) \on{d}\tau- \sum_{i} K_{i}(\zz|\tau)\on{d}z_{i}$ 
is flat. This can be interpreted as follows.

Let $N_{+}\subset \on{SL}_{2}(\CC)$ be the connected subgroup with Lie algebra
$\CC\Delta_{0}$. Set $\tilde{\bold N}_{n} 
:= \on{exp}((\t_{1,n}\rtimes\d_{+})^{\wedge})\rtimes N_{+}$, 
${\bold N}_{n} 
:= \on{exp}((\bar\t_{1,n}\rtimes\d_{+})^{\wedge})\rtimes N_{+}$ and 
$\tilde{\bold G}_{n} := \on{exp}((\t_{1,n}\rtimes\d_{+})^{\wedge})\rtimes 
\on{SL}_{2}(\CC)$. Then we have a diagram of groups
$$
\begin{matrix}
\tilde{\bold N}_{n} & \to & {\bold N}_{n}\\
\downarrow &     & \downarrow\\
\tilde{\bold G}_{n} & \to & {\bold G}_{n} 
\end{matrix}$$
The trivial ${\bold N}_{n}$-bundle on $(\HH\times\CC^{n}) - \on{Diag}_{n}$
with flat connection $\on{d} - \bar\Delta(\zz|\tau)\on{d}\tau 
- \sum_{i} \bar K_{i}(\zz|\tau)
\on{d}z_{i}$ admits a reduction to $\tilde{\bold N}_{n}$, 
where the bundle is again trivial and the connection is $\on{d} - \Delta(\zz|\tau)
\on{d}\tau - \sum_{i} K_{i}(\zz|\tau) \on{d}z_{i}$.

$((\ZZ^{2})^{2}\times \CC)\rtimes \on{SL}_{2}(\ZZ)$ contains the
  subgroups $(\ZZ^n)^2$, $(\ZZ^n)^2\times 
\CC$, $(\ZZ^n)^2\rtimes \on{SL}_2(\ZZ)$. We denote the corresponding 
quotients of $(\CC^n\times\HH) - \on{Diag}_{n}$ by $C(n)$, $\bar C(n)$, 
$\tilde{\cal M}_{1,n}$. These fit in the diagram 
$$
\begin{matrix}
\tilde C(n) & \to & C(n)\\
\downarrow & & \downarrow \\
\tilde{\cal M}_{1,n} & \to & {\cal M}_{1,n}
\end{matrix}
$$
The pair $({\cal P}_n,\nabla_{{\cal P}_n})$ can be pulled back to ${\bold
G}_n$-bundles over these covers of ${\cal M}_{1,n}$. These  pull-backs 
admit $G$-structures, where $G$ is the corresponding group in the 
above diagram of groups.

We have natural projections $C(n)\to \HH$, 
$\bar C(n) \to \HH$. The fibers of $\tau\in\HH$ are respectively 
$C(E_\tau,n)$ and $\bar C(E_\tau,n)$. The pair $({\cal P}_n,\nabla_n)$ 
can be pulled back to $C(E_\tau,n)$ and $\bar C(E_\tau,n)$; these 
pull-backs admit $G$-structures, where $G = \on{exp}(\t_{1,n})$
and $\on{exp}(\bar\t_{1,n})$, which coincide with
$(P_{n,\tau},\nabla_{n,\tau})$ and $(\bar P_{n,\tau},\bar\nabla_{n,\tau})$. 
\end{remark}

\subsection{Bundle with flat connection over ${\cal M}_{1,[n]}$}

The semidirect product $((\ZZ^{n})^{2}\times \CC)\rtimes 
(\on{SL}_{2}(\ZZ)\times S_{n})$ 
acts on $(\CC^{n}\times\HH) - \on{Diag}_{n}$ as follows: the action of 
$((\ZZ^{n})^{2}\times \CC)\rtimes \on{SL}_{2}(\CC)$ is as above and the 
action of $S_{n}$ is $\sigma * (z_{1},...,z_{n},\tau) := (z_{\sigma^{-1}(1)},...,
z_{\sigma^{-1}(n)},\tau)$. The quotient then identifies with ${\cal M}_{1,[n]}$.

We will define a principal ${\bold G}_{n}\rtimes S_{n}$-bundle with a flat connection 
$({\cal P}_{[n]},\nabla_{{\cal P}_{[n]}})$ over ${\cal M}_{1,[n]}$.

\begin{proposition}
There exists a unique principal ${\bold G}_{n}\rtimes S_{n}$-bundle 
${\cal P}_{[n]}$ over ${\cal M}_{1,[n]}$, such that a section of 
$U\subset {\cal M}_{1,[n]}$ is a function $f : \tilde\pi^{-1}(U) \to 
{\bold G}_{n}\rtimes S_{n}$, satisfying the conditions 
of Proposition \ref{prop:bundle} as well as 
$f(\sigma\zz|\tau) = \sigma f(\zz|\tau)$ for $\sigma\in S_n$
(here $\tilde\pi : (\CC^{n}\times \HH) - \on{Diag}_{n} \to {\cal M}_{1,[n]}$ 
is the canonical projection). 
\end{proposition}

{\em Proof.} One checks that $\sigma c_{\tilde g}(\zz|\tau) \sigma^{-1} = c_{\sigma \tilde g
\sigma^{-1}}(\sigma^{-1}\zz)$, where $\tilde g\in ((\ZZ^{n})^{2}\times \CC)
\rtimes \on{SL}_{2}(\ZZ)$, $\sigma\in S_{n}$. It follows that 
there is a unique cocycle $c_{(\tilde g,\sigma)} : 
\CC^{n}\times\HH\to {\bold G}_{n}\rtimes S_{n}$ 
such that $c_{(\tilde g,1)} = c_{\tilde g}$ and 
$c_{(1,\sigma)}(\zz|\tau)=\sigma$. \hfill \qed\medskip

\begin{theorem}
There is a unique flat connection $\nabla_{{\cal P}_{[n]}}$ on ${\cal P}_{[n]}$, 
whose pull-back to $(\CC^{n}\times \HH) - \on{Diag}_{n}$ is the connection 
$\on d - \bar\Delta(\zz|\tau) \on{d}\tau - \sum_{i}\bar K_{i}(\zz|\tau)\on{d}z_{i}$
on the trivial ${\bold G}_{n}\rtimes S_{n}$-bundle. 
\end{theorem}

{\em Proof.} Taking into account Theorem \ref{thm:nabla}, 
it remains to show that this connection is 
$S_{n}$-equivariant. We have already mentioned that 
$\sum_{i}\bar K_{i}(\zz|\tau)\on{d}z_{i}$ is equivariant; $\bar\Delta(\zz|\tau)$
is also checked to be equivariant. \hfill \qed \medskip

\section{The monodromy morphisms $\Gamma_{1,[n]}\to {\bold G}_{n}\rtimes S_n$}
\label{sect:5}

Let $\Gamma_{1,[n]}$ be the mapping class group of genus 1 surfaces 
with $n$ unordered marked points. It can be viewed as the fundamental group 
$\pi_1({\cal M}_{1,[n]},*)$, where $*$ is a base point at infinity 
which will be specified later. The flat connection on ${\cal M}_{1,[n]}$
introduced above gives rise to morphisms $\gamma_n : \Gamma_{1,[n]} \to 
{\bold G}_n\rtimes S_n$, which we now study. This study in divided in two 
parts: in the first, analytic part, we show that $\gamma_n$ can be obtained 
from $\gamma_1$ and $\gamma_2$, and show that the restriction of 
$\gamma_n$ to $\overline{\on{B}}_{1,n}$ can be expressed 
in terms of the KZ associator only. In the second part, we show that 
morphisms $\overline{\on{B}}_{1,n} \to \on{exp}(\hat{\bar{\t}}_{1,n})
\rtimes S_n$ can be constructed algebraically using an arbitrary associator. 
Finally, we introduce the notion of an elliptic structure over a 
quasi-bialgebra.

\subsection{The solution $F^{(n)}(\zz|\tau)$}

The elliptic KZB system is now 
$$
(\partial/\partial z_i)F(\zz|\tau) = \bar K_i(\zz|\tau)F(\zz|\tau), \quad 
(\partial/\partial\tau)F(\zz|\tau) = \bar\Delta(\zz|\tau)F(\zz|\tau), 
$$
where $F(\zz|\tau)$ is a function $(\CC^n \times \HH)-\on{Diag}_n 
\supset U \to {\bold G}_n\rtimes S_{n}$ invariant under translation by 
$\CC(\sum_i\delta_i)$. 
Let $D_n := \{(\zz,\tau) \in\CC^n\times \HH | z_i = a_i + b_i\tau, a_i,b_i\in
\RR, a_1<a_2<...<a_n<a_1+1, 
b_1<b_2<...<b_n<b_1+1\}$. Then $D_n\subset (\CC^n \times \HH)-\on{Diag}_n$
is simply connected and invariant under $\CC(\sum_i\delta_i)$. 
A solution of the elliptic KZB system on this domain is then unique, 
up to right multiplication by a constant. We now determine a 
particular solution $F^{(n)}(\zz|\tau)$.

Let us study the elliptic KZB system in the region 
$z_{ij}\ll 1$, $\tau\to \i\infty$.
Then $\bar K_i(\zz|\tau) = \sum_{j|j\neq i} \bar t_{ij}/(z_i - z_j)
+ O(1)$.

We now compute the expansion of $\bar\Delta(\zz|\tau)$. 
The heat equation for $\vartheta$ implies the 
expansion $\vartheta(x|\tau) = \eta(\tau)^3 \big(x + 2\pi\i \partial_\tau
\on{log}\eta(\tau)x^3 + O(x^5)\big)$, so 
$\theta(x|\tau) = x + 2\pi\i \partial_\tau
\on{log}\eta(\tau)x^3 + O(x^5)$, hence 
$$
g(0,x|\tau) = ({\theta'\over\theta})'(x|\tau) + {1\over{x^2}}
= 4\pi\i \partial_\tau\on{log}\eta(\tau) + O(x) = - (\pi^2/3)E_2(\tau)+O(x)
$$ 
since $E_2(\tau) = {{24}\over{2\pi\i}}\partial_\tau
\on{log}\eta(\tau)$. 
We have 
$g(0,x|\tau) = g(0,0|\tau) - \varphi(x|\tau)$, so 
$$
g(0,x|\tau) = -\sum_{k\geq 0} a_{2k}x^{2k}E_{2k+2}(\tau), 
$$ 
where $a_0 = \pi^2/3$. 
Then 
$$
\bar\Delta(\zz|\tau) = -{1\over{2\pi\i}}
\Big( \Delta_0 + \sum_{k\geq 0} a_{2k}E_{2k+2}(\tau) \big( \delta_{2k} + 
\sum_{i,j|i<j} (\on{ad}\bar x_i)^{2k}(\bar t_{ij})\big) \Big) + o(1)
$$
for $z_{ij}\ll 1$ and any $\tau\in\HH$. Since we have an expansion 
$E_{2k}(\tau) = 1 + \sum_{l>0} a_{kl}e^{2\pi\i l\tau}$ 
as $\tau\to\i\infty$, and using Proposition \ref{prop:app:3}
with $u_{n} = z_{n1}$, $u_{n-1} = z_{n-1,1}/z_{n1}$,..., $u_{2}=z_{21}/z_{31}$, 
$u_{1} = q = e^{2\pi\i\tau}$, 
there is a unique solution $F^{(n)}(\zz|\tau)$ with the expansion 
$$
F^{(n)}(\zz|\tau) \simeq z_{21}^{\bar t_{12}}z_{31}^{\bar t_{13} + \bar t_{23}}
... z_{n1}^{\bar t_{1n} + ... + \bar t_{n-1,n}}
\on{exp}\Big(-{\tau\over{2\pi\i}}\Big(\Delta_0 + \sum_{k\geq 0}
a_{2k}\big( \delta_{2k} + \sum_{i<j}(\on{ad}\bar x_i)^{2k}(\bar t_{ij})\big)
\Big)\Big) 
$$
in the region $z_{21}\ll z_{31}\ll... \ll z_{n1}\ll 1$, $\tau\to
\i\infty$, $(\zz,\tau)\in D_n$ (here $z_{ij} = z_i - z_j$);  here 
the sign $\simeq$ 
means that any of the ratios of both sides has the form $1+\sum_{k>0} 
\sum_{i,a_{1},...,a_{n}}r_{k}^{i,a_{1},...,a_{n}}(u_{1},...,u_{n})$, where 
the second sum is finite with $a_{i}\geq 0$, $i\in \{1,...,n\}$, 
$r_{k}^{i,a_{1},...,a_{n}}(u_{1},...,u_{n})$ has degree $k$, 
and is $O(u_{i}(\on{log}u_{1})^{a_{1}}...(\on{log}u_{n})^{a_{n}})$.

\subsection{Presentation of $\Gamma_{1,[n]}$}

According to \cite{Bi2}, $\Gamma_{1,[n]} = \{\overline{\on{B}}_{1,n}\rtimes
\widetilde{\on{SL}_2(\ZZ)}\}/\ZZ$, where $\widetilde{\on{SL}_2(\ZZ)}$ 
is a central 
extension $1\to \ZZ \to \widetilde{\on{SL}_2(\ZZ)} \to \on{SL}_2(\ZZ)\to 1$; 
the action $\alpha : \widetilde{\on{SL}_2(\ZZ)} \to 
\on{Aut}(\overline{\on{B}}_{1,n})$ is such that for $Z$ the 
central element $1\in \ZZ\subset \widetilde{\on{SL}_2(\ZZ)}$, 
$\alpha_Z(x) = Z'x(Z')^{-1}$, where $Z'$ is the image of 
a generator of the center of ${\on{PB}}_n$ (the pure braid group of $n$ points 
on the plane) under the natural morphism 
${\on{PB}}_n \to \overline{\on{B}}_{1,n}$; 
$\overline{\on{B}}_{1,n}\rtimes \widetilde{\on{SL}_2(\ZZ)}$ is then 
$\overline{\on{B}}_{1,n}\times \widetilde{\on{SL}_2(\ZZ)}$ 
with the product $(p,A)(p',A') = (p\alpha_A(p'),AA')$; this 
semidirect product is then factored by its central subgroup 
(isomorphic to $\ZZ$) generated by $((Z')^{-1},Z)$.

$\Gamma_{1,[n]}$ is presented explicitly as follows. Generators are 
$\sigma_i$ ($i=1,...,n-1$), $A_i,B_i$ ($i=1,...,n$), $C_{jk}$ 
($1\leq j<k\leq n$), $\Theta$ and $\Psi$, and relations are: 
$$
\sigma_i\sigma_{i+1}\sigma_i = \sigma_{i+1}\sigma_{i}\sigma_{i+1}\; 
(i=1,...,n-2), 
\quad \sigma_i \sigma_j = \sigma_j \sigma_i \; 
(1\leq i<j\leq n), 
$$
$$
\sigma_i^{-1} X_i \sigma_i^{-1} = X_{i+1}, \quad
\sigma_i Y_i\sigma_i=Y_{i+1} \; (i=1,...,n-1), 
$$
$$
(\sigma_i,X_j) = (\sigma_i,Y_j)=1 \; 
(i\in \{1,...,n-1\}, j\in\{1,...,n\},j\neq i,i+1), 
$$
$$
\quad \sigma_i^2 = C_{i,i+1}C_{i+1,i+2}C_{i,i+2}^{-1}\; (i=1,...,n-1), 
$$
$$
(A_i,A_j) = (B_i,B_j)=1 \on{(any\ }i,j), \quad A_1 = B_1 = 1, 
$$ 
$$
(B_k,A_kA_j^{-1}) = (B_kB_j^{-1},A_k) = C_{jk} \ (1\leq j<k\leq n), 
$$
$$
(A_i,C_{jk}) = (B_i,C_{jk})=1 \ (1\leq i\leq j<k\leq n), 
$$
$$
\Theta A_i\Theta^{-1} = B_i^{-1}, \quad 
\Theta B_i \Theta^{-1} = B_iA_iB_i^{-1}, 
$$
$$
\Psi A_i\Psi^{-1} = A_i, \quad 
\Psi B_i\Psi^{-1} = B_iA_i, \quad 
(\Theta,\sigma_i) = (\Psi,\sigma_i)=1,
$$
$$
(\Psi,\Theta^2)=1, \quad (\Theta\Psi)^3 =
\Theta^4 = C_{12}...C_{n-1,n}. 
$$
Here $X_i = A_i A_{i+1}^{-1}$, $Y_i = B_i B_{i+1}^{-1}$ for $i=1,...,n$
(with the convention $A_{n+1}=B_{n+1}=C_{i,n+1}=1$). 
The relations imply 
$$
C_{jk} = \sigma_{j,j+1...k}...\sigma_{j+n-k,j+n-k+1...n}
\sigma_{j,j+1...n-k+j+1}...\sigma_{k-1,k...n}, 
$$
where $\sigma_{i,i+1...j} = \sigma_{j-1}...\sigma_i$. 
Observe that 
$C_{12},...,C_{n-1,n}$ commute with each other.

The group $\widetilde{\on{SL}_2(\ZZ)}$ is presented by generators 
$\Theta,\Psi$ and $Z$, and relations: $Z$ is central, 
$\Theta^4 = (\Theta\Psi)^3 = Z$ and $(\Psi,\Theta^2)=1$. 
The morphism $\widetilde{\on{SL}_2(\ZZ)}\to \on{SL}_2(\ZZ)$
is $\Theta\mapsto \big( \begin{smallmatrix} 0 & 1 \\ -1 & 0
\end{smallmatrix} \big)$, 
$\Psi\mapsto \big( \begin{smallmatrix} 1 & 1 \\ 0 & 1
\end{smallmatrix} \big)$, and the morphism $\Gamma_{1,[n]}\to 
\on{SL}_2(\ZZ)$ is given by the same formulas and 
$A_i,B_i,\sigma_i\mapsto 1$.

The elliptic braid group $\overline{\on{B}}_{1,n}$
is the kernel of $\Gamma_{1,[n]}\to \on{SL}_2(\ZZ)$; it has the same
presentation as $\Gamma_{1,[n]}$, except for the omission of the generators 
$\Theta,\Psi$ and the relations involving them. The ``pure'' mapping class 
group $\Gamma_{1,n}$ is the kernel of 
$\Gamma_{1,[n]}\to S_n$, $A_i,B_i,C_{jk}\mapsto 1$, 
$\sigma_i\mapsto \sigma_i$; it has the same presentation as $\Gamma_{1,[n]}$, 
except for the omission of the $\sigma_i$. Finally, recall that 
$\overline{\on{PB}}_{1,n}$ is the kernel of 
$\Gamma_{1,[n]}\to \on{SL}_2(\ZZ)\times S_n$.

\begin{remark}
The extended mapping class group $\tilde\Gamma_{1,n}$
of classes of non necessarily orientation-preserving 
self-homeomorphisms of a surface of type $(1,n)$ fits
in a split exact sequence $1\to\Gamma_{1,n}\to\tilde\Gamma_{1,n}
\to\ZZ/2\ZZ\to 1$; it may be viewed as $\{\overline{\on{PB}}_{1,n}
\rtimes \widetilde{\on{GL}_2(\ZZ)}\}/\ZZ$; it has the 
same presentation as $\Gamma_{1,n}$ with the additional generator
$\Sigma$ subject to 
$$
\Sigma^2 = 1, \quad \Sigma\Theta\Sigma^{-1} = \Theta^{-1}, \quad 
\Sigma\Psi\Sigma^{-1} = \Psi^{-1}, \quad 
\Sigma A_i \Sigma^{-1} = A_i^{-1}, \quad 
\Sigma B_i \Sigma^{-1} = A_i B_i A_i^{-1}. 
$$
\end{remark}

\subsection{The monodromy morphisms $\gamma_n : \Gamma_{1,[n]}
\to{\bold G}_n\rtimes S_{n}$}

Let $F(\zz|\tau)$ be a solution of the elliptic KZB system defined on 
$D_n$.

Recall that $D_n := \{(\zz,\tau) \in\CC^n\times \HH | 
z_i = a_i + b_i\tau, a_i,b_i\in\RR, a_1<a_2<...<a_n<a_1+1, 
b_1<b_2<...<b_n<b_1+1\}$. 
The domains $H_n:= \{(\zz,\tau)\in\CC^n\times\HH | z_i = a_i + b_i \tau, 
a_i,b_i\in\RR, a_1<a_2<...<a_n<a_1+1\}$ and $D_n := \{(\zz,\tau)\in
\CC^n\times\HH | z_i = a_i + b_i \tau, a_i,b_i\in\RR, 
b_1<b_2<...<b_n<b_1+1\}$ are also simply connected and invariant, 
and we denote by $F^H(\zz|\tau)$ and $F^V(\zz|\tau)$ the prolongations
of $F(\zz|\tau)$ to these domains.

Then $(\zz,\tau)\mapsto F^H(\zz + \sum_{j=i}^n \delta_i|\tau)$
and $(\zz,\tau)\mapsto e^{2\pi\i(\bar x_i + ... +\bar x_n)}F^V(\zz + 
\tau(\sum_{j=i}^n \delta_i)|\tau)$ are solutions of the elliptic 
KZB system on $H_n$ and $D_n$ respectively. We define $A_i^F,B_i^F
\in {\bold G}_n$ by 
$$
F^H(\zz + \sum_{j=i}^n \delta_i|\tau) = F^H(\zz|\tau)A_i^F, 
\quad 
e^{2\pi\i(\bar x_i + ... +\bar x_n)}F^V(\zz + 
\tau(\sum_{j=i}^n \delta_i)|\tau) = F^V(\zz|\tau)B_i^F. 
$$

The action of $T^{-1} = \big(\begin{smallmatrix} 0 & 1 \\ -1 & 0 
\end{smallmatrix}\big)$ is $(\zz,\tau)\mapsto (-\zz/\tau,-1/\tau)$; 
this transformation takes $H_n$ to $V_n$. Then 
$(\zz,\tau)\mapsto c_{T^{-1}}(\zz|\tau)^{-1} F^V(-\zz/\tau|-1/\tau)$
is a solution of the elliptic KZB system on $H_n$
(recall that $c_{T^{-1}}(\zz|\tau)^{-1}= 
e^{2\pi\i(-\sum_i z_i\bar x_i + \tau X)}(-\tau)^d
= (-\tau)^d e^{(2\pi\i/\tau)(\sum_i z_i\bar x_i + X)}$). 
We define $\Theta^F$ by 
$$
c_{T^{-1}}(\zz|\tau)^{-1} F^V(-\zz/\tau|-1/\tau)
= F^H(\zz|\tau)\Theta^F. 
$$

The action of $S = \big(\begin{smallmatrix} 1 & 1 \\ 0 & 1 
\end{smallmatrix}\big)$ is $(\zz,\tau)\mapsto (\zz,\tau+1)$. 
This transformation takes $H_n$ to itself. Since $c_S(\zz|\tau)=1$, 
the function $(\zz,\tau)\mapsto F^H(\zz,\tau+1)$ is a solution of 
the elliptic KZB system on $H_n$. We define $\Psi^F$ by 
$$
F^H(\zz|\tau+1) = F^H(\zz|\tau)\Psi^F. 
$$
Finally, define $\sigma_{i}^{F}$ by 
$$
\sigma_{i}F(\sigma_{i}^{-1}\zz|\tau) = F(\zz|\tau)\sigma_{i}^{F},
$$ 
where on the l.h.s. $F$ is extended to the universal cover of 
$(\CC^{n}\times\HH) - 
\on{Diag}_{n}$ ($\sigma_{i}$ exchanges $z_{i}$ and $z_{i+1}$, $z_{i+1}$ 
passing to the right of $z_{i}$).

\begin{lemma} \label{lemma:holonomy}
There is a unique morphism 
$\Gamma_{1,[n]}\to {\bold G}_{1,n}\rtimes S_{n}$, taking $X$ to $X^F$, 
where $X = A_i,B_i,\Theta$ or $\Psi$. 
\end{lemma}

{\em Proof.} This follows from the geometric description of 
generators of $\Gamma_{1,[n]}$: if $(\zz_0,\tau_0)\in D_n$, 
then $A_i$ is the class of the projection of the path 
$[0,1]\ni t\mapsto (\zz_0 + t\sum_{j=i}^n \delta_j,\tau_0)$, 
$B_i$ is the class of the projection of $[0,1]\ni t \mapsto 
(\zz_0 + t\tau \sum_{j=i}^n \delta_j,\tau_0)$, 
$\Theta$ is the class of the projection of any path 
connecting $(\zz_0,\tau_0)$ to $(-\zz_0/\tau_0,-1/\tau_0)$
contained in $H_n$, and $\Psi$ is the class of the projection of 
any path connecting $(\zz_0,\tau_0)$ to $(\zz_0,\tau_0+1)$
contained in $H_n$. \hfill \qed \medskip

We will denote by $\gamma_n : \Gamma_{1,[n]} \to {\bold G}_n\rtimes S_{n}$ 
the morphism induced by the solution $F^{(n)}(\zz|\tau)$.

\subsection{Expression of $\gamma_n : \Gamma_{1,[n]} \to 
{\bold G}_n\rtimes S_{n}$ using $\gamma_1$ and $\gamma_2$} \label{5:4}

\begin{lemma}
There exists a unique Lie algebra morphism $\d
\to \bar\t_{1,n}\rtimes\d$, $x\mapsto [x]$, 
such that $[\delta_{2n}] = \delta_{2n} + \sum_{i<j}
(\on{ad}\bar x_i)^{2n}(\bar t_{ij})$, $[X] = X$, $[\Delta_0] = \Delta_0$, 
$[d]=d$.

It induces a group morphism ${\bold G}_1 \to {\bold G}_n$, also denoted
$g\mapsto [g]$. 
\end{lemma}

\begin{lemma}
For each map $\phi : \{1,...,m\} \to \{1,...,n\}$, 
there exists a Lie algebra morphism $\bar\t_{1,n}\to \bar\t_{1,m}$, 
$x\mapsto x^\phi$, defined by $(\bar x_i)^\phi:= \sum_{i'\in\phi^{-1}(i)}
\bar x_{i'}$, $(\bar y_i)^\phi:= \sum_{i'\in\phi^{-1}(i)}
\bar y_{i'}$, $(\bar t_{ij})^\phi:= \sum_{i'\in\phi^{-1}(i),j'\in \phi^{-1}(j)}
\bar t_{i'j'}$.

It induces a group morphism $\on{exp}(\hat{\bar\t}_{1,n}) \to 
\on{exp}(\hat{\bar\t}_{1,m})$, also denoted $g\mapsto g^\phi$. 
\end{lemma}

The proofs are immediate. We now recall the definition and 
properties of the KZ associator (\cite{Dr:Gal}).

If $\kk$ is a field with $\on{char}(\kk)=0$, we let $\t_n^\kk$ be the 
$\kk$-Lie algebra generated by $t_{ij}$, where 
$i\neq j\in \{1,...,n\}$, with relations 
$$
t_{ji} = t_{ij}, \quad [t_{ij}+t_{ik},t_{jk}]=0, \quad
[t_{ij},t_{kl}]=0
$$
for $i,j,k,l$ distinct (in this section, we set $\t_n := \t_n^\CC$). 
For each partially defined map 
$\{1,...,m\}\supset D_\phi \stackrel{\phi}{\to} \{1,...,n\}$, 
we have a Lie algebra morphism $\t_n \to \t_m$, $x\mapsto x^\phi$, 
defined by\footnote{We will also use the notation $x^{I_1,...,I_n}$ 
for $x^{\phi}$, where $I_i = \phi^{-1}(i)$.} 
$(t_{ij})^\phi := \sum_{i'\in \phi^{-1}(i),j'\in\phi^{-1}(j)}
t_{i'j'}$. We also have morphisms $\t_n\to \t_{1,n}$, $t_{ij}\mapsto 
\bar t_{ij}$, compatible with the maps $x\mapsto x^\phi$ on both sides.

The KZ associator $\Phi = \Phi(t_{12},t_{23})\in \on{exp}(\hat\t_3)$ is defined by 
$G_0(z) = G_1(z) \Phi$, where $G_i : ]0,1[\to \on{exp}(\hat\t_3)$ 
are the solutions of $G'(z)G(z)^{-1} = t_{12}/z + t_{23}/(z-1)$
with $G_0(z) \sim z^{t_{12}}$ as $z\to 0$ and $G_1(z) \sim 
(1-z)^{t_{23}}$ as $z\to 1$. The KZ associator satisfies the duality, 
hexagon and pentagon equation (\ref{assoc:1}), (\ref{assoc:2}) below 
(where $\lambda = 2\pi\i$).

\begin{lemma}
$\gamma_2(A_2)$ and $\gamma_2(B_2)$ belong to 
$\on{exp}(\hat{\bar{\t}}_{1,2}) \subset {\bold G}_2$. 
\end{lemma}

{\em Proof.} If $F(\zz|\tau) : H_2 \to {\bold G}_2$ is a solution of the 
KZB equation for $n=2$, then $A_2^F = F^H(\zz + \delta_2|\tau)
F^H(\zz|\tau)^{-1}$ is expressed as the iterated integral, from 
$\zz_0\in D_n$ to $\zz_0 + \delta_2$, of $\bar K_2(\zz|\tau)\in
\hat{\bar\t}_{1,2}$, hence $A_2^F\in \on{exp}(\hat{\bar\t}_{1,2})$. 
Since $\gamma_2(A_2)$ is a conjugate of $A_2^F$, it belongs to 
$\on{exp}(\hat{\bar\t}_{1,2})$ as $\on{exp}(\hat{\bar\t}_{1,2})\subset
{\bold G}_2\rtimes S_{2}$ is normal. One proves similarly that $\gamma_2(B_2) 
\in \on{exp}(\hat{\bar{\t}}_{1,2})$. 
\hfill \qed \medskip

Set 
$$
\Phi_i:= \Phi^{1...i-1,i,i+1...n}...\Phi^{1...n-2,n-1,n}\in \on{exp}(\hat\t_{n}). 
$$
We denote by $x\mapsto \{x\}$ the morphism $\on{exp}(\hat\t_{n}) \to 
\on{exp}(\hat{\bar\t}_{1,n})$ induced by $t_{ij}\mapsto \bar t_{ij}$.

\begin{proposition}
If $n\geq 2$, then 
$$
\gamma_n(\Theta) = 
[\gamma_1(\Theta)] e^{\i{\pi\over 2}\sum_{i<j}\bar t_{ij}}, 
\quad 
\gamma_n(\Psi) = [\gamma_1(\Psi)] e^{\i{\pi\over 6}\sum_{i<j}\bar t_{ij}}, 
$$
and if $n\geq 3$, then 
$$
\gamma_n(A_i) = \{\Phi_i\}^{-1}
\gamma_2(A_2)^{1...i-1,i...n} \{\Phi_i\}, \quad 
\gamma_n(B_i) = \{\Phi_i\}^{-1}
\gamma_2(B_2)^{1...i-1,i...n}\{\Phi_i\}, \; (i=1,...,n), 
$$
$$
\gamma_{n}(\sigma_{i}) = \{\Phi^{1...i-1,i,i+1}\}^{-1}
e^{\i\pi\bar t_{i,i+1}} \{\Phi^{1...i-1,i,i+1}\},\; 
(i=1,...,n-1). 
$$
\end{proposition}

{\em Proof.} In the region $z_{21}\ll z_{31}\ll...\ll z_{n1}\ll 1$, 
$(\zz,\tau)\in D_{n}$, we have 
$$
F^{(n)}(\zz|\tau)\simeq z_{21}^{\bar t_{12}}... 
z_{n1}^{\bar t_{1n}+...+\bar t_{n-1,n}}
\on{exp}(-{{a_{0}}\over{2\pi\i}}(\int_{\i}^{\tau}E_{2}+C)
(\sum_{i<j}\bar t_{ij}))[F(\tau)], 
$$
where $F(\tau) = F^{(1)}(z|\tau)$ for any $z$. Here $C$ is the constant
such that $\int_{\i}^{\tau}E_{2} + C = \tau + o(1)$ as $\tau\to\i\infty$.

We have $F(\tau+1)=F(\tau)\gamma_{1}(\Psi)$, $F(-1/\tau)=F(\tau)
\gamma_{1}(\Theta)$. Since $\sum_{i<j}\bar t_{ij}$ commutes with the image of
$x\mapsto [x]$, we get 
$F^{(n)}(\zz|\tau+1) = F^{(n)}(\zz|\tau)\on{exp}(-{{a_{0}}\over{2\pi\i}}
(\sum_{i<j}\bar t_{ij}))[\gamma_{1}(\Psi)]$, so 
$$
\gamma_{n}(\Psi) = \on{exp}(\i{\pi\over 6}\sum_{i<j}\bar t_{ij})[\gamma_{1}(\Psi)]. 
$$

In the same region, 
\begin{align*}
c_{T^{-1}}(\zz|\tau)^{-1}F^{(n)V}(-{\zz\over\tau}|-{1\over\tau})
\simeq & (-\tau)^{d} e^{{{2\pi\i}\over\tau}(\sum_{i}z_{i}\bar x_{i}+X)}
(-z_{21}/\tau)^{\bar t_{12}}...(-z_{n1}/\tau)^{\bar t_{1n}+...+\bar t_{n-1,n}}
\\ & \on{exp}(-{{a_{0}}\over{2\pi\i}}(\int_{\i}^{-1/\tau}E_{2}+C)
(\sum_{i<j}\bar t_{ij}))[F(-1/\tau)]. 
\end{align*}

Now $E_{2}(-1/\tau) = \tau^{2}E_{2}(\tau) + (6\i/\pi)\tau$, so 
$\int_{\i}^{-1/\tau}E_{2} - \int_{\i}^{\tau}E_{2} = (6\i/\pi)
[\on{log}(-1/\tau) - \on{log}\i]$ (where $\on{log}(re^{i\theta}) = \on{log}r
+\i\theta$ for $\theta\in]-\pi,\pi[$).

It follows that 
\begin{align*}
& c_{T^{-1}}(\zz|\tau)^{-1}F^{(n)V}(-{\zz\over\tau}|-{1\over\tau})
\simeq e^{2\pi\i(\sum_{i}z_{i}\bar x_{i})} z_{21}^{\bar t_{12}}...
z_{n1}^{\bar t_{1n}+...+\bar t_{n-1,n}}\on{exp}( - {{a_{0}}\over{2\pi\i}}
(\int_{\i}^{\tau}E_{2}+C)(\sum_{i<j}\bar t_{ij}))\\
& (\on{exp} - {{a_{0}}\over{2\pi\i}}{{-6\i}\over\pi}
(\on{log}\i) (\sum_{i<j}\bar t_{ij}))
[(-\tau)^{d} e^{(2\pi\i/\tau)X}F(-1/\tau)]\\
& \simeq 
z_{21}^{\bar t_{12}}...
z_{n1}^{\bar t_{1n}+...+\bar t_{n-1,n}}\on{exp}( - {{a_{0}}\over{2\pi\i}}
(\int_{\i}^{\tau}E_{2}+C)(\sum_{i<j}\bar t_{ij})) 
[F(\tau)\gamma_{1}(\Theta)] \on{exp}({{\i\pi}\over 2}
\sum_{i<j}\bar t_{ij})
\\ & 
\simeq F^{(n)H}(\zz|\tau)[\gamma_{1}(\Theta)] \on{exp}({{\i\pi}\over 2}
\sum_{i<j}\bar t_{ij})
\end{align*}
(the second $\simeq$ follows from $\sum_{i}z_{i}\bar x_{i} = \sum_{i>1}z_{i1}
\bar x_{i}$ and $z_{i1}\to 0$), 
so 
$$
\gamma_{n}(\Theta) = [\gamma_{1}(\Theta)]\on{exp}(\i{\pi\over 2}\sum_{i<j}
\bar t_{ij}).
$$

Let $G_i(\zz|\tau)$ be the solution of the elliptic KZB 
system, such that 
\begin{align*}
& G_i(\zz|\tau) \\ & = z_{21}^{\bar t_{12}}...
z_{i-1,1}^{\bar t_{12}+...+\bar t_{1,i-1}}
z_{n,i}^{\bar t_{i,n}+...+\bar t_{n-1,n}}
...z_{n,n-1}^{\bar t_{n-1,n}}
\on{exp}\Big(-{\tau\over{2\pi\i}}\Big(\Delta_0 + \sum_{n\geq 0}
a_{2n}\big(
\delta_{2n} + \sum_{i<j} (\on{ad}\bar x_i)^{2n}(\bar t_{ij})\big)\Big)\Big) 
\end{align*}
when $z_{21}\ll...\ll z_{i-1,1}\ll 1$, $z_{n,n-1}\ll...\ll z_{n,i}\ll 1$, 
$\tau\to\i\infty$ and $(\zz,\tau)\in D_n$. 
Then $G_i(\zz+\sum_{j=i}^n \delta_i|\tau) = G_i(\zz|\tau)
\gamma_2(A_2)^{1...i-1,i...n}$, because in the domain considered
$\bar K_i(\zz|\tau)$ is close to $\bar K_2(z_1,z_n|\tau)^{1...i-1,i...n}$
(where $\bar K_2(...)$ corresponds to the 2-point system); on the other hand, 
$F(\zz|\tau) = G_i(\zz|\tau)\{\Phi_i\}$, which implies the 
formula for $\gamma_n(A_i)$. The formula for $\gamma_n(B_i)$
is proved in the same way. Finally, the behavior of $F^{(n)}(\zz|\tau)$
for $z_{21}\ll ... \ll z_{n1}\ll 1$ is similar to that of a solution of the KZ
equations, which implies the formula for $\gamma_{n}(\sigma_{i})$. 
\hfill \qed \medskip

\begin{remark} One checks that the composition 
$\on{SL}_{2}(\ZZ) \simeq \Gamma_{1,1}\to {\bold G}_{1}\to \on{SL}_{2}(\CC)$
is a conjugation of the canonical inclusion. It follows that the 
composition $\widetilde{\on{SL}_{2}(\ZZ)} \subset \Gamma_{1,n}\to 
{\bold G}_{1}\to \on{SL}_{2}(\CC)$ is a conjugation of the 
canonical projection for any $n\geq 1$. \hfill \qed \medskip 
\end{remark}

Let us set $\tilde A := \gamma_2(A_2)$, $\tilde B := \gamma_2(B_2)$.
The image of $A_{2}A_{3}^{-1} = \sigma_{1}^{-1}A_{2}^{-1}\sigma_{1}^{-1}$ 
by $\gamma_{3}$ yields 
\begin{equation} \label{A:ids}
\tilde A^{12,3} = e^{\i\pi\bar t_{12}} \{\Phi\}^{3,1,2}\tilde A^{2,13}
\{\Phi\}^{2,1,3}e^{\i\pi\bar t_{12}} \cdot 
\{\Phi\}^{3,2,1}\tilde A^{1,23}\{\Phi\}^{1,2,3}
\end{equation}
and the image of $B_{2}B_{3}^{-1} = \sigma_{1}B_{2}^{-1}\sigma_{1}$ 
yields 
\begin{equation} \label{B:ids}
\tilde B^{12,3} = e^{-\i\pi\bar t_{12}} \{\Phi\}^{3,1,2}\tilde B^{2,13}
\{\Phi\}^{2,1,3}e^{-\i\pi\bar t_{12}} \cdot \{\Phi\}^{3,2,1}
\tilde B^{1,23}\{\Phi\}^{1,2,3}. 
\end{equation}

Since $(\gamma_3(A_2),\gamma_3(A_3)) = (\gamma_3(B_2),\gamma_3(B_3))=1$, 
we get 
\begin{equation} \label{gamma13}
(\{\Phi\}^{3,2,1}\tilde A^{1,23}\{\Phi\},\tilde A^{12,3}) = 
(\{\Phi\}^{3,2,1}\tilde B^{1,23}\{\Phi\},\tilde B^{12,3}) = 1 
\end{equation}
(this equation can also be directly derived from (\ref{A:ids}) and 
(\ref{B:ids}) by noting that the l.h.s. is invariant $x\mapsto x^{2,1,3}$ and 
commutes with $e^{\pm\i\pi\bar t_{12}}$). 
We have for $n=2$, $C_{12} = (B_{2},A_{2})$, so 
$(\tilde A,\tilde B) = \gamma_{2}(C_{12})^{-1}$. 
Also $\gamma_{1}(\Theta)^{4}=1$, so $\gamma_{2}(C_{12})
= \gamma_{2}(\Theta)^{4} = (e^{\i\pi\bar t_{12}/2}[\gamma_{1}(\Theta)])^{4}
= e^{2\pi\i\bar t_{12}}[\gamma_{1}(\Theta)^{4}] = e^{2\pi\i\bar t_{12}}$, 
so 
\begin{equation} \label{gamma12''}
(\tilde A,\tilde B) = e^{-2\pi\i\bar t_{12}}. 
\end{equation}

For $n=3$, we have $\gamma_{3}(\Theta)^{4} = e^{2\pi\i(\bar t_{12} + \bar t_{13} 
+ \bar t_{23})} = \gamma_{3}(C_{12}C_{23})$; 
since $\gamma_{3}(C_{12}) = (\gamma_{3}(B_{2}),\gamma_{3}(A_{2}))
= \{\Phi\}^{-1}(\tilde B,\tilde A)^{1,23}\{\Phi\} 
= \{\Phi\}^{-1}e^{2\pi\i(\bar t_{12}+\bar t_{13})}\{\Phi\}$, 
we get $\gamma_{3}(C_{23}) = \{\Phi\}^{-1}e^{2\pi\i\bar t_{23}}\{\Phi\}$. 
The image by $\gamma_{3}$
of $(B_{3},A_{3}A_{2}^{-1}) = (B_{3}B_{2}^{-1},A_{3}) = C_{23}$ then gives 
\begin{equation} \label{gamma13'}
(\tilde B^{12,3},\tilde A^{12,3}\{\Phi\}^{-1}(\tilde A^{1,23})^{-1}\{\Phi\}) =
(\tilde B^{12,3}\{\Phi\}^{-1}(\tilde B^{1,23})^{-1}\{\Phi\},\tilde A^{12,3}) 
= \{\Phi\}^{-1}e^{2\pi\i\bar t_{23}}\{\Phi\} 
\end{equation}
(applying $x\mapsto x^{\emptyset,1,2}$, this identity implies 
(\ref{gamma12''})).

Let us set $\tilde\Theta := \gamma_{1}(\Theta)$, $\tilde\Psi:= \gamma_{1}(\Theta)$. 
Since $\gamma_{1},\gamma_{2}$ are group morphisms, we have
\begin{equation} \label{gamma11}
\tilde\Theta^{4} = (\tilde\Theta\tilde\Psi)^{3}= (\tilde\Theta^{2},
\tilde\Psi)=1, 
\end{equation}
\begin{equation} \label{gamma12}
[\tilde\Theta]e^{\i{\pi\over 2}\bar t_{12}}
\tilde A
([\tilde\Theta]e^{\i{\pi\over 2}\bar t_{12}})^{-1} = \tilde B^{-1}, \quad 
[\tilde\Theta]e^{\i{\pi\over 2}\bar t_{12}}
\tilde B
([\tilde\Theta]e^{\i{\pi\over 2}\bar t_{12}})^{-1} = \tilde B\tilde A
\tilde B^{-1}, 
\end{equation}
\begin{equation} \label{gamma12'}
[\tilde\Psi]e^{\i{\pi\over 6}\bar t_{12}}
\tilde A
([\tilde\Psi]e^{\i{\pi\over 6}\bar t_{12}})^{-1} = \tilde A, \quad 
[\tilde\Psi]e^{\i{\pi\over 6}\bar t_{12}}
\tilde B
([\tilde\Psi]e^{\i{\pi\over 6}\bar t_{12}})^{-1} = \tilde B\tilde A. 
\end{equation}
(\ref{gamma11}) (resp., (\ref{gamma12}), (\ref{gamma12'})) 
are identities in ${\bold G}_{1}$ (resp., ${\bold G}_{2}$); in 
(\ref{gamma12}), (\ref{gamma12'}), 
$x\mapsto [x]$ is induced by the map $\d\to\d\rtimes\bar\t_{1,2}$ defined above.

\subsection{Expression of $\tilde\Psi$ and of $\tilde A$ and $\tilde B$ 
in terms of $\Phi$}

In this section, we compute 
$\tilde A$ and $\tilde B$ in terms of the KZ associator $\Phi$. We also 
compute $\tilde\Psi$.

Recall the definition of $\tilde\Psi$. The elliptic KZB system for $n=1$ is 
$$
2\pi\i\partial_{\tau}F(\tau) + \big(\Delta_{0} + 
\sum_{k\geq 1}a_{2k}E_{2k+2}(\tau)\delta_{2k} \big) F(\tau)=0. 
$$
The solution $F(\tau) := F^{(1)}(z|\tau)$ (for any $z$) is determined 
by $F(\tau) \simeq \on{exp}( - {\tau\over{2\pi\i}}(\Delta_{0} + \sum_{k\geq 1}
a_{2k}\delta_{2k}))$. Then $\tilde\Psi$ is determined by 
$F(\tau+1) = F(\tau)\tilde\Psi$. We have therefore:

\begin{lemma}
$\tilde\Psi = \on{exp}(- {1\over{2\pi\i}}(\Delta_{0} + \sum_{k\geq 1} 
a_{2k}\delta_{2k}))$. 
\end{lemma}

Recall the definition of $\tilde A$ and $\tilde B$. The elliptic KZB 
system for $n=2$ is 
\begin{equation} \label{ell:KZ}
\partial_z F(z|\tau) = -\big( {{\theta(z+\on{ad}x|\tau)\on{ad}x}\over
{\theta(z|\tau)\theta(\on{ad}x|\tau)}} \big)(y) \cdot  F(z|\tau),
\end{equation}
\begin{equation} \label{ell:KZ2}
2\pi\i\partial_\tau F(z|\tau) + \big(\Delta_0 
+ \sum_{k\geq 1} a_{2k} E_{2k+2}(\tau)
\delta_{2k} - g(z,\on{ad}x|\tau)(t) \big)F(z|\tau) = 0, 
\end{equation}
where $z = z_{21}$, $x = \bar x_2 = -\bar x_1$, 
$y = \bar y_2 = -\bar y_1$, $t = \bar t_{12} = -[x,y]$.

The solution $F(z|\tau):= F^{(2)}(z_1,z_2|\tau)$ is determined by its behavior 
$F(z|\tau) \simeq z^t \on{exp}( - {\tau\over{2\pi\i}} \big( \Delta_0 
+ \sum_{k\geq 0} a_{2k}(\delta_{2k} + (\on{ad}x)^{2k})(t)\big))$ 
when $z\to 0^+$, $\tau\to\i\infty$. We then have 
$F^H(z+1|\tau) = F^H(z|\tau)\tilde A$, $e^{2\pi\i x}
F^V(z+\tau|\tau) = F^V(z|\tau)\tilde B$.

\begin{proposition} We have\footnote{By convention, if 
$z\in \CC \setminus \RR_{-}$ and $x\in \n$, where $\n$ is a pronilpotent 
Lie algebra, then $z^{x}$ is $\exp(x\log z)\in \on{exp}(\n)$, where 
$\log z$ is chosen with imaginary part in $]-\pi,\pi[$.}
$$ 
\tilde A =  (2\pi/\i)^{t}\Phi(\tilde y,t)e^{2\pi\i \tilde y}
\Phi(\tilde y,t)^{-1}(\i/2\pi)^{t} = (2\pi)^t \i^{-3t}
\Phi(-\tilde y-t,t) e^{2\pi\i(\tilde y+t)}\Phi(-\tilde y-t,t)^{-1}
(2\pi\i)^{-t},
$$ 
where $\tilde y = - {{\on{ad}x}\over{e^{2\pi\i\on{ad}x}-1}}(y)$. 
\end{proposition}

{\em Proof.} $\tilde A = F^H(z|\tau)^{-1}F^H(z+1|\tau)$, which we will 
compute in the limit $\tau\to\i\infty$. For this, we will 
compute $F(z|\tau)$ in the limit $\tau\to\i\infty$. 
In this limit, $\theta(z|\tau) =  (1/\pi)\on{sin}(\pi z)[1 + O(e^{2\pi\i\tau})]$  so 
the system becomes
\begin{equation} \label{trigo:KZ}
\partial_z F(z|\tau) = \big(
\pi\on{cotg}(\pi z)t  - \pi\on{cotg}(\pi\on{ad}x)\on{ad}x(y) 
+ O(e^{2\pi\i\tau}) \big) F(z|\tau)
\end{equation}
$$
2\pi\i\partial_\tau F(z|\tau) + \big( \Delta_0 + \sum_{k\geq 1} 
a_{2k}\delta_{2k}
+ ( {{\pi^2} \over {\on{sin}^2(\pi\on{ad}x)}} - {1 \over 
{(\on{ad}x)^2}})(t) + O(e^{2\pi\i\tau}) \big) F(z|\tau) = 0
$$
where the last equation is 
$$ 
2\pi\i\partial_\tau F(z|\tau) + \big( \Delta_0 + a_{0}t + \sum_{k\geq 1} a_{2k}
(\delta_{2k} + (\on{ad}x)^{2k}(t)) + O(e^{2\pi\i\tau}) \big)F(z|\tau) = 0. 
$$
We set 
$$
\Delta := \Delta_0 + \sum_{k\geq 1} a_{2k}\delta_{2k}, \quad \on{so} \quad 
\Delta_0 + a_{0}t + \sum_{k\geq 1} a_{2k}
(\delta_{2k} + (\on{ad}x)^{2k}(t)) = [\Delta] + a_0 t. 
$$
The compatibility of this system implies that 
$[\Delta] + a_0 t$ commutes with $t$ and 
$(\pi\on{ad}x)\on{cotg}(\pi\on{ad}x)(y) = \i\pi(-t-2\tilde y)$, hence with 
$t$ and $\tilde y$; actually $t$ commutes with each 
$[\delta_{2k}] = \delta_{2k} + (\on{ad}x)^{2k}(t)$.

Equation (\ref{ell:KZ}) can be written $\partial_z F(z|\tau) 
= (t/z + O(1))F(z|\tau)$. We then let $F_0(z|\tau)$ be the solution of 
(\ref{ell:KZ}) in $V:= \{(z,\tau)|\tau\in\HH,z = a + b \tau, a\in ]0,1[, b\in\RR\}$
such that $F_0(z|\tau) \simeq z^t$ when 
$z\to 0^+$, for any $\tau$. This means that the left (equivalently, right)
ratio of these quantities has the form $1 + \sum_{k>0}(\on{degree}\ k)
O(z(\on{log}z)^{f(k)})$ where $f(k)\geq 0$.

We now relate $F(z|\tau)$ and $F_{0}(z|\tau)$. Let $F(\tau) = F^{(1)}(z|\tau)$
for any $z$ be the solution of the KZB system for $n=1$, such that 
$F(\tau) \simeq \on{exp}(- {\tau\over{2\pi\i}}\Delta)$ as $\tau\to\i\infty$
(meaning that the left, or equivalently right, ratio of these quantities 
has the form $1 + \sum_{k>0}(\on{degree}\ k)O(\tau^{f(k)}
e^{2\pi\i\tau})$, where $f(k)\geq 0$).

\begin{lemma}
We have $F(z|\tau) = F_{0}(z|\tau) 
\on{exp}(-{a_{0}\over{2\pi\i}}(\int_{\i}^{\tau}
E_{2}+C)t)[F(\tau)]$, where $C$ is such that $\int_{\i}^{\tau}E_{2}+C = \tau 
+ O(e^{2\pi\i\tau})$. 
\end{lemma}

{\em Proof of Lemma.}
$F(z|\tau) = F_0(z|\tau)X(\tau)$, where 
$X : \HH\to {\bold G}_2$ is a map. We have $g(z,\on{ad}x|\tau)(t) 
= a_{0}E_2(\tau)t + \sum_{k>0}a_{2k}E_{2k+2}(\tau)
(\on{ad}x)^{2k}(t) + O(z)$ when $z\to 0^{+}$ and for any $\tau$, so 
(\ref{ell:KZ2}) is written as 
$$
2\pi\i\partial_{\tau}F(z|\tau) + \big(\Delta_{0} + a_{0}E_2(\tau)t
+\sum_{k>0}a_{2k}E_{2k+2}(\tau)[\delta_{2k}] + O(z)
\big)F(z|\tau) = 0
$$
where $O(z)$ has degree $>0$. 
Since $\Delta_{0},t$ and the $[\delta_{2k}]$ all commute with $t$, 
the ratio $F_{0}(z|\tau)^{-1}F(z|\tau)$ satisfies 
$$
2\pi\i\partial_{\tau}(F_{0}^{-1}F(z|\tau))
+ \big( \Delta_{0}+ a_{0}E_2(\tau)t 
+ \sum_{k>0} a_{2k}E_{2k+2}(\tau)[\delta_{2k}] 
+  \sum_{k>0} (\on{degree}\ k)O(z(\on{log}z)^{h(k)})\big) 
(F_{0}^{-1}F(z|\tau))=0
$$
where $h(k)\geq 0$. Since $F_{0}(z|\tau)^{-1}F(z|\tau) = X(\tau)$ 
is in fact independent on $z$, we have 
$$
2\pi\i\partial_{\tau}(X(\tau)) + \big( \Delta_{0}+ a_{0}E_2(\tau)t + 
\sum_{k>0} a_{2k}E_{2k+2}(\tau)[\delta_{2k}] \big) (X(\tau))=0, 
$$
which implies that $X(\tau) = \on{exp}( 
- {a_{0}\over{2\pi\i}}(\int_{\i}^{\tau}
E_{2}+C)t)[F(\tau)] X_{0}$, where $X_{0}$ is a suitable element in 
${\bold G}_{2}$. 
The asymptotic behavior of $F(z|\tau)$ when $\tau\to\i\infty$ 
and $z\to 0^{+}$ then implies $X_{0}=1$. 
\hfill \qed \medskip

{\em End of proof of Proposition.}
We then have $F(z|\tau) = F_0(z|\tau)X(\tau)$, where 
$X(\tau) \simeq \on{exp}( - {\tau\over{2\pi\i}}([\Delta] + a_0 t))$
as $\tau\to\i\infty$, where this means that the left ratio (equivalently, the right ratio)
of these quantities has the form $1 + \sum_{k>0} (\on{degree}\ k)
O(\tau^{x(k)}e^{2\pi\i\tau})$, where $x(k)\geq 0$.

If we set $u := e^{2\pi\i z}$, then (\ref{ell:KZ}) is rewritten as 
\begin{equation} \label{KZ'}
\partial_{u} \bar F (u|\tau)= (\tilde y/u + t/(u-1) + 
O(e^{2\pi\i\tau})) \bar F(u|\tau), 
\end{equation}
where $\bar F(u|\tau) = F(z|\tau)$.

Let $D' := \{u||u|\leq 1\} - [0,1]$ be the complement of the unit interval 
in the unit disc. Then we have a bijection 
$\{(z,\tau) | \tau\in \i\RR_+^\times, z = a + \tau b, a\in [0,1],b\geq 0\} 
\to D' \times \i\RR_+^\times$, given by $(z,\tau)\mapsto (u,\tau) 
:= (e^{2\pi\i z},\tau)$.

Let $\bar F_{a},\bar F_{f}$ be the solutions of (\ref{KZ'}) in 
$D' \times \i\RR_+$, such that $\bar F_{a}(u|\tau) \simeq 
((u-1)/(2\pi\i))^{t}$ when $u = 1 +\i 0^{+}$, and for any $\tau$, and 
$\bar F_{f}(u|\tau) \simeq e^{\i\pi t}((1-u)/(2\pi\i))^{t}$ 
when $u=1-\i 0^{+}$, for any $\tau$.

Then  one checks that $F_0(z|\tau) = \bar F_{a}(e^{2\pi\i z}|\tau)$, 
$F_0(z-1|\tau) = \bar F_{f}(e^{2\pi\i z}|\tau)$ when 
$(z,\tau)\in \{(z,\tau)|\tau\in\i\RR_+^\times, 
z = a+\tau b|a\in [0,1],b\geq 0\}$.

We then define $\bar F_{b},...,\bar F_{e}$ as the solutions of (\ref{KZ'}) 
in $D' \times \i\RR_+^\times$, such that: 
$\bar F_{b}(u|\tau) \simeq (1-u)^{t}$ as $u = 1 -0^{+}$, $\Im(u)>0$ 
for any $\tau$, 
$\bar F_{c}(u|\tau) \simeq u^{\tilde y}$ as $u\to 0^{+}$, $\Im(u)>0$ for 
any $\tau$, 
$\bar F_{d}(u|\tau) \simeq u^{\tilde y}$ as $u\to 0^{+}$, $\Im(u)<0$ for 
any $\tau$, 
$\bar F_{e}(u|\tau) \simeq (1-u)^{t}$ as $u = 1 -0^{+}$, $\Im(u)<0$ 
for any $\tau$.

Then $\bar F_{b} = \bar F_{a}(-2\pi\i)^{t}$, $\bar F_{c}(-|\tau) = 
\bar F_{b}(-|\tau)[\Phi(\tilde y,t) + O(e^{2\pi\i\tau})]$, 
$\bar F_{d}(-|\tau) = \bar F_{c}(-|\tau) e^{-2\pi\i \tilde y}$, 
$\bar F_{e}(-|\tau) 
= \bar F_{d}(-|\tau)[\Phi(\tilde y,t)^{-1} + O(e^{2\pi\i\tau})]$, 
$\bar F_{f} = \bar F_{e} (\i/2\pi)^{t}$.

So $\bar F_{f}(-|\tau) = \bar F_{a}(-|\tau) \big( 
(-2\pi\i)^{t}\Phi(\tilde y,t)e^{-2\pi\i \tilde y}
\Phi(\tilde y,t)^{-1}(\i/2\pi)^{t} + O(e^{2\pi\i\tau})\big)$. 
It follows that $F_0(z+1|\tau) = F_0(z|\tau)A(\tau)$, where
$$
A(\tau)  =  (-2\pi\i)^{t}\Phi(\tilde y,t)e^{2\pi\i \tilde y}
\Phi(\tilde y,t)^{-1}(\i/2\pi)^{t} + O(e^{2\pi\i\tau}).
$$

Now 
\begin{align*}
& \tilde A = F(z|\tau)^{-1}F(z+1|\tau) = X(\tau)^{-1} A(\tau)X(\tau)
= \big( 1 + \sum_{k>0}(\on{degree}\ k) O(\tau^{x(k)}e^{2\pi\i\tau}) \big)^{-1} 
\\ & 
\on{exp}({\tau\over{2\pi\i}}([\Delta] + a_0t))
\big(  (-2\pi\i)^{t}\Phi(\tilde y,t)e^{2\pi\i \tilde y}
  \Phi(\tilde y,t)^{-1}(\i/2\pi)^{t} + O(e^{2\pi\i\tau})\big)
\\ & \on{exp}(-{\tau\over{2\pi\i}}([\Delta] + a_0t))
\big( 1 + \sum_{k>0}(\on{degree}\ k) O(\tau^{x(k)}e^{2\pi\i\tau})\big).
\end{align*} 
As we have seen, $[\Delta]+a_0t$ commutes with $\tilde y$ and $t$; 
on the other hand, 
\begin{align*}
& \on{exp}({\tau\over{2\pi\i}}([\Delta] + a_0t))
O(e^{2\pi\i\tau})
\on{exp}(-{\tau\over{2\pi\i}}([\Delta] + a_0t))
\\ & = \on{exp}(\tau \on{ad}({{[\Delta] + a_0t}\over{2\pi\i}}))
(O(e^{2\pi\i\tau}))
  = \sum_{k\geq 0} (\on{degree}\ k) O(\tau^{n_1(k))}e^{2\pi\i\tau})
  \end{align*}
where $n_{1}(k)\geq 0$, as $[\Delta]+ a_0t$ is a sum of terms of positive 
degree and of $\Delta_0$, which is locally ad-nilpotent.

Then 
\begin{align*}
& \tilde A = 
\big( 1 + \sum_{k>0}(\on{degree}\ k) O(\tau^{x(k)}e^{2\pi\i\tau}) \big)^{-1} 
\big(  (-2\pi\i)^{t}\Phi(\tilde y,t)e^{2\pi\i \tilde y}
  \Phi(\tilde y,t)^{-1}(\i/2\pi)^{t} \\ & + \sum_{k\geq 0}(\on{degree\ }k) 
O(\tau^{n_1(k)}e^{2\pi\i\tau})\big)
  \big( 1 + \sum_{k>0}(\on{degree}\ k) O(\tau^{x(k)}e^{2\pi\i\tau})\big).
\end{align*}
It follows that 
$$
\tilde A =  (-2\pi\i)^{t}\Phi(\tilde y,t)e^{2\pi\i \tilde y}
  \Phi(\tilde y,t)^{-1}(\i/2\pi)^{t} + \sum_{k\geq 0}(\on{degree\ }k) 
O(\tau^{n_2(k)}e^{2\pi\i\tau}), 
$$
where $n_{2}(k)\geq 0$, which implies the first formula for $\tilde A$. 
The second formula either follows from the first one by using the hexagon 
identity, or can be obtained repeating the above argument using a path 
$1 \to +\infty \to 1$, winding around $1$ and $\infty$. \hfill \qed \medskip

We now prove:

\begin{theorem} \label{thm:B}
$$
\tilde B =  (2\pi\i)^{t} \Phi(-\tilde y-t,t)  e^{2\pi\i x} 
\Phi(\tilde y,t)^{-1}(2\pi/\i)^{-t}. 
$$
\end{theorem}

{\em Proof.} We first define $F_{0}(z|\tau)$ as the solution 
in $V:= \{a+b\tau | a\in]0,1[,b\in\RR\}$ of (\ref{ell:KZ}) such that 
$F_{0}(z|\tau)\sim z^{t}$ as $z\to 0^{+}$. Then there exists $B(\tau)$ such 
that $e^{2\pi\i x}F_{0}(z+\tau|\tau) 
= F_{0}(z|\tau) B(\tau)$. We compute the asymptotics of $B(\tau)$ as 
$\tau\to\i\infty$.

We define four asymptotic zones ($z$ is assumed to remain on the segment 
$[0,\tau]$, and $\tau$ in the line $\i\RR_{+}$): 
(1) $z\ll 1 \ll \tau$, (2) $1\ll z \ll \tau$, (3) $1\ll\tau-z \ll \tau$, 
(4) $\tau-z\ll1 \ll \tau$.

In the transition (1)-(2), the system takes the form (\ref{trigo:KZ}), or 
if we set $u := e^{2\pi\i z}$, (\ref{KZ'}).

In the transition (3)-(4), $G(z'|\tau) := e^{2\pi\i x}F(\tau+z'|\tau)$
satisfies (\ref{ell:KZ}), so $\bar G(u'|\tau) = 
e^{2\pi\i x}F(\tau+z'|\tau)$ satisfies (\ref{KZ'}), where 
$u' = e^{2\pi\i z'}$.

We now compute the form of the system in the transition (2)-(3). 
We first prove:

\begin{lemma}
Set $u := e^{2\pi\i z}$, $v := e^{2\pi\i(\tau-z)}$. When 
$0<\Im(z)<\Im(\tau)$, we have $|u|<1$, $|v|<1$. 
When $k\geq 0$, $(\theta^{(k)}/\theta)(z|\tau) = (-\i\pi)^{k} 
+ \sum_{s,t\geq 0,s+t>0} a_{st}^{(k)}u^{s}v^{t}$, where the sum in the
r.h.s. is convergent in the domain $|u|<1$, $|v|<1$. 
\end{lemma}

{\em Proof.} This is clear if $k=0$. Set $q = uv = e^{2\pi\i\tau}$. 
We have 
$\theta(z|\tau) = u^{1/2}\prod_{s>0}(1-q^s u) \prod_{s\geq 0}(1-q^s u^{-1})
\cdot (2\pi\i)^{-1}\prod_{s>0}(1-q^s)^{-2}$, so 
\begin{align*}
(\theta'/\theta)(z|\tau) & = \i\pi - 2\pi\i \sum_{s>0} q^s u/(1-q^s u)
+ 2\pi\i \sum_{s\geq 0} q^s u^{-1}/(1-q^s u^{-1})
\\
  & = -\i\pi - 2\pi\i\sum_{s\geq 0} {{u^{s+1}v^{s}}\over{1-u^{s+1}v^{s}}}
  + 2\pi\i \sum_{s\geq 0} {{u^{s}v^{s+1}}\over {1-u^{s}v^{s+1}}}
  = -\i\pi + \sum_{s+t>0} a_{st} u^{s}v^{t}, 
\end{align*} 
where $a_{st} = 2\pi\i$ if $(s,t) = k(r,r+1)$, $k>0$, $r\geq 0$, and 
$a_{st}=-2\pi\i$ if $(s,t) = k(r+1,r)$, $k>0$, $r\geq 0$. One checks 
that this series is convergent in the domain $|u|<1$, $|v|<1$. 
This proves the lemma for $k=1$.

We then prove the remaining cases by induction, using 
$$
{{\theta^{(k+1)}} \over \theta}(z|\tau) 
=  {{\theta^{(k)}} \over \theta}(z|\tau) {\theta' \over \theta}(z|\tau) 
+ {\partial\over{\partial z}} {\theta^{(k)} \over \theta}(z|\tau).
$$
  \hfill \qed \medskip

Using the expansion 
\begin{align*}
& {{\theta(z+x|\tau)x}\over{\theta(z|\tau)\theta(x|\tau)}} 
= {x\over{\theta(x|\tau)}} \sum_{k\geq 0} 
(\theta^{(k)}/\theta)(z|\tau){{x^{k}}\over{k!}} \\
  & = 
{{\pi x}\over{\on{sin}(\pi x)}} (1 + \sum_{n>0} q^{n} P_{n}(x))
\Big(\sum_{k\geq 0} ((-\i\pi)^{k} + \sum_{s+t>0} a_{st}^{(k)}u^{s}v^{t})
{{x^{k}}\over{k!}}\Big) 
\\ & = {{\pi x}\over {\on{sin}(\pi x)}}e^{-\i\pi x} + 
\sum_{s+t>0} a_{st}(x)u^{s}v^{t}
= {{2\i\pi x}\over{e^{2\i\pi x}-1}} + \sum_{s+t>0} a_{st}(x)u^{s}v^{t} ,
\end{align*} 
the form of the system in the transition (2)-(3) is 
\begin{align} \label{transition}
\partial_{z}F(z|\tau) & = \big(
- {{2\i\pi\on{ad}x}\over{e^{2\i\pi\on{ad}x}-1}}(y) 
+ \sum_{s,t|s+t>0} a_{st} u^{s}v^{t}\big) 
F(z|\tau) \nonumber \\
& = \big( 2\i\pi\tilde y + \sum_{s,t|s+t>0} a_{st} u^{s}v^{t}\big)F(z|\tau), 
\end{align}
where each homogeneous part of $\sum_{s,t} a_{st} u^{s}v^{t}$
converges for $|u|<1$, $|v|<1$.

\begin{lemma}
There exists a solution $F_{c}(z|\tau)$ of (\ref{transition}) defined for 
$0<\Im(z)<\Im(\tau)$, such that 
$$
F_{c}(z|\tau) = u^{\tilde y}(1 + \sum_{k>0} \sum_{s\leq s(k)}
\on{log}(u)^{s} f_{ks}(u,v)) 
$$
($\on{log}u = i\pi z$, $u^{\tilde y} = e^{2\pi\i z \tilde y}$),
where $f_{ks}(u,v)$ is an analytic function taking its values in the 
homogeneous part of the algebra of degree $k$, convergent for $|u|<1$ and 
$|v|<1$, and vanishing at $(0,0)$. 
This function is uniquely defined up to right multiplication by 
an analytic function of the form $1 + \sum_{k>0} a_{k}(q)$
(recall that $q = uv$), where $a_{k}(q)$ is an analytic function on 
$\{q||q|<1\}$, vanishing at $q=0$, with values in the degree $k$ part 
of the algebra. 
\end{lemma}

{\em Proof of Lemma.} We set $G(z|\tau) := u^{-\tilde y}F(z|\tau)$, 
so $G(z|\tau)$ should satisfy 
$$
\partial_{z}G(z|\tau) = \on{exp}(- \on{ad}(\tilde y)\on{log}u)
\{\sum_{s+t>0} a_{st}u^{s}v^{t}\} G(z|\tau), 
$$
which has the general form 
$$
\partial_{z}G(z|\tau) = \Big(
\sum_{k>0} \sum_{s\leq a(k)} \on{log}(u)^{s} a_{ks}(u,v)\Big) G(z|\tau), 
$$
where $a_{ks}(u,v)$ is analytic in $|u|<1$, $|v|<1$ and vanishes at $(0,0)$. 
We show that this system 
admits a solution of the form $1 + \sum_{k>0} \sum_{s \leq s(k)} 
\on{log}(u)^{s} f_{ks}(u,v)$, with $f_{ks}(u,v)$ analytic in $|u|<1$, $|v|<1$, 
in the degree $k$ part of the algebra, vanishing at $(0,0)$ for $s\neq 0$.
For this, we solve inductively  (in $k$) the system of equations
\begin{equation} \label{inductive}
\partial_{z}\big(\sum_{s} (\on{log}u)^{s} f_{ks}(u,v)\big) = 
\sum_{s',s'',k',k''|k'+k'' = k} (\on{log}u)^{s'+s''} 
a_{k's'}(u,v)f_{k''s''}(u,v). 
\end{equation}
Let $\cO$ be the ring of analytic functions on $\{(u,v)| |u|<1 , |v|<1\}$ 
(with values in a finite dimensional vector space) and $\m\subset \cO$
be the subset of functions vanishing at $(0,0)$. We have an injection 
$\cO[X] \to \{$analytic functions in $(u,v)$, $|u|<1$, $|v|<1$, $u\notin
\RR_{-}\}$, given by $f(u,v)X^{k}\mapsto (\on{log}u)^{k}f(u,v)$. The 
endomorphism ${\partial \over {\partial z}} 
= 2\pi\i(u{\partial\over{\partial u}} - v {\partial\over{\partial v}})$ 
then corresponds to the endomorphism of  $\cO[X]$ given by 
$2\pi\i ( {\partial\over{\partial X}} + u{\partial\over{\partial u}}
- v {\partial\over{\partial v}})$. It is surjective, and restricts to a 
surjective endomorphism of $\m[X]$. The latter surjectivity implies that 
equation (\ref{inductive}) can be solved.

Let us show that the solution $G(z|\tau)$ is unique up to right 
multiplication by functions of $q$ like in the lemma. The ratio of two 
solutions is of the form $1 + \sum_{k>0} \sum_{s\leq s(k)}
\on{log}(u)^{s} f_{ks}(u,v)$ and is killed by $\partial_{z}$. Now 
the kernel of the endomorphism of $\m[X]$ given by 
$2\pi\i ({\partial\over{\partial X}} + u{\partial\over{\partial u}}
- v {\partial\over{\partial v}})$ is $m^{*}(\m_1)$, where 
$m^{*}(\m_{1}) \subset \m$ is the set of all functions of the form $a(uv)$, 
where $a$ is an analytic function on $\{q||q|<1\}$ vanishing at $0$. 
This implies that the ratio of two solutions is as above. 
\hfill \qed \medskip

{\em End of proof of Theorem.}
Similarly, there exists a solution $F_{d}(z|\tau)$ of (\ref{transition})
defined in the same domain, such that 
$$
F_{d}(z|\tau) = v^{-\tilde y}(1+\sum_{k>0}\sum_{s\leq t(k)}
\on{log}(v)^{t} g_{ks}(u,v)), 
$$
where $b_{ks}(u,v)$ is as above (and $\on{log}v = \i\pi(\tau-z)$, 
$v^{-\tilde y} = \on{exp}(2\pi\i(z-\tau)\tilde y)$). $F_{d}(z|\tau)$ 
is defined up to right multiplication by a function of $q$ as above.

We now study the ratio $F_{c}(z|\tau)^{-1}F_{d}(z|\tau)$. This is a 
function of $\tau$ only, and it has the form 
$$
q^{-\tilde y} \Big( 1 + \sum_{k>0} \sum_{s\leq s(k), t\leq t(k)}
(\on{log}u)^{s}(\on{log}v)^{t} a_{kst}(u,v) \Big) 
$$
where $a_{kst}(u,v)\in \m$ (as $v^{-\tilde y}(1 + \sum_{k>0} \sum_{s\leq s(k)}
(\on{log}u)^{s} c_{ks}(u,v))v^{\tilde y}$ has the form 
$1 + \sum_{k>0}\sum_{s,t\leq t(k)}$ $ (\on{log}u)^{s}(\on{log}v)^{t} d_{ks}(u,v)$, 
where $d_{ks}(u,v)\in \m$ if $c_{ks}(u,v)\in\m$). 
Set $\on{log}q:= \on{log}u+\on{log}v = 2\pi\i\tau$, 
then this ratio can be rewritten 
$q^{-\tilde y}\{
1 + \sum_{k>0} \sum_{s\leq s(k), t\leq t(k)}(\on{log}u)^{s}(\on{log}q)^{t} 
b_{kst}(u,v)\}$ where $b_{kst}(u,v)\in\m$, and since the product of this 
ratio with $q^{\tilde y}$ is killed by $\partial_{z}$
(which identifies with the endomorphism $2\pi\i ({\partial\over{\partial X}} 
+ u{\partial\over{\partial u}} - v {\partial\over{\partial v}})$ of $\cO[X]$), 
the ratio is in fact of the form 
$$
F_{c}^{-1}F_{d}(z|\tau) = 
q^{\-\tilde y}(1 + \sum_{k>0}\sum_{s\leq s(k)} (\on{log}q)^{s} a_{ks}(q)), 
$$
where $a_{ks}$ is analytic in $\{q||q|<1\}$, vanishing at $q=0$.

It follows that 
\begin{equation} \label{ratio:cd}
F_{c}^{-1}F_{d}(z|\tau) 
= e^{-2\pi\i\tau\tilde y}(1 + \sum_{k>0} (\on{degree}\ k)
O(\tau^{k}e^{-2\pi\i\tau})). 
\end{equation}

In addition to $F_{c}$ and $F_{d}$, which have prescribed behaviors in zones
(2) and (3), we define solutions of (\ref{ell:KZ}) in $V$ by prescribing 
behaviors in the remaining asymptotic zones: 
$F_{a}(z|\tau) \simeq z^{t}$ when $z\to 0^{+}$ for any $\tau$; 
$F_{b}(z|\tau) \simeq (2\pi z/\i)^{t}$ when $z\to \i 0^{+}$ for any $\tau$
(in particular in zone (1)); 
$e^{2\pi\i x}F_{e}(z|\tau) \simeq (2\pi(\tau-z)/\i)^{t}$ when 
$z = \tau - \i 0^{+}$ for any $\tau$; 
$e^{2\pi\i x}F_{f}(z|\tau) \simeq (z-\tau)^{t}$ when 
$z = \tau + 0^{+}$ for any $\tau$ (in particular in zone (4)).

Then $F_{0}(z|\tau) = F_{a}(z|\tau)$, and 
$e^{-2\pi\i x}F_{0}(z-\tau|\tau) = F_{f}(z|\tau)$. 
We have $F_{b} = F_{a}(2\pi/\i)^{t}$,    $F_{f} = F_{e} (2\pi\i)^{-t}$.

Let us now compute the ratio between $F_{b}$ and $F_{c}$. 
Recall that $u = e^{2\pi\i z}$, $v = e^{2\pi\i(\tau-z)}$. 
Set $\bar F(u,v) := F(z|\tau)$. Using the expansion of $\theta(z|\tau)$, one shows that 
(\ref{ell:KZ}) has the form 
$$
\partial_{u}\bar F(u,v) = ({{A(u,v)}\over u} + {{B(u,v)}\over {u-1}})\bar F(u,v), 
$$
where $A(u,v)$ is holomorphic in the region $|v|<1/2$, $|u|<2$, and 
$A(u,0) = \tilde y$, $B(u,0) = t$. We have $\bar F_{b}(u,v)
= (1-u)^{t}(1 + \sum_{k} \sum_{s\leq s(k)}\on{log}(1-u)^{k} b_{ks}(u,v))$
and $\bar F_{b}(u,v)
= u^{\tilde t}(1 + \sum_{k} \sum_{s\leq s(k)}\on{log}(u)^{k} a_{ks}(u,v))$, 
with $a_{ks},b_{ks}$ analytic, and $a_{ks}(0,v) = b_{ks}(1,v)=0$. The 
ratio $\bar F_{b}^{-1}\bar F_{c}$ is an analytic function of $q$ only, 
which coincides with $\Phi(\tilde y,t)$ for $q=0$, so it has the form 
$\Phi(\tilde y,t) + \sum_{k>0} a_{k}(q)$, where $a_{k}(q)$ has degree $k$, is 
analytic in the neighborhood of $q=0$ and vanishes at $q=0$. Therefore
$$
F_{c}(z|\tau) = F_{b}(z|\tau) \big( \Phi(\tilde y,t) + O(e^{2\pi\i\tau})\big). 
$$

In the same way, one proves that 
$$
F_{e}(z|\tau)=F_{d}\big(e^{-2\pi\i x}\Phi(-\tilde y-t,t)^{-1} 
+ O(e^{2\pi\i\tau})\big).
$$ 
Indeed, let us set $\bar G_{d}(u',v') := e^{2\pi\i x} F_{d}(\tau+z'|\tau)$,
$\bar G_{e}(u',v') := e^{2\pi\i x} F_{e}(\tau+z'|\tau)$, 
where $u' = e^{2\pi\i(\tau+z')}$, $v' = e^{-2\pi\i z'}$, then 
$\bar G_{d}(u',v') \simeq (v')^{-\tilde y - t}e^{2\pi\i x}$ as $(u',v')\to (0^{+},0^{+})$ 
and 
$\bar G_{e}(u',v') \simeq (1-v')^{t}$ as $v'\to 1^{-}$ for any $u'$, and 
both $\bar G_{d}$ and $\bar G_{e}$ are solutions of $\partial_{v'}\bar G(u',v') 
= [-(\tilde y + t)/v' + t/(v'-1) + O(u')]\bar G(v')$. 
Therefore $\bar G_{d} = \bar G_{e}[\Phi(-\tilde y - t, t)e^{2\pi\i x} + O(u')]$.

Combining these results, we get:

\begin{lemma} 
$$
B(\tau) \simeq(2\pi\i)^{t} \Phi(-\tilde y-t,t) e^{2\pi\i x} 
e^{2\i\pi\tau\tilde y}
\Phi(\tilde y,t)^{-1}(2\pi/\i)^{-t}, 
$$
in the sense that the left (equivalently, right) ratio of these quantities
has the form $1 + \sum_{k>0} (\on{degree}\ k)O(\tau^{n(k)} e^{2\pi\i\tau})$
for $n(k)\geq 0$. 
\end{lemma}

Recall that we have proved: 
$$
F(z|\tau) = F_{0}(z|\tau)\exp( - {{a_0}\over{2\pi\i}}
(\int_{\i}^{\tau} E_{2} + C)t)[F(\tau)],
$$ 
where $C$ is such that 
$\int_{\i}^\tau E_2 + C = \tau + O(e^{2\pi\i\tau})$.

Set $X(\tau) := \exp( - {{a_0}\over{2\pi\i}} (\int_{\i}^{\tau} E_{2} + C)t)
[F(\tau)]$.

When $\tau\to\i\infty$, 
$X(\tau) = \on{exp}(- {\tau\over{2\pi\i}}
([\Delta]+a_0 t))(1+\sum_{k>0}(\on{degree}\ k)O(\tau^{f(k)}e^{2\pi\i\tau}))$.

Then 
\begin{align*}
& \tilde B = F(z|\tau)^{-1}e^{2\pi\i x}F(z+\tau|\tau)
= X(\tau)^{-1} B(\tau) X(\tau)
\\ & 
= \on{Ad}\Big( (1+\sum_{k>0}(\on{degree}\ k)O(\tau^{f(k)}e^{2\pi\i\tau}))^{-1}
\on{exp}({\tau\over{2\pi\i}}([\Delta]+a_{0}t)) \Big)
\\ & 
\Big(
\big( (2\pi\i)^{t}\Phi(-\tilde y-t,t)e^{2\pi\i x}e^{2\pi\i\tau\tilde y}
\Phi(\tilde y,t)^{-1}(2\pi/\i)^{-t}\big)
\big(1 + \sum_{k>0} (\on{degree}\ k)O(\tau^{n(k)}e^{2\pi\i\tau})\big) \Big), 
\end{align*}
where $\on{Ad}(u)(x) = uxu^{-1}$.

$[\Delta]+a_{0}t$ commutes with $\tilde y$ and $t$; assume for a moment that 
$\on{Ad}(\on{exp}({\tau\over{2\pi\i}}([\Delta]+a_{0}t)))(e^{2\pi\i x}e^{2\pi\i\tau
\tilde y}) = e^{2\pi\i x}$ (Lemma \ref{adjoint} below), then 
\begin{align*}
& \on{Ad}(\on{exp}({\tau\over{2\pi\i}}([\Delta]+a_{0}t)))
\Big( (2\pi\i)^{t}\Phi(-\tilde y-t,t)e^{2\pi\i x}e^{2\pi\i\tau\tilde y}
\Phi(\tilde y,t)^{-1}(2\pi/\i)^{-t}\Big)
\\ & = (2\pi\i)^{t}\Phi(-\tilde y-t,t)e^{2\pi\i x}
\Phi(\tilde y,t)^{-1}(2\pi/\i)^{-t}.
\end{align*} 
On the other hand, 
$\on{Ad}(\on{exp}({\tau\over{2\pi\i}}([\Delta]+a_{0}t)))
(1 + \sum_{k>0}(\on{degree}\ k)O(\tau^{n(k)}e^{2\pi\i\tau}))$
has the form $1 + \sum_{k>0}(\on{degree}\ k)O(\tau^{n'(k)}e^{2\pi\i\tau})$, 
where $n'(k)\geq 0$. It follows that 
\begin{align*}
\tilde B = & \on{Ad}\big( 1 + \sum_{k>0}(\on{degree}\ k)
O(\tau^{f(k)}e^{2\pi\i\tau})\big)
\\ & \Big( \big((2\pi\i)^{t}\Phi(-\tilde y-t,t)e^{2\pi\i x}
\Phi(\tilde y,t)^{-1}(2\pi/\i)^{-t} \big)
\big(1 + \sum_{k>0}(\on{degree}\ k)O(\tau^{n'(k)}e^{2\pi\i\tau})\big)\Big); 
\end{align*}
now 
\begin{align*}
& \on{Ad} \Big((2\pi\i)^{t}\Phi(-\tilde y-t,t)e^{2\pi\i x}
\Phi(\tilde y,t)^{-1}(2\pi/\i)^{-t} \Big)^{-1}
(1 + \sum_{k>0}(\on{degree}\ k)O(\tau^{f(k)}e^{2\pi\i\tau})) \\
  & = 
1 + \sum_{k>0}(\on{degree}\ k)O(\tau^{f(k)}e^{2\pi\i\tau}), 
\end{align*}
so 
\begin{align*}
\tilde B = & \Big((2\pi\i)^{t}\Phi(-\tilde y-t,t)e^{2\pi\i x}
\Phi(\tilde y,t)^{-1}(2\pi/\i)^{-t} \Big)
(1 + \sum_{k>0}(\on{degree}\ k)O(\tau^{f(k)}e^{2\pi\i\tau}))
\\ & (1 + \sum_{k>0}(\on{degree}\ k)O(\tau^{n'(k)}e^{2\pi\i\tau}))
\\ & 
= \Big((2\pi\i)^{t}\Phi(-\tilde y-t,t)e^{2\pi\i x}
\Phi(\tilde y,t)^{-1}(2\pi/\i)^{-t} \Big)
(1 + \sum_{k>0}(\on{degree}\ k)O(\tau^{n''(k)}e^{2\pi\i\tau}))
\end{align*}
for $n''(k)\geq 0$. Since $\tilde B$ is constant w.r.t. $\tau$, this implies 
$$
\tilde B  = 
(2\pi\i)^{t}\Phi(-\tilde y-t,t)e^{2\pi\i x}
\Phi(\tilde y,t)^{-1}(2\pi/\i)^{-t},
$$
as claimed.

We now prove the conjugation used above.

\begin{lemma} \label{adjoint}
For any $\tau\in\CC$, we have 
$$
e^{ {\tau\over{2\pi\i}}([\Delta]+a_0 t)} 
e^{2\pi\i x}e^{- {\tau\over{2\pi\i}}([\Delta]+a_0 t)} 
e^{2\i\pi\tau\tilde y} = e^{2\pi\i x}. 
$$
\end{lemma}

{\em Proof.} We have $[\Delta] + a_0 t = \Delta_0 
+ \sum_{k\geq 0} a_{2k}(\delta_{2k}
+ (\on{ad}x)^{2k}(t))$ (where $\delta_0=0$), so 
$[[\Delta]+a_0t,x] = y - \sum_{k\geq 0}
a_{2k} (\on{ad}x)^{2k+1}(t)$. Recall that 
$$
\sum_{k\geq 0} a_{2k}u^{2k} = {{\pi^2}\over{\on{sin}^2(\pi u)}} 
- {1\over{u^2}}, 
$$
then $[[\Delta]+a_0t,x] 
= y - (\on{ad}x) ({{\pi^2}\over{\on{sin}^2(\pi \on{ad}x)}} 
- {1\over{(\on{ad}x)^2}})(t)$. 
So 
\begin{align*}
& e^{-2\pi\i x} ({1\over{2\pi\i}}([\Delta]+a_0 t)) 
e^{2\pi\i x}  = 
{1\over{2\pi\i}}([\Delta]+a_0 t) 
+ {{e^{-2\pi\i \on{ad}x}-1}\over{\on{ad}x}} 
([x,{1\over{2\pi\i}}([\Delta]+a_0 t)])
\\ & 
= {1\over{2\pi\i}}([\Delta]+a_0 t) 
-{1\over{2\pi\i}} {{e^{-2\pi\i \on{ad}x}-1}\over{\on{ad}x}} 
\big( y - (\on{ad}x) ({{\pi^2}\over{\on{sin}^2(\pi \on{ad}x)}} 
- {1\over{(\on{ad}x)^2}})(t)\big). 
\end{align*}
We have 
$$
-{1\over{2\pi\i}} {{e^{-2\pi\i \on{ad}x}-1}\over{\on{ad}x}} 
\big(y - (\on{ad}x) ({{\pi^2}\over{\on{sin}^2(\pi \on{ad}x)}} 
- {1\over{(\on{ad}x)^2}})(t)\big) = -2\pi\i \tilde y, 
$$
therefore we get 
$$
e^{-2\pi\i x} ({1\over{2\pi\i}}([\Delta]+a_0 t)) 
e^{2\pi\i x} = 
{1\over{2\pi\i}}([\Delta]+a_0 t) - 2\pi\i\tilde y. 
$$
Multiplying by $\tau$, taking the exponential, and using the 
fact that $[\Delta]+a_0t$ commutes with $\tilde y$, we get 
$e^{-2\pi\i x} e^{ {\tau\over{2\pi\i}}([\Delta]+a_0 t)} 
e^{2\pi\i x} = 
e^{ {\tau\over{2\pi\i}}([\Delta]+a_0 t) } e^{- 2\pi\i\tau\tilde y}$, 
which proves the lemma. 
\hfill \qed \medskip

This ends the proof of Theorem \ref{thm:B}. \hfill \qed \medskip

\section{Construction of morphisms $\Gamma_{1,[n]}\to 
{\bold G}_{n}\rtimes S_n$} \label{sect:5bis}

In this section, we fix a field $\kk$ of characteristic zero. 
We denote the algebras $\bar\t_{1,n}^\kk$, $\t_n^\kk$
simply by  $\bar\t_{1,n}$, $\t_n$. The above group 
${\bold G}_n$ is the set of $\CC$-points of a group scheme
defined over $\QQ$, and we now again denote by ${\bold G}_n$
the set of its $\kk$-points.

\subsection{Construction of morphisms $\Gamma_{1,[n]}\to 
{\bold G}_{n}\rtimes S_n$
from a $5$-uple $(\Phi_\lambda,\tilde A,\tilde B,\tilde\Theta,\tilde\Psi)$}

Let $\Phi_\lambda$ be a $\lambda$-associator defined over $\kk$. 
This means that $\Phi_{\lambda}\in \on{exp}(\hat\t_{3})$ 
(the Lie algebras are now over $\kk$), 
\begin{equation} \label{assoc:1}
\Phi_{\lambda}^{3,2,1} = \Phi_{\lambda}^{-1}, \quad 
\Phi_{\lambda}^{2,3,4}\Phi_{\lambda}^{1,23,4}\Phi_{\lambda}^{1,2,3}
= \Phi_{\lambda}^{1,2,34}\Phi_{\lambda}^{12,3,4}, 
\end{equation}
\begin{equation} \label{assoc:2}
e^{\lambda t_{31}/2}
\Phi_\lambda^{2,3,1}
e^{\lambda t_{23}/2}\Phi_\lambda e^{\lambda t_{12}/2}\Phi_\lambda^{3,1,2}
= e^{\lambda(t_{12}+t_{23}+t_{13})/2}. 
\end{equation}
E.g., the KZ associator is a $2\pi\i$-associator over $\CC$.

\begin{proposition} \label{construction}
If  $\tilde\Theta,\tilde\Psi\in {\bold G}_{1}$
and $\tilde A,\tilde B\in \on{exp}(\hat{\bar\t}_{1,2})$
satisfy: the ``$\Gamma_{1,1}$ identities'' (\ref{gamma11}), 
the ``$\Gamma_{1,2}$ identities'' (\ref{gamma12}), (\ref{gamma12'}), 
and the ``$\Gamma_{1,[3]}$ identities'' 
(\ref{B:ids}), (\ref{A:ids}), (\ref{gamma13'}) (with $2\pi\i$ 
replaced by $\lambda$), as well as $\tilde A^{\emptyset,1} = \tilde
A^{1,\emptyset} = \tilde B^{\emptyset,1} = \tilde
B^{1,\emptyset}=1$, then one defines a morphism 
$\Gamma_{1,[n]}\to {\bold G}_{n}\rtimes S_{n}$ by 
$$ 
\Theta\mapsto [\tilde\Theta]e^{\i{\pi\over 2}\sum_{i<j}\bar t_{ij}}, 
\quad 
\Psi\mapsto [\tilde\Psi]e^{\i{\pi\over 6}\sum_{i<j}\bar t_{ij}}, \quad
\sigma_{i}\mapsto \{\Phi_{\lambda}^{1...i-1,i,i+1}\}^{-1} 
e^{\lambda\bar t_{i,i+1}/2}
(i,i+1) \{\Phi_{\lambda}^{1...i-1,i,i+1}\}, 
$$
$$
C_{jk} \mapsto \{\Phi_{\lambda,j}^{-1}\Phi_\lambda^{j,j+1,...n}...
\Phi_\lambda^{j...,k-1,...n} (e^{\lambda t_{12}})^{j...k-1,k...n}
(\Phi_\lambda^{j,j+1,...n}...\Phi_\lambda^{j...,k-1,...n})^{-1}
\Phi_{\lambda,j}\}, 
$$
$$
A_{i}\mapsto \{\Phi_{\lambda,i}\}^{-1}\tilde A^{1...i-1,i...n}
\{\Phi_{\lambda,i}\}, \quad 
B_{i}\mapsto \{\Phi_{\lambda,i}\}^{-1}\tilde B^{1...i-1,i...n}
\{\Phi_{\lambda,i}\}, 
$$
where $\Phi_{\lambda,i} =
\Phi_\lambda^{1...i-1,i,i+1...n}...\Phi_\lambda^{1...n-2,n-1,n}$. 
\end{proposition}

According to Section \ref{5:4}, the representations $\gamma_{n}$ are 
obtained by the procedure described in this proposition from 
the KZ associator, $\tilde\Theta,\tilde\Psi$ arising from $\gamma_{1}$,
and $\tilde A,\tilde B$ arising from $\gamma_{2}$.

Note also that the analogue of (\ref{A:ids}) is equivalent to the pair of
equations 
$$
e^{\lambda \bar t_{12}/2}\tilde A^{2,1}
e^{\lambda \bar t_{12}/2}\tilde A = 1, \quad 
(e^{\lambda \bar t_{12}/2}\tilde A)^{3,12}\Phi_\lambda^{3,1,2}
(e^{\lambda \bar t_{12}/2}\tilde A)^{2,31}\Phi_\lambda^{2,3,1}
(e^{\lambda \bar t_{12}/2}\tilde A)^{1,23}\Phi_\lambda^{1,2,3}=1, 
$$
and similarly (\ref{B:ids}) is equivalent to the same equations, 
with $\tilde A,\lambda$ replaced by $\tilde B,-\lambda$.

\begin{remark}
One can prove that it $\Phi_{\lambda}$ satisfies only the pentagon equation and 
$\tilde\Theta,\tilde\Psi,\tilde A,\tilde B$ satisfy the 
the ``$\Gamma_{1,1}$
identities'' (\ref{gamma11}), the ``$\Gamma_{1,2}$ identities'' 
(\ref{gamma12}), (\ref{gamma12'}), 
and the ``$\Gamma_{1,3}$ identities'' (\ref{gamma13}), (\ref{gamma13'}),
then the above formulas (removing $\sigma_i$) 
define a morphism  $\Gamma_{1,n}\to {\bold G}_{n}$. In the same way, if 
$\Phi_\lambda$ satisfies all the associator conditions and 
$\tilde A,\tilde B$ satisfy the $\Gamma_{1,[3]}$ identities (\ref{A:ids}), 
(\ref{B:ids}), (\ref{gamma13'}), then 
the above formulas (removing $\Theta,\Psi$) define a morphism 
$\overline{\on{B}}_{1,n} \to \on{exp}(\hat{\bar\t}_{1,n})\rtimes S_n$. 
\end{remark}

\medskip 
{\em Proof.} Let us prove that the identity $(A_{i},A_{j})=1$ ($i<j$) 
is preserved. 
Applying $x\mapsto x^{1...i-1,i...j-1,j...n}$ to the first identity 
of (\ref{gamma13}), we get 
$$
(\tilde A^{1...i-1,i...n},\Phi_\lambda^{1...,i...j-1,...n}
\tilde A^{1...j-1,j...n} (\Phi_\lambda^{-1})^{1...,i...j-1,...n})=1. 
$$
The pentagon identity implies 
\begin{equation} \label{gen:pentagon}
\Phi_\lambda^{1...,i,...n}...\Phi_\lambda^{1...,j-1,...n} = 
(\Phi_\lambda^{i,i+1,...n}...\Phi_\lambda^{i...,j-1,..,n}) 
\Phi_\lambda^{1...,i...j-1,...n} (\Phi_\lambda^{1...,i,...j-1}...
\Phi_\lambda^{1...,j-2,j-1}), 
\end{equation}
so the above identity is rewritten 
\begin{align*}
& \big(\Phi_\lambda^{i,i+1,...n}...\Phi_\lambda^{i...,j-1,..,n}
\tilde A^{1...i-1,i...n}
(\Phi_\lambda^{i,i+1,...n}...\Phi_\lambda^{i...,j-1,..,n})^{-1},
\Phi_\lambda^{1...,i,...n}...\Phi_\lambda^{1...,j-1,...n} \\ & 
(\Phi_\lambda^{1...,i,...j-1}...\Phi_\lambda^{1...,j-2,...j-1})^{-1}
\tilde A^{1...j-1,j...n}
\Phi_\lambda^{1...,i,...j-1}...\Phi_\lambda^{1...,j-2,...j-1}
(\Phi_\lambda^{1...,i,...n}...\Phi_\lambda^{1...,j-1,...n})^{-1}\big) = 1. 
\end{align*}
Now $\Phi_\lambda^{i,i+1,...n},...,\Phi_\lambda^{i...,j-1,..,n}$
commute with $\tilde A^{1...i-1,i...n}$, and 
$\Phi_\lambda^{1...,i,...j-1},...,\Phi_\lambda^{1...,j-2,...j-1}$
commute with $\Phi_\lambda^{1...,i,...j-1}...\Phi_\lambda^{1...,j-2,...j-1}$, 
which implies 
$$
(\tilde A^{1...i-1,i...n},
  \Phi_\lambda^{1...,i,...n}...\Phi_\lambda^{1...,j-1,...n}
\tilde A^{1...j-1,j...n}
(\Phi_\lambda^{1...,i,...n}...\Phi_\lambda^{1...,j-1,...n})^{-1}) = 1, 
$$
so that $(A_i,A_j)=1$ is preserved. In the same way, one shows that 
$(B_i,B_j)=1$ is preserved.

Let us show that $(B_k,A_k A_j^{-1}) = C_{jk}$ is preserved
(if $j\leq k$). 
\begin{align*}
& (\Phi_{\lambda,k}^{-1}\tilde B^{1...k-1,k...n}\Phi_{\lambda,k},
\Phi_{\lambda,k}^{-1}\tilde A^{1...k-1,k...n}\Phi_{\lambda,k}
\Phi_{\lambda,j}^{-1}(\tilde A^{1...j-1,j...n})^{-1}\Phi_{\lambda,j})
\\ & = 
\Phi_{\lambda,j}^{-1}\big((\Phi_\lambda^{1...,j,...n}...
\Phi_\lambda^{1...,k-1,...n})
\tilde B^{1...k-1,k...n}(\Phi_\lambda^{1...,j,...n}...
\Phi_\lambda^{1...,k-1,...n})^{-1},
\\ & (\Phi_\lambda^{1...,j,...n}...\Phi_\lambda^{1...,k-1,...n})
\tilde A^{1...k-1,k...n}
(\Phi_\lambda^{1...,j,...n}...\Phi_\lambda^{1...,k-1,...n})^{-1}
(\tilde A^{1...j-1,j...n})^{-1}\big)\Phi_{\lambda,j} 
\\ & = 
\Phi_{\lambda,j}^{-1}\big(\Phi_\lambda^{j,j+1,...n}...
\Phi_\lambda^{j...,k-1,...n}\Phi_\lambda^{1...,j...k-1,...n}
\tilde B^{1...k-1,k...n}
(\Phi_\lambda^{j,j+1,...n}...\Phi_\lambda^{j...,k-1,...n}
\Phi_\lambda^{1...,j...k-1,...n})^{-1},
\\ & 
\Phi_\lambda^{j,j+1,...n}...\Phi_\lambda^{j...,k-1,...n}
\Phi_\lambda^{1...,j...k-1,...n}
\tilde A^{1...k-1,k...n}
(\Phi_\lambda^{j,j+1,...n}...\Phi_\lambda^{j...,k-1,...n}
\Phi_\lambda^{1...,j...k-1,...n})^{-1}
\\ & (\tilde A^{1...j-1,j...n})^{-1}\big)\Phi_{\lambda,j} 
= \Phi_{\lambda,j}^{-1}\Phi_\lambda^{j,j+1,...n}...\Phi_\lambda^{j...,k-1,...n}
\big(\Phi_\lambda^{1...,j...k-1,...n}
\tilde B^{1...k-1,k...n}
(\Phi_\lambda^{1...,j...k-1,...n})^{-1}, \\ & 
\Phi_\lambda^{1...,j...k-1,...n}
\tilde A^{1...k-1,k...n}
(\Phi_\lambda^{1...,j...k-1,...n})^{-1}
(\tilde A^{1...j-1,j...n})^{-1}\big)
(\Phi_\lambda^{j,j+1,...n}...\Phi_\lambda^{j...,k-1,...n})^{-1}\Phi_{\lambda,j} 
\\ & 
= \Phi_{\lambda,j}^{-1}\Phi_\lambda^{j,j+1,...n}...\Phi_\lambda^{j...,k-1,...n}
\{\Phi (\tilde B^{12,3},\tilde A^{12,3}\Phi_\lambda^{-1}
(\tilde A^{1,23})^{-1} \Phi_\lambda)\Phi_\lambda^{-1}\}^{1...,j...k-1,...n} 
\\ & (\Phi_\lambda^{j,j+1,...n}...\Phi_\lambda^{j...,k-1,...n})^{-1}
\Phi_{\lambda,j}
\\ & = \Phi_{\lambda,j}^{-1}\Phi_\lambda^{j,j+1,...n}...
\Phi_\lambda^{j...,k-1,...n}
(e^{2\pi\i\bar t_{12}})^{j...k-1,k...n}
(\Phi_\lambda^{j,j+1,...n}...\Phi_\lambda^{j...,k-1,...n})^{-1}\Phi_{\lambda,j}, 
\end{align*}
where the second identity uses (\ref{gen:pentagon}) and the invariance of 
$\Phi_\lambda$, the third identity uses the fact that 
$\Phi_\lambda^{j,j+1,...n},...,\Phi_\lambda^{j...,k-1,...n}$ commute with 
$\tilde A^{1...j-1,j..n}$ (again by the invariance of $\Phi_\lambda$), 
and the last
identity uses (\ref{gamma13'}). So $(B_k,A_k A_j^{-1}) = C_{jk}$ is preserved. 
One shows similarly that 
\begin{align*}
& \big(\Phi_{\lambda,k}^{-1}\tilde B^{1...k-1,k...n}\Phi_{\lambda,k}
\Phi_{\lambda,j}^{-1}(\tilde B^{1...j-1,j...n})^{-1}\Phi_{\lambda,j}, 
\Phi_{\lambda,k}^{-1}\tilde A^{1...k-1,k...n}\Phi_{\lambda,k}\big)
\\ & = \Phi_j^{-1}\Phi^{j,j+1,...n}...\Phi^{j...,k-1,...n}
(e^{2\pi\i\bar t_{12}})^{j...k-1,k...n}
(\Phi_\lambda^{j,j+1,...n}...\Phi_\lambda^{j...,k-1,...n})^{-1}\Phi_{\lambda,j}, 
\end{align*}
so that $(B_kB_j^{-1},A_k) = C_{jk}$ is preserved.

Let us show that $(A_i,C_{jk})=1$ ($i\leq j\leq k$) is preserved. 
We have 
\begin{align*}
& \big(\Phi_{\lambda,i}^{-1}\tilde A^{1...i-1,i...n}\Phi_{\lambda,i},
\Phi_{\lambda,j}^{-1}\Phi_\lambda^{j,j+1,...n}...\Phi_\lambda^{j...,k-1,...n}
(e^{2\pi\i\bar t_{12}})^{j...k-1,k...n}
(\Phi_\lambda^{j,j+1,...n}...\Phi_\lambda^{j...,k-1,...n})^{-1}
\Phi_{\lambda,j}\big) 
\\ & 
= \Phi_{\lambda,i}^{-1}\big(\tilde A^{1...i-1,i...n},
\Phi_\lambda^{1...,i,...n}...\Phi_\lambda^{1...,j-1,...n}
\Phi_\lambda^{j,j+1,...n}...\Phi_\lambda^{j...,k-1,...n}
(e^{2\pi\i\bar t_{12}})^{j...k-1,k...n}
\\ & (\Phi_\lambda^{1...,i,...n}...\Phi_\lambda^{1...,j-1,...n}
\Phi_\lambda^{j,j+1,...n}...\Phi_\lambda^{j...,k-1,...n})^{-1} \big)
\Phi_{\lambda,i}
\\ & 
= \Phi_{\lambda,i}^{-1}\big(\tilde A^{1...i-1,i...n},
\Phi_\lambda^{i,i+1,...n}...\Phi_\lambda^{i...,j-1,...n}
\Phi_\lambda^{1...,i...j-1,...n}\Phi_\lambda^{1...,i,...j-1}...
\Phi_\lambda^{1...,j-2,j-1}
\\ & \Phi_\lambda^{j,j+1,...n}...\Phi_\lambda^{j...,k-1,...n}
(e^{2\pi\i\bar t_{12}})^{j...k-1,k...n}
\\ & 
(\Phi_\lambda^{i,i+1,...n}...\Phi_\lambda^{i...,j-1,...n}
\Phi_\lambda^{1...,i...j-1,...n}\Phi_\lambda^{1...,i,...j-1}...
\Phi_\lambda^{1...,j-2,j-1}
\Phi_\lambda^{j,j+1,...n}...\Phi_\lambda^{j...,k-1,...n})^{-1}\big)
\Phi_{\lambda,i}
\\ & = \Phi_{\lambda,i}^{-1}\big(\tilde A^{1...i-1,i...n},
\Phi_\lambda^{i,i+1,...n}...\Phi_\lambda^{i...,j-1,...n}
\Phi_\lambda^{1...,i...j-1,...n}
\Phi_\lambda^{j,j+1,...n}...\Phi_\lambda^{j...,k-1,...n}
\\ & (e^{2\pi\i\bar t_{12}})^{j...k-1,k...n}
(\Phi_\lambda^{i,i+1,...n}...\Phi_\lambda^{i...,j-1,...n}
\Phi_\lambda^{1...,i...j-1,...n}
\Phi_\lambda^{j,j+1,...n}...\Phi_\lambda^{j...,k-1,...n})^{-1}\big)
\Phi_{\lambda,i}
\\ & = \Phi_{\lambda,i}^{-1}\Phi_\lambda^{i,i+1,...n}...
\Phi_\lambda^{i...,j-1,...n}
\big(\tilde A^{1...i-1,i...n},
\Phi_\lambda^{1...,i...j-1,...n}
\Phi_\lambda^{j,j+1,...n}...\Phi_\lambda^{j...,k-1,...n}
\\ & (e^{2\pi\i\bar t_{12}})^{j...k-1,k...n}
(\Phi_\lambda^{1...,i...j-1,...n}
\Phi_\lambda^{j,j+1,...n}...\Phi_\lambda^{j...,k-1,...n})^{-1}\big)
(\Phi_\lambda^{i,i+1,...n}...\Phi_\lambda^{i...,j-1,...n})^{-1}
\Phi_{\lambda,i}
\\ & 
= \Phi_{\lambda,i}^{-1}\Phi_\lambda^{i,i+1,...n}...\Phi_\lambda^{i...,j-1,...n}
(\tilde A^{1...i-1,i...n},
\Phi_\lambda^{j,j+1,...n}...\Phi_\lambda^{j...,k-1,...n}
\Phi_\lambda^{1...,i...j-1,...n}
\\ & (e^{2\pi\i\bar t_{12}})^{j...k-1,k...n}
(\Phi_\lambda^{j,j+1,...n}...\Phi_\lambda^{j...,k-1,...n}
\Phi_\lambda^{1...,i...j-1,...n})^{-1})
(\Phi_\lambda^{i,i+1,...n}...\Phi_\lambda^{i...,j-1,...n})^{-1}\Phi_{\lambda,i}
\\ & = 1, 
\end{align*}
where the second equality follows from the generalized pentagon identity
(\ref{gen:pentagon}), the third equality follows from the fact that 
$\Phi_\lambda^{1...,i,...j-1}$, ..., $\Phi_\lambda^{1...,j-2,j-1}$ commute with 
$(e^{2\pi\i\bar t_{12}})^{j...k-1,k...n}$, $\Phi_\lambda^{j,j+1,...n}$, 
..., $\Phi_\lambda^{j...,k-1,...n}$, the fourth equality follows 
from the fact that $\Phi_\lambda^{i,i+1,...n}$, ..., 
$\Phi_\lambda^{i...,j-1,...n}$ commute with 
$\tilde A^{1...i-1,i...n}$ (as $\Phi_\lambda$ is invariant), the last equality 
follows from the fact that $\Phi_\lambda^{1...,i...j-1,j...n}$ commutes with 
$\Phi_\lambda^{j,j+1,...n}$, ..., $\Phi_\lambda^{j...,k-1,...n}$ 
(again as $\Phi_\lambda$ is invariant) and with 
$(e^{2\pi\i\bar t_{12}})^{j...k-1,k...n}$
(as $t_{34}$ commutes with the image of $\t_{3}\to\t_{4}$, 
$x\mapsto x^{1,2,34}$). 
Therefore $(A_{i},C_{jk})=1$ is preserved. One shows similarly that 
$(B_{i},C_{jk})=1$ ($i\leq j\leq k$), 
$X_{i+1} = \sigma_i X_i \sigma_i$
and $Y_{i+1} = \sigma_i^{-1} Y_i \sigma_i^{-1}$ are preserved.

The fact that the relations 
$\Theta A_{i}\Theta^{-1}=B_{i}^{-1}$, $\Theta B_{i}\Theta^{-1}
=B_{i}A_{i}B_{i}^{-1}$, 
$\Psi A_{i}\Psi^{-1}=A_{i}$, $\Psi B_{i}\Psi^{-1}=B_{i}A_{i}$, are preserved 
follows from the identities (\ref{gamma12}), (\ref{gamma12'}) and that if we 
denote by $x\mapsto [x]_{n}$ the morphism $\d\to\d\rtimes\bar\t_{1,n}$ defined 
above, then: 
(a) $\Phi_{i}$ commutes with $\sum_{i,j|i<j}\bar t_{ij}$
and with the image of $\d\to\d\rtimes\bar\t_{1,n}$, $x\mapsto [x]_n$; 
(b) for $x\in \d$, $y\in \bar\t_{1,2}$, 
we have $[[x]_{n},y^{1...i-1,i...n}] = [[x]_{2},y]^{1...i-1,i...n}$. 
Let us prove (a): the first part follows from the fact that $\Phi$
commutes with $t_{12}+t_{13}+t_{23}$; the second part follows from the 
fact that $X,d,\Delta_{0}$ and 
$\delta_{2n}+\sum_{k<l}(\on{ad}\bar x_{k})^{2n}(\bar t_{kl})$
commute with $\bar t_{ij}$ for any $i<j$. Let us prove (b): the identity 
holds for $[x,x']$ 
whenever it holds for $x$ and for $x'$, so it suffices to check it for 
$x$ a generator of $\d$; 
$x$ being such a generator, both sides are (as functions of $y$) 
derivations $\bar\t_{1,2}\to\bar\t_{1,n}$ w.r.t. the morphism 
$\bar\t_{1,2}\to\bar\t_{1,n}$, 
$y\mapsto y^{1...i-1,i...n}$, so it suffices to check the identity for 
$y$ a generator of $\bar\t_{1,2}$. 
The identity is obvious if $x\in\{\Delta_{0},d,X\}$ and $y\in \{
\bar x_{1},\bar y_{1}, \bar x_{2},\bar y_{2}\}$. If $x = \delta_{2s}$ 
and $y=\bar x_{1}$, then the identity holds because we have 
\begin{align*}
& [\delta_{2s} + (\on{ad}\bar x_1)^{2s}(\bar t_{12}),\bar x_1]^{1...i-1,i...n}
= - \big( (\on{ad}\bar x_1)^{2s+1}(\bar t_{12})\big)^{1...i-1,i...n}
\\ & = - (\on{ad}(\sum_{u'=1}^{i-1}\bar x_{u'}))^{2s+1}
(\sum_{1\leq u<i\leq v\leq n}\bar t_{uv})
= -\sum_{1\leq u<i\leq v\leq n} (\on{ad}\bar x_u)^{2s+1}(\bar t_{uv}), 
\end{align*}
while 
\begin{align*}
& [\delta_{2s}+\sum_{1\leq u<v\leq n}(\on{ad}\bar x_u)^{2s}(\bar t_{uv}),
\sum_{u'=1}^{i-1}\bar x_{u'}] = 
[\sum_{1\leq u<i\leq v\leq n}(\on{ad}\bar x_u)^{2s}(\bar t_{uv}),
\sum_{u'=1}^{i-1}\bar x_{u'}] 
\\ & = 
-\sum_{1\leq u<i\leq v\leq n}(\on{ad}\bar x_u)^{2s+1}(\bar t_{uv}) 
\end{align*}
where the first equality follows from the fact that 
$(\on{ad}\bar x_u)^{2s}(\bar t_{uv})$ commutes with 
$\sum_{u'=1}^{i-1}\bar x_{u'}$ whenever $u<v<i$ or $i\leq u<v$. 
If $x = \delta_{2s}$ 
and $y=\bar x_{2}$, then the identity follows because 
$[\delta_{2s} + (\on{ad}\bar x_1)^{2s}(\bar t_{12}),\bar x_1+\bar x_2]=0$
and $[\delta_{2s}+\sum_{1\leq u<v\leq n}(\on{ad}\bar x_u)^{2s}(\bar t_{uv}),
\sum_{u'=1}^{n}\bar x_{u'}] = 0$.

If $x = \delta_{2s}$ and $y=\bar y_{1}$, then 
\begin{align*}
& [\delta_{2s}+(\on{ad}\bar x_1)^{2s}(\bar t_{12}),\bar y_1]^{1...i-1,i...n}
\\ & 
= \{{1\over 2}\sum_{p+q=2s-1}[(\on{ad}\bar x_1)^p(\bar t_{12}),
(-\on{ad}\bar x_1)^q(\bar t_{12})] 
+ [(\on{ad}\bar x_1)^{2s}(\bar t_{12}),\bar y_1]\}^{1...i-1,i...n}
\\ & 
= {1\over 2} \sum_{p+q=2s-1} 
[\sum_{1\leq u<i\leq v\leq n}(\on{ad}\bar x_u)^p(\bar t_{uv}), 
\sum_{1\leq u'<i\leq v'\leq n}(\on{ad}\bar x_{u'})^q(\bar t_{u'v'})] 
\\ & + [\sum_{1\leq u<i\leq v\leq n}(\on{ad}\bar x_u)^{2s}(\bar t_{uv}),
\bar y_1+...+\bar y_{i-1}]; 
\end{align*}
on the other hand, 
\begin{align*}
& [\delta_{2s} + \sum_{1\leq u<v\leq n}(\on{ad}\bar x_u)^{2s}(\bar t_{uv}),
\bar y_1+...+\bar y_{i-1}] 
\\ & 
= \sum_{1\leq u<v\leq n} [(\on{ad}\bar x_u)^{2s}(\bar t_{uv}),
\bar y_1 + ... + \bar y_{i-1}] 
+ \sum_{u=1}^{i-1}\sum_{v|v\neq u}\sum_{p+q=2s-1}{1\over 2}[(\on{ad}\bar x_u)^p
(\bar t_{uv}),(-\on{ad}\bar x_u)^q(\bar t_{uv})]
\\ & 
= \sum_{1\leq u<v\leq n} [(\on{ad}\bar x_u)^{2s}(\bar t_{uv}),
\bar y_1 + ... + \bar y_{i-1}] 
+ \sum_{1\leq u<i\leq v\leq n}\sum_{p+q=2s-1}{1\over 2}[(\on{ad}\bar x_u)^p
(\bar t_{uv}),(-\on{ad}\bar x_u)^q(\bar t_{uv})], 
\end{align*}
where the second equality follows from the fact that 
$[(\on{ad}\bar x_u)^p(\bar t_{uv}),(-\bar x_u)^q(\bar t_{uv})] 
+ [(\on{ad}\bar x_v)^p(\bar t_{uv}),$ $(-\on{ad}\bar x_v)^q(\bar t_{uv})] =0$
as $p+q$ is odd.

Then 
\begin{align*}
& [\delta_{2s}+(\on{ad}\bar x_1)^{2s}(\bar t_{12}),\bar y_1]^{1...i-1,i...n}
- [\delta_{2s} + \sum_{1\leq u<v\leq n}(\on{ad}\bar x_u)^{2s}(\bar t_{uv}),
\bar y_1+...+\bar y_{i-1}] 
\\ & 
= - \sum_{1\leq u<v<i} [(\on{ad}\bar x_u)^{2s}(\bar t_{uv}),
\bar y_1 + ... + \bar y_{i-1}]
- \sum_{i\leq u<v \leq n} [(\on{ad}\bar x_u)^{2s}(\bar t_{uv}),
\bar y_1 + ... + \bar y_{i-1}]
\\ &
 + {1\over 2} \sum_{p+q=2s-1}\sum_{\matrix 
\scriptstyle{1\leq u<i\leq v\leq n} \\ 
\scriptstyle{1\leq u'<i\leq v'\leq n},  \scriptstyle{(u,v)\neq (u',v')}
\endmatrix} [(\on{ad}\bar x_u)^p(\bar t_{uv}),
(-\on{ad}\bar x_{u'})^q(\bar t_{u'v'})]
\\ & 
= \sum_{1\leq u<v<i} [(\on{ad}\bar x_u)^{2s}(\bar t_{uv}),
\bar y_i + ... + \bar y_n]
- \sum_{i\leq u<v \leq n} [(\on{ad}\bar x_u)^{2s}(\bar t_{uv}),
\bar y_1 + ... + \bar y_{i-1}]
\\ & 
+ {1\over 2} \sum_{p+q=2s-1}\sum_{\matrix 
\scriptstyle{1\leq u<i\leq v\leq n} \\ 
\scriptstyle{1\leq u<i\leq v'\leq n},  \scriptstyle{v\neq v'}
\endmatrix} [(\on{ad}\bar x_u)^p(\bar t_{uv}),
(-\on{ad}\bar x_u)^q(\bar t_{uv'})]
\\ & 
+ {1\over 2} \sum_{p+q=2s-1}\sum_{\matrix 
\scriptstyle{1\leq u<i\leq v\leq n} \\ 
\scriptstyle{1\leq u'<i\leq v\leq n},  \scriptstyle{u\neq u'}
\endmatrix} [(\on{ad}\bar x_u)^p(\bar t_{uv}),
(-\on{ad}\bar x_{u'})^q(\bar t_{u'v})]
\end{align*}
where the second equality follows from the centrality of $\bar y_1+...+\bar
y_n$, the last equality follows for the fact that $(\on{ad}\bar x_u)^p
(\bar t_{uv})$ and $(-\on{ad}\bar x_{u'})^q(\bar t_{u'v'})$
commute for $u,v,u',v'$ all distinct. Since $p+q$ is odd, it follows that 
\begin{align*}
& [\delta_{2s}+(\on{ad}\bar x_1)^{2s}(\bar t_{12}),\bar y_1]^{1...i-1,i...n}
- [\delta_{2s} + \sum_{1\leq u<v\leq n}(\on{ad}\bar x_u)^{2s}(\bar t_{uv}),
\bar y_1+...+\bar y_{i-1}] 
\\ & 
= \sum_{1\leq u<v<i} [(\on{ad}\bar x_u)^{2s}(\bar t_{uv}),
\bar y_i + ... + \bar y_n]
- \sum_{i\leq u<v \leq n} [(\on{ad}\bar x_u)^{2s}(\bar t_{uv}),
\bar y_1 + ... + \bar y_{i-1}]
\\ & 
+ \sum_{p+q=2s-1}\sum_{1\leq u<i\leq v<v'\leq n}
[(\on{ad}\bar x_u)^p(\bar t_{uv}),
(-\on{ad}\bar x_u)^q(\bar t_{uv'})]
\\ & 
+  \sum_{p+q=2s-1}\sum_{1\leq u<u'<i\leq v\leq n} 
[(\on{ad}\bar x_u)^p(\bar t_{uv}),
(-\on{ad}\bar x_{u'})^q(\bar t_{u'v})]. 
\end{align*}
Now if $1\leq u<v<i$, we have 
\begin{align*}
& [(\on{ad}\bar x_u)^{2s}(\bar t_{uv}),
\bar y_i + ... + \bar y_n]
  = \sum_{p+q = 2s-1} (\on{ad}\bar x_u)^p 
\on{ad}(\bar t_{ui}+...+\bar t_{un})
(\on{ad}\bar x_u)^q(\bar t_{uv}) 
\\ & = 
\sum_{w=i}^n \sum_{p+q = 2s-1} (\on{ad}\bar x_u)^p 
[\bar t_{uw},(-\on{ad}\bar x_v)^q(\bar t_{uv})]
= 
\sum_{w=i}^n  \sum_{p+q = 2s-1} 
(\on{ad}\bar x_u)^p (-\on{ad}\bar x_v)^q
([\bar t_{uw},\bar t_{uv}])
\\ & 
= 
-\sum_{w=i}^n  \sum_{p+q = 2s-1} 
(\on{ad}\bar x_u)^p (-\on{ad}\bar x_v)^q
([\bar t_{uw},\bar t_{vw}])
= 
-\sum_{w=i}^n  \sum_{p+q = 2s-1} 
[(\on{ad}\bar x_u)^p (\bar t_{uw}),
(-\on{ad}\bar x_v)^q(\bar t_{vw})]; 
\end{align*}
one shows in the same way that 
if $i\leq u<v\leq n$, then 
$[(\on{ad}\bar x_u)^{2s}(\bar t_{uv}),
\bar y_1 + ... + \bar y_{i-1}]
= \sum_{w=1}^{i-1}\sum_{p+q=2s-1}
[(\on{ad}\bar x_u)^p (\bar t_{uw}),
(-\on{ad}\bar x_v)^q(\bar t_{vw})]$; all this implies that 
$$
[\delta_{2s}+(\on{ad}\bar x_1)^{2s}(\bar t_{12}),\bar y_1]^{1...i-1,i...n}
- [\delta_{2s} + \sum_{1\leq u<v\leq n}(\on{ad}\bar x_u)^{2s}(\bar t_{uv}),
(\bar y_1)^{1...i-1}]. 
$$
Since $[\delta_{2s}+(\on{ad}\bar x_1)^{2s}(\bar t_{12}),
\bar y_1 + \bar y_2]=0$ and $[\delta_{2s} + 
\sum_{1\leq u<v\leq n}(\on{ad}\bar x_u)^{2s}(\bar t_{uv}),
\bar y_1 + ... + \bar y_n]=0$, this equality implies 
$$
[\delta_{2s}+(\on{ad}\bar x_1)^{2s}(\bar t_{12}),\bar y_2]^{1...i-1,i...n}
- [\delta_{2s} + \sum_{1\leq u<v\leq n}(\on{ad}\bar x_u)^{2s}(\bar t_{uv}),
(\bar y_2)^{1...i-1}], 
$$
which ends the proof of (b) above, and therefore of the fact that
the identities $\Theta A_i\Theta^{-1} = B_i^{-1}$, ..., 
$\Psi B_i\Psi^{-1} = B_i A_i$ are preserved.

The relation $(\Theta,\Psi^{2})=1$ is preserved because 
$$([\tilde\Theta]e^{\i{\pi\over 2}\sum_{i<j}\bar t_{ij}},
([\tilde\Psi]e^{\i{\pi\over 6}\sum_{i<j}\bar t_{ij}})^{2}) = 
([\tilde\Theta]e^{\i{\pi\over 2}\sum_{i<j}\bar t_{ij}},
[\tilde\Psi]^{2}e^{\i{\pi\over 3}\sum_{i<j}\bar t_{ij}}) = ([\tilde\Theta],
[\tilde\Psi]^{2}) = [(\tilde\Theta,\tilde\Psi^{2})]=1,
$$ 
where the two first identities follow from the fact that $\sum_{i<j}\bar t_{ij}$
commutes with the image of $\d\to \d\rtimes\bar\t_{1,n}$, $x\mapsto [x]$, the 
third identity follows from the fact that ${\bold G}_{1} \to {\bold G}_{n}$, 
$g\mapsto [g]$
is a group morphism, and the last identity follows from (\ref{gamma11}).

The image of $C_{i,i+1}$ is 
$\Phi_{\lambda,i}^{-1}(e^{2\pi\i\bar t_{12}})^{i,i+1...n}\Phi_{\lambda,i}$, 
to the product of the images of $C_{12},..., C_{n-1,n}$ is 
\begin{align*}
& \Phi_{\lambda,1}^{-1} (e^{2\pi\i\bar t_{12}})^{1,2...n}
(\Phi_{\lambda,1}\Phi_{\lambda,2}^{-1}) (e^{2\pi\i\bar t_{12}})^{2,3...n}
(\Phi_{\lambda,2}\Phi_{\lambda,3}^{-1}) (e^{2\pi\i\bar t_{12}})^{3,4...n}
... (\Phi_{\lambda,n-1}\Phi_{\lambda,n}^{-1}) e^{2\pi\i\bar t_{n-1,n}}
\Phi_{\lambda,n} 
\\ & 
= 
\Phi_{\lambda,1}^{-1} (e^{2\pi\i\bar t_{12}})^{1,2...n}
(e^{2\pi\i\bar t_{12}})^{2,3...n}
\Phi_\lambda^{1,2,3...n} (e^{2\pi\i\bar t_{12}})^{3,4...n} ... 
\Phi_\lambda^{1...,i-1,...n} (e^{2\pi\i\bar t_{12}})^{i,i+1...n}
\\ & ... 
\Phi_\lambda^{1...,n-2,n-1\ n} e^{2\pi\i\bar t_{n-1,n}} 
\\ & 
= \Phi_{\lambda,1}^{-1}
(e^{2\pi\i\bar t_{12}})^{1,2...n}
(e^{2\pi\i\bar t_{12}})^{2,3...n}
(e^{2\pi\i\bar t_{12}})^{3,4...n}
... (e^{2\pi\i\bar t_{12}})^{i,i+1...n}
...  e^{2\pi\i\bar t_{n-1,n}} 
\\ & \Phi_\lambda^{1,2,3...n}  ...
\Phi_\lambda^{1...,i-1,...n} ...\Phi_\lambda^{1...,n-2,n-1\ n}
= \Phi_{\lambda,1}^{-1} e^{2\pi\i\sum_{i<j}\bar t_{ij}}\Phi_{\lambda,1}
= e^{2\pi\i\sum_{i<j}\bar t_{ij}}, 
\end{align*}
where the second equality follows from the fact that $\Phi^{1...,i,...n}$
commutes with $(e^{2\pi\i\bar t_{12}})^{j,j+1...n}$ whenever $j>i$, 
and the last equality follows from the fact that $\sum_{i<j}t_{ij}$
is central is $\t_n$.

So the product of the images of $C_{12}...C_{n-1,n}$ is 
$e^{2\pi\i\sum_{i<j}\bar t_{ij}}$.

The relation $(\Theta\Psi)^{3}=C_{12}...C_{n-1,n}$ is then preserved
because 
$([\tilde\Theta]e^{\i{\pi\over 2}\sum_{i<j}\bar t_{ij}}
[\tilde\Psi]e^{\i{\pi\over 6}\sum_{i<j}\bar t_{ij}})^3 = 
([\tilde\Theta][\tilde\Psi])^3 e^{2\pi\i\sum_{i<j}\bar t_{ij}}
= [(\tilde\Theta\tilde\Psi)^3] e^{2\pi\i\sum_{i<j}\bar t_{ij}} = 
e^{2\pi\i\sum_{i<j}\bar t_{ij}}$, where the first equality follows 
from the fact that $\sum_{i<j}\bar t_{ij}$ commutes with the 
image of ${\bold G}_{1} \to {\bold G}_{n}$, 
$g\mapsto [g]$, the second equality follows from the fact that 
$g\mapsto [g]$ is a group morphism and the last equality follows from 
(\ref{gamma11}). In the same way, one proves that 
$\Theta^{4}=C_{12}...C_{n-1,n}$, $\sigma_i^2 = C_{i,i+1} 
C_{i+1,i+2} C_{i,i+1}^{-1}$ and $(\Theta,\sigma_{i})=(\Psi,\sigma_{i})=1$ 
are preserved. \hfill \qed \medskip

\subsection{Construction of morphisms $\overline{\on{B}}_{1,n}
\to \on{exp}(\widehat{{\bar\t}_{1,n}^\kk})\rtimes S_{n}$ using an 
associator $\Phi_\lambda$}

Let us keep the notation of the previous section. 
Set $a_{2n}(\lambda):= -(2n+1)B_{2n+2}\lambda^{2n+2}/(2n+2)!$, 
$\tilde y_\lambda:= - {{\on{ad}x}\over{e^{\lambda\on{ad}x}-1}}(y)$, 
$$
\tilde A_\lambda := \Phi_\lambda(\tilde y_\lambda,t)
e^{\lambda \tilde y_\lambda} \Phi_\lambda(\tilde y_\lambda,t)^{-1} = 
e^{-\lambda t/2} \Phi_\lambda(-\tilde y_\lambda-t,t) 
e^{\lambda(\tilde y_\lambda+t)}
\Phi_\lambda(-\tilde y_\lambda-t,t)^{-1}e^{-\lambda t/2},
$$
$$
\tilde B_\lambda := e^{\lambda t/2} \Phi_\lambda(-\tilde y_\lambda-t,t) 
e^{\lambda x} \Phi_\lambda(\tilde y_\lambda,t)^{-1} 
$$
(the identity in the definition of $A_{\lambda}$ follows from 
the hexagon relation).

\begin{proposition} \label{construction:bis}
We have 
$$
\tilde A_\lambda^{12,3} =  e^{\lambda \bar t_{12}/2} \{\Phi_\lambda\}^{3,1,2} 
\tilde A_\lambda^{2,13} \{\Phi_\lambda\}^{2,1,3} e^{\lambda
\bar t_{12}/2}\cdot \{\Phi_\lambda\}^{3,2,1} 
\tilde A_\lambda^{1,23}\{\Phi_\lambda\}^{1,2,3}, 
$$
$$
\tilde B_\lambda^{12,3} =  e^{-\lambda\bar t_{12}/2} \{\Phi_\lambda\}^{3,1,2} 
\tilde B_\lambda^{2,13} \{\Phi_\lambda\}^{2,1,3} e^{-\lambda
\bar t_{12}/2}\cdot \{\Phi_\lambda\}^{3,2,1} 
\tilde B_\lambda^{1,23}\{\Phi_\lambda\}^{1,2,3}, 
$$
\begin{align*}
(\tilde B_\lambda^{12,3}, 
e^{\lambda \bar t_{12}/2} \{\Phi_\lambda\}^{3,1,2} 
\tilde A_\lambda^{2,13} \{\Phi_\lambda\}^{2,1,3} e^{\lambda
\bar t_{12}/2})  & = 
( e^{-\lambda\bar t_{12}/2} \{\Phi_\lambda\}^{3,1,2} 
\tilde B_\lambda^{2,13} \{\Phi_\lambda\}^{2,1,3} e^{-\lambda
\bar t_{12}/2} , \tilde A_\lambda^{12,3}) 
\\ & = 
\{\Phi_\lambda\}^{3,2,1} e^{\lambda \bar t_{23}} \{\Phi_\lambda\}^{1,2,3}, 
\end{align*}
so the formulas of Proposition \ref{construction} (restricted to 
the generators $A_i,B_i,\sigma_i,C_{jk}$) induce a morphism 
$\overline{\on{B}}_{1,n} \to \on{exp}(\widehat{{\bar\t}^\kk_{1,n}})
\rtimes S_n$ (here $\widehat{{\bar\t}^\kk_{1,n}}$ is the degree 
completion of ${\bar\t}_{1,n}^\kk$). 
\end{proposition}

{\em Proof.} In this proof, we shift the indices of the generators of 
$\t_{n+1}$ by $1$, so these generators are now $t_{ij}$, $i\neq j
\in \{0,...,n\}$ (recall that $\t_{n+1} = \t_{n+1}^\kk$, $\bar t_{1,n}
= \bar\t_{1,n}^\kk$).

We have a morphism $\alpha_{n} : \t_{n+1}\to\bar\t_{1,n}$, defined
by $t_{ij}\mapsto \bar t_{ij}$ if $1\leq i<j\leq n$ and 
$t_{0i}\mapsto \tilde y_i := - {{\on{ad}\bar x_i}\over
{e^{\lambda\on{ad}\bar x_i} -1}}(\bar y_i)$ if $1\leq i\leq n$
(it takes the central element $\sum_{0\leq i<j\leq n}t_{ij}$ to $0$).

Let $\phi : \{1,...,m\}\to \{1,...,n\}$ be a map and 
$\phi' : \{0,...,m\}\to \{0,...,n\}$ be given by 
$\phi'(1)=1$, $\phi'(i) = \phi(i)$ for $i=1,...,m$. 
The diagram 
$$
\begin{matrix}
\t_{n+1} & \stackrel{x\mapsto x^{\phi'}}{\to} & \t_{m+1} \\ 
{\scriptstyle{\alpha_{n}}}\downarrow & & 
\downarrow\scriptstyle{\alpha_{m}} \\
\bar\t_{1,n} & \stackrel{x\mapsto x^\phi}{\to} & \bar\t_{1,m}
\end{matrix}
$$
is not commutative, we have instead the identity
$$
\alpha_{m}(x^{\phi'}) = \alpha_{n}(x)^{\phi} - \sum_{i=1}^{n}\xi_{i}(x)
(\sum_{i',j' \in \phi^{-1}(i)| i'<j'} \bar t_{i'j'}), 
$$
where $\xi_{i} : \bar\t_{1,n} \to \kk$ is the linear form defined by 
$\xi_{i}(t_{0i})=1$, $\xi_{i}($any other homogeneous Lie polynomial in the 
$t_{kl}) = 0$.

Since the various $\sum_{i',j' \in \phi^{-1}(i)| i'<j'} \bar t_{i'j'}$ 
commute with each other and with the image of 
$x\mapsto x^{\phi}$, this implies 
$$
\alpha_{m}(g^{\phi'}) = \alpha_{n}(g)^{\phi}
\prod_{i=1}^{n} e^{- \xi_{i}(\on{log}g)
(\sum_{i',j'\in\phi^{-1}(i),i'<j'}\bar t_{i'j'})}
$$
for $g\in \on{exp}(\hat\t_{n+1})$.

Set $\bar A_\lambda:= \Phi_\lambda^{0,1,2} e^{\lambda t_{01}}
(\Phi_\lambda^{0,1,2})^{-1}\in \on{exp}(\hat\t_3)$. 
One proves that 
$$
\bar A_\lambda^{0,12,3} e^{\lambda t_{12}} = 
e^{\lambda t_{12}/2} \Phi_\lambda^{3,1,2} 
\bar A_\lambda^{0,2,13} \Phi_\lambda^{2,1,3}
e^{\lambda t_{12}/2} \cdot \Phi_\lambda^{3,2,1} \bar A_\lambda^{0,1,23}
\Phi_\lambda^{1,2,3} 
$$
(relation in $\on{exp}(\hat\t_4)$). 
We then have $\alpha_{2}(\bar A_{\lambda}) = \tilde A_{\lambda}$, 
$\alpha_{3}(\Phi_{\lambda}^{1,2,3}) = \Phi_{\lambda}^{1,2,3}$, 
and the relation between the $\alpha_i$ and coproducts
implies $\alpha_{3}(\bar A_{\lambda}^{0,1,23}) = 
\tilde A_{\lambda}^{1,23}$ and $\alpha_{3}(\bar A_{\lambda}^{0,12,3}
e^{\lambda t_{12}}) = \tilde A_{\lambda}^{12,3}$. 
Taking the image by $\alpha_3$, we get the first identity.

As we have already mentioned, this identity implies 
$(\Phi_{\lambda}^{-1}\tilde A_{\lambda}^{1,23}\Phi_{\lambda},
\tilde A_{\lambda}^{12,3})=1$.



Let $\on{exp}(\hat\t_{n+1}) * \ZZ^n / I_n$ be the quotient 
of the free product of $\on{exp}(\hat\t_{n+1})$ with 
$\ZZ^n = \oplus_{i=1}^n \ZZ X_i$ by the normal subgroup generated by 
the rations of the exponentials of the sides of each of the equations 
$$
X_i t_{0i} X_i^{-1} = \sum_{0\leq \alpha\leq n,
\alpha\neq i} t_{\alpha i}, \quad 
X_i(t_{0j}+t_{ij})X_i^{-1} = t_{0j}, \; 
X_i t_{jk} X_i^{-1} = t_{jk}, \; X_jX_k t_{jk}(X_jX_k)^{-1} = t_{jk}
$$ 
where $i,j,k$ are distinct in $\{1,...,n\}$. Then the morphism 
$\alpha_n : \t_{n+1}\to \bar\t_{1,n}$ extends to 
$\tilde \alpha_n : \on{exp}(\hat\t_{n+1}) * \ZZ^n / I_n \to 
\on{exp}(\hat{\bar\t}_{1,n})$ by $X_i\mapsto e^{\lambda x_i}$.

If $\phi : \{1,...,m\}\to \{1,...,n\}$ is a map, then the Lie algebra
morphism $\t_{n+1} \to \t_{m+1}$, $x\mapsto x^{\phi'}$ extends to 
a group morphism $\on{exp}(\hat\t_n) * \ZZ^n/I_n \to 
\on{exp}(\hat\t_m) * \ZZ^m/I_m$
by $X_i\mapsto \prod_{i'\in\phi^{-1}(i)}X_{i'}$.

Let 
$$
\bar B_\lambda:= e^{\lambda t_{12}/2} \Phi_\lambda^{0,2,1} X_1 
\Phi_\lambda^{2,1,0} \in \on{exp}(\hat\t_{3}) * \ZZ^2 / I_2, 
$$
then $\alpha_2(\bar B_\lambda) = \tilde B_\lambda$.

We will prove that 
\begin{equation} \label{rel:B}
\bar B_{\lambda}^{0,12,3} = e^{-\lambda t_{12}/2} \Phi_\lambda^{3,1,2} 
\bar B_{\lambda}^{0,2,13} \Phi_\lambda^{2,1,3} e^{-\lambda t_{12}/2} \cdot 
\Phi^{3,2,1} \bar B_{\lambda}^{0,1,23} \Phi_\lambda^{1,2,3}. 
\end{equation}

The l.h.s. is 
$$
\bar B_{\lambda}^{0,12,3} = e^{\lambda t_{3,12}/2} \Phi_\lambda^{0,3,12} 
X_{1}X_{2}\Phi_\lambda^{3,21,0}
$$
and the r.h.s. is 
$$
e^{-\lambda t_{12}/2} \Phi_\lambda^{3,1,2} e^{\lambda t_{31,2}/2} 
\Phi_\lambda^{0,13,2} X_{2} \Phi_\lambda^{13,2,0} \Phi_\lambda^{2,1,3} 
e^{-\lambda t_{12}/2} \Phi_\lambda^{3,2,1} e^{\lambda t_{23,1}/2} 
\Phi_\lambda^{0,23,1} X_{1}\Phi_\lambda^{32,1,0}\Phi_\lambda^{1,2,3}. 
$$
The equality between these terms is rewritten as 
$$
X_{1}X_{2} = \Phi_\lambda^{03,1,2}\Phi_\lambda^{1,3,0} e^{-\lambda t_{13}/2} 
X_{2}\Phi_\lambda^{13,2,0}e^{\lambda t_{13}/2} \Phi_\lambda^{2,3,1} 
\Phi_\lambda^{0,23,1} X_{1} \Phi_\lambda^{01,2,3}\Phi_\lambda^{2,1,0}, 
$$
or, using the fact that $X_{i}$ commutes with $t_{jk}$ ($i,j,k$ distinct), as 
$$
X_{1}X_{2} = \Phi_\lambda^{03,1,2}\Phi_\lambda^{1,3,0} X_{2}
\Phi_\lambda^{02,3,1} \Phi_\lambda^{3,2,0} X_{1} \Phi_\lambda^{01,2,3}
\Phi_\lambda^{2,1,0}. 
$$
Now $X_{2}\Phi_\lambda^{02,3,1} = \Phi_\lambda^{0,3,1}X_{2}$, $X_{1}
\Phi_\lambda^{01,2,3}= \Phi_\lambda^{0,2,3}X_{1}$ and 
$X_{1}X_{2}\Phi_\lambda^{2,1,0} = \Phi_\lambda^{2,1,03}X_{1}X_{2}$, 
so the r.h.s. is rewritten as
$\Phi_\lambda^{03,1,2}\Phi_\lambda^{1,3,0}\Phi_\lambda^{0,3,1} X_{2}
\Phi_\lambda^{3,2,0} \Phi_\lambda^{0,2,3} X_{1}\Phi_\lambda^{2,1,0} 
= X_{1}X_{2}$. This ends the proof of (\ref{rel:B}). Taking the image by 
$\alpha_4$, we then get the second identity of the Proposition.


Let us prove the next identity. We have 
\begin{align*}
& (\bar B_\lambda^{0,12,3},e^{\lambda\bar t_{12}/2} \Phi_\lambda^{3,1,2}
\bar A_\lambda^{0,2,13}\Phi_\lambda^{2,1,3}e^{\lambda \bar t_{12}/2})
\\ & = 
e^{\lambda t_{12,3}/2}\Phi_\lambda^{0,3,12}X_1X_2 \Phi_\lambda^{3,12,0}
e^{\lambda\bar t_{12}/2} \Phi_\lambda^{3,1,2}
\Phi_\lambda^{0,2,13}e^{\lambda t_{0,2}}\Phi_\lambda^{13,2,0}
\Phi_\lambda^{2,1,3}e^{\lambda \bar t_{12}/2}
\Phi_\lambda^{0,12,3}(X_1X_2)^{-1}
\\ & \Phi_\lambda^{12,3,0}
e^{-\lambda t_{12,3}/2}
e^{-\lambda\bar t_{12}/2} \Phi_\lambda^{3,1,2}
\Phi_\lambda^{0,2,13}e^{-\lambda t_{0,2}}\Phi_\lambda^{13,2,0}
\Phi_\lambda^{2,1,3}e^{-\lambda \bar t_{12}/2}. 
\end{align*}
Now 
\begin{align*}
& X_1X_2 \Phi_\lambda^{3,12,0}
e^{\lambda\bar t_{12}/2} \Phi_\lambda^{3,1,2}
\Phi_\lambda^{0,2,13}e^{\lambda t_{0,2}}\Phi_\lambda^{13,2,0}
\Phi_\lambda^{2,1,3}e^{\lambda \bar t_{12}/2}
\Phi_\lambda^{0,12,3}(X_1X_2)^{-1}
\\ & = 
e^{\lambda\bar t_{12}/2} 
X_1X_2 \Phi_\lambda^{3,12,0}
\Phi_\lambda^{3,1,2}
\Phi_\lambda^{0,2,13}e^{\lambda t_{0,2}}\Phi_\lambda^{13,2,0}
\Phi_\lambda^{2,1,3}\Phi_\lambda^{0,12,3}(X_1X_2)^{-1}
e^{\lambda \bar t_{12}/2}
\\ & = 
e^{\lambda\bar t_{12}/2} 
X_1X_2 \Phi_\lambda^{0,2,1}
\Phi_\lambda^{3,1,02}
e^{\lambda t_{0,2}}\Phi_\lambda^{02,1,3}
\Phi_\lambda^{0,2,1}(X_1X_2)^{-1}
e^{\lambda \bar t_{12}/2}
\\ & = 
e^{\lambda\bar t_{12}/2} 
X_1X_2 \Phi_\lambda^{0,2,1}
e^{\lambda t_{0,2}}
\Phi_\lambda^{0,2,1}(X_1X_2)^{-1}
e^{\lambda \bar t_{12}/2}
\\ & = 
e^{\lambda\bar t_{12}/2}
  \Phi_\lambda^{03,2,1}
X_1X_2 e^{\lambda t_{0,2}} (X_1X_2)^{-1}
\Phi_\lambda^{03,2,1}e^{\lambda \bar t_{12}/2}
\\ & = 
e^{\lambda\bar t_{12}/2}
  \Phi_\lambda^{03,2,1}
  e^{\lambda t_{03,2}} \Phi_\lambda^{03,2,1}e^{\lambda \bar t_{12}/2}. 
\end{align*}
Plugging this in the above expression for 
$(\bar B_\lambda^{0,12,3},e^{\lambda\bar t_{12}/2} \Phi_\lambda^{3,1,2}
\bar A_\lambda^{0,2,13}\Phi_\lambda^{2,1,3}e^{\lambda t_{12}/2})$, one then 
finds
$(\bar B_\lambda^{0,12,3},e^{\lambda\bar t_{12}/2} \Phi_\lambda^{3,1,2}
\bar A_\lambda^{0,2,13}\Phi_\lambda^{2,1,3}e^{\lambda \bar t_{12}/2}) =
= \Phi_\lambda^{3,2,1} e^{\lambda t_{23}} \Phi_\lambda^{1,2,3}$. 
Taking the image by $\alpha_4$, we then obtain
$(\tilde B_\lambda^{12,3},e^{\lambda\bar t_{12}/2} \Phi_\lambda^{3,1,2}
\tilde A_\lambda^{2,13}\Phi_\lambda^{2,1,3}e^{\lambda \bar t_{12}/2})
= \Phi_\lambda^{3,2,1} e^{\lambda \bar t_{23}} \Phi_\lambda^{1,2,3}$.

Let us prove that last identity. For this, we will show 
$$
(e^{-\lambda t_{12}/2}\Phi_\lambda^{3,1,2}\bar
B_\lambda^{0,2,13}\Phi_\lambda^{2,1,3}e^{-\lambda t_{12}/2},
\bar A_\lambda^{0,12,3}e^{\lambda t_{12}}) 
= \Phi_\lambda^{3,2,1}e^{\lambda t_{23}}\Phi_\lambda^{1,2,3} 
$$
and take the image by $\alpha_4$.

We have 
\begin{align*}
& (e^{-\lambda t_{12}/2}\Phi_\lambda^{3,1,2}\bar
B_\lambda^{0,2,13}\Phi_\lambda^{2,1,3}e^{-\lambda t_{12}/2},
\bar A_\lambda^{0,12,3}e^{\lambda t_{12}}) 
\\ & 
=
e^{-\lambda t_{12}/2}\Phi_\lambda^{3,1,2}
e^{\lambda t_{2,13}/2}\Phi_\lambda^{0,13,2} X_2 \Phi_\lambda^{13,2,0}
\Phi_\lambda^{2,1,3}e^{-\lambda t_{12}/2}
\Phi_\lambda^{0,12,3}e^{\lambda t_{0,12}}\Phi_\lambda^{3,12,0}
e^{\lambda t_{12}}
e^{\lambda t_{12}/2}\Phi_\lambda^{3,1,2}
\Phi_\lambda^{0,2,13} X_2^{-1} 
\\ & 
\Phi_\lambda^{2,13,0} e^{-\lambda t_{2,13}/2}
\Phi_\lambda^{2,1,3}e^{\lambda t_{12}/2}
\Phi_\lambda^{0,12,3}e^{-\lambda t_{0,12}}\Phi_\lambda^{3,12,0}
e^{-\lambda t_{12}}
\\ & = 
e^{-\lambda t_{12}/2}\Phi_\lambda^{3,1,2}
e^{\lambda t_{2,13}/2}\Phi_\lambda^{0,13,2} X_2 \Phi_\lambda^{13,2,0}
\Phi_\lambda^{2,1,3}
\Phi_\lambda^{0,12,3}e^{\lambda t_{0,12} + \lambda t_{12}}
\Phi_\lambda^{3,12,0}\Phi_\lambda^{3,1,2}
\Phi_\lambda^{0,2,13} X_2^{-1} 
\\ & 
\Phi_\lambda^{2,13,0} e^{-\lambda t_{2,13}/2}
\Phi_\lambda^{2,1,3}e^{-\lambda t_{12}/2}
\Phi_\lambda^{0,12,3}e^{-\lambda t_{0,12}}\Phi_\lambda^{3,12,0}. 
\end{align*}
Now 
\begin{align*}
& X_2 \Phi_\lambda^{13,2,0}
\Phi_\lambda^{2,1,3}
\Phi_\lambda^{0,12,3}e^{\lambda t_{0,12} + \lambda t_{12}}
\Phi_\lambda^{3,12,0}\Phi_\lambda^{3,1,2}
\Phi_\lambda^{0,2,13} X_2^{-1}
\\ & 
= X_2 \Phi_\lambda^{02,1,3}
\Phi_\lambda^{1,2,0}
e^{\lambda t_{0,12} + \lambda t_{12}}
\Phi_\lambda^{0,2,1}
\Phi_\lambda^{3,1,02} X_2^{-1}
= \Phi_\lambda^{0,1,3} X_2 
\Phi_\lambda^{1,2,0}
e^{\lambda t_{0,12} + \lambda t_{12}}
\Phi_\lambda^{0,2,1}
X_2^{-1} \Phi_\lambda^{3,1,0}
\\ & 
= \Phi_\lambda^{0,1,3} X_2 
e^{\lambda(t_{01}+t_{02}+t_{12})}
X_2^{-1} \Phi_\lambda^{3,1,0}
= \Phi_\lambda^{0,1,3} 
e^{\lambda(t_{01}+t_{02}+t_{12} + t_{23})}
\Phi_\lambda^{3,1,0}. 
\end{align*}
So 
\begin{align*}
& (e^{-\lambda t_{12}/2}\Phi_\lambda^{3,1,2}\bar
B_\lambda^{0,2,13}\Phi_\lambda^{2,1,3}e^{-\lambda t_{12}/2},
\bar A_\lambda^{0,12,3}e^{\lambda t_{12}}) 
\\ & 
=
e^{-\lambda t_{12}/2}\Phi_\lambda^{3,1,2}
e^{\lambda t_{2,13}/2}\Phi_\lambda^{0,13,2} 
\Phi_\lambda^{0,1,3} 
e^{\lambda(t_{01}+t_{02}+t_{12} + t_{23})}
\\ & 
\Phi_\lambda^{3,1,0} 
\Phi_\lambda^{2,13,0} e^{-\lambda t_{2,13}/2}
\Phi_\lambda^{2,1,3}e^{-\lambda t_{12}/2}
\Phi_\lambda^{0,12,3}e^{-\lambda t_{0,12}}\Phi_\lambda^{3,12,0}; 
\end{align*}
after some computation, we find that this equals 
$\Phi_{\lambda}^{3,2,1}e^{\lambda t_{23}}\Phi_{\lambda}^{1,2,3}$. 
\hfill \qed \medskip

In particular, $(\Phi_\lambda,\tilde A_\lambda,\tilde B_\lambda)$
give rise to a morphism $\overline{\on{B}}_{1,n} \to 
\exp(\widehat{{\bar\t}_{1,n}^\kk})\rtimes S_n$; one proves as in Section 
\ref{sect:2} that it induces an isomorphism of filtered Lie algebras
$\on{Lie}(\overline{\on{PB}}_{1,n})_\kk \simeq \widehat{{\bar\t}_{1,n}^\kk}$. 
Taking $\Phi_\lambda$ to be a rational associator (\cite{Dr:Gal}), 
we then obtain:

\begin{corollary} We have a filtered isomorphism
$\on{Lie}(\overline{\on{PB}}_{1,n})_\QQ \simeq \widehat{{\bar\t}_{1,n}^\QQ}$, 
which can be extended to an isomorphism 
$\overline{\on{B}}_{1,n}(\QQ) \simeq \on{exp}(\widehat{{\bar\t}_{1,n}^\QQ})
\rtimes S_n$. 
\end{corollary}

\subsection{Construction of morphisms $\Gamma_{1,[n]}\to {\bold G}_{1,n}
\rtimes S_{n}$ using a pair $(\Phi_\lambda,\tilde\Theta_\lambda)$}

Keep the notation of the previous section and set 
$$
\tilde\Psi_{\lambda}:= \on{exp}( -{1\over\lambda}(\Delta_0 + \sum_{k\geq 1}
a_{2k}(\lambda)\delta_{2k})). 
$$

\begin{proposition}
We have 
$$
[\tilde\Psi_\lambda]e^{\lambda\bar t_{12}/12} \tilde A_\lambda
([\tilde\Psi_\lambda]e^{\lambda\bar t_{12}/12})^{-1} = \tilde A_\lambda,
\quad 
[\tilde\Psi_\lambda]e^{\lambda\bar t_{12}/12} \tilde B_\lambda
([\tilde\Psi_\lambda]e^{\lambda\bar t_{12}/12})^{-1} = \tilde B_\lambda
\tilde A_\lambda. 
$$ 
\end{proposition}

{\em Proof.}
The first identity follows from the fact that 
$\Delta_0 + \sum_{k\geq 1} a_{2k}(\lambda)[\delta_{2k}] - \lambda^2 t/12$ 
commutes with $t$ and $\tilde y_\lambda$; the second identity 
follows from these facts and the analogue of Lemma \ref{adjoint}, where
$2\pi\i$ is replaced by $\lambda$. 
\hfill \qed \medskip

Assume that $\tilde\Theta_\lambda\in{\bold G}_1$ satisfies 
$$
\tilde\Theta_\lambda^4 = (\tilde\Theta_\lambda\tilde\Psi_\lambda)^3 =
(\tilde\Theta_\lambda^2,\tilde\Psi_\lambda)=1, 
$$
$$
[\tilde\Theta_\lambda] e^{\lambda\bar t_{12}/4} \tilde A_\lambda 
([\tilde\Theta_\lambda] e^{\lambda\bar t_{12}/4})^{-1} = 
\tilde B_\lambda^{-1}, \quad 
[\tilde\Theta_\lambda] e^{\lambda\bar t_{12}/4} \tilde B_\lambda 
([\tilde\Theta_\lambda] e^{\lambda\bar t_{12}/4})^{-1} = 
\tilde B_\lambda \tilde A_\lambda \tilde B_\lambda^{-1} 
$$
(one can show that the two last equations are equivalent), 
then $\Theta \mapsto [\tilde\Theta_\lambda]e^{\lambda(\sum_{i<j}\bar
t_{ij})/4}$, $\Psi \mapsto [\tilde\Psi_\lambda]e^{\lambda(\sum_{i<j}
\bar t_{ij})/12}$ extends the morphism defined in Proposition
\ref{construction:bis} to a morphism 
$\Gamma_{1,[n]}\to {\bold G}_n\rtimes S_{n}$.

We do not know whether for each $\Phi_\lambda$ defined over $\kk$, 
there exists a $\tilde\Theta_\lambda$ defined over $\kk$, satisfying 
the above conditions.

\subsection{Elliptic structures over QTQBA's} \label{5:6}

Let $(H,\Delta_{H},R_{H},\Phi_{H})$ be a quasitriangular 
quasibialgebra (QTQBA). Recall that this means that (\cite{Dr:QH}): 
$(H,m_{H})$ is an algebra, $\Delta_{H} : H \to H^{\otimes 2}$
is an algebra morphism, $R_{H}\in H^{\otimes 2}$ and $\Phi_{H}
\in H^{\otimes 3}$ are invertible, and 
$$
\Delta_{H}(x)^{2,1} = R_{H}\Delta_{H}(x)R_{H}^{-1}, 
\quad
(\on{id}\otimes \Delta_{H})\circ \Delta_{H}(x) = 
\Phi_{H} (\Delta_{H}\otimes \on{id})\circ \Delta_{H}(x) \Phi_{H}^{-1}, 
$$
$$
R_{H}^{12,3} = \Phi_{H}^{3,1,2}R_{H}^{1,3}(\Phi_{H}^{1,3,2})^{-1}R_{H}^{2,3}
\Phi_{H}^{1,2,3}, \quad 
R_{H}^{1,23} = (\Phi_{H}^{2,3,1})^{-1}R_{H}^{1,3}\Phi_{H}^{2,1,3}
R_{H}^{1,2} (\Phi_{H}^{1,2,3})^{-1}, 
$$
$$ 
\Phi_{H}^{1,2,34}\Phi_{H}^{12,3,4} = \Phi_{H}^{2,3,4}\Phi_{H}^{1,23,4}
\Phi_{H}^{1,2,3}. 
$$
One also assumes the existence of a unit $1_{H}$ and a counit $\varepsilon_{H}$.

If ${\bf A}$ is an algebra and $J_{1},J_{2}\subset {\bf A}$ are left ideals, 
define the Hecke bimodule  ${\cal H}({\bf A}|J_{1},J_{2})$
or ${\cal H}(J_{1},J_{2})$ as $\on{Hom}_{{\bf A}}({\bf A}/J_1,{\bf A}/J_2) 
= ({\bf A}/J_{2})^{J_{1}}$ where $J_{1}$ acts on the quotient from the 
left; we have thus ${\cal H}(J_{1},J_{2})
= \{x\in {\bf A} | J_{1}x \subset J_{2}\} / J_{2}$. The product of ${\bf A}$
induces a product ${\cal H}(J_{1},J_{2}) \otimes {\cal H}(J_{2},J_{3}) \to 
{\cal H}(J_{1},J_{3})$. When $J_{1} = J_{2} = J$, ${\cal H}(J):= {\cal H}(J,J)$
is the usual Hecke algebra, and ${\cal H}(J_{1},J_{2})$ is a $({\cal H}(J_{1}),
{\cal H}(J_{2}))$-bimodule. Recall that we have a functor 
${\bf A}\on{-mod}\to {\cal H}(J)\on{-mod}$, 
$V\mapsto V^{J}:= \{v\in V | Jv = 0\}$.

If $H$ is an algebra with unit equipped with a morphism 
$\Delta_{H} : H \to H^{\otimes 2}$
and $a: H \to D$ is a morphism of algebras with unit, we define 
for each $n\geq 1$ and each pair of words $w,w'$ in the free magma generated by 
$1,...,n$ containing $1,...,n$ exactly once (recall that a magma is a set with a 
non-necessarily associative binary operation) the Hecke bimodule
$$
{\cal H}^{w,w'}(D,H):= {\cal H}(D\otimes H^{\otimes n} | J_{w},J_{w'}), 
$$
(or simply ${\cal H}^{w,w'}$)
where $J_{w} \subset D\otimes H^{\otimes n}$ is the left ideal generated by 
the image of $(a\otimes \Delta_{H}^{w})\circ \Delta_{H} : H_{+}\to 
D\otimes H^{\otimes n}$. Here $H_{+} = \on{Ker}(H
\stackrel{\varepsilon_{H}}{\to}\kk)$ and for example 
$\Delta_{H}^{(21)3} = (213) \circ 
(\Delta_{H}\otimes \on{id}_{H}) \circ \Delta_{H}$, etc. 
We have products ${\cal H}^{w,w'}\otimes {\cal H}^{w',w''}\to {\cal H}^{w,w''}$. 
We denote the Hecke algebra ${\cal H}^{w,w}$ by ${\cal H}^{w}(D,H)$
or ${\cal H}^{w}$; 
we denote by $1_{w}$ its unit. We denote by $({\cal H}^{w,w'})^{\times}$
the set of invertible elements of ${\cal H}^{w,w'}$, i.e., the set of elements 
$X$ such that for some $X'\in {\cal H}^{w',w}$, $X'X = 1_{w'}$, 
$XX' = 1_{w}$. 
The symmetric group $S_{n}$ acts on the system of bimodules 
${\cal H}^{w,w'}$ by permuting the factors, so we get 
maps $\on{Ad}(\sigma) : {\cal H}^{w,w'}\to 
{\cal H}^{\sigma(w),\sigma(w')}$
(where $\sigma(w)$ is the word $w$, where $i$ is replaced by $\sigma(i)$). 
If $w_{0} = ((12)...)n$, we define an algebra structure on 
$\oplus_{\sigma\in S_{n}} {\cal H}^{w_{0},\sigma(w_{0})}\sigma$
by $(\sum_{\sigma\in S_{n}} h_{\sigma}\sigma)(\sum_{\tau\in S_{n}} 
h'_{\tau}\tau):= \sum_{\sigma,\tau\in S_{n}} h_{\sigma}
\on{Ad}(\sigma)(h'_{\tau})\sigma\tau$. 
Then $\sqcup_{\sigma\in S_{n}} 
({\cal H}^{w_{0},\sigma(w_{0})})^{\times}\sigma \subset 
\oplus_{\sigma\in S_{n}} {\cal H}^{w_{0},\sigma(w_{0})}\sigma$
is a group with unit $1_{w_{0}}$. We have an exact sequence 
$1\to ({\cal H}^{w_{0}})^{\times} \to \sqcup_{\sigma\in S_{n}}
({\cal H}^{w_{0},\sigma(w_{0})})^{\times}\sigma \to S_{n}$, 
but the last map is not necessarily surjective (and if it is, does not necessarily split).

If $H$ is a quasibialgebra, then $\Phi_{H}$ gives rise to an element 
of ${\cal H}^{1(23),(12)3}(D,H)$, which we also denote $\Phi_{H}$; 
similarly $\Phi_{H}^{-1}$ gives rise to the inverse (w.r.t. composition 
of Hecke bimodules) element $\Phi_{H}^{-1}\in {\cal H}^{(12)3,1(23)}(D,H)$. 
We have algebra morphisms ${\cal H}^{12}(D,H) \to {\cal H}^{(12)3}(D,H)$
induced by $X\mapsto X^{0,12,3}:= (\on{id}_{D}\otimes 
\Delta_{H}\otimes \on{id}_{H})(X)$ ($0$ is the 
index of $D$) and similarly morphisms 
${\cal H}^{12}(D,H)\to {\cal H}^{2(13)}(D,H)$, 
$X\mapsto X^{0,2,13}$, ${\cal H}^{12}(D,H)\to {\cal H}^{1}(D,H)$, 
$X\mapsto X^{0,1,\emptyset}$ and $X^{0,\emptyset,1}$, etc. 
If moreover $H$ is quasitriangular, then $R_{H}\in {\cal H}^{21,12}(D,H)$, 
$R_{H}^{-1}\in {\cal H}^{12,21}(D,H)$, so in that case 
$\sqcup_{\sigma\in S_n} {\cal H}^{w_0,\sigma(w_0)}\sigma\to S_n$ 
is surjective, and we have a morphism $\on{B}_n \to 
\sqcup_{\sigma\in S_n} {\cal H}^{w_0,\sigma(w_0)}\sigma$ such that 
the composition $\on{B}_n \to \sqcup_{\sigma\in S_n} 
{\cal H}^{w_0,\sigma(w_0)}\sigma \to S_n$ is the canonical projection.

\begin{definition}
If $H$ is a QTQBA, an elliptic structure on $H$ is a triple $(D,A,B)$, 
where $D$ is an algebra with unit, 
equipped with an algebra morphism  $a : H\to D$, 
and $A,B\in {\cal H}^{12}(D,H)$ are invertible such that 
$A^{0,1,\emptyset} = A^{0,\emptyset,1} = B^{0,1,\emptyset} = 
B^{0,\emptyset,1} = 1_{D}\otimes 1_{H}$, 
\begin{equation} \label{A:id}
A^{0,12,3} = R_{H}^{2,1}(\Phi_{H}^{2,1,3})^{-1} A^{0,2,13}
\Phi_{H}^{2,1,3} R_{H}^{1,2}(\Phi_{H}^{1,2,3})^{-1}A^{0,1,23}\Phi_{H}^{1,2,3}, 
\end{equation}
\begin{equation} \label{B:id}
B^{0,12,3} = (R_{H}^{1,2})^{-1} (\Phi_{H}^{2,1,3})^{-1} B^{0,2,13} 
\Phi_{H}^{2,1,3}
(R_{H}^{2,1})^{-1} (\Phi_{H}^{1,2,3})^{-1} B^{0,1,23} \Phi_{H}^{1,2,3} 
\end{equation}
and 
\begin{align*}
& (B^{0,12,3},R_{H}^{2,1}(\Phi_{H}^{2,1,3})^{-1}
A^{0,2,13}\Phi_{H}^{2,1,3}R_{H}^{1,2}) 
\\ & = 
((R_{H}^{1,2})^{-1}(\Phi_{H}^{2,1,3})^{-1}B^{0,2,13}\Phi_{H}^{2,1,3}
(R_{H}^{2,1})^{-1},A^{0,12,3})
= (\Phi_{H}^{1,2,3})^{-1}R_{H}^{3,2}R_{H}^{2,3}\Phi_{H}^{1,2,3} 
\end{align*}
(identities in ${\cal H}^{(12)3}(D,H)$).
\end{definition}

The pair of identities (\ref{A:id}), (\ref{B:id}) is equivalent to 
$$
R_{H}^{2,1}A^{0,2,1}R_{H}^{1,2}A^{0,1,2} =1, \quad 
R_{H}^{3,12}A^{0,3,12}\Phi_{H}^{3,1,2}R_{H}^{2,31}A^{0,2,31}\Phi_{H}^{2,3,1}
R_{H}^{1,23}A^{0,1,23}\Phi_{H}^{1,2,3}=1, 
$$
and
$$
(R_{H}^{1,2})^{-1}B^{0,2,1}(R_{H}^{2,1})^{-1}B^{0,1,2}=1, \quad 
(R_{H}^{-1})^{12,3}B^{0,3,12}\Phi_{H}^{3,1,2}
(R_{H}^{-1})^{31,2}B^{0,2,31}\Phi_{H}^{2,3,1}
(R_{H}^{-1})^{23,1}B^{0,1,23}\Phi_{H}^{1,2,3}=1, 
$$
so the invertibility conditions on $A,B$ follow from (\ref{A:id}), (\ref{B:id}).

If $F\in H^{\otimes 2}$ is invertible with 
$(\varepsilon_{H}\otimes \on{id}_{H})(F)
= (\on{id}_{H}\otimes \varepsilon_{H})(F)=1_{H}$, then the twist of $H$ by $F$ is
the quasi-Hopf algebra ${}^{F}\!H$ with product $m_{H}$, coproduct 
$\tilde \Delta_{H}(x) = F \Delta_{H}(x)F^{-1}$, $R$-matrix $\tilde R_{H} = 
F^{2,1}R_{H}F^{-1}$ and associator $\tilde \Phi_{H} = F^{2,3}F^{1,23}\Phi_{H}
(F^{1,2}F^{12,3})^{-1}$. If $a : H\to D$ is an algebra morphism, it can be viewed
as a morphism ${}^{F}\!H\to D$, and we have an algebra isomorphism 
${\cal H}^{(12)3}(D,H)\to {\cal H}^{(12)3}(D,{}^{F}\!H)$, induced by 
$X\mapsto F^{1,2}F^{0,12}X(F^{1,2}F^{0,12})^{-1}$ (more generally, 
we have an isomorphism of the systems of bimodules ${\cal H}^{w,w'}(D,H)
\to {\cal H}^{w,w'}(D,{}^{F}\!H)$ induced by $X\mapsto F_{w}X F_{w'}^{-1}$
for suitable $F_{w}$).

If $(D,A,B)$ is an elliptic structure on $H$, then an elliptic structure ${}^{F}\!H$ is 
$(D,\tilde A,\tilde B)$, where $\tilde A = F^{1,2}F^{0,12}A
(F^{1,2}F^{0,12})^{-1}$
and $\tilde B = F^{1,2} F^{0,12} B (F^{1,2}F^{0,12})^{-1}$.

An elliptic structure $(D,A,B)$ over $H$ gives rise to a unique group morphism 
$$
\overline{\on{B}}_{1,n} \to \sqcup_{\sigma\in S_{n}}
{\cal H}^{w_{0},\sigma(w_{0})}(D,H)^{\times}\sigma, 
$$
such that 
$$
\sigma_{i}\mapsto \Big(\Phi_{H}^{(((12)3)...i-1),i,i+1}
\Big)^{-1} R_{H}^{i,i+1}(i,i+1) \Phi_{H}^{(((12)3)...i-1),i,i+1}, 
$$
$$
A_{i}\mapsto \Phi_{H,i}^{-1}A^{0,(((12)3)...i-1),(i...(n-1,n))}\Phi_{H,i},
\quad B_{i}\mapsto\Phi_{H,i}^{-1}B^{0,(((12)3)...i-1),(i...(n-1,n))}\Phi_{H,i}, 
$$
where 
$$
\Phi_{H,i} = \Phi_{H}^{((12)...i-1),i,(i+1(...(n-1,n)))}...
\Phi_{H}^{((12)...n-2),n-1,n}; 
$$
here we have for example $x^{((12)3)} = (\Delta_{H}\otimes \on{id}_{H})
\circ \Delta_{H}(x)$ for $x\in H$.

If $\g$ is a Lie algebra and $t_{\g}\in S^{2}(\g)^{\g}$ is nondegenerate, 
then $H = U(\g)[[\hbar]]$ is a QTQBA, with $m_{H},\Delta_{H}$ are the 
undeformed product and coproduct, $R_{H} = e^{\hbar t_{\g}/2}$ and $\Phi_{H} 
= \Phi(\hbar t_{\g}^{1,2},\hbar t_{\g}^{2,3})$, where $\Phi$ is an 
$1$-associator. The results of next Section then imply that $(D,A,B)$ is 
an elliptic structure over $H$, where $D = D(\g)[[\hbar]]$ ($D(\g)$ is the 
algebra of algebraic differential operators on $\g$) and $A,B$ are given 
by the formulas for $\tilde A_{\lambda},\tilde B_{\lambda}$ with $t$ 
replaced by $\hbar t_{\g}^{1,2}$, $x$ replaced by 
$\hbar\sum_{\alpha}\on{x}_{\alpha} \otimes (e_{\alpha}^{1}
+ e_{\alpha}^{2})$, $y$ replaced by $\hbar\sum_{\alpha}\partial_{\alpha}
\otimes (e_{\alpha}^{1} + e_{\alpha}^{2})$.


\begin{remark}
If $H$ is a Hopf algebra,  we have an isomorphism 
$$
{\cal H}^{w_{0}}(D,H) \simeq (D\otimes H^{\otimes n-1})^{H}, 
$$
where the right side is the commutant of the diagonal map 
$H\to D\otimes H^{\otimes n-1}$, $h\mapsto (a\otimes \on{id}_H^{\otimes n-1})
\circ \Delta_H^{(n)}(h)$. 
This map takes the class of $d\otimes h_{1}\otimes ... \otimes h_{n}$
to $d a(S_{H}(h_{n}^{(n)})) \otimes h_{1}S_{H}(h_{n}^{(n-1)}) \otimes ...
\otimes h_{n-1}S_{H}(h_{n}^{(1)})$ ($S_{H}$ is the antipode of $H$). 
So $A,B$ identify with elements ${\cal A},{\cal B}\in 
(D\otimes H)^{H}$; the conditions 
are then 
$$
{\cal A}^{0,12} = R_H^{2,1} {\cal A}^{0,2} R_H^{1,2}{\cal A}^{0,1}, \quad 
\cB^{0,12} = (R_H^{1,2})^{-1} \cB^{0,2} (R_H^{2,1})^{-1} \cB^{0,1}, 
$$
$$
(\cB^{0,12},R_H^{2,1}\cA^{0,2}R_H^{1,2}) = ((R_H^{1,2})^{-1}\cB^{0,2}
(R_H^{2,1})^{-1},\cA^{0,12}) = 
(R_H^{3,2} R_H^{1,2} R_H^{0,2}R_H^{2,0}R_H^{2,1}R_H^{2,3}
)^{\tilde 0, \tilde 1, 2\cdot \tilde 3}
$$
(conditions in $(D\otimes H^{\otimes 2})^{H}$), where the superscript 
$\on{B}'_n\rtimes \ZZ^{n-1} \to \on{B}_{n-1}\rtimes \ZZ^{n-1}$
is the map $x_0 \otimes ...\otimes x_3 \mapsto 
S_H(x_0) \otimes S_H(x_1) \otimes x_2 S_H(x_3)$.

Moreover, the morphism $\on{PB}_n \to ({\cal H}^{w_0})^\times \simeq 
(D\otimes H^{\otimes n-1})^H$ factors through 
$\on{PB}_n \to \on{PB}_{n-1} \times \ZZ^{n-1}
\to (D\otimes H^{\otimes n-1})^H$, 
where: (a) the first morphism is induced by $\ZZ^{n-1} \rtimes 
\on{B}'_n \to \ZZ^{n-1} \rtimes \on{B}_{n-1}$
(where $\on{B}'_n = \on{B}_n \times_{S_n} S_{n-1}$ is the group of 
braids leaving the last strand fixed), constructed as follows: 
we have a composition 
$\on{B}'_{n+1} \to \pi_1( (\PP^1)^{n+1} - \on{diagonals}/S_n) \to
\pi_1(\CC^n-\on{diagonals}/S_n)=\on{B}_n$, where the first 
map is induced by $\CC \subset \PP^1$, and the middle map comes from 
the fibration $\CC^n-\on{diagonals} \to (\PP^1)^{n+1}-\on{diagonals}
\to \PP^1$, 
$(z_1,...,z_n) \to (z_1,...,z_n,\infty)$ and $(z_1,...,z_{n+1}) \to z_{n+1}$
[the second projection has a section so the map between $\pi_1$'s 
is an isomorphism]; viewing $\ZZ^{n-1} \rtimes \on{B}'_n$, 
$\ZZ^{n-1} \rtimes \on{B}_{n-1}$ as fundamental groups of 
configuration spaces of points equipped with a
nonzero tangent vector, we then get the morphism 
$\ZZ^{n-1} \rtimes \on{B}'_n \to \ZZ^{n-1}\rtimes \on{B}_{n-1}$ 
(which does not restrict to a morphism $\on{B}'_n\to \on{B}_{n-1}$); 
(b) the second map is induced by the standard map 
$\on{PB}_{n-1}\times \ZZ^{n-1}\to (H^{\otimes n-1})^\times$ 
induced by $R_H = \sum_\alpha r'_\alpha 
\otimes r''_\alpha$ and the map taking the $i$th generator of $\ZZ^{n-1}$ to 
$1\otimes ... \otimes u S_H(u) \otimes ... \otimes 1$, where 
$u = \sum_i S_H(r''_\alpha)r'_\alpha$ (see \cite{Dr:coco}). The morphism 
$\on{B}_n \to \on{Aut}(({\cal H}^{w_0})^\times) 
= \on{Aut}((D\otimes H^{\otimes n-1})^H)$
extends the inner action of $\on{PB}_n$ by 
$$
\sigma_{n-1}\cdot X := 
\{R_H^{n-1,n...2n-1}X^{0,1,...,n-2,n...2n-1}R_H^{n...2n-1,n-1}\}^{
0\cdot \widetilde{2n-1},...,n-1\cdot \tilde n}
$$
(where the superscript means 
that $x_0\otimes ... \otimes x_{2n-1}$ maps to $x_0 S_H(x_{2n-1})\otimes ... 
\otimes x_{n-1}S_H(x_n)$).

We have then $\sqcup_{\sigma\in S_n} ({\cal H}^{w_0,\sigma(w_0)})^\times 
\sigma \simeq ((D\otimes H^{\otimes n-1})^\times)^H \rtimes_{\on{PB}_n} 
\on{B}_n$ (the index means that $\on{PB}_n \subset \on{B}_n$ is identified
with its image in $((D\otimes H^{\otimes n-1})^\times)^H$).

Then if $(\cA,\cB)$ is an elliptic structure over $a:H\to D$, the morphism 
$\on{B}_n \to ((D\otimes H^{\otimes n-1})^\times)^H \rtimes_{\on{PB}_n} 
\on{B}_n$ extends to a morphism 
$$
\overline{\on{B}}_{1,n} \to 
((D\otimes H^{\otimes n-1})^\times)^H \rtimes_{\on{PB}_n} 
\on{B}_n 
$$ 
via $A_i\mapsto \cA^{0,1...i-1}$, $B_i\mapsto \cB^{0,1...i-1}$.

This interpretation of 
${\cal H}^{w_{0}}$ and of the relations between $\cA,\cB$ can be extended to 
the case when $H$ is a quasi-Hopf algebra. 
\end{remark}

\begin{remark}
Let $\cal C$ be a rigid braided monoidal category. We define an elliptic 
structure 
on $\cal C$ as a quadruple $({\cal E},A,B,F)$, where $\cal E$ is a category, 
$F: {\cal E} \to {\cal C}$ is a functor,
and $A,B$ are functorial automorphisms of $F(?)\otimes ?$,
which reduce to the identity if the second factor is the neutral
object ${\bf 1}$, and such that the following equalities of automorphisms of 
$F(M)\otimes (X\otimes Y)$ hold (we write them omitting 
associativity maps, as they can be put in automatically): 
$$
A_{M,X\otimes Y}=\beta_{Y,X}A_{M,Y}\beta_{X,Y}A_{M,X}, 
$$
$$
B_{M,X\otimes Y}=\beta_{X,Y}^{-1}B_{M,Y}\beta_{Y,X}^{-1} B_{M,X}, 
$$
\begin{align*}
& (B_{M,X\otimes Y}, \beta_{Y,X}A_{M,Y}\beta_{X,Y})=
(\beta_{Y,X}^{-1}B_{M,Y}\beta_{X,Y}^{-1},A_{M,X\otimes Y})
\\ & = \beta_{(M\otimes X \otimes Y)^*,Y}\beta_{Y,(M\otimes X \otimes Y)^*}
\circ \on{can}_{M\otimes X \otimes Y}, 
\end{align*} 
where $\on{can}_X \in \on{Hom}_{\cC}({\bf 1},X \otimes X^*)$ is the 
canonical map and the r.h.s. of the last identity is viewed as 
an element of $\on{End}_{\cC}(M\otimes X \otimes Y)$ using its identification 
with $\on{Hom}_{\cC}({\bf 1},(M\otimes X \otimes Y) \otimes 
(M\otimes X \otimes Y)^*)$. An elliptic structure on a quasitriangular 
quasi-Hopf algebra $H$ gives rise to an elliptic structure on 
$H$-mod. An elliptic structure over a rigid braided monoidal category 
${\cal C}$ gives rise to representations of $\overline{\on{B}}_{1,n}$ by 
${\cal C}$-automorphisms of $F(M) \otimes X^{\otimes n-1}$. 
\end{remark}

\section{The KZB connection as a realization of the universal KZB connection}
\label{sect:6}

\subsection{Realizations of $\bar\t_{1,n}$}\label{sect:real:1n}

Let $\g$ be a Lie algebra and $t_\g\in S^2(\g)^\g$ be nondegenerate. 
We denote by $(a,b)\mapsto \langle a,b \rangle$ 
the corresponding invariant pairing.

Let $D(\g)$ be the algebra of algebraic differential operators on 
$\g$. It has generators 
$\x_a$, $\partial_a$, $a\in \g$, and relations: $a\mapsto \x_a$, 
$a\mapsto \partial_a$ are linear, $[\x_a,\x_b] = [\partial_a,\partial_b] = 0$, 
$[\partial_a,\x_b] = \langle a,b \rangle$.

There is a unique Lie algebra morphism $\g\to D(\g)$, 
$a\mapsto X_a$, where $X_a := \sum_{\alpha} \x_{[a,e_\alpha]}
\partial_{e_\alpha}$, and $t_\g = \sum_\alpha e_\alpha\otimes e_\alpha$
(it is the infinitesimal of the adjoint action). 
We also have a Lie algebra morphism $\g\to A_{n} := 
D(\g) \otimes U(\g)^{\otimes n}$, $a\mapsto Y_a := X_a \otimes 1
+ 1\otimes (\sum_{i=1}^n a^{(i)})$. We denote by $\g^{\on{diag}}$ the image 
of this morphism. We denote by ${\cal H}_{n}(\g)$ the Hecke algebra of
$(A_{n},\g^{\on{diag}})$. It is defined as the quotient 
$\{x\in A_{n} | \forall a\in\g, Y_a x \in
A_{n}\g^{\on{diag}}\}/A_{n}\g^{\on{diag}}$. 
We have a natural action of $S_{n}$ on $A_{n}$, which induces
an action of $S_{n}$ on ${\cal H}_{n}(\g)$.

If $(V_i)_{i=1,\ldots,n}$ are $\g$-modules, then 
$(S(\g) \otimes (\otimes_{i=1}^n V_i))^\g$ is a module over 
${\cal H}_{n}(\g)$. If moreover $V_{1}= ... = V_{n}$, this is a module 
over ${\cal H}_{n}(\g)\rtimes S_{n}$.

\begin{proposition}\label{prop:realization}
There is a unique Lie algebra morphism $\rho_\g : \bar\t_{1,n} \to 
{\cal H}_{n}(\g)$, $\bar x_i\mapsto \sum_\alpha \x_{\alpha} \otimes 
e_\alpha^{(i)}$, $\bar y_i \mapsto -\sum_\alpha \partial_{\alpha} \otimes 
e_\alpha^{(i)}$, $\bar t_{ij} \mapsto 1\otimes t_\g^{(ij)}$ (we set 
$\x_\alpha:= \x_{e_\alpha}$, $\partial_\alpha:= \partial_{e_\alpha}$).
\end{proposition}

{\em Proof.} The images of all the generators of $\bar\t_{1,n}$
are contained in the commutant of $\g^{\on{diag}}$ in $A_n$, therefore also 
in its normalizer. According to Lemma \ref{lemma:pres}, we will use the 
following presentation of $\bar\t_{1,n}$. Generators are
$\bar x_{i}, \bar y_{i},\bar t_{ij}$, relations are 
$[\bar x_{i},\bar x_{j}] = [\bar y_{i},\bar y_{j}]=0$, 
$[\bar x_{i},\bar y_{j}] = \bar t_{ij}$ ($i\neq j$), 
$\bar t_{ij} = \bar t_{ji}$, 
$\sum_{i}\bar x_{i} = \sum_{i}\bar y_{i}=0$, $[\bar x_{i},\bar t_{jk}] = 
[\bar y_{i},\bar t_{jk}]=0$ ($i,j,k$ distinct).

The relations $[\bar x_{i},\bar x_{j}] = [\bar y_{i},\bar y_{j}]=0$, 
$[\bar x_{i},\bar y_{j}] = \bar t_{ij}$ ($i\neq j$), 
$\bar t_{ij} = \bar t_{ji}$ and $[\bar x_{i},\bar t_{jk}] = 
[\bar y_{i},\bar t_{jk}]=0$ are obviously preserved. 
Let us check that $\sum_{i}\bar x_{i} = \sum_{i}\bar y_{i}=0$ are preserved.

We have
\begin{align*}
& \sum_{i}\rho_{\g}(\bar x_{i}) = \sum_{\alpha} \on{x}_{\alpha} \otimes 
(\sum_{i} e_{\alpha}^{(i)}) = \sum_{\alpha} (\on{x}_{\alpha} \otimes 1)
(Y_{\alpha} - X_{\alpha}\otimes 1)  \\
  &  \equiv -\sum_{\alpha} \on{x}_{\alpha}X_{\alpha} \otimes 1
  = \sum_{\alpha,\beta} \on{x}_{e_\alpha} \on{x}_{[e_{\alpha},e_{\beta}]}
  \partial_{e_{\beta}} \otimes 1 = 0
\end{align*}
since  $\on{x}_{\alpha}$ commutes with $\on{x}_{[e_{\alpha},e_{\beta}]}$
and $\sum_{\beta} e_{\beta} \otimes e_{\beta} = t_{\g}$ is invariant. 
We also have 
\begin{align*}
& \sum_{i} \rho_{\g}(\bar y_{i}) = - \sum_{\alpha} \partial_{\alpha}
\otimes (\sum_{i} e_{\alpha}^{(i)}) = -\sum_{\alpha} (\partial_{\alpha}
\otimes 1) (Y_{\alpha} - X_{\alpha} \otimes 1) \equiv \sum_{\alpha} 
\partial_{\alpha} X_{\alpha} \otimes 1
\\
& = - \sum_{\alpha,\beta} \partial_{e_{\alpha}} \on{x}_{[e_{\alpha},e_{\beta}]}
\partial_{e_{\beta}} = - \sum_{\alpha,\beta} \langle e_{\alpha},
[e_{\alpha},e_{\beta}]\rangle \partial_{e_{\beta}} - \sum_{\alpha,\beta}
\on{x}_{[e_{\alpha},e_{\beta}]} \partial_{e_{\alpha}} \partial_{e_{\beta}}; 
\end{align*}
since $t_{\g}$ is invariant and $\langle-,-\rangle$ is symmetric, we have 
$\sum_{\alpha} \langle e_{\alpha},[e_{\alpha},e_{\beta}] \rangle = 0$ for 
any $\beta$, and since $[\partial_{e_{\alpha}},\partial_{e_{\beta}}]=0$, 
we have $\sum_{\alpha,\beta}  \on{x}_{[e_{\alpha},e_{\beta}]} 
\partial_{e_{\alpha}} \partial_{e_{\beta}}$, so 
$\sum_{i}\rho_{\g}(\bar y_{i})=0$. \hfill \qed \medskip

\subsection{Realizations of $\bar\t_{1,n}\rtimes\d$}

Let $(\g,t_\g)$ be as in Subsection \ref{sect:real:1n}. We keep the 
same notations.

\begin{proposition} \label{real:der}
The Lie algebra morphism $\rho_\g : \bar\t_{1,n} \to {\cal H}_{n}(\g)$ of 
Proposition \ref{prop:realization} extends to a Lie algebra
morphism $\bar\t_{1,n} \rtimes \d \to {\cal H}_{n}(\g)$, defined by 
$\Delta_{0}\mapsto - {1\over 2}(\sum_\alpha \partial_\alpha^2) \otimes 1$, 
$X\mapsto {1\over 2}(\sum_{\alpha} \x_{\alpha}^{2}) \otimes 1$, $d\mapsto 
{1\over 2}(\sum_{\alpha}\x_{\alpha}\partial_{\alpha} 
+ \partial_{\alpha}\x_{\alpha}) \otimes 1$, 
and 
$$
\delta_{2m} \to {1\over 2}\sum_{\alpha_1,\ldots,\alpha_{2m},\alpha} 
\x_{\alpha_1} \cdots \x_{\alpha_{2m}}
\otimes (\sum_{i=1}^n (\ad(e_{\alpha_1}) \cdots
\ad(e_{\alpha_{2m}})(e_\alpha) \cdot e_\alpha)^{(i)} )
$$ 
for $m\geq 1$. This morphism further extends to a morphism 
$U(\bar\t_{1,n}\rtimes\d) \rtimes S_{n}\to {\cal H}_{n}(\g)\rtimes S_{n}$
by $\sigma\mapsto\sigma$. 
\end{proposition}

{\em Proof.} We have 
\begin{align*}
& [\rho_\g(\delta_{2m}),\rho_\g(\bar x_i)] =
{1\over 2}\sum_{\alpha_1,\ldots,\alpha_{2m},\alpha,\beta}
\x_{\alpha_1}\cdots\x_{\alpha_{2m}}\x_\beta \otimes 
[e_\beta,\ad(e_{\alpha_1}) \cdots 
\ad(e_{\alpha_{2m}})(e_\alpha)e_\alpha]^{(i)}
\\ & 
=
{1\over 2}\sum_{\alpha_1,\ldots,\alpha_{2m},\alpha,\beta}
\x_{\alpha_1}\cdots\x_{\alpha_{2m}}\x_\beta \otimes 
\sum_{\ell = 1}^{2m}
\big( 
\ad(e_{\alpha_1}) \cdots
\ad([e_\beta,e_{\alpha_\ell}])\cdots 
\ad(e_{\alpha_{2m}})(e_\alpha)e_\alpha \big) ^{(i)} = 0
\end{align*}
the second equality follows from the invariance of $t_\g$, and the 
last equality follows from the fact that the first factor is
symmetric in $(\beta,\alpha_\ell)$ while the second is antisymmetric
in $(\beta,\alpha_l)$.

$\rho_\g$ preserves the relation $[\delta_{2m},\bar t_{ij}] = [\bar t_{ij},
\ad(\bar x_i)^{2m}(\bar t_{ij})]$, because $\rho_\g(\delta_{2m} + \sum_{i<j}
\ad(\bar x_i)^{2m}(\bar t_{ij}))$ belongs to $D(\g) \otimes 
\on{Im}(\Delta^{(n)} : U(\g) \to U(\g)^{\otimes n})$, where $\Delta^{(n)}$ 
is the $n$-fold coproduct and $U(\g)$ is equipped with its standard 
bialgebra structure.

Now 
\begin{align*}
& [\rho_\g(\delta_{2m}),\rho_\g(\bar y_i)] = 
{1\over 2}\sum_{\alpha_1,\ldots,\alpha_{2m},\alpha,\beta}
\Big( \sum_j[\partial_\beta,\x_{\alpha_1}\cdots\x_{\alpha_{2m}}]
\otimes e_\beta^{(i)}
\ad(e_{\alpha_1})\cdots \ad(e_{\alpha_{2m}})(e_\alpha)^{(j)}
e_\alpha^{(j)}
\\ & 
+ \x_{\alpha_1} \cdots \x_{\alpha_{2m}} \partial_\beta
\otimes [e_\beta, \ad(e_{\alpha_1})\cdots 
\ad(e_{\alpha_{2m}})(e_\alpha)\cdot e_\alpha]^{(i)} \Big) 
\\ & = 
{1\over 2}\sum_{l=1}^{2m}
\sum_{\alpha_1,\ldots,\alpha_{2m},\alpha}
\Big( \sum_j \x_{\alpha_1}\cdots \check{\x}_{\alpha_l}\cdots
\x_{\alpha_{2m}}
\otimes e_{\alpha_l}^{(i)}
\ad(e_{\alpha_1})\cdots 
\ad(e_{\alpha_{2m}})(e_\alpha)^{(j)} e_\alpha^{(j)}
\\ & 
+ \x_{\alpha_1} \cdots \x_{\alpha_{2m}} \partial_\beta
\otimes 
\ad(e_{\alpha_1})\cdots \ad([e_\beta,e_{\alpha_l}])
\cdots \ad(e_{\alpha_{2m}})(e_\alpha)^{(i)}e_\alpha^{(i)} \Big) 
\\ & 
\equiv 
{1\over 2}\sum_{l=1}^{2m}
\sum_{\alpha_1,\ldots,\alpha_{2m},\alpha} \sum_j
  \Big( \x_{\alpha_1}\cdots \check{\x}_{\alpha_l}\cdots
\x_{\alpha_{2m}}
\otimes e_{\alpha_l}^{(i)}
\ad(e_{\alpha_1})\cdots 
\ad(e_{\alpha_{2m}})(e_\alpha)^{(j)} e_\alpha^{(j)}
\\  &
- \x_{\alpha_1}\cdots \check{\x}_{\alpha_l}\cdots
\x_{\alpha_{2m}}
\otimes \ad(e_{\alpha_1})\cdots 
\ad(e_{\alpha_{2m}})(e_\alpha)^{(i)} e_\alpha^{(i)} e_{\alpha_l}^{(j)}
\Big). 
\end{align*}
The term corresponding to $j=i$ is 
$$
{1\over 2}\sum_{l=1}^{2m}
\sum_{\alpha_1,\ldots,\alpha_{2m},\alpha}
   \x_{\alpha_1}\cdots \check{\x}_{\alpha_l}\cdots \x_{\alpha_{2m}}
\otimes 
[e_{\alpha_l}, \ad(e_{\alpha_1})\cdots 
\ad(e_{\alpha_{2m}})(e_\alpha)\cdot e_\alpha]^{(i)} 
$$
It corresponds to the linear map $S^{2m-1}(\g) \to U(\g)$, such that for 
$x\in\g$, 
\begin{align*}
& x^{2m-1} \mapsto {1\over 2}\sum_{p+q = 2m-1} \sum_{\alpha,\beta} [e_\beta,
\ad(x)^p \ad(e_\beta) \ad(x)^q(e_\alpha)\cdot e_\alpha] 
\\ & = {1\over 2}\sum_{\alpha,\beta}
\sum_{p+q+r = 2m-2} \ad(x)^p \ad([e_\beta,x]) \ad(x)^q \ad(e_\beta)
\ad(x)^r(e_\alpha)\cdot e_\alpha 
\\ & + \ad(x)^p \ad(e_{\beta}) \ad(x)^q \ad([e_\beta,x])
\ad(x)^r(e_\alpha)\cdot e_\alpha 
\end{align*}
since $\mu(t_\g) = 0$ ($\mu : \g^{\otimes 2}\to\g$ is the Lie bracket)
and $t_\g$ is $\g$-invariant. Now this is zero since $t_\g = \sum_\beta 
e_\beta \otimes e_\beta$ is invariant.

The term corresponding to $j\neq i$ corresponds to the map 
$S^{2m-1}(\g) \to U(\g)^{\otimes n}$, such that for $x\in\g$
\begin{align*}
& x^{2m-1}\mapsto -{1\over 2} \sum_{l=1}^{2m} \sum_{\alpha,\beta}
\big( (\ad x)^{l-1} (\ad e_\beta) (\ad x)^{2m-l}(e_\alpha)\cdot e_\alpha 
\big)^{(i)} e_\beta^{(j)} - (i\leftrightarrow j)
\\ & = {1\over 2} \sum_{l=1}^{2m} (-1)^{l+1} \sum_{\alpha,\beta}
\big((\ad x)^{l-1}([e_\beta,e_\alpha]) \cdot 
(\ad x)^{2m-l}(e_\alpha)\big)^{(i)} e_\beta^{(j)} - (i\leftrightarrow j) 
\\ & = {1\over 2} \sum_{l=1}^{2m} (-1)^{l-1} \sum_{\alpha,\beta}
\big( (\ad x)^{l-1}(e_\beta) \cdot (\ad x)^{2m-l}(e_\alpha) \big)^{(i)} 
[e_\alpha,e_\beta]^{(j)} - (i\leftrightarrow j) 
\\ & = {1\over 2} \sum_{l=1}^{2m} (-1)^{l} 
\big[ \sum_\alpha \big( (\ad x)^{l-1}(e_\alpha)\big)^{(i)} e_\alpha^{(j)}, 
\sum_\beta \big( (\ad x)^{2m-l}(e_\beta)\big)^{(i)} e_\beta^{(j)} \big], 
\end{align*}
which coincides with the image of ${1\over 2}\sum_{p+q=2m-1} (-1)^q
[(\ad \bar x_i)^p(\bar t_{ij}),(\ad \bar x_i)^q(\bar t_{ij})]$.

It is then clear that $\rho_\g$ preserves the commutation relations of 
$\Delta_{0},X$ and $d$ with $\delta_{2m}$. \hfill \qed \medskip

\subsection{Reductions} 
Assume that $\g$ is finite dimensional and 
we have a reductive decomposition 
$\g = \h\oplus\n$, i.e., $\h\subset\g$ is a Lie subalgebra and $\n\subset\g$
is a vector subspace such that $[\h,\n]\subset \n$; assume also that 
$t_{\g} = t_{\h} + t_{\n}$, 
where $t_{\h}\in S^{2}(\h)^{\h}$ and $t_{\n}\in S^{2}(\n)^{\h}$.

We assume that for a generic $h\in\h$, $\on{ad}(h)_{|\n}\in\on{End}(\n)$
is invertible. This condition is equivalent to the nonvanishing of 
$P(\lambda) := \on{det}(\on{ad}(\lambda^{\vee})_{|\n})\in 
S^{\on{dim}\n}(\h)$, where $\lambda\mapsto \lambda^{\vee}$ 
is the map $\h^{*}\to \h$, with 
$\lambda^{\vee} := (\lambda\otimes \on{id})(t_{\h})$. If $G$ is a Lie group 
with Lie algebra $\g$, an equivalent condition is that a generic
element of $\g^*$ is conjugate to some element in $\h^*$ (see \cite{EE}).

Let us set, for $\lambda\in\h^{*}$, 
$$
r(\lambda) := (\on{id}\otimes (\on{ad}\lambda^{\vee})_{|\n}^{-1})(t_{\n}), 
$$
Then $r : \h^*_{\on{reg}} \to \wedge^2(\n)$ is an $\h$-equivariant map (here 
$\h^*_{\on{reg}} = \{\lambda\in \h^* | P(\lambda)\neq 0\}$), 
satisfying the classical dynamical Yang-Baxter (CDYB) equation 
$$
\on{CYB}(r) - \on{Alt}(\on{d}r)=0
$$ (see \cite{EE}). Here for
$r = \sum_\alpha a_\alpha \otimes b_\alpha \otimes \ell_\alpha
\in (\n^{\otimes 2} \otimes 
S(\h)[1/P])^\h$, we set $\on{CYB}(r) = \sum_{\alpha,\alpha'} 
([a_\alpha,a_{\alpha'}] \otimes b_\alpha \otimes b_{\alpha'} + a_\alpha
  \otimes [b_\alpha,a_{\alpha'}] \otimes b_{\alpha'}
+ a_\alpha \otimes a_{\alpha'} \otimes [b_\alpha,b_{\alpha'}])\otimes
  \ell_\alpha\ell_{\alpha'}$, 
$\on{d}r := \sum_\alpha a_\alpha \otimes b_\alpha \otimes \on{d}\ell_\alpha$, 
where 
$\on{d}$ extends $S(\h)\to \h\otimes S(\h)$, $x^k \mapsto k x \otimes x^{k-1}$
and $\on{Alt}(X \otimes \ell) = (X + X^{2,3,1} + X^{3,1,2})\otimes \ell$.

We also set 
$$
\psi(\lambda) := (\on{id} \otimes (\on{ad}\lambda^\vee)_{|\n}^{-2})(t_\n). 
$$
We write $\psi(\lambda) = \sum_\alpha A_\alpha \otimes B_\alpha \otimes
L_\alpha$.

Let $D(\h)[1/P]$ be the localization at $P$ of the algebra $D(\h)$ of 
differential operators on $\h$; the latter algebra is generated by 
$\bar{\on{x}}_{h}$, $\bar\partial_{h}$, $h\in \h$, with relations 
$h\mapsto \bar{\on{x}}_{h}$, 
$h\mapsto \bar\partial_{h}$ linear, $[\bar{\on{x}}_{h},\bar{\on{x}}_{h'}] 
= [\bar\partial_{h},\bar\partial_{h'}]=0$, and
$[\bar\partial_{h},\bar{\on{x}}_{h'}] = \langle h,h'\rangle$.

Set $B_{n}:= D(\h)[1/P] \otimes U(\g)^{\otimes n}$. 
For $h\in \h$, we define $\bar X_{h} := \sum_{\nu} \bar{\on{x}}_{[h,h_{\nu}]} 
\bar\partial_{h_{\nu}}\in D(\h)$, where 
$t_{\h} = \sum_{\nu} h_{\nu}\otimes h_{\nu}$. 
We then set $\bar Y_{h}:= \bar X_{h} + \sum_{i=1}^{n} h^{(i)}$. The map 
$\h\to B_{n}$ is a Lie algebra morphism; we denote by $\h^{\on{diag}}$ 
its image.

We denote by ${\cal H}_{n}(\g,\h)$ the Hecke algebra of $B_{n}$ relative to 
$\h^{\on{diag}}$. Explicitly, ${\cal H}_{n}(\g,\h) = \{x\in B_{n} | 
\forall h\in \h, \bar Y_{h}x \in B_{n}\h^{\on{diag}}\} / B_{n}\h^{\on{diag}}$.

\begin{proposition} \label{prop:red:1}
There is a unique Lie algebra morphism 
$$
\rho_{\g,\h} : \bar\t_{1,n} 
\to {\cal H}_{n}(\g,\h), 
$$
such that $\bar x_{i}\mapsto \sum_{\nu}\bar{\on{x}}_{\nu} \otimes 
h_{\nu}^{(i)}$, $\bar y_{i}\mapsto -\sum_{\nu} \bar\partial_{\nu} 
\otimes h_{\nu}^{(i)} + \sum_{j} \sum_{\alpha} \ell_{\alpha} \otimes 
a_{\alpha}^{(i)}b_{\alpha}^{(j)}$, $\bar t_{ij}\mapsto t_{\g}^{(ij)}$. 
Here $r(\lambda) = \sum_{\alpha} 
\ell_{\alpha}(\lambda)(a_{\alpha} \otimes b_{\alpha})$. 
\end{proposition}

If $V_{1},...,V_{n}$ are $\g$-modules, then $S(\h)[1/P] \otimes 
(\otimes_{i}V_{i})$ is a module over $D(\h)[1/P] \otimes 
U(\g)^{\otimes n}$, and $(S(\h)[1/P] \otimes (\otimes_{i}V_{i}))^{\h}$ 
is a module over $H_{n}(\g,\h)$.

Moreover, we have a restriction morphism 
$(S(\g)\otimes (\otimes_{i}V_{i}))^{\g} 
\to (S(\h)[1/P]\otimes (\otimes V_{i}))^{\h}$. Note that 
$(S(\g)\otimes (\otimes_{i}V_{i}))^{\g}$ is a $\bar\t_{1,n} 
$-module
using the morphism $\bar\t_{1,n} 
\to {\cal H}_{n}(\g)$, while 
$(S(\h)[1/P]\otimes (\otimes V_{i}))^{\h}$ is a $\bar\t_{1,n} 
$-module
using the morphism $\bar\t_{1,n} 
\to {\cal H}_{n}(\g,\h)$. Then one checks that 
the restriction morphism $(S(\g)\otimes (\otimes_{i}V_{i}))^{\g} 
\to (S(\h)[1/P]\otimes (\otimes V_{i}))^{\h}$ is a $\bar\t_{1,n}
$-modules morphism.

\medskip

{\em Proof.} The images of the above elements are all $\h$-invariant. 
To lighten the notation, we will imply summation over repeated indices and 
denote elements of $B_n$ as follows: $\bar\partial_\nu\otimes 1$ 
by $\bar\partial_\nu$, $\bar{\on{x}}_\nu\otimes 1$ by 
$\langle \lambda, h_\nu\rangle$, 
$1\otimes x^{(i)}$ by $x^i$. Then $\rho_{\g,\h}(\bar x_i) = (\lambda^\vee)^i$, 
$\rho_{\g,\h}(\bar y_{i}) = - h_{\nu}^{i}\bar\partial_{\nu}
+ \sum_{j=1}^n r(\lambda)^{ij}$ (here for $x\otimes y\in \g^{\otimes 2}$, 
$(x\otimes y)^{ii}:= x^i y^i$).

We will use the same 
presentation of $\bar\t_{1,n}$ as in Proposition \ref{prop:realization}. 
The relations $[\bar x_{i},\bar x_{j}] = 0$ and $\bar t_{ij}=\bar t_{ji}$
are obviously preserved.

Let us check that $[\bar x_{i},\bar y_{j}]=\bar t_{ij}$ is preserved. 
We have for $i\neq j$, $[\rho_{\g,\h}(\bar x_{i}),\rho_{\g,\h}(\bar y_{j})] = 
[\bar x_{\nu}h_{\nu}^{i}, - h_{\nu}^{j}\bar\partial_{\nu}
+ \sum_{k} r(\lambda)^{jk}] = t_{\h}^{ij} + [\lambda^{i},r(\lambda)^{ji}] = 
t_{\h}^{ij} + t_{\n}^{ij} = t_{\g}^{ij} = \rho_{\g,\h}(\bar t_{ij})$.

Let us check that $\sum_{i}\bar x_{i} = \sum_{i}\bar y_{i}=0$ are preserved. 
We have $\sum_{i}\rho_{\g,\h}(\bar x_{i}) = 0$ by the same argument 
as above and $\sum_{i}\rho_{\g,\h}(\bar y_{i}) = \sum_{i} (\lambda^\vee)^{i}$
(by the antisymmetry of $r(\lambda)$), which vanishes by the same argument 
as above.

Let us check that $[\bar y_{i},\bar y_{j}]=0$ is preserved, 
for $i\neq j$. We have 
\begin{align*}
& [\rho_{\g,\h}(\bar y_{i}),\rho_{\g,\h}(\bar y_{j})] \\
& = \sum_{k|k\neq i,j} 
\big( -h_{\nu}^{i} (\partial_{\nu}r(\lambda))^{jk}
+h_{\nu}^{j} (\partial_{\nu}r(\lambda))^{ik}
+ [r(\lambda)^{ij},r(\lambda)^{jk}] + [r(\lambda)^{ik},r(\lambda)^{jk}]
+ [r(\lambda)^{ik},r(\lambda)^{ji}] \big) 
\\
& + [(h_{\nu}^{i} + h_{\nu}^{j})\bar\partial_{\nu},r(\lambda)^{ij}]
- [h_{\nu}^{i}\bar\partial_{\nu},r(\lambda)^{jj}] 
+ [h_{\nu}^{j}\bar\partial_{\nu},r(\lambda)^{ii}] 
+ [r(\lambda)^{ij},r(\lambda)^{ii} + r(\lambda)^{jj}]
\\
  &  = \sum_{k|k\neq i,j} 
h_{\nu}^{k} (\partial_{\nu}r(\lambda))^{ij}  
+ [(h_{\nu}^{i} + h_{\nu}^{j})\bar\partial_{\nu},r(\lambda)^{ij}]
- [h_{\nu}^{i}\bar\partial_{\nu},r(\lambda)^{jj}] 
+ [h_{\nu}^{j}\bar\partial_{\nu},r(\lambda)^{ii}] 
+ [r(\lambda)^{ij},r(\lambda)^{ii} +  r(\lambda)^{jj}]
\\
  & 
\equiv (\partial_{\nu}r(\lambda))^{ij}
(-h_{\nu}^{i} - h_{\nu}^{j} - \bar X_{\nu}) 
+ [(h_{\nu}^{i} + h_{\nu}^{j})\bar\partial_{\nu},r(\lambda)^{ij}]
- h_{\nu}^{i} (\partial_{\nu}r(\lambda))^{jj} 
+ h_{\nu}^{j} (\partial_{\nu}r(\lambda))^{ii} \\
  & + [r(\lambda)^{ij},r(\lambda)^{ii}  +  r(\lambda)^{jj}]
  = [h_{\nu}^{i}+h_{\nu}^{j},r(\lambda)^{ij}]\bar\partial_{\nu} 
- (\partial_{\nu}r^{ij}(\lambda))\bar X_{\nu}
\\
& + [h_{\nu}^{i}+h_{\nu}^{j},\partial_{\nu}r(\lambda)^{ij}]
- h^{i}_{\nu}(\partial_{\nu}r(\lambda))^{jj} 
+ h^{j}_{\nu}(\partial_{\nu}r(\lambda))^{ii} 
+ [r(\lambda)^{ij},r(\lambda)^{ii}+r(\lambda)^{jj}]. 
\end{align*}
The second equality follows from the CDYBE and the antisymmetry on 
$r(\lambda)$. Then 
$$
[h_{\nu}^{i}+h_{\nu}^{j},r(\lambda)^{ij}]\bar\partial_{\nu} 
- (\partial_{\nu}r^{ij}(\lambda))\bar X_{\nu}
= \big( [h_{\nu'}^{i}+h_{\nu'}^{j},r(\lambda)^{ij}]
- \partial_{\nu}r^{ij}(\lambda) \langle \lambda , [h_{\nu},h_{\nu'}] \rangle
  \big)\bar\partial_{\nu'} = 0
$$
using the $\h$-invariance of $r(\lambda)$. Applying 
$x^{i}y^{j}z^{k}\mapsto x^{i}(yz)^{i}$ to the CDYB identity 
$$
[r(\lambda)^{ij},r(\lambda)^{ik}]+[r(\lambda)^{ij},r(\lambda)^{jk}]
+ [r(\lambda)^{ik},r(\lambda)^{jk}]-h_{\nu}^{i}\partial_{\nu}r(\lambda)^{jk}
+ h_{\nu}^{j}\partial_{\nu}r(\lambda)^{ik} 
- h_{\nu}^{j}\partial_{\nu}r(\lambda)^{ij}=0, 
$$
we get 
$$
(1/2)\sum_{\alpha,\beta} \ell_{\alpha}\ell'_{\beta}(\lambda) 
[a_{\alpha},a_{\beta}]^{i} [b_{\alpha},b_{\beta}]^{j} 
+ [r(\lambda)^{ij},r(\lambda)^{ii}] - h_{\nu}^{i}(\partial_{\nu}r(\lambda))^{jj}
+ [h_{\nu}^{j},\partial_{\nu}r(\lambda)^{ij}]=0. 
$$
Since $r(\lambda)$ is antisymmetric, the sum $(1/2)\sum_{\alpha,\beta}...$ 
is symmetric in $(i,j)$; antisymmetrizing in $(i,j)$, we get 
$$
[h_{\nu}^{i}+h_{\nu}^{j},\partial_{\nu}r(\lambda)^{ij}]
- h^{i}_{\nu}(\partial_{\nu}r(\lambda))^{jj} 
+ h^{j}_{\nu}(\partial_{\nu}r(\lambda))^{ii} 
+ [r(\lambda)^{ij},r(\lambda)^{ii}+r(\lambda)^{jj}] = 0. 
$$
All this implies that $[\rho_{\g,\h}(\bar y_{i}),\rho_{\g,\h}(\bar y_{j})] = 0$.

Let us check that $[\bar x_{i},\bar t_{jk}]=0$ is preserved ($i,j,k$ distinct). 
We have $[\rho_{\g,\h}(\bar x_{i}),\rho_{\g,\h}(\bar t_{jk})] 
= [(\lambda^\vee)^i,t_{\g}^{jk}]=0$.

Let us prove that $[\bar y_{i},\bar t_{jk}] = 0$ is preserved 
($i,j,k$ distinct). 
We have $[\rho_{\g,\h}(\bar y_{i}),\rho_{\g,\h}(\bar t_{jk})] = 
[-h_{\nu}^{i}\bar\partial_{\nu}+\sum_{l}r(\lambda)^{il},t_{\g}^{jk}] 
= [r(\lambda)^{ij}
+r(\lambda)^{ik},t_{\g}^{jk}]=0$ because $t_{\g}$ is $\g$-invariant. 
\hfill \qed \medskip

\begin{proposition}
If $V_1,...,V_n$ are $\g$-modules, then 
$(S(\h)[1/P]\otimes (\otimes_i V_i))^\h$ is a $\bar\t_{1,n}\rtimes \d$-module. 
The $\bar\t_{1,n}$-module structure is induced by the morphism 
$\bar\t_{1,n}\to {\cal H}_n(\g,\h)$ of Proposition \ref{prop:red:1}, so 
$$
\rho_{(V_i)}(\bar x_i)(f(\lambda)\otimes (\otimes_i v_i)) = 
(\lambda^\vee)^i (f(\lambda) \otimes (\otimes_i v_i)), 
$$ 
$$
\rho_{(V_i)}(\bar y_i)(f(\lambda)\otimes (\otimes_i v_i)) = 
(-h_\nu^i\partial_\nu + \sum_{j}r(\lambda)^{ij})
(f(\lambda) \otimes (\otimes_i v_i)), 
$$ 
$$
\rho_{(V_i)}(\bar t_{ij})(f(\lambda)\otimes (\otimes_i v_i)) = 
t_\g^{ij}(f(\lambda) \otimes (\otimes_i v_i)), 
$$
and the $\d$-module structure is given by 
$$
\rho_{(V_i)}(\delta_{2m})(f(\lambda)\otimes (\otimes_i v_i))
= {1\over 2}(\sum_i \{(\on{ad}\lambda^\vee)^{2m}(e_\alpha)\cdot e_\alpha\}^i)
(f(\lambda) \otimes (\otimes_i v_i)), 
$$
\begin{align*}
& \rho_{(V_i)}(\Delta_0)(f(\lambda)\otimes (\otimes_i v_i))
\\
  & = \Big( - {1\over 2}\partial_\nu^2  + {1\over 2} \langle \mu(r(\lambda)), 
h_{\nu}\rangle \partial_{\nu}
+ \{ {1\over 2}\psi(\lambda)^{11} 
- {1\over 2} (\on{ad}\lambda^{\vee})_{|\n}^{-1}
(\mu(r(\lambda))_{\n}) \}^{12...n} \Big)
(f(\lambda) \otimes (\otimes_i v_i)),
\end{align*}
$$
\rho_{(V_i)}(d)(f(\lambda)\otimes (\otimes_i v_i))
=  {1\over 2}(\langle \lambda, h_\nu\rangle \partial_{\nu} +  \partial_{\nu}
\langle \lambda, h_{\nu}\rangle + \langle \mu(r(\lambda)), 
\lambda^{\vee}\rangle) (f(\lambda) \otimes (\otimes_i v_i)), 
$$
$$
\rho_{(V_i)}(X)(f(\lambda)\otimes (\otimes_i v_i))
= (1/2)\langle \lambda^\vee, \lambda^\vee \rangle 
(f(\lambda) \otimes (\otimes_i v_i)). 
$$
Here $x_{\n}$ is the projection of $x\in\g$ on $\n$ along $\h$. 
\end{proposition}

To summarize, we have a diagram 
$$
\begin{matrix}
\bar\t_{1,n} & \to & {\cal H}_n(\g,\h) & \to & \on{End}(
(S(\h)[1/P]\otimes (\otimes_i V_i))^\h) \\
  & \scriptstyle{\subset}\searrow
  & \scriptstyle{(1)}\uparrow & \nearrow  & \\
  & & \bar\t_{1,n}\rtimes \d & & 
\end{matrix}
$$
As before, the restriction morphism $(S(\g)\otimes(\otimes_i V_i))^\g\to 
(S(\h)[1/P]\otimes (\otimes_i V_i))^\h$ 
extends to a $\bar\t_{1,n} \rtimes \d$-modules morphism.

The action of $\bar\t_{1,n}\rtimes \d$ factors through a morphism 
$\tilde\rho_{\g,\h} : \bar\t_{1,n}\rtimes\d\to {\cal H}_n(\g,\h)$
extending $\rho_{\g,\h} : \bar\t_{1,n}\to {\cal H}_n(\g,\h)$ 
(denoted by (1) in the diagram). 

\medskip

{\em Proof.} Let $\lambda\in \h^*_{\on{reg}}$. Then if $V$ is a $\g$-module, 
we have 
$(\hat\cO_{\g^*,\lambda}\otimes V)^\g = (\hat\cO_{\h^*,\lambda} \otimes V)^\h$
(where $\hat\cO_{X,x}$ is the completed local ring of a variety $X$ at 
the point $x$). We then have a morphism 
$\bar\t_{1,n}\rtimes \d \to {\cal H}_n(\g) \to 
\on{End}((\hat\cO_{\g^*,\lambda}\otimes 
(\otimes_i V_i))^\g)$ for any $\lambda\in\g^*$, so when $\lambda\in
\h^*_{\on{reg}}$
we get a morphism $\bar\t_{1,n}\rtimes \d \to 
\on{End}((\hat\cO_{\h^*,\lambda}\otimes 
(\otimes_i V_i))^\h)$.

Let show that the images of the generators of $\bar\t_{1,n}\rtimes\d$
under this morphism are given by the above formulas.

Since the actions of $\bar x_{i}$, $\bar t_{ij}$ and $X$ on 
$(\hat\cO_{\g^{*},\lambda}\otimes (\otimes_{i}V_{i}))^{\g}$
are given by multiplication by elements of 
$(\hat\cO_{\g^{*},\lambda}\otimes U(\g)^{\otimes n})^{\g}$, 
their actions on 
$(\hat\cO_{\h^{*},\lambda}\otimes (\otimes_{i}V_{i}))^{\h}$
are given by multiplication by restrictions of these elements to $\h^{*}$.

Let us compute the action of $\bar y_{i}$. Let $\tilde f(\lambda)
\in (\hat\cO_{\h^{*},\lambda}\otimes (\otimes_{i}V_{i}))^{\h}$ and 
$\tilde F(\lambda)
\in (\hat\cO_{\g^{*},\lambda}\otimes (\otimes_{i}V_{i}))^{\g}$
be its equivariant extension to a formal map $\g^{*}\to \otimes_{i}V_{i}$. 
Then for $x\in\n$, we have $(\partial_{x^{\wedge}} 
+ \sum_{i} (\on{ad}\lambda^{\vee})^{-1}(x)^{i})
(\tilde F(\lambda))_{|\h^{*}} = 0$ (the map $x\mapsto x^{\wedge}$ is the 
inverse of $\g^*\to\g$, $\lambda\mapsto \lambda^{\vee}$). 
Then $\rho_{(V_{i})}(\bar y_{i})(\tilde f(\lambda)) = 
\Big(-h^{i}_{\nu}\partial_{\nu} + \sum_{j} e^{i}_{\beta}
\big( (\on{ad}\lambda^{\vee})^{-1}(e_{\beta})\big)^{j}\Big)\tilde f(\lambda)
= (-h^{i}_{\nu}\partial_{\nu} + \sum_{j}r(\lambda)^{ij})(\tilde f(\lambda))$.

Let us now compute the action of $\Delta_{0}$. Let $\lambda_{0}\in \h^{*}$
be such that $\lambda_{0}^{\vee}\in U$ and $\lambda\in \g^{*}$ be close to 
$\lambda_{0}$. We set $\delta\lambda := \lambda - \lambda_{0}$. 
We then have $\lambda = e^{\on{ad}x}(\lambda_{0} + h^{\wedge})$, where 
$x\in\n$ and $h\in \h$ are close to $0$. 
We have the expansions 
$$
h = (\delta\lambda)_{\h}^{\vee} + {1\over 2}
[(\on{ad}\lambda_{0}^{\vee})_{|\n}^{-1}((\delta\lambda)_{\n}^{\vee}),
(\delta\lambda)_{\n}^{\vee}]_{\h}, 
$$
$$
x = - (\on{ad}\lambda_{0}^{\vee})_{|\n}^{-1}\Big(
(\delta\lambda)_{\n}^{\vee} 
+ [(\on{ad}\lambda_{0}^{\vee})_{|\n}^{-1}((\delta\lambda)_{\n}^{\vee}),
(\delta\lambda)_{\h}^{\vee}] 
+ {1\over 2}[(\on{ad}\lambda_{0}^{\vee})_{|\n}^{-1}
((\delta\lambda)_{\n}^{\vee}), (\delta\lambda)_{\n}^{\vee}]_{\n}\Big)
$$ 
up to terms of order $>2$; here the indices $u_{\n}$ and $u_{\h}$ mean the 
projections of $u\in\g$ to $\n$ and $\h$. If now $\tilde f(\lambda) : \h^{*}
\supset V(\lambda_{0},\h^{*})\to \otimes_{i}V_{i}$ is an $\h$-equivariant 
function defined at the vicinity of 
$\lambda_{0}$ and $\tilde F(\lambda) : \g^{*}\supset V(\lambda_{0},\g^{*})
\to \otimes_{i}V_{i}$ it its $\g$-equivariant extension to a neighborhood 
of $\lambda_{0}$ in $\g^{*}$, then 
$\tilde F(\lambda) = (e^{x})^{1...n}\tilde f(\lambda_{0} + h)$, 
which implies the expansion 
\begin{align*}
& \tilde F(\lambda) = \tilde f(\lambda_{0}) 
+  \Big( (\delta\lambda)_{\nu} + {1\over 2} \langle 
[(\on{ad}\lambda_{0}^{\vee})_{|\n}^{-1}(e_{\beta}),
e_{\beta'}], h_{\nu}\rangle (\delta\lambda)_{\beta}
(\delta\lambda)_{\beta'}\Big) \partial_{\nu}\tilde f(\lambda_{0}) 
+ {1\over 2} (\delta\lambda)_{\nu}(\delta\lambda)_{\nu'}
\partial^{2}_{\nu\nu'}\tilde f(\lambda_{0})
\\
  & + \Big( - (\on{ad}\lambda_{0}^{\vee})_{|\n}^{-1}(e_{\beta})
  (\delta\lambda)_{\beta}
  - (\on{ad}\lambda_{0}^{\vee})^{-1}([(\on{ad}\lambda_{0}^{\vee})_{|\n}^{-1}
  (e_{\beta}),h_{\nu}])(\delta\lambda)_{\nu} (\delta\lambda)_{\beta}
  \\
   & - {1\over 2}(\on{ad}\lambda_{0}^{\vee})_{|\n}^{-1}
   ([(\on{ad}\lambda_{0}^{\vee})_{|\n}^{-1}
  (e_{\beta}),e_{\beta'}]_{\n})(\delta\lambda)_{\beta}
  (\delta\lambda)_{\beta'}
  + {1\over 2} (\on{ad}\lambda_{0}^{\vee})_{|\n}^{-1}(e_{\beta})
  (\on{ad}\lambda_{0}^{\vee})_{|\n}^{-1}(e_{\beta'}) (\delta\lambda)_{\beta}
(\delta\lambda)_{\beta'}
   \Big)^{1...n}\tilde f(\lambda_{0}) \\
  & 
-(\on{ad}\lambda_{0}^{\vee})_{|\n}^{-1}(e_{\beta})^{1...n}
  (\delta\lambda)_{\beta}(\delta\lambda)_{\nu}
  \partial_{\nu}\tilde f(\lambda_{0})
\end{align*}
up to terms of order $>2$.

Then 
\begin{align*}
& (\partial_{\alpha}^{2}F)(\lambda_{0}) = 
(\partial_{\nu}^{2}\tilde f)(\lambda_{0})
+ \langle [(\on{ad}\lambda_{0}^{\vee})_{|\n}^{-1}
(e_{\beta}),e_{\beta}], h_{\nu}\rangle \partial_{\nu}\tilde f(\lambda_{0})
\\
& + \Big( - (\on{ad}\lambda_{0}^{\vee})_{|\n}^{-1}
([(\on{ad}\lambda_{0}^{\vee})_{|\n}^{-1} (e_{\beta}),e_{\beta}]_{\n}) 
+ ((\on{ad}\lambda_{0}^{\vee})_{|\n}^{-1}(e_{\beta}))^{2}
\Big)^{1...n}\tilde f(\lambda_{0}), 
\end{align*}
which implies the formula for the action of $\Delta_{0}$.

Then $(S(\h)[1/P] \otimes (\otimes_i V_i))^\h \subset 
\prod_{\lambda\in \h^*_{\on{reg}}} (\hat\cO_{\h^*,\lambda}\otimes 
(\otimes_i V_i))^\h$
is preserved by the action of the generators of 
$\bar\t_{1,n}\rtimes\d$-module, hence it is a
sub-$(\bar\t_{1,n}\rtimes\d)$-module, with action given by the 
above formulas. 
\hfill \qed \medskip

\subsection{Realization of the universal KZB system}

The realization of the flat connection $\on{d} - \sum_i \bar K_i(\zz|\tau)
\on{d}z_i - \bar\Delta(\zz|\tau)\on{d}\tau$ on $(\HH\times \CC^n)-\on{Diag}_n$
is a flat connection on the trivial bundle with fiber $(\cO_{\h^*_{\on{reg}}}
\otimes (\otimes_i V_i))^\h$.

We now compute this realization, under the assumption that 
$\h\subset \g$ is a maximal abelian subalgebra. In this case, 
two simplifications occur:

(a) $(\on{ad}\lambda^\vee)(h_\nu) = 0$ since $\h$ is abelian,

(b) $[(\on{ad}\lambda^\vee)_{|\n}^{-1}(e_\beta),e_\beta]_{\n} = 0$
since $[(\on{ad}\lambda^\vee)_{|\n}^{-1}(e_\beta),e_\beta]$
commutes with any element in $\h$, so that it belongs to $\h$.

The image of $\bar K_{i}(\zz|\tau)$ is then the operator 
\begin{align*}
& K_i^{(V_i)}(\zz|\tau) = h^{i}_{\nu}\partial_{\nu} - \sum_{j} r(\lambda)^{ij} 
+ \sum_{j|j\neq i} k(z_{ij},(\on{ad}\lambda^{\vee})^{i}|\tau)(t_{\n}^{ij}
+t_{\h}^{ij})
\\  & = 
h^{i}_{\nu}\partial_{\nu} - r(\lambda)^{ii} 
+ \sum_{j|j\neq i} {{\theta(z_{ij} + (\on{ad}\lambda^{\vee})^{i}|\tau)}
\over{\theta(z_{ij}|\tau)\theta((\on{ad}\lambda^{\vee})^{i}|\tau)}}(t_{\n}^{ij})
+ \sum_{j|j\neq i} {{\theta'}\over\theta}(z_{ij}|\tau)t_{\h}^{ij}
\end{align*}

The image of $2\pi\i\bar\Delta(\zz|\tau)$ is the operator 
\begin{align*}
& 2\pi\i\Delta^{(V_i)}(\zz|\tau) 
=  {1\over 2}\partial_\nu^2 + {1\over 2} \langle [(\on{ad}\lambda^\vee)^{-1}
(e_\beta),e_\beta], h_\nu\rangle \partial_\nu 
- g(0,0|\tau)\sum_i{1\over 2}t_\g^{ii}
\\ & 
+ \sum_{i,j} {1\over 2} \big( [g(z_{ij},\on{ad}\lambda^{\vee}|\tau) -
(\on{ad}\lambda^{\vee})^{-2}](e_\beta)\big)^i e_\beta^j 
+ \sum_{i,j} {1\over 2} g(z_{ij},0|\tau)h_\nu^i h_\nu^j 
\end{align*}
and the connection is now
$$
\nabla^{(V_i)} = \on{d} - \sum_i K_i^{(V_i)}(\zz|\tau) \on{d}z_i 
- \Delta^{(V_i)}(\zz|\tau) \on{d}\tau. 
$$
Recall that $P(\lambda) = \on{det}((\on{ad}\lambda^\vee)_{|\n})$. 
We compute the conjugation $P^{1/2}\nabla^{(V_i)} P^{-1/2}$, where 
$P^{\pm 1/2}$ is the operator of multiplication by (inverse branches of) 
$P^{\pm 1/2}$ on $\cO_{\h^*_{\on{reg}}} \otimes (\otimes_i V_i)^\h$.

\begin{lemma}
$\partial_\nu\on{log}P(\lambda) = - \langle h_\nu, \mu(r(\lambda))\rangle$, 
$P^{1/2}[h^i_\nu\partial_\nu - r(\lambda)^{ii}]P^{-1/2} = 
h^i_\nu\partial_\nu$, 
$P^{1/2}[\partial_\nu^2 + \langle [(\on{ad}\lambda^\vee)^{-1}_{|\n}
(e_\beta),e_\beta], h_\nu \rangle \partial_\nu] P^{-1/2} = \partial_\nu^2 
+ \partial_\nu\big(\langle h_\nu, {1\over 2}\mu(r(\lambda))\rangle \big) 
- \langle h_\nu,{1\over 2}\mu(r(\lambda)) \rangle^2$. 
\end{lemma}

{\em Proof.} $\partial_\nu \on{log}P(\lambda) = 
(d/dt)_{|t=0}\on{det}[(\on{ad}(\lambda^\vee + t h_\nu)_{|\n})
(\on{ad}\lambda^\vee)_{|\n}^{-1}] = \on{tr} 
[(\on{ad}h_\nu)_{|\n} \circ (\on{ad}\lambda^\vee)^{-1}_{|\n}] = 
\langle e_\beta, (\on{ad}h_\nu) \circ (\on{ad}\lambda^\vee)_{|\n}^{-1}
(e_\beta)\rangle  = \langle [(\on{ad}\lambda^\vee)_{|\n}^{-1}(e_\beta),e_\beta]
, h_\nu\rangle  = - \langle h_\nu, \mu(r(\lambda))\rangle$. 
The next equality follows from $\mu(r(\lambda))^i = 2r(\lambda)^{ii}$. 
The last equality is a direct consequence. 
\hfill \qed \medskip

We then get:

\begin{proposition}
$P^{1/2}\nabla^{(V_i)}P^{-1/2} = \on{d} - \sum_i \tilde
K_i(\zz|\tau)\on{d}z_i - \tilde\Delta(\zz|\tau)\on{d}\tau$, where 
$$
\tilde K_i(\zz|\tau) = h^i_\nu\partial_\nu 
+ \sum_{j|j\neq i} {{\theta(z_{ij} + (\on{ad}\lambda^{\vee})^{i}|\tau)}
\over{\theta(z_{ij}|\tau)\theta((\on{ad}\lambda^{\vee})^{i}|\tau)}}
(t_{\n}^{ij}) + 
\sum_{j|j\neq i} {{\theta'}\over\theta}(z_{ij}|\tau)t_{\h}^{ij}
$$
\begin{align*}
& 2\pi\i\tilde\Delta(\zz|\tau) = {1\over 2}\partial_\nu^2 
+ \partial_\nu\big( \langle h_\nu, {1\over 2}\mu(r(\lambda))\rangle \big) 
- \langle h_\nu,{1\over 2}\mu(r(\lambda)) \rangle^2
- g(0,0|\tau)\sum_i{1\over 2}t_\g^{ii}
\\ & 
+ \sum_{i,j} {1\over 2} \Big( \big( g(z_{ij},\on{ad}\lambda^{\vee}|\tau) -
(\on{ad}\lambda^{\vee})^{-2} \big) (e_\beta)\Big)^i e_\beta^j 
+ \sum_{i,j} {1\over 2} g(z_{ij},0|\tau)h_\nu^i h_\nu^j, 
\end{align*}
where 
$$
g(z,0|\tau) = {1\over 2}{{\theta''}\over\theta}(z|\tau) - 2\pi\i
{{\partial_\tau\eta}\over\eta}(\tau) 
$$
and 
$$
g(z,\alpha|\tau) - \alpha^{-2} = {1\over 2}
{{\theta(z+\alpha|\tau)}\over{\theta(x|\tau)\theta(\alpha|\tau)}}
({\theta'\over\theta}(z+\alpha|\tau) - {\theta'\over\theta}(\alpha|\tau))
$$
\end{proposition}

The term in $\sum_i (1/2)t_\g^{ii}$ is central and can be absorbed 
by a suitable further conjugation. 
Rescaling $t_{\g}$ into $\kappa^{-1} t_{\g}$, where 
$\kappa\in\CC^{\times}$, $\tilde K_i(\zz|\tau)$ and 
$\tilde\Delta(\zz|\tau)$ get multiplied by $\kappa$. 
Moreover, we have:

\begin{lemma}
When $\g$ is simple and $\h\subset\g$ is the 
Cartan subalgebra,
$ \partial_\nu\{\langle h_\nu, {1\over 2}\mu(r(\lambda))\rangle \} 
= \langle h_\nu,{1\over 2}\mu(r(\lambda)) \rangle^2$. 
\end{lemma}

{\em Proof.} Let $D(\lambda) := \prod_{\alpha\in\Delta^+}(\alpha,\lambda)$, 
where $\Delta^+$ is the set of positive roots of $\g$. Then 
$D(\lambda)$ is $W$-antiinvariant, where $W$ is the Weyl group. 
Therefore $\partial_\nu^2 D(\lambda)$ is also $W$-antiinvariant, 
so it is divisible (as a polynomial on $\h^*$) by all the 
$(\alpha,\lambda)$, where $\alpha\in\Delta^+$, so it is divisible by 
$D(\lambda)$; since $\partial_\nu^2 D(\lambda)$ has degree strictly 
lower than $D(\lambda)$, we get $\partial_\nu^2 D(\lambda)=0$.

Now if $(e_\alpha,f_\alpha,h_\alpha)$ is a basis of the $\sl_2$-triple
associated with $\alpha$, we have $r(\lambda) = \sum_{\alpha\in\Delta^+}
-(e_\alpha \otimes f_\alpha - f_\alpha \otimes e_\alpha)/(\alpha,\lambda)$, 
so ${1\over 2}\mu(r(\lambda)) = -\sum_{\alpha\in\Delta^+}
h_\alpha/(\alpha,\lambda)$. Therefore ${1\over 2}\mu(r(\lambda))
= -\partial_\nu \on{log}D(\lambda) h_\nu$. 
Then $\partial_\nu^2 D(\lambda)=0$ implies that 
$\partial_\nu^2 \on{log}D + (\partial_\nu\on{log}D)^2 = 0$, which 
implies the lemma. 
\hfill \qed \medskip

The resulting flat connection
then coincides with that of \cite{Be1,FW}.

\section{The universal KZB connection and representations of
Cherednik algebras}

\subsection{The rational Cherednik algebra of type $A_{n-1}$}

Let $k$ be a complex number, and $n\ge 1$ an integer.
The rational Cherednik algebra $H_n(k)$ of type $A_{n-1}$ 
is the quotient of the algebra ${\CC}[S_n]\ltimes 
{\CC}[{\rm x}_1,...,{\rm x}_n,{\rm y}_1,...,{\rm y}_n]$ 
by the relations
$$
\sum_i {\rm x}_i=0,\ \sum_i {\rm y}_i=0,\ [{\rm x}_i,{\rm x}_j]=0=
[{\rm y}_i,{\rm y}_j],\
$$
$$
[{\rm x}_i,{\rm y}_j]=\frac{1}{n}-ks_{ij},\ i\ne j,
$$
where $s_{ij}\in S_n$ is the permutation of $i$ and $j$
(see e.g. \cite{EG}). 
\footnote{The generators ${\rm x}_\alpha,\partial_\alpha$ of Section 
\ref{sect:real:1n} will be henceforth renamed $q_\alpha,p_\alpha$.}

Let $e := {1\over {n!}}\sum_{\sigma\in S_n} \sigma\in {\CC}[S_n]$ 
be the Young symmetrizer.
The spherical subalgebra $B_n(k)$ (often called the spherical
Cherednik algebra) is defined to be the algebra 
$e H_n(k)e$.

We define an important element 
$$
{\bold h}:=\frac{1}{2}\sum_i ({\rm x}_i{\rm y}_i+{\rm y}_i{\rm
x}_i).
$$
We recall that category
${\mathcal O}$ is the category of $H_n(k)$-modules which are
locally nilpotent under the action of the operators ${\rm y}_i$ and decompose
into a direct sum of finite dimensional generalized eigenspaces of
${\bold h}$. Similarly, one defines category ${\mathcal O}$ 
over $B_n(k)$ to be the category of $B_n(k)$-modules which are
locally nilpotent under the action of ${\CC}[{\rm y}_1,...,{\rm y}_n]^{S_n}$
and decompose into a direct sum of finite dimensional 
generalized eigenspaces of ${\bold h}$.

\subsection{The homorphism from $\bar{\mathfrak t}_{1,n}$
to the rational Cherednik algebra}

\begin{proposition}\label{mapxi}
For each $k,a,b\in {\CC}$, we have a homomorphism of Lie algebras 
$\xi_{a,b}: \bar{\mathfrak t}_{1,n}\to H_n(k)$,
defined by the formula 
$$
\bar x_i\mapsto a{\rm x}_i, \quad \bar y_i\mapsto b{\rm y}_i, \quad 
\bar t_{ij}\mapsto ab\left(\frac{1}{n}-ks_{ij}\right).
$$ 
\end{proposition}

{\em Proof.} 
Straightforward. 
\hfill \qed \medskip

\begin{remark} Obviously, $a,b$ can be rescaled independently,
by rescaling the generators $\bar x_i$ and $\bar y_i$ of the source 
algebra
$\bar{\mathfrak t}_{1,n}$. On the other hand, if we are only allowed to 
apply
automorphisms of the target algebra $H_n(k)$, then 
$a,b$ can only be rescaled in such a way that the product $ab$ is
preserved. \hfill \qed \medskip 
\end{remark}

This shows that any representation $V$ of the rational Cherednik
algebra $H_n(k)$ yields a family of realizations for $\bar{\mathfrak
t}_{1,n}$ parametrized by $a,b\in {\CC}$, 
and gives rise to a family of flat connections $\nabla_{a,b}$ 
over the configuration
space $\bar C(E_\tau,n)$.

\subsection{Monodromy representations of double affine Hecke algebras}

Let ${\mathcal H}_n(q,t)$ be Cherednik's double affine Hecke algebra of 
type
$A_{n-1}$. By definition, ${\mathcal H}_n(q,t)$ 
is the quotient of the group algebra of the 
orbifold fundamental group $\overline{{\rm B}}_{1,n}$ 
of $\bar C(E_\tau,n)/S_n$ by the additional relations
$$
(T-q^{-1}t)(T+q^{-1}t^{-1})=0, 
$$
where $T$ is any element of $\overline{{\rm B}}_{1,n}$ homotopic (as a free
loop) to a small loop around the divisor of diagonals 
in the counterclockwise direction.

Let $V$ be a representation of $H_n(k)$, and let 
$\nabla_{a,b}(V)$ be the universal connection 
$\nabla_{a,b}$ evaluated in $V$. 
In some cases, for example if $a,b$ are formal, 
or if $V$ is finite dimensional, we can
consider the monodromy of this connection, which obviously gives 
a representation of ${\mathcal H}_n(q,t)$ on $V$, with 
$$
q=e^{-2\pi {\rm i} ab/n},\ t=e^{-2\pi {\rm i}kab}.
$$ 
In particular, taking $a=b$, $V=H_n(k)$, this monodromy
representation defines an homomorphism 
$\theta_a: {\mathcal H}_n(q,t)\to H_n(k)[[a]]$, 
where 
$$
q=e^{-2\pi {\rm i} a^2/n},\ t=e^{-2\pi {\rm i}ka^2}.
$$ 
It is easy to check that this homomorphism becomes an isomorphism
upon inverting $a$. 
The existence of such an isomorphism was pointed out 
by Cherednik (see \cite{Ch2}, end of Section 6, and 
the end of \cite{Ch1}), but his proof is different.

\begin{example} Let $k=r/n$, where $r$ is an integer relatively
prime to $n$. In this case, it is known (see e.g. \cite{BEG1}) 
that the algebra $H_n(k)$ admits an irreducible finite
dimensional representation $Y(r,n)$ of dimension $r^{n-1}$. 
By virtue of the above construction, the space $Y(r,n)$ carries 
an action of ${\mathcal H}_n(q,t)$ with any nonzero $q,t$ such
that $q^r=t$. This finite dimensional representation of 
${\mathcal H}_n(q,t)$ is irreducible for generic $q$, 
and is called a perfect representation; it was first constructed
in \cite{E}, p. 500, and later in \cite{Ch2}, Theorem 6.5,
in a greater generality. 
\end{example}

\subsection{The modular extension of $\xi_{a,b}$.}

Assume that $a,b\ne 0$.

\begin{proposition} The homomorphism $\xi_{a,b}$ can be extended 
to the algebra $U(\bar {\mathfrak t}_{1,n}\rtimes
{\mathfrak d})\rtimes S_n$ by the formulas
$$
\xi_{a,b}(s_{ij})=s_{ij},
$$ 
$$
\xi_{a,b}(d)={\bold h}=\frac{1}{2}\sum_i ({\rm x}_i{\rm y}_i+{\rm y}_i{\rm 
x}_i), \quad 
\xi_{a,b}(X)=-\frac{1}{2}ab^{-1}\sum_i {\rm x}_i^2,
$$
$$
\xi_{a,b}(\Delta_0)=\frac{1}{2}ba^{-1}\sum_i {\rm y}_i^2,
\quad 
\xi_{a,b}(\delta_{2m})=-\frac{1}{2}
a^{2m-1}b^{-1}\sum_{i<j}({\rm x}_i-{\rm x}_j)^{2m}.
$$
\end{proposition}

{\em Proof.} 
Direct computation.
\hfill \qed \medskip

Thus, the flat connections $\nabla_{a,b}$ extend to flat connections
on ${\mathcal M}_{1,[n]}$.

This shows that the monodromy representation 
of the connection $\nabla_{a,b}(V)$, 
when it can be defined, is a representation
of the double affine Hecke algebra 
${\mathcal H}_n(q,t)$ with a compatible action 
of the extended modular group 
$\widetilde{{\rm SL}_2(\ZZ)}$. In particular, this is the case 
if $V=Y(r,n)$. Such representations of 
$\widetilde{{\rm SL}_2(\ZZ)}$ were considered by
Cherednik, \cite{Ch2}. The element $T$ of 
$\widetilde{{\rm SL}_2(\ZZ)}$ acts in
this representation by ``the Gaussian'', and the element $S$ by 
the ``Fourier-Cherednik transform''. They are generalizations
of the $\widetilde{{\rm SL}_2(\ZZ)}$-action on Verlinde algebras.

\section{Explicit realizations of certain highest 
weight representations of the rational
Cherednik algebra of type $A_{n-1}$}

\subsection{The representation $V_N$.}

Let $N$ be a divisor of $n$, and $\g={\mathfrak {sl}}_N({\CC})$,
$G={\rm SL}_N({\CC})$. Let $V_N=({\CC}[\g]\otimes({\CC}^N)^{\otimes n})^\g$
(the divisor condition is needed for this space to be nonzero). 
It turns out that $V_N$ has a natural structure of a
representation of $H_n(k)$ for $k=N/n$.

\begin{proposition}\label{cheract}
We have a homomorphism 
$\zeta_N: H_n(N/n)\to {\rm End}(V_N)$, defined by the formulas
$$
\zeta_N(s_{ij})=s_{ij}, \quad \zeta_N({\rm x}_i)=X_i, 
\quad \zeta_N({\rm y}_i)=Y_i, \quad (i=1,...,n)
$$
where for $f\in V_N, A\in \g$ we have
$$
(X_if)(A)=A_i f(A),
$$
$$
(Y_if)(A)=\frac{N}{n}\sum_p (b_p)_i \frac{\partial f}{\partial b_p}(A),
$$
where $\lbrace{ b_p\rbrace}$ is an orthonormal basis of $\g$ with respect 
to the trace form. 
\end{proposition}

{\em Proof.} 
Straightforward verification.
\hfill \qed \medskip

The relationship of the representation $V_N$ to other results in
this paper is described by the following proposition.

\begin{proposition}\label{classkzb} The connection 
$\nabla_{a,1}(V_N)$ corresponding to the representation $V_N$ 
is the usual KZB connection 
for the $n$-point correlation functions on the elliptic curve for 
the Lie algebra ${\mathfrak {sl}}_N$ and $n$ copies of the vector
representation ${\CC}^N$, at level $K=-\frac{n}{aN}-N$. 
\end{proposition}

{\em Proof.} 
We have a sequence of maps 
$$
U(\bar {\mathfrak t}_{1,n}\rtimes
{\mathfrak d})\rtimes S_n \to H_n(N/n) \to  {\cal H}_n(\g)\rtimes S_n
\to {\rm End}(V_N),
$$
where the first map is $\xi_{a,b}$, the second map sends $s_{ij}$
to $s_{ij}$, ${\rm x}_i$ to the class of $\sum_\alpha q_\alpha 
\otimes e_\alpha^i$,
and ${\rm y}_i$ to the class of $\sum_\alpha p_\alpha \otimes 
e_\alpha^i$ (recall that the ${\rm x}_a, \partial_a$ of Section 
\ref{sect:real:1n} have been renamed $q_a,p_a$),
and the last map is explained in Section \ref{sect:real:1n}.
The composition of the two first maps is
then that of  Proposition \ref{real:der}, and 
the composition of the two last maps is the map $\zeta_N$ 
of Proposition \ref{cheract}. This implies the statement. 
\hfill \qed \medskip

\begin{remark} Suppose that $K$ is a nonnegative integer, i.e. 
$a=-\frac{n}{N(K+N)}$, where $K\in \ZZ_+$. Then the connection
$\nabla_{a,1}$ on the infinite dimensional vector bundle with
fiber $V_N$  preserves a finite dimensional subbundle of conformal
blocks for the WZW model at level $K$. Th subbundle 
gives rise to a finite dimensional monodromy representation 
$V_N^K$ of the Cherednik algebra ${\mathcal H}_n(q,t)$ 
with 
$$
q=e^{\frac{2\pi {\rm i}}{N(K+N)}}, t=q^N,
$$
(so both parameters are roots of unity). 
The dimension of $V_N^K$ is given by the Verlinde formula, 
and it carries a compatible action of $\widetilde{{\rm SL}_2(\ZZ)}$ 
to the action of the Cherednik algebra. 
Representations of this type were studied by Cherednik in \cite{Ch2}. 
\end{remark}

\subsection{The spherical part of $V_N$.}

Note that 
\begin{equation}\label{xi}
((\sum_{i=1}^n X_i^p)f)(A)=\frac{n}{N}(\tr A^p)f(A),
\end{equation}
\begin{equation}\label{yi}
((\sum_{i=1}^n Y_i^p)f)(A)=
\left(\frac{N}{n}\right)^{p-1}(\tr \partial_A^p) f(A)
\end{equation}

Consider the space $U_N=eV_N=({\CC}[\g]\otimes S^n{\CC}^N)^\g$
as a module over the spherical subalgebra $B_n(k)$. 
It is known (see e.g. \cite{BEG2}) that 
the spherical subalgebra is generated by 
the elements $(\sum {\rm x}_i^p)e$ and 
$(\sum {\rm y}_i^p)e$. Thus formulas (\ref{xi},\ref{yi})
determine the action of $B_n(k)$ on $U_N$.

We note that by restriction to the set $\h$ of diagonal matrices
$\diag(\lambda_1,...,\lambda_N)$, and dividing by $\Delta^{n/N}$, where
$\Delta=\prod_{i<j}(\lambda_i-\lambda_j)$, one identifies 
$U_N$ with ${\CC}[\h]^{S_N}$. Moreover, it follows from 
\cite{EG} that formulas (\ref{xi},\ref{yi}) can 
be viewed as defining an action of another spherical Cherednik algebra, 
namely $B_N(1/k)$, on ${\CC}[\h]^{S_N}$. Moreover, 
this representation is the symmetric part $W$ of the standard 
polynomial representation of $H_N(1/k)$, which is faithful and irreducible
since $1/k=n/N$ is an integer (\cite{GGOR}). 
In other words, we have the following proposition.

\begin{proposition}\label{surje}
There exists a surjective homomorphism 
$\phi: B_n(N/n)\to B_N(n/N)$, such that 
$\phi^*W=U_N$. In particular, $U_N$ is an irreducible
representation of $B_n(N/n)$. 
\end{proposition}

Proposition \ref{surje} can be generalized as follows.
Let $0\le p\le n/N$ be an integer. 
Consider the partition $\mu(p)=(n-p(N-1),p,...,p)$ of $n$. 
The representation of $\g$ attached to $\mu(p)$ 
is $S^{n-pN}{\CC}^N$.

Let $e(p)$ be a primitive idempotent of the representation 
of $S_n$ attached to $\mu(p)$. Let $U_N^p=e(p)V_N
=(\CC[\g]\otimes S^{n-pN}{\CC}^N)^\g$. Then the algebra 
$e(p)H_n(N/n)e(p)$ acts on $U_N^p$, and the above situation 
of $U_N$ is the special case $p=0$.

\begin{proposition}
There exists a surjective homomorphism 
$\phi_p: e(p)H_n(N/n)e(p)\to B_N(n/N-p)$, such that 
$\phi_p^*W=U_N^p$. In particular, $U_N^p$ is an irreducible
representation of $B_n(N/n-p)$. 
\end{proposition}

{\em Proof.} 
Similar to the proof of Proposition \ref{surje}. 
\hfill \qed \medskip

\begin{example}
$p=1$, $n=N$. In this case $e(p)=e_- = {1\over {n!}}
\sum_{\sigma\in S_n} \varepsilon(\sigma)\sigma$, the antisymmetrizer, 
and the map $\phi_p$ is the shift isomorphism 
$e_-H_N(1)e_-\to eH_N(0)e$. 
\end{example}

\subsection{Coincidence of the two ${\mathfrak{sl}}_2$ actions}

As before, let $\lbrace{b_p\rbrace}$ be an orthonormal basis of
$\g$ (under some invariant inner product). 
Consider the ${\mathfrak{sl}}_2$-triple 
\begin{equation}\label{seo} 
H=\sum b_p\frac{\partial}{\partial b_p}+\frac{\dim \g}{2}
\end{equation}
(the shifted Euler field), 
\begin{equation}\label{ef}
F=\frac{1}{2}\sum_p b_p^2, \quad E=\frac{1}{2}\Delta_\g,
\end{equation} 
where $\Delta_\g$ is the Laplace operator on $\g$. 
Recall also (see e.g. \cite{BEG2}) that the rational Cherednik 
algebra contains the ${\mathfrak{sl}}_2$-triple 
${\bold h}=\frac{1}{2}\sum_i ({\rm x}_i{\rm y}_i+{\rm y}_i{\rm
x}_i)$, ${\bold e}=\frac{1}{2}\sum_i {\rm y}_i^2$, 
${\bold f}=\frac{1}{2}\sum_i {\rm x}_i^2$.

The following proposition shows that the actions of these two 
${\mathfrak{sl}}_2$ algebras on $V_N$ essentially coincide.

\begin{proposition}\label{eulh}
On $V_N$, one has 
$$ 
{\bold h}=H, \quad {\bold e}=\frac{N}{n}E, \quad {\bold f}=\frac{n}{N}F.
$$ 
\end{proposition}

{\em Proof.} 
The last two equations follow from formulas (\ref{xi},\ref{yi}), and the
first one follows from the last two by taking commutators. 
\hfill \qed \medskip

\subsection{The irreducibility of $V_N$.}

Let $\Delta(n,N)$ be the representation of the symmetric group
$S_n$ corresponding to the rectangular Young diagram with $N$ rows
(and correspondingly $n/N$ columns), i.e. to the partition 
$(\frac{n}{N},...,\frac{n}{N})$; e.g., 
$\Delta(n,1)$ is the trivial representation.

For a representation $\pi$ of $S_n$, let 
$L(\pi)$ denote the irreducible 
lowest weight representation 
of $H_n(k)$ with lowest weight
$\pi$.

\begin{theorem}\label{irre} 
The representation $V_N$ is 
isomorphic to $L(\Delta(n,N))$. 
\end{theorem}

{\em Proof.} 
The representation $V_N$ is graded by the degree of polynomials, 
and in degree zero we have $V_N[0]=(({\CC}^N)^{\otimes 
n})^\g=\Delta(n,N)$ by the Weyl duality.

Let us show that the module $V_N$ is semisimple. 
It is sufficient to show that $V_N$ is a unitary representation, 
i.e. admits a positive definite
contravariant Hermitian form. Such a form can be 
defined by the formula 
$$
(f,g)= \langle f(\partial_A),g(A) \rangle|_{A=0},
$$
where $\langle -,-\rangle$ is the Hermitian form on $({\CC}^N)^{\otimes n}$ 
obtained by tensoring the standard forms on the factors. 
This form is obviously positive definite, and satisfies the
contravariance properties: 
$$
(Y_if,g)=\frac{N}{n}(f,X_ig),\ (f,Y_ig)=\frac{N}{n}(X_if,g).
$$
The existence of the form $(-,-)$ implies 
the semisimplicity of $V_N$. 
In particular, we have a natural inclusion 
$L(\Delta(n,N))\subset V_N$.

Next, formula (\ref{xi}) implies that $V_N$ is a torsion-free
module over $R:={\CC}[{\rm x}_1,...,{\rm x}_N]^{S_N} = 
\CC[\sum_{i=1}^N  {\rm x}_i^p, 2\le p\le N]$.
Since $V_N$ is semisimple, this implies that $V_N/L(\Delta(n,N))$
is torsion-free as well.

On the other hand, we will now show that the quotient 
$V_N/L(\Delta(n,N))$ is a torsion module over $R$. 
This will imply that the quotient is zero, as desired.

Let $v_1,...,v_N$ be the standard basis of ${\CC}^N$, 
and for each sequence $J=(j_1,...,j_n)$, $j_i\in \{1,...,N\}$, 
let $v_J:=v_{j_1}\otimes...\otimes v_{j_n}$. Let 
us say that a sequence $J$ is balanced if it contains 
each of its members exactly $n/N$ times. 
Let $B$ be the set of balanced sequences. 
The set $B$ has commuting left and right actions 
$S_N$ and $S_n$, $\sigma *(j_1,...,j_n)*\tau = (\sigma(j_{\tau(1)}),...,
\sigma(j_{\tau(n)}))$. Let $J_0 = (1...1,2...2,...,N...N)$, then 
any $J\in B$ has the form $J = J_0 *\tau$ for some $\tau\in S_n$.

Let $f\in V_N$. Then $f$ is a function $\h \to ((\CC^N)^{\otimes n})^{\h}$, 
equivariant under the action of $S_N$ (here $\h\subset\g$ is the Cartan 
subalgebra, so $\h = \{(\lambda_1,...,\lambda_N) | \sum_i \lambda_i = 0\}$), 
so 
\begin{equation}
f(\lambda)=\sum_{J\in B} f_J(\lambda)v_J,
\end{equation} 
where $\lambda=(\lambda_1,...,\lambda_N)$, and $f_J$ 
are scalar functions (the summation is over $B$ since
$f(\lambda)$ must have zero weight). By the $S_N$-invariance, 
we have $f_{\sigma * J}(\sigma(\lambda)) = f_J(\lambda)$. 
We then decompose $f(\lambda) = \sum_{o\in S_N\setminus B} f_o(\lambda)$, 
where $f_o(\lambda) = \sum_{J\in o} f_J(\lambda)v_J$.

For each $o\in S_N\setminus B$, we construct a nonzero 
$\phi_o\in \CC[{\rm x}_{1},....,{\rm x}_{n}]$
such that $\phi_o \cdot f_o(\lambda)\in L(\Delta(n,N))$. Then 
$\phi := \prod_{o\in S_N \setminus B}\prod_{\sigma\in S_{N}}\sigma(\phi_{o})
\in R$ is nonzero and 
such that $\phi \cdot f(\lambda)\in L(\Delta(n,N))$.

We first construct $\phi_o$ when $o = o_0$, the class of $J_0$. 
By $S_N$-invariance, $f_{o_0}(\lambda)$ has the form 
$$
f_{o_0}(\lambda) = \sum_{\sigma\in S_N}
g(\lambda_{\sigma(1)},...,\lambda_{\sigma(N)}) v_{\sigma(1)}^{\otimes n/N}
\otimes ... \otimes v_{\sigma(N)}^{\otimes n/N}, \quad 
{\rm where} \quad g(\lambda,...,\lambda_N)\in
\CC[\lambda_1,...,\lambda_N].
$$ 
For $\phi_{o_0}\in \CC[{\rm x}_1,...,{\rm x}_N]$, we have 
\begin{equation} \label{ffirst}
\phi_{o_0}\cdot f_{o_{0}}(\lambda) = \sum_{\sigma\in S_N} 
(\phi_{o_0}g)
(\lambda_{\sigma(1)},...,\lambda_{\sigma(N)})
v_{\sigma(1)}^{\otimes n/N}
\otimes ... \otimes v_{\sigma(N)}^{\otimes n/N}.
\end{equation}
On the other hand, let $v\in \Delta(n,N)$; expand $v = \sum_{J\in B} c_J v_J$. 
One checks that $v$ can be chosen such that $c_{J_0}\neq 0$ (one starts with 
a nonzero vector $v'$ and $J'\in B$ such that the coordinate of $v'$ along 
$J'$ is nonzero, and then acts on $v'$ by an element of $S_{n}$ bringing
$J'$ to $J_{0}$). Then since $v$ is $\g$-invariant (and therefore 
$S_N$-invariant), we have 
\begin{equation} \label{invce:c}
c_{\sigma(1)...\sigma(1)...\sigma(N)...\sigma(N)} = c_{J_0}
\end{equation}
for any $\sigma\in S_N$.

If $Q\in \CC[{\rm x}_1,...,{\rm x}_n]$, then 
\begin{equation} \label{ssecond}
(Q \cdot v)(\lambda) = \sum_{(j_1,...,j_n)\in B} 
c_{j_1...j_n} Q(\lambda_{j_1},...,\lambda_{j_n}) v_{j_1}\otimes ... \otimes
v_{j_n}\in L(\Delta(n,N)). 
\end{equation}

Set $Q_0(\lambda_{1},...,\lambda_{n})
:= \prod_{1\leq a< b \leq n, j^0_a \neq j^0_b} (\lambda_a - \lambda_b)$, 
where $J_0 = (1...1,...,N...N) = (j^0_1,...,j^0_n)$, 
$q_0(\lambda_{1},...,\lambda_{N}) := Q_{0}(\lambda_{1}...\lambda_{1},...,
\lambda_{N}...\lambda_{N})$, so $q_0(\lambda_{1},...,\lambda_{N}) = 
\big( \prod_{1\leq i<j\leq N}(\lambda_i - \lambda_j) \big)^{(n/N)^2}$.

Set $\phi_{o_{0}}(\lambda_{1},...,\lambda_{N}):= q_{0}(\lambda_{1},....,\lambda_{N})$
and 
$$
Q(\lambda_{1},...,\lambda_{n}):= Q_{0}(\lambda_{1},...,\lambda_{n})
q(\lambda_{1},\lambda_{(n/N)+1},...,\lambda_{(N-1){n\over N}+1}).
$$ 
Then (\ref{ffirst}) and (\ref{ssecond}) coincide, as: (a) for $J\notin o_{0}$, 
$Q_{0}(\lambda_{j_{1}},...,\lambda_{j_{n}}) = 0$ so the coefficient of $v_{J}$
in both expressions is zero, (b) the coefficients of $v_{J_{0}}$ in both expressions
coincide, (c) for $J\in o_{0}$, the coefficients of $v_{J}$ coincide because of (b)
and of (\ref{invce:c}).  The functions $\phi_{o}$ are constructed in the same 
way for a general $o\in S_{N}\setminus B$. This ends the proof of the theorem. 
%
%
\hfill \qed \medskip

\begin{remark} Theorem \ref{irre} 
is a special case of a much more general (but much less
elementary) Theorem \ref{funconirr}, which is proved below. 
\end{remark}

\subsection{The character formula for $V_N$.}

For each partition $\mu$ of $n$, 
let $V(\mu)$ be the representation of $\g$,
  and $\pi(\mu)$ the representation of $S_n$ 
corresponding to $\mu$.

Let $P_\mu(q)$ be the $q$-analogue of the weight multiplicity 
of the zero weight in $V(\mu)$. Namely, we have a filtration
$F^\bullet$ on $V(\mu)[0]$ such that $F^i$ is the space of
vectors in $V(\mu)[0]$ killed by the $i+1$-th power 
of the principal nilpotent element $\sum e_i$ of $\g$.
Then $P_\mu(q)=\sum_{j\ge 0} \dim(F^j/F^{j-1})q^j$. 
The coefficients of $P_\mu(q)$ are called the generalized
exponents of $V(\mu)$ (see \cite{K,He,Lu1} for more details).

We have $V_N=\oplus_\mu \pi(\mu) \otimes ({\CC}[\g]\otimes
V(\mu))^\g$. This together with Theorem \ref{irre}
implies the following.

\begin{corollary}\label{chara} The character of $L(\Delta(n,N))$
is given by the formula 
$$
{\rm Tr}|_{L(\Delta(n,N))}(w\cdot q^{\bold h})=
q^{(N^2-1)/2}\frac{\sum_\mu 
\chi_{\pi(\mu)}(w)P_\mu(q)}{(1-q^2)...(1-q^N)},
$$
where $w\in S_n$, and $\chi_{\pi(\mu)}$ is the character of $\pi(\mu)$.
Here the summation is over partitions $\mu$ of $n$ with at most $N$ parts. 
\end{corollary}

{\em Proof.}  The formula follows, using Proposition \ref{eulh}, 
from Kostant's result (\cite{K})
that $({\CC}[\g]\otimes V(\mu))^\g$ is a free 
module over ${\CC}[\g]^\g$, and the fact that the Hilbert polynomial 
of the space of generators for this module 
is the $q$-weight multiplicity of the zero weight, $P_\mu(q)$
(\cite{K,Lu1,He}).
\hfill \qed \medskip

\begin{remark} It would be interesting to compare this formula with the
character formula of \cite{Ro} for the same module.
\end{remark}

\section{Equivariant $D$-modules and representations of the
rational Cherednik algebra}

\subsection{The category of equivariant $D$-modules
on the nilpotent cone}

The theory of equivariant $D$-modules on the nilpotent cone 
arose from Harish-Chandra's work on invariant distributions on
nilpotent orbits of real groups, and was developed further 
in many papers, see e.g. \cite{HK,LS,L,Mi} and references therein. 
Let us recall some of the basics of this theory.

Let $G$ be a simply connected simple algebraic group over ${\CC}$, and 
$\g$ its Lie algebra. Let ${\mathcal N}\subset \g$ be the
nilpotent cone of $\g$. 
We denote by ${\mathcal D}(\g)$ the category of finitely
generated $D$-modules on $\g$, by ${\mathcal D}_G(\g)$ the
subcategory of $G$-equivariant $D$-modules, and 
by ${\mathcal D}_G({\mathcal N})$ the category of $G$-equivariant
$D$-modules which are set-theoretically supported 
on ${\mathcal N}$ (here we do not make a 
distinction between a $D$-module on an affine
space and the space of its global sections). 
Since $G$ acts on ${\mathcal N}$ with finitely
many orbits, it is well known that any object 
in ${\mathcal D}_G({\mathcal N})$ is regular and holonomic.

Moreover, the category ${\mathcal D}_G({\mathcal N})$ has finitely many
simple objects, and every object 
of this category has finite length (so this category is
equivalent to the category of modules over a finite dimensional
algebra).

\subsection{Simple objects in ${\mathcal D}_G({\mathcal N})$}

Recall (see e.g. \cite{Mi} and references) 
that irreducible objects in the category 
${\mathcal D}_G({\mathcal N})$ are parametrized by pairs 
$(O,\chi)$, where $O$ is a nilpotent orbit of $G$ in $\g$, 
and $\chi$ is an irreducible representation of the fundamental group 
$\pi_1(O)$, 
which is clearly isomorphic to the component group $A(O)$ 
of the centralizer $G_x$ of a point $x\in O$. Namely, $\chi$ defines a 
local system $L_\chi$ on $O$, and 
the simple object $M(O,\chi)\in {\mathcal D}_G({\mathcal N})$
is the direct image of the Goresky-Macpherson extension of $L_\chi$ to the 
closure
$\bar O$ of $O$, under the inclusion of $\bar O$ into $\g$.

\subsection{Semisimplicity of ${\mathcal D}_G({\mathcal N})$.}

The proof of the following theorem was explained to us by
G. Lusztig.

\begin{theorem}\label{semisi}
The category ${\mathcal D}_G({\mathcal N})$ is semisimple. 
\end{theorem}

{\em Proof.}  We may replace the category ${\mathcal
D}_G({\mathcal N})$ by the category of $G$-equivariant perverse 
sheaves (of complex vector spaces) on $\g$ supported on
${\mathcal N}$, ${\rm Perv}_G({\mathcal N})$, 
as these two categories are known to be equivalent. 
We must show that ${\rm Ext}^1(P,Q)=0$ for every two simple
objects $P,Q\in {\rm Perv}_G({\mathcal N})$.

Let $P',Q'$ be the Fourier transforms of $P,Q$. 
Then $P',Q'$ are character sheaves on $\g$, and it suffices 
to show that ${\rm Ext}^1(P',Q')=0$.

Recall that to each character sheaf $S$
one can naturally attach a conjugacy class of pairs $(L,\theta)$, where 
$L$ is a Levi
subgroup of $G$, and $\theta$ is a cuspidal local system on 
a nilpotent orbit for $L$. It is shown by arguments parallel to 
those in \cite{Lu3} (which treats the more difficult case of
character sheaves on the group) that if $(L_i,\theta_i)$ corresponds to
$S_i$, $i=1,2$, and $(L_1,\theta_1)$ is not conjugate to $(L_2,\theta_2)$ 
then 
${\rm Ext}^*(S_1,S_2)=0$. Thus it is sufficient to assume that
the pair $(L,\theta)$ attached to $P'$ and $Q'$ is the same.

Using standard properties of constructible sheaves (in
particular, Poincar\'e duality), we have 
$$
{\rm Ext}^1(P',Q')=H^1(\g,\underline{\rm Hom}(P',Q'))=
$$
$$
H^{2\dim \g-1}_c(\g,\underline{\rm Hom}(P',Q')^\ast)^*
=H^{2\dim \g-1}_c(\g,(Q')^\ast\otimes P')^*,
$$
where $\ast$ for sheaves denotes the Verdier duality functor.

Recall that to each character sheaf one can attach 
an irreducible representation of a certain Weyl group, 
via the generalized Springer correspondence. Let $R$ be the direct sum of 
all character sheaves corresponding to a given pair $(L,\theta)$ 
with multiplicities given by the dimensions 
of the corresponding representations. Then it is sufficient to
show that $H^{2\dim \g-1}_c(\g,(R')^\ast\otimes R')=0$.

This fact is essentially proved in \cite{Lu2}.
Namely, it follows from the computations of \cite{Lu2} that 
$H^i_c(\g,(R')^\ast\otimes R')$ is the cohomology with compact
support of a certain generalized Steinberg variety with
twisted coefficients, and it is shown that this cohomology 
is concentrated in even degrees.\footnote{More precisely, in the arguments 
of \cite{Lu2} 
the vanishing of odd cohomology is proved for $G$-equivariant
cohomology with compact supports, 
and in the non-equivariant case one should use parallel
arguments, rather than exactly
the same arguments.} 
The theorem is proved. 
\hfill \qed \medskip

\subsection{Monodromicity}
We will need the following lemma.

\begin{lemma} \label{monodromic}
Let $Q\in {\mathcal D}_G({\mathcal N})$. 
Then for any finite dimensional representation $U$ of $\g$, 
the action of the shifted Euler operator $H$ defined by (\ref{seo})
on $(Q\otimes U)^\g$ is locally finite (so $Q$ is a monodromic $D$-module), 
and has finite dimensional generalized eigenspaces. 
Moreover, the eigenvalues of $H$ on 
$(Q\otimes U)^\g$ are bounded from above. In particular, 
$(Q\otimes U)^\g$ belongs to category ${\mathcal O}$ 
for the ${\mathfrak{sl}}_2$-algebra spanned by $H$ and the
elements $E,F$ given by (\ref{ef}). 
\end{lemma}

{\em Proof.}  Since $Q$ has finite length, it is sufficient to 
assume that $Q$ is irreducible. We may further assume that $Q$ is 
generated by an irreducible $G$-submodule 
$V$, annihilated by multiplication by any invariant polynomial on $\g$
of positive degree. Indeed, let $V_0$ be an irreductible $G$-submodule
of $Q$, let $J_{V_0}:= \{f\in \CC[\g]^\g | fV_0 = 0\}$ and for any 
$v\in V_0$, let $J_v := \{f\in \CC[\g]^\g | fv=0\}$. Then if $v\in V_0$
is nonzero, $J_v = J_{V_0}$ as $Gv = V_0$. Moreover, the support 
condition implies that $J_v \subset \m^k$ for some $k\geq 0$, where 
$\m = \CC[\g]_+^\g$. So $J_{V_0} \subset \m^k$ and is an ideal of $\CC[\g]^\g$. 
Let $f\in \CC[\g]^\g$ be such that $f\notin J_{V_0}$ and $f\m \subset 
J_{V_0}$; we set $V:= fV_0$.

Then $Q$ is a quotient of 
the $D$-module $\tilde Q\otimes V$ by a $G$-stable submodule, where
$$
\tilde Q:=D(\g)/(D(\g){\rm ad}({\rm Ann}(V))+D(\g)I),
$$
${\rm Ann}(V)$ is the annihilator of $V$ in $U(\g)$, 
and $I$ is the ideal in ${\CC}[\g]$ generated by 
invariant polynomials on $\g$ of positive degree. 
Thus, it suffices to show that the lemma holds 
for the module $\tilde Q$ (which is only weakly $G$-equivariant,
i.e. the group action and the Lie algebra action coming from 
differential operators do not agree, in general).

The algebra $D(\g)$ has a grading in which $\deg(\g^*)=-1$,
$\deg(\g)=1$. This grading descends to a grading 
on $\tilde Q$. We will show that 
for each $U$, this grading on $(\tilde Q\otimes U)^\g$ 
has finite dimensional pieces, and is bounded from above. 
This implies the lemma, since the Euler operator 
preserves the grading.

Consider the associated graded module $\tilde Q_0$ of $\tilde Q$ 
under the Bernstein filtration. This is a bigraded module over 
${\CC}[\g\oplus \g]$ (where we identify $\g$ and $\g^*$ using 
the trace form). We have to show that the homogeneous subspaces 
of $(\tilde Q_0\otimes U)^\g$ under the grading defined by 
$\deg(\g\oplus 0)=-1$, $\deg(0\oplus \g)=1$ are finite dimensional.

The associated graded of the ideal ${\rm Ann}(V) \subset U(\g)$ is
such that $\CC[\g]_+^k \subset {\rm gr}{\rm Ann}(V) \subset \CC[\g]_+$
for some $k\geq 1$, therefore 
$$
\tilde Q_0={\CC}[\g\oplus \g]/J,
$$
where $J$ is a (not necessarily radical) ideal whose zero set 
is the variety ${\mathcal Z}$ of pairs 
$(u,v)\in {\mathcal N}\times \g$ such that $[u,v]=0$.
Let 
$$
Q_0'={\CC}[\g\oplus \g]/\sqrt{J}.
$$
Because of the Hilbert basis theorem, it suffices to 
prove that the homogeneous subspaces of 
$(Q_0'\otimes U)^\g$ are finite dimensional, and the degree is bounded 
above. 
But $Q_0'$ is the algebra of regular functions 
on $\mathcal Z$. By the result of \cite{J}, one has ${\CC}[{\mathcal 
Z}]^\g=
{\CC}[\g]^\g$, the algebra of invariant polynomials of $Y$. But
it follows from the Hilbert's theorem on invariants that every isotypic 
component of $\CC[{\mathcal Z}]$ is a finitely generated module over 
$\CC[{\mathcal Z}]^\g$. This implies the result. 
\hfill \qed \medskip

\subsection{Characters}

Lemma \ref{monodromic} allows one to define the character of 
an object $M\in {\mathcal D}_G({\mathcal N})$. Namely, 
let $\mu=(\mu_1,...,\mu_N)$ be a dominant integral weight for
$\g$, and $V(\mu)$ the irreducible representation 
of $\g$ with highest weight $\mu$.
Let $K_M(\mu)=(M\otimes V(\mu))^\g$.
Then the character of $M$ is defined by the formula 
$$
{\rm Ch}_M(t,g)=\Tr_{|M}(gt^{-H})=\sum_\mu
\Tr_{|K_M(\mu)}(t^{-H})\chi_\mu(g),\ g\in G,
$$
where $\chi_\mu$ denotes the character of $\mu$. 
It can be viewed as a linear functional from ${\CC}[G]^G$ 
to ${\Bbb F}:=\oplus_{\beta\in {\CC}}t^\beta{\CC}[[t]]$, via the
integration pairing.

In other words, the multiplicity spaces $K_M(\mu)$ are
representations from category ${\mathcal O}$ of the Lie 
algebra ${\mathfrak{sl}}_2$ spanned by $E,F,H$, 
and the character of $M$ carries the
information about the characters of these representations.

The problem of computing characters of simple objects 
in ${\mathcal D}_G({\mathcal N})$ is interesting and, to our
knowledge, open. Below we will show how these characters for
$G={\rm SL}_N(\CC)$ can be expressed via characters of irreducible
representations of the rational Cherednik algebra.

\begin{example} Recall (see e.g. \cite{Mi}) that an object 
$M\in {\mathcal D}_G({\mathcal N})$ is cuspidal iff 
${\mathcal F}(M)\in {\mathcal D}_G({\mathcal N})$, where ${\mathcal F}$
is the Fourier transform (Lusztig's criterion). If follows that in the
case of cuspidal objects $M$, the spaces $K_M(\mu)$ are also in the 
category ${\mathcal O}$ for the opposite Borel subalgebra of 
${\mathfrak sl}_2$, hence are finite dimensional representations of 
${\mathfrak{sl}}_2$, and, in particular, their dimensions are of interest. 
\end{example}

\subsection{The functors $F_n$, $F_n^*$}

The representation $V_N$ is a special case of 
representations of the rational Cherednik algebra 
which can be constructed via a functor 
similar to the one defined in \cite{GG1}. 
Namely, the construction of $V_N$ can be generalized as follows.

Let $n$ and $N$ be positive integers (we no longer assume that 
$N$ is a divisor of $n$), and $k=N/n$. 
We again consider the special case $G={\rm SL}_N(\CC)$,
$\g={\mathfrak{sl}}_N(\CC)$. Then we have a functor 
$F_n: {\mathcal D}(\g)\to H_n(k)$-mod 
defined by the formula 
$$
F_n(M)=(M\otimes ({\CC}^N)^{\otimes n})^\g,
$$ 
where $\g$ acts on $M$ by adjoint vector fields.
The action of $H_n(k)$ on $F_n(M)$ is defined by the same formulas
as in Proposition \ref{cheract}, and Proposition \ref{eulh} 
remains valid.

Note that $F_n(M)=F_n(M_{\rm fin})$, where $M_{\rm fin}$ is the set of
$\g$-finite vectors in $M$. Clearly $M_{\rm fin}$ is a $G$-equivariant 
$D$-module. Thus, it is sufficient to consider 
the restriction of $F_n$ to the subcategory 
${\mathcal D}_G(\g)$, which we will do from now on.

In general, $F_n(M)$ does not belong to category ${\mathcal
O}$. However, we have the following lemma.

\begin{lemma}\label{cato} If the Fourier transform ${\mathcal
F}(M)$ of $M$ is set-theoretically supported on the 
nilpotent cone ${\mathcal N}$
of $\g$, then $F_n(M)$ belongs to the category ${\mathcal O}$. 
\end{lemma}

{\em Proof.}  Since ${\mathcal F}(M)$ is supported on ${\mathcal N}$, 
invariant polynomials on $\g$ act locally nilpotently on
${\mathcal F}(M)$. Hence invariant differential operators on $\g$ 
with constant coefficients act locally nilpotently 
on $M$. Thus, it follows from formula (\ref{yi}) that 
the algebra ${\CC}[{\rm y}_1,...,{\rm y}_n]^{S_n}$ acts locally 
nilpotently 
on $F_n(M)$. Also, by Lemma \ref{monodromic}, the operator 
${\bold h}$ acts with finite dimensional generalized eigenspaces
on $F_n(M)$. This implies the statement.
\hfill \qed \medskip

Thus we obtain an exact 
functor $F_n^*=F_n\circ {\mathcal F}: {\mathcal D}_G({\mathcal
N})\to {\mathcal O}(H_n(k))$.

\subsection{The symmetric part of $F_n$}

Consider the symmetric part $eF_n(M)$
of $F_n(M)$. We have $eF_n(M)=(M\otimes S^n{\CC}^N)^\g$, 
and we have an action of the spherical subalgebra $B_n(k)$ 
on $eF_n(M)$, given by formulas (\ref{xi},\ref{yi}).

This allows us to relate the functor $F_n$ with the functor 
defined in \cite{GG1}. Namely, recall from \cite{GG1} 
that for any $c\in \ZZ$, one may define the category 
${\mathcal D}_c(\g\times \PP^{N-1})$ of coherent $D$-modules 
on $\g\times \PP^{N-1}$ which are twisted by the $c$-th power 
of the tautological line bundle on the second factor (this makes
sense for all complex $c$ even though the $c$-th power is defined
only for integer $c$). Then the paper \cite{GG1}\footnote{There
seems to be a misprint in \cite{GG1}: in the definition of $\Bbb H$, 
$c$ should be replaced by $c/N$.} defines a functor
$$
{\Bbb H}: {\mathcal D}_c(\g\times {\PP}^{N-1})\to B_N(c/N)\operatorname{-mod},
$$
given by $\HH(M)=M^\g$.

\begin{proposition}\label{gg} (i) If $n$ is divisible by $N$ 
then one has a functorial isomorphism 
$eF_n(M)\simeq \phi^*{\Bbb H}(M\otimes S^n{\CC}^N)$,
where $S^n{\CC}^N$ is regarded as a twisted $D$-module on
${\Bbb P}^{N-1}$ (with $c=n$).

(ii) For any $n$, the actions of $B_n(N/n)$ and $B_N(n/N)$ 
on the space $eF_n(M)={\Bbb H}(M\otimes S^n{\CC}^N)$ 
have the same image in the algebra of endomorphisms 
of this space. 
\end{proposition}

{\em Proof.} 
This follows from the definition of $\Bbb H$ and 
formulas (\ref{xi},\ref{yi}).
\hfill \qed \medskip

\begin{corollary}\label{irred}
The functor $eF_n^*$ on the category
${\mathcal D}_G({\mathcal N})$ maps irreducible 
objects into irreducible ones. 
\end{corollary}

{\em Proof.} 
This follows from Proposition \ref{gg}, (ii) and 
Proposition 7.4.3 of \cite{GG1}, which states that the functor
$\Bbb H$ maps irreducible objects to irreducible ones. 
\hfill \qed \medskip

Formulas \ref{xi},\ref{yi} can also be used to 
study the support of $F_n^*(M)$ for $M\in {\mathcal
D}_G({\mathcal N})$, as a ${\CC}[{\rm x}_1,...,{\rm
x}_n]$-module. Namely, we have the following proposition.

\begin{proposition}\label{support} 
Let $q=GCD(n,N)$ be the greatest common divisor of $n$ and $N$. 
Then the support $S$ of $F_n^*(M)$ is contained in the union of the 
$S_n$-translates of the subspace $E_q$ of ${\CC}^n$ defined by the 
equations
$\sum_{i=1}^n x_i=0$ and $x_i=x_j$ if $\frac{n}{q}(l-1)+1\le i,j\le
\frac{nl}{q}$ for some $1\le l\le q$. 
\end{proposition}

{\em Proof.} 
It follows from equation (\ref{yi}) that for any $(x_1,...,x_n)\in
S$ there exists a point $(z_1,...,z_N)\in {\CC}^N$ such that one has 
$$
\frac{1}{n}\sum_{i=1}^n x_i^p=\frac{1}{N}\sum_{j=1}^Nz_j^p
$$
for all positive integer $p$. In particular, writing generating
functions, we find that
$$
N\sum_{i=1}^n
\frac{1}{1-tx_i}=n\sum_{j=1}^N\frac{1}{1-tz_j}. 
$$
In particular, every fraction occurs on both sides at least $LCM(n,N)$ 
times, 
and hence the numbers $x_i$ fall into $n/q$-tuples of equal
numbers (and the numbers $z_j$ into $N/q$-tuples of equal
numbers). The proposition is proved. 
\hfill \qed \medskip

\subsection{Irreducible equivariant $D$-modules on the nilpotent
cone for $G = {\rm SL}_N(\CC)$}

Nilpotent orbits for ${\rm SL}_N(\CC)$ are labelled by Young diagrams, or
partitions. Namely, if $x\in {\mathfrak{sl}}_N(\CC)$ is a nilpotent element, 
then we let $\mu_i$ be the sizes of its Jordan blocks enumerated in 
the decreasing order. The partition $\mu=(\mu_1,...,\mu_m)$ and the 
corresponding Young diagram whose rows have lengths $\mu_i$ are attached 
to $x$. If $O$ is the orbit of $x$ then we will denote $\mu$ by $\mu(O)$.
For instance, if $O=\lbrace{0\rbrace}$ then 
$\mu=(1^N)$ and if $O$ is the open orbit then $\mu=(N)$.

It is known (and easy to show) that the group $A(O)$ is
naturally isomorphic to ${\Bbb Z}/d{\Bbb Z}$, where $d$ 
is the greatest common divisor of the $\mu_i$. Namely, 
let $Z={\Bbb Z}/N{\Bbb Z}$ be the center of $G$
(we identify ${\Bbb Z}/N{\Bbb Z}$ with $Z$
by $p\to e^{2\pi ip/N}{\rm Id}$). Then
we have a natural surjective homomorphism 
$\theta: Z\to A(O)$ induced by the inclusion 
$Z\to G_x$, $x\in O$. This homomorphism sends 
$d$ to $0$, and thus $A(O)$ gets identified with 
${\Bbb Z}/d{\Bbb Z}$.

Thus, any character $\chi: A(O)\to {\CC}^*$ 
is defined by the formula $\chi(p)=e^{-2\pi i ps/d}$, where $0\le
s<d$. We will denote this character by $\chi_s$.

\subsection{The action of $F_n^*$ on irreducible objects}

Obviously, the center $Z$ of $G$ acts on $F_n^*(M)$ by 
$z\to z^{-sN/d}$. Thus, a necessary condition for 
$F_n^*(M(O,\chi_s))$ to be nonzero is 
\begin{equation}\label{nescond}
n=N\left(p+\frac{s}{d}\right),
\end{equation}
where $p$ is a nonnegative integer.

Our main result in this section is the following theorem.

\begin{theorem}\label{funconirr} The functor $F_n^*$ maps irreducible 
objects into irreducible ones or zero. Specifically, if condition 
(\ref{nescond}) holds, then we have 
$$
F_n^*(M(O,\chi_s))=L(\pi(n\mu(O)/N)),
$$
the irreducible representation of $H_n(k)$ whose lowest weight 
is the representation of $S_n$ corresponding to the partition 
$n\mu(O)/N$.
\end{theorem}

\begin{remark}
Here if $\mu$ is a partition and $c\in \Bbb Q$ is a rational
number, then we denote by $c\mu$ the partition whose parts are
$c\mu_i$, provided that these numbers are all integers.
In our case, this integrality condition holds since all parts of
$\mu(O)$ are divisible by $d$. \hfill \qed \medskip 
\end{remark}

\begin{corollary} Let $\lambda$ be a partition of $n$ into at most $N$ parts. 
Let $M=M(O_\mu,\chi_s)$, and assume that condition (\ref{nescond})
is satisfied. Then 
$$
(M\otimes V(\lambda))^\g={\rm Hom}_{S_n}(\pi(\lambda),L(\pi(n\mu/N)))
$$
as graded vector spaces.
\end{corollary}

This corollary allows us to express the characters 
of the irreducible $D$-modules $M(O,\chi)$ in terms of characters 
of certain special lowest weight irreducible representations 
of $H_n(k)$. We note that characters of lowest weight irreducible 
representations of rational Cherednik algebras of type $A$ 
have been computed by Rouquier, \cite{Ro}.

\begin{remark} Note that Theorem \ref{irre} is the special case of 
Theorem \ref{funconirr} for $O=\lbrace{0\rbrace}$.
\end{remark}

\subsection{Proof of Theorem \ref{funconirr}}

Our proof of Theorem \ref{funconirr} is based on the following result 
(see \cite{GS} for $c$ not a half-integer, and \cite{BE}
in general).

\begin{theorem}\label{gst}
Let $k>0$. Then the functor $V\mapsto eV$ is an equivalence of categories
between $H_n(k)$-modules and $B_n(k)$-modules. 
\end{theorem}


Theorem \ref{gst} implies the first statement of the theorem, i.e. that if (\ref{nescond}) holds then 
$F_n^*(M(O_\mu,\chi_s))$ is irreducible. Indeed, it follows from Corollary \ref{irred}
that $eF_n^*(M(O_\mu,\chi_s))$ is irreducible over $B_n(k)$. Thus, it remains to find the lowest weight 
of $F_n^*(M(O_\mu,\chi_s))$.

Let $\mu=(\mu_1,...,\mu_N)$ be a partition of $N$ ($\mu_i\ge
0$). Let $O_\mu$ be the nilpotent orbit of $\g$ corresponding to the
partition $\mu$. Denote by $d$ the greatest common divisor 
of $\mu_i$, and by $m$ a divisor of $d$. 
Define the following function $f$ on $O_\mu$ 
with values in $\otimes_{i=1}^N S^{\mu_i}{\CC}^N$:
$$
f(X,\xi_1,...,\xi_N)=\bigwedge_{i=1}^N\bigwedge_{j=0}^{\mu_i-1} \xi_i X^j,
$$
$\xi_i\in ({\CC}^N)^*$ (here $X^j \in {\rm M}_N(\CC)$ is the $j$th 
power of $X$, so $\xi_i X^j\in \CC^N$).

\begin{lemma}\label{nilorb} (i) For any $X\in O_\mu$, 
$f(X,\dots)^{1/m}$ is a polynomial in
$\xi_1,...,\xi_N$. Thus, $f^{1/m}$ is a regular function on the 
universal cover $\tilde O_\mu$ of $O_\mu$ with values in 
$\otimes_{i=1}^N S^{\mu_i/m}{\CC}^N$.

(ii) For any $X\in O_\mu$, the function $f(X,\dots)^{1/m}$ 
generates a copy of the representation 
$V(\mu/m)$ inside $\otimes_{i=1}^N S^{\mu_i/m}{\CC}^N$.

(iii) Specifically, let the standard basis 
$u_1,...,u_N$ of $({\CC}^N)^*$ be filled 
into the squares of the Young diagram of $\mu$ 
(filling the first column top to bottom, then the second one,
etc.), and let $X$ be the matrix $J$ acting by the 
horizontal shift to the right on this basis. 
Then $f(J,\dots)^{1/m}$ is a highest weight vector 
of the representation $V(\mu/m)$. 
\end{lemma}

{\em Proof.}  It is sufficient to prove (iii).
Let $\mu^*=(\mu_1^*,...,\mu_N^*)$ be the conjugate partition.
Let $p_j$ be the number of times the part $j$ occurs in this
partition. Clearly, $p_j$ is divisible by $m$. 
By looking at the matrix whose determinant is $f$, 
we see that we have, up to sign: 
$$
f(J,\xi_1,...,\xi_N)=\prod_j \Delta_j(\xi_1,...,\xi_N)^{p_j},
$$
where $\Delta_j$ is the left upper $j$-by-$j$ minor of the matrix
$(\xi_1,...,\xi_N)$. Thus $f^{1/m}=\prod_j \Delta_j^{p_j/m}$ 
is clearly a highest weight vector 
of weight $\sum_j p_j\varpi_j/m$, where $\varpi_j$ 
are the fundamental weights. But $\sum p_j\varpi_j=\mu$, so we
are done. 
\hfill \qed \medskip

\begin{corollary}\label{incl} 
The function $f$ gives rise to a $G$-equivariant regular map 
$f: \tilde O_\mu\to V(\mu/d)$, whose image is the orbit of
the highest weight vector. In particular, we have a
$G$-equivariant inclusion of commutative algebras 
$$
f^*: \oplus_{\ell\ge 0}V(\ell \mu/d)^*\to
{\CC}[\tilde O_\mu].
$$ 
\end{corollary}

Now let $0\le s\le d-1$, and denote by 
${\CC}[\tilde O_\mu]_s$ the subspace of 
${\CC}[\tilde O_\mu]$, on which central elements $z\in G$
act by $z\to z^{-s}$. Then we have an inclusion 
$$
f^*: \oplus_{\ell: d^{-1}(\ell-s)\in \Bbb Z}V(\ell \mu/d)^*\to 
{\CC}[\tilde O_\mu]_s.
$$

Now recall that by construction, ${\CC}[\tilde O_\mu]_s$ 
sits inside $M=M(O_\mu,\chi_s)$ as a ${\CC}[O_\mu]$-submodule. 
In particular, the operators $X_i$ act on the space 
$({\CC}[\tilde O_\mu]_s\otimes ({\CC}^N)^{\otimes n})^\g$.

Let $\pi(\mu)$ be the representation of $S_n$ corresponding to
$\mu$, and regard $V(\lambda)\otimes \pi(\lambda)$, for any
partition $\lambda$ of $n$, as a subspace of 
$({\CC}^N)^{\otimes n}$ using the Weyl duality. 
Then for any $u\in \pi(n\mu/N)$, we can define 
the element $a(u)\in F_n^*(M)$ by $a(u)=f_n^*\otimes u$, where 
$f_n^*\in {\CC}[\tilde O_\mu]_s\otimes V(n\mu/N)$ is the homogeneous 
part of $f^*$ of degree $n$.

\begin{lemma}\label{annih}
$a(u)$ is annihilated by the elements ${\rm y}_i$
of $H_n(k)$. 
\end{lemma}

{\em Proof.} 
We need to show that the operators $X_i$ 
(or, equivalently, the elements ${\rm x}_i\in H_n(k)$) annihilate
$a(u)\in F_n(M)$. Since $a(u)$ is $G$-invariant, 
it is sufficient to prove the statement at the point $X=J$. 
This boils down to showing that for any $j$ not exceeding the
number of parts of $\mu$ (i.e. $j\le \mu_1^*$), 
the application of $J$ in any component annihilates
the element $\Delta_j(\xi_1,...,\xi_N)\in \wedge^j \CC^N\subset 
({\CC}^N)^{\otimes j}$. This is clear, since 
the first $\mu_1^*$ columns of $J$ are zero. 
\hfill \qed \medskip

This implies that the lowest weight of $F_n^*(M(O_\mu,\chi_s))$
is $\pi(n\mu/N)$, as desired. The theorem is proved.

\begin{remark} 
Here is another, short proof of Theorem \ref{funconirr} for $n=N$. 
We have 
$$
e_-F_N^*(M(O,1))={\mathcal F}(M(O,1))^G.
$$
According to \cite{L,LS}, 
$$
{\mathcal F}(M(O,1))^G=
({\CC}[\h]\otimes \pi(\mu(O)))^{S_N}
$$
as a module over $D(\h)^W=e_-H_N(1)e_-$. 
Thus, $e_-F_N^*(M(O,1))=e_-L(\pi(\mu(O)))$ as 
$e_-H_N(1)e_-$-modules. But the functor $V\to e_-V$ 
is an equivalence of categories $H_N(1)$-mod $\to 
e_-H_N(1)e_-$-mod (see \cite{BEG2}). 
Thus, $F_N^*(M(O,1))=L(\pi(\mu(O)))$ as $H_N(1)$-modules, as
desired. 
\end{remark}

\subsection{The support of $L(\pi(n\mu/N))$}

\begin{corollary}\label{superr} 
Let $\mu$ be a partition of $N$ such that $n\mu_i/N$ are
integers. Then the support of the representation
$L(\pi(n\mu/N))$ of $H_n(N/n)$ as a module over ${\Bbb
C}[{\rm x}_1,...,{\rm x}_n]$ is contained in the union of
$S_n$-translates of $E_q$, $q=GCD(n,N)$. 
\end{corollary}

{\em Proof.}  This follows from Theorem \ref{funconirr}
and Proposition \ref{support}. 
\hfill \qed \medskip

We note that in the case when $\mu=(N)$, Corollary \ref{superr} 
follows from Theorem 3.2 from \cite{CE}.

\subsection{The cuspidal case} 
An interesting special case of Theorem \ref{funconirr} 
is the cuspidal case. In this case $N$ and $n$ are relatively
prime, $d=N$ (i.e., $O$ is the open orbit), 
and $s$ is relatively prime to $N$.

Here is a short proof of Theorem \ref{funconirr}
in the cuspidal case.

Since the Fourier transform of $M(O,\chi_s)$ in the cuspidal case
is supported on the nilpotent cone, $F_n^*(M(O,\chi_s))$ belongs 
not only to the category ${\mathcal O}$ generated by
lowest weight modules, but also to the ``dual'' category 
${\mathcal O}_-$ generated by highest weight modules over $H_n(k)$. 
Thus, by the results of \cite{BEG1}, 
$F_n^*(M(O,\chi_s))$ is a multiple of the unique finite
dimensional irreducible $H_n(k)$-module $L({\CC})=Y(N,n)$, 
of dimension $N^{n-1}$. But this multiple must be a single 
copy by Corollary \ref{irred}, so the theorem is proved.

Theorem \ref{funconirr} implies the following 
formula for the characters of the cuspidal $D$-modules $M(O,\chi_s)$.

Let $\mu$ be a dominant integral weight for $\g$, such that
the center $Z$ of $G$ acts on $V(\mu)$ via 
$z\to z^s=z^n$. Let $\rho$ be the half-sum of positive
roots of $\g$. Let $K_s(\mu)=(M(O,\chi_s)\otimes V(\mu))^\g$
be the isotypic components of $M(O,\chi_s)$.

\begin{theorem}\label{deco} We have 
$$
{\rm Tr}_{|K_s(\mu)}(q^{2H})=\frac{q-q^{-1}}{q^N-q^{-N}}
\varphi_\mu(q),
$$
where 
$$
\varphi_\mu(q):=\prod_{1\le p<r\le
N}\frac{q^{\mu_r-\mu_p+r-p}-q^{\mu_p-\mu_r+p-r}}{q^{r-p}-q^{p-r}}=
\chi_{V(\mu)}(q^{2\rho}),
$$
where $\chi_{V(\mu)}$ is the character of $V(\mu)$. 
In particular, 
$$
\dim K_s(\mu)=\frac{1}{N}\prod_{1\le p<r\le
N}\frac{\mu_r-\mu_p+r-p}{r-p}= 
\frac{1}{N}\dim V(\mu). 
$$
\end{theorem}

{\em Proof.} 
We extend the representation $V(\mu)$ to ${\rm GL}_N(\CC)$ by
setting $z\to z^n$ for all scalar matrices $z$, 
so that its ${\rm GL}_N(\CC)$-highest weight is
$$
\tilde\mu:=(\mu_1+n/N,...,\mu_N+n/N).
$$
Note that we automatically have $\mu_i+n/N\in \Bbb Z$. 
Assume that $n$ is so big that $\tilde\mu$ is a partition
of $n$ (i.e., $\mu_i+n/N\ge 0$).

It follows from the results of \cite{BEG1} that the character 
of the irreducible representation $L({\CC})$ 
of the rational Cherednik algebra $H_n(k)$, $k=N/n$, 
is given by the formula

\begin{equation}\label{lc}
\Tr_{|L({\CC})}(g q^{2\bold h})=\frac{q-q^{-1}}{q^N-q^{-N}}
\frac{\det(q^{-N}-q^Ng)}{\det(q^{-1}-qg)}, \quad g\in S_n,
\end{equation}
where the determinants are taken in ${\CC}^n$.

Let us equip ${\CC}^N$ with the structure of an irreducible
representation of ${\mathfrak{sl}}_2$ with basis $e,f,h$. Let $g\in S_n$. 
Then 
$$
\Tr_{|{\rm Hom}_{S_n}(\pi(\tilde\mu),({\CC}^N)^{\otimes n})}(q^h)=
\Tr_{|V(\mu)}(q^{2\rho})=\varphi_\mu(q), 
$$
by the Weyl character formula. 
On the other hand, it is easy to show that 
$$
\Tr_{|({\CC}^N)^{\otimes n}}(g q^h)=
\frac{\det(q^{-N}-q^Ng)}{\det(q^{-1}-qg)}.
$$
Thus, 
\begin{align*}
\Tr_{|{\rm Hom}_{S_n}(\pi(\tilde\mu),L({\CC}))}(q^{2\bold h}) & =
\frac{q-q^{-1}}{q^N-q^{-N}}
\Tr_{|{\rm Hom}_{S_n}(\pi(\tilde\mu),({\CC}^N)^{\otimes n})}(q^h)
\\ & =
\frac{q-q^{-1}}{q^N-q^{-N}}\varphi_\mu(q). 
\end{align*}
By Theorem \ref{funconirr} and Weyl duality, this implies that 
$$
\Tr_{|(M(O,\chi_s)\otimes V(\mu))^\g}(q^{2H})=
\frac{q-q^{-1}}{q^N-q^{-N}}\varphi_\mu(q),
$$
as desired. 
\hfill \qed \medskip

\begin{example}
Let $N=2$, $s=1$. In this case Theorem \ref{deco} 
gives us the following decomposition of $M(O,\chi_s)$: 
$$
M(O,\chi_s)=\oplus_{j\ge 1} N_j\otimes V_{2j-1},
$$
where $V_j$ is the irreducible representation of ${\mathfrak{sl}}_2$ 
of dimension $j+1$, and the spaces $N_j$ 
satisfy the equation 
$$
{\rm Tr}_{|N_j}(q^{2H})=\frac{q^{2j}-q^{-2j}}{q^2-q^{-2}}.
$$
This shows that $N_j=V_{j-1}$ as a representation of the
${\mathfrak{sl}}_2$-subalgebra spanned by $E,F,H$, which commutes with $\g$. 
\end{example}

\subsection{The case of general orbits}

Let $W=S_N$ the Weyl group of $G$, 
$\lambda\in \h/W$, and ${\mathcal N}_\lambda$ be the closure
in $\g$ of the adjoint orbit of a regular element of $\g$ whose
semisimple part is $\lambda$. 
Denote by ${\mathcal D}_G({\mathcal N}_\lambda)$ the category of
$G$-equivariant $D$-modules on $G$ which are concentrated on
${\mathcal N}_\lambda$. We also let ${\mathcal O}_\lambda$
be the category of finitely generated $H_n(k)$-modules in which
the subalgebra ${\Bbb C}[{\rm y}_1,...,{\rm y}_n]^{S_n}$ acts through the
character $\lambda$. Then one can show, similarly to the above,
that the functor $F_n^*$ restricts to a functor $F^*_{n,\lambda}:
{\mathcal D}_G({\mathcal N}_\lambda)\to {\mathcal O}_\lambda$. 
The functor considered above is $F_{n,0}^*$. We plan to study 
the functor $F_{n,\lambda}^*$ for general $\lambda$ in a 
future work.

\subsection{The trigonometric case}

Our results about rational Cherednik algebras can be extended 
to the trigonometric case. For this purpose, $D$-modules 
on the Lie algebra $\g$ should be replaced with $D$-modules on the group 
$G$. Let us describe this generalization.

First, let us introduce some notation. As above, we let $G={\rm SL}_N(\CC)$. 
For $b\in \g$, let $L_b$ be the right invariant vector field 
on $G$ equal to $b$ at the identity element; that is, 
$L_b$ generates the group of left translations by $e^{tb}$. 
As before, we let $k=N/n$.

Now let $M$ be a $D$-module on $G$. Similarly to the above, 
we define $F_n(M)$ to be the space 
$$
F_n(M)=(M\otimes ({\Bbb C}^N)^{\otimes n})^G,
$$
where $G$ acts on itself by conjugation.

Consider the operators $X_i,Y_i$, $i=1,...,n$, on $F_n(M)$, defined 
by the formulas similar to the rational case:
$$
X_i=\sum_{j,l}A_{jl}\otimes (E_{lj})_i,\ Y_i=
\frac{N}{n}\sum_p L_{b_p}\otimes (b_p)_i,
$$
where $A_{jl}$ is the $jl$-th matrix element of $A\in G$ regarded as 
the multiplication operator in $M$ by a regular function on $G$.

\begin{proposition}\label{rela}
The operators $X_i,Y_i$ satisfy the following relations:
$$
\prod_i X_i=1,\ \sum_i Y_i+k\sum_{i<j}s_{ij}=0,
$$
$$
s_{ij}X_i=X_js_{ij},\ s_{ij}Y_i=Y_js_{ij},\ [s_{ij},X_l]=[s_{ij},Y_l]=0,
$$
$$
[X_i,X_j]=0,\ [Y_i,Y_j]=ks_{ij}(Y_i-Y_j),
$$
$$
[Y_i,X_j]=\left(ks_{ij}-\frac{1}{n}\right)X_j,
$$
where $i,j,l$ denote distinct indices. 
\end{proposition}

{\em Proof.}
Straightforward computation. 
\hfill \qed \medskip

\begin{corollary}
The operators $\bar Y_i=Y_i+k\sum_{j<i}s_{ij}$ 
pairwise commute.
\end{corollary}

The relations of Proposition \ref{rela} are nothing but 
the defining relations of the {\it degenerate double affine Hecke 
algebra} of type $A_{n-1}$, which we will denote $H_n^{\rm tr}(k)$ 
(where ``tr'' stands for trigonometric, to illustrate the fact that 
this algebra is a trigonometric deformation of the rational 
Cherednik algebra $H_n(k)$). Thus we have defined an exact functor 
$F_n: {\mathcal D}(G)\to H_n^{\rm tr}(k)$-mod. 
As before, it is sufficient to consider the restriction of this functor 
to the category of equivariant finitely generated $D$-modules, 
${\mathcal D}_G(G)$.

This allows us to generalize much of our story for rational
Cherednik algebras to the trigonometric case. In particular, let
${\mathcal U}$ be the unipotent variety on $G$, and ${\mathcal
D}_G({\mathcal U})$ be the category of finitely generated
$G$-equivariant $D$-modules on $G$ concentrated on ${\mathcal
U}$. If we restrict the functor $F_n$ to this category, we get a
situation identical to that in the rational case. Indeed, one can
show that for any $M$ in this category, $F_n(M)$ belongs to the
category ${\mathcal O}_-^{\rm tr}$ of finitely generated modules
over $H_n^{\rm tr}(k)$ which are locally unipotent with respect
to the action of $X_i$. The latter category is equivalent to
the category ${\mathcal O}_-$ over the rational Cherednik algebra
$H_n(k)$, because the completion of $H_n^{\rm tr}(k)$ with
respect to the ideal generated by $X_i-1$ is isomorphic to the
completion of $H_n(k)$ with respect to the ideal generated by
${\rm x}_i$. On the other hand, the exponential map identifies
the categories ${\mathcal D}_G({\mathcal U})$ and ${\mathcal
D}_G({\mathcal N})$. It is clear that after we make these two
identifications, the functor $F_n$ becomes the functor $F_n$ in
the rational case that we considered above.

On the other hand, because of the absence of Fourier transform on the 
group
(as opposed to Lie algebra), the trigonometric story is richer than the 
rational one.
Namely, we can consider another subcategory of ${\mathcal D}_G(G)$, 
the category of character sheaves. By definition, a character sheaf on $G$ 
is 
an object $M$ in ${\mathcal D}_G(G)$ which is locally finite with respect 
to the action 
of the algebra of biinvariant differential operators, $U(\g)^G$. 
This category is denoted by ${\rm Char}(G)$. It is known that one has a 
decomposition 
$$
{\rm Char}(G)=\oplus_{\lambda\in T^\vee/W}{\rm Char}_\lambda(G),
$$
where $T^\vee$ is dual torus, and 
${\rm Char}_\lambda(G)$ the category of those $M\in {\mathcal D}_G(G)$ 
for which the generalized eigenvalues of $U(\g)^G$ (which we identify
with $U(\h)^W$ via the Harish-Chandra homomorphism) project to $\lambda$ 
under the natural projection $\h^*\to T^\vee$.

On the other hand, one can define the category ${\rm
Rep}_{Y{\rm -fin}}(H_n^{\rm tr}(k))$ of modules over $H_n^{\rm
tr}(k)$ on which the commuting elements $\bar Y_i$ act in
a locally finite manner. We have a similar decomposition 
$$ 
{\rm Rep}_{Y{\rm -fin}}(H_n^{\rm tr}(k))=\oplus_{\lambda\in T^\vee/W} {\rm
Rep}_{Y{\rm -fin}}(H_n^{\rm tr}(k))_\lambda, 
$$ 
where ${\rm Rep}_{Y{\rm -fin}}(H_n^{\rm tr}(k))_\lambda$ is the subcategory of
all objects where the generalized eigenvalues of $\bar Y_i$ project to 
$\lambda\in T^\vee/W$. Then one can show, similarly to the rational case,
that the functor $F_n$ gives rise to the functors 
$$
F_{n,\lambda}: {\rm Char}_\lambda(G)\to 
{\rm Rep}_{Y-fin}(H_n^{\rm tr}(k))_\lambda 
$$ 
for each $\lambda\in T^\vee/W$. The most interesting case is $\lambda=0$ 
(unipotent character sheaves). We plan to study these functors in 
subsequent works.

\subsection{Relation with the Arakawa-Suzuki functor}

Note that the elements $Y_i$ and $s_{ij}$ generate the degenerate 
affine Hecke algebra ${\mathcal H}_n$ of Drinfeld and Lusztig (of type 
$A_{n-1}$). To define the action of this algebra on 
$F_n(M)=(M\otimes ({\Bbb C}^N)^{\otimes n})^\g$ 
by the formula of Proposition \ref{rela}, we only need the action
of the operators $L_b$, $b\in \g$ in $M$. So $M$ can be taken to be
an arbitrary $\g$-bimodule which is locally finite with respect to
the diagonal action of $\g$ (in this case,  $\sum_i Y_i+\sum_{i<j}s_{ij}$ 
is a central element which does not necessarily act by zero, so we get 
a representation of a central extension 
$\widetilde{\mathcal H}_n$ of ${\mathcal H}_n$). In particular, we
have an exact functor $F_n: {\rm HC}(\g)\to \widetilde{\mathcal H}_n$-mod 
from the category of Harish-Chandra bimodules over $\g$ to the
category of finite dimensional representations of the degenerate
affine Hecke algebra $\widetilde{\mathcal H}_n$. 
This functor was essentially considered in \cite{AS} (where it
was applied to the Harish-Chandra modules of the form $M={\rm
Hom}_{{\mathfrak{g}}{\rm -finite}}(M_1,M_2)$, where $M_1$ and
$M_2$ are modules from category ${\mathcal O}$ over $\g$).
We note that the paper \cite{AST} describes the extension of this
construction to affine Lie algebras, which yields representations
of degenerate double affine Hecke algebras.

\subsection{Directions of further study}

In conclusion we would like to discuss (in a fairly speculative manner) 
several directions of further study and generalizations
(we note that these generalizations can be combined with each other).

1. The $q$-case: the group $G$ is replaced with the corresponding
quantum group, $D$-modules with $q$-$D$-modules, and degenerate double affine
Hecke algebras with the usual double affine Hecke algebras
(defined by Cherednik). It is especially interesting to consider
this generalization if $q$ is a root of unity.

2. The quiver case. This generalization was suggested by Ginzburg,
and will be studied in his subsequent work with the third author.
In this case, one has a finite subgroup $\Gamma\subset {\rm SL}_2(\CC)$, 
and one should consider equivariant $D$-modules on the
representation space of the affine quiver attached to $\Gamma$
(with some orientation). Then there should exist an analog
of the functor $F_n$, which takes values in the category of
representations of an appropriate symplectic reflection algebra
for the wreath product $S_n\ltimes \Gamma^n$, \cite{EG} (or, equivalently,
the Gan-Ginzburg algebra, \cite{GG2}). This generalization should
be especially nice in the case when $\Gamma$ is a cyclic group,
when the symplectic reflection algebra is a Cherednik algebra 
for a complex reflection group, and 
one has the notion of category ${\mathcal O}$ for it.

3. The symmetric space case. This is the trigonometric version 
of the previous generalization for $\Gamma={\Bbb Z}/2$. In this
generalization one considers (monodromic) equivariant $D$-modules on the
symmetric space ${\rm GL}_{p+q}(\CC)/({\rm GL}_p\times {\rm GL}_q)(\CC)$ 
(see \cite{Gin}), and one expects a functor from this category to the 
category of representations of an appropriate degenerate double affine 
Hecke algebra of type $C^\vee C_n$. This functor should be related, 
similarly to the previous subsection, to an analog of the Arakawa-Suzuki
functor, which would attach to a Harish-Chandra module 
for the pair $({\rm GL}_{p+q}(\CC),{\rm GL}_p(\CC)\times {\rm GL}_q(\CC))$, 
a finite dimensional representation of the degenerate double affine 
Hecke algebra of type $BC_n$.

\begin{appendix}

\section{}

Let $\cO$ be the ring $\CC[[u_{1},...,u_{n}]][\ell_{1},...,\ell_{n}]$. 
Define commuting derivations $D_{i}$ of $\cO$ by 
$D_{i}(u_{j}) = \delta_{ij}u_{i}$, $D_{i}(\ell_{j}) = \delta_{ij}$
(we will later think of $\ell_{i}$ and $D_{i}$ as $\on{log}u_{i}$ and 
$u_{i}{\partial\over{\partial u_{i}}}$).

We set $\cO_{+}:= \m[\ell_{1},...,\ell_{n}]$, where 
$\m = \on{Ker}(\CC[[u_{1},...,u_{n}]]\to\CC)$ is the augmentation
ideal. Let $A = \oplus_{k\geq 0} A_{k}$ be a graded ring with 
finite dimensional homogeneous components.

\begin{proposition} \label{prop:app:1}
Let $X_{i}(u_{1},...,\ell_{n}) \in \hat\oplus_{k>0} (A_{k}\otimes \cO_{+})$
be such that $D_{i}(X_{j}) = D_{j}(X_{i})$. Then there exists a unique 
$F(u_{1},...,\ell_{n})\in \hat\oplus_{k>0} (A_{k}\otimes \cO_{+})$
such that $D_{i}(F) = X_{i}$ for $i=1,...,n$.

Let us say that $f\in \cO$ has radius of convergence $R>0$
if $f = \sum_{k_{1},...,k_{n}\geq 0} f_{k_{1},...,k_{n}}(u_{1},...,u_{n})
\ell_{1}^{k_{1}}...\ell_{n}^{k_{n}}$, where each $f_{k_{1},...,k_{n}}
(u_{1},...,u_{n})$ converges for $|u_{1}|,...,|u_{n}|\leq R$. Then if 
$X_{1},...,X_{n}$ have radius of convergence $R$, so does $F$. 
\end{proposition}

{\em Proof.} For each $i$, $D_{i}$ restricts to an endomorphism of $\cO_{+}$; 
one checks that $\cap_{i=1}^{n}\on{Ker}(D_{i} : \cO_{+}\to \cO_{+}) = 0$
which implies the uniqueness. To prove the existence, we work by induction. 
One proves that $D_{n} : \cO_{+}\to \cO_{+}$ is surjective, and its kernel 
is $\m_{n-1}[\ell_{1},...,\ell_{n-1}]$, where $\m_{n-1} = 
\on{Ker}(\CC[[u_{1},...,u_{n-1}]]\to\CC)$. Let $G$ be a solution of 
$D_{n}(G) = X_{n}$, then the system $D_{i}(F') = X_{i} - D_{i}(G)$
($i=1,...,n$)
is compatible, which implies $D_{n}(X'_{i}) = 0$, where 
$X'_{i}:= X_{i}-D_{i}(G)$, so $X'_{i}\in \oplus_{k>0}(A_{k}\otimes \cO_{+}^{(n-1)})$, 
where $\cO_{+}^{(n-1)}$ is the analogue of $\cO_{+}$ at order $n-1$. 
Hence the system $D_{i}(F') = X_{i} - D_{i}(G)$
($i=1,...,n-1$) is compatible and we may apply to it the result at order 
$n-1$ to obtain a solution $F'$. Then a solution of $D_{i}(F)=X_{i}$
is $F'+G$.

Let $D : u\CC[[u]]\to u\CC[[u]]$ be the map $u{\partial\over{\partial u}}$
and let $I:= D^{-1}$. The map $D_{1} : u\CC[[u]][\ell] \to u\CC[[u]][\ell]$
is bijective and its inverse is given by 
$D_{1}^{-1}(F(u)\ell^{a}) = \sum_{k=0}^{a} (-1)^{k} a(a-1)...(a-k+1)
(I^{k+1}(F))(u)\ell^{a-k}$.

We have $\cO_{+} = \cO^{(n-1)}\hat\otimes u_{n}\CC[[u_{n}]][\ell_{n}]
\oplus \m^{(n-1)}\hat\otimes \CC[\ell_{n}]$
(where $\cO^{(n-1)},\m^{(n-1)}$ are the analogues of $\cO,\m$ at order $n-1$, 
$\hat\otimes$ is the completed tensor product). The endomorphism $D_{n}$
preserves this decomposition and a section of $D_{n}$ is given by 
$(\on{id}\otimes D_{1}^{-1}) \oplus (\on{id}\otimes J)$, where 
$J\in\on{End}(\CC[\ell])$ is a section of $\partial/\partial\ell$.

It follows from the fact that $I$ preserves the radius of convergence 
of a series that the same holds for the section of $D_{n}$ defined above. 
One then follows the above construction of a solution $X$ of 
$D_{i}(X)=X_{i}$ and uses the fact that $D_{i}$ also preserves 
the radius of convergence to show by induction that $X$ has radius $R$
if the $X_{i}$ do.  \hfill \qed \medskip

\begin{proposition} \label{prop:app:2}
Let $X_{i}(u_{1},...,\ell_{n})\in \hat\oplus_{k>0}(A_{k}\otimes\cO_{+})$ 
be such that $D_{i}(X_{j})-D_{j}(X_{i})=[X_{i},X_{j}]$. Then there exists a
unique $F(u_{1},...,\ell_{n}) \in 1 + \hat\oplus_{k>0} (A_{k}\otimes\cO_{+})$
such that $D_{i}(F) = X_{i}F$ for $i=1,...,n$. 
If the $X_{i}$ have radius $R$, then so does $F$. 
\end{proposition}

{\em Proof.} Let us prove the uniqueness. If $F,F'$ are two solutions, then 
$F^{-1}F'$ is a constant (as $\cap_{i=0}^{n}\on{Ker}(D_{i} : \cO\to \cO)=0$), 
and it also belongs to $1+\hat\oplus_{k>0}(A_{k}\otimes \cO_{+})$, 
which implies that $F = F'$. To prove the existence, one sets 
$F = 1 + f_{1}+f_{2}+...$, $X_{i} = x^{(i)}_{1}+...$, 
where $f_{k},x^{(i)}_{k}\in A_{k}\otimes \cO_{+}$ and solves by induction the 
system $D_{i}(f_{k}) = x_{1}^{(i)}f_{k-1}+...+x_{k}^{(i)}$ using 
Proposition \ref{prop:app:1}. 
\hfill \qed \medskip

\begin{proposition} \label{prop:app:3}
Let $C_{i}(u_{1},...,u_{n})\in \hat\oplus_{k>0} A_{k}[[u_{1},...,u_{n}]]$
($i=1,...,n$) be such that 
$u_{i}\partial_{u_{i}}(C_{j}) - u_{j}\partial_{u_{j}}(C_{i}) = [C_{i},C_{j}]$
for any $i,j$. Assume that the series $C_{i}$ have radius $R$.

Then there exists a unique solution of the system 
$u_{i}\partial_{u_{i}}(X) = C_{i}X$, analytic in the domain
$\{u||u|\leq R, u\notin\RR_{-}\}^{n}$, such that 
the ratio 
$(u_{1}^{C_{0}^{1}}...u_{n}^{C_{0}^{n}})^{-1}X(u_{1},...,u_{n})$
(we set $C_{0}^{i}:= C_{i}(0,...,0)$)
has the form $1 + \sum_{k>0} \sum_{a_{1},...,a_{n},i}
r_{k}^{a_{1},...,a_{n},i}(u_{1},...,u_{n})$ 
(the second sum is finite for any $k$), $r_{k}^{a_{1},...,a_{n},i}$ 
has degree $k$, $a_{i}\geq 0$, $i\in \{1,...,n\}$, and 
$r_{k}^{a_{1},...,a_{n},i}(u_{1},...,u_{n}) = O(u_{i} 
(\on{log}u_{1})^{a_{1}}... (\on{log}u_{n})^{a_{n}})$. 
The same is then true of the ratio 
$X(u_{1},...,u_{n}) (u_{1}^{C_{0}^{1}}...u_{n}^{C_{0}^{n}})^{-1}$; 
we write $X(u_{1},...,u_{n}) \simeq u_{1}^{C_{0}^{1}}...u_{n}^{C_{0}^{n}}$.
  \end{proposition}

{\em Proof.} Let us show the existence of $X$. 
The compatibility condition implies that $[C_{0}^{i},C_{0}^{j}]=0$. 
If we set $Y(u_{1},...,u_{n}):= 
(u_{1}^{C_{0}^{1}}...u_{n}^{C_{0}^{n}})^{-1}X(u_{1},...,u_{n})$, 
then $X$ is a solution iff $Y$ is a solution of 
$u_{i}\partial_{u_{i}}(Y) = \on{exp}(-\sum_{j=1}^{n} 
(\on{log}u_{j}) C_{j}^{0})(C_{i} - C_{i}^{0}) \cdot Y$.

Let us set $X_{i}(u_{1},...,\ell_{n}) := \on{exp}(-\sum_{j=1}^{n}
\ell_{j}C_{j}^{0})(C_{i}(u_{1},...,u_{n}) - C_{i}(0,...,0))$, then 
$X_{i}(u_{1},...,\ell_{n})\in \hat\oplus_{k>0}(A_{k}\otimes \cO_{+})$. 
We then apply Proposition \ref{prop:app:2} and find a solution 
$Y\in 1 + \hat\oplus_{k>0} A_{k}\otimes \cO_{+}$ of 
$D_{i}(Y) = X_{i}Y$. Let $Y_{k}$ be the component of $Y$ of degree $k$. 
Since $Y$ has radius $R$, the replacement 
$\ell_{i} = \on{log}u_{i}$ in $Y_{k}$
for $u_{i}\in \{u||u|\leq R, u\notin\RR_{-}\}$
gives an analytic function on $\{u||u|\leq R, u\notin \RR_{-}\}^{n}$.
Moreover, $\cO_{+} = \sum_{i=1}^{n} 
u_{i}\CC[[u_{1},...,u_{n}]][\ell_{1},...,\ell_{n}]$, 
which gives a decomposition $Y_{k} = \sum_{i,a_{1},...,a_{n}} 
u_{i} \ell_{1}^{a_{1}}...
\ell_{n}^{a_{n}} y^{k}_{i,a_{1},...,a_{n}}(u_{1},...,u_{n})$ and leads 
(after substitution $\ell_{i} = \on{log}u_{i}$) to the above estimates.

The ratio $X(u_{1},...,u_{n}) (u_{1}^{C_{0}^{1}}...u_{n}^{C_{0}^{n}})^{-1}$
is then $1 + \on{exp}(\sum_{j} C^{j}_{0}\on{log}u_{j})(Y(u_{1},...,u_{n})-1)$; 
the term of degree $k$ has finitely many contributions to which we apply the 
above estimates.

Let us prove the uniqueness of $X$. Any other solution has the form 
$X = X(1+c_{k}+...)$ where $c_j\in A_{j}$, and $c_{k}\neq 0$. Then 
the degree $k$ term is transformed by the addition of $c_{k}$, which 
cannot be split as a sum of terms in the various 
$O(u_{i}(\on{log}u_{1})^{a_{1}}... (\on{log}u_{n})^{a_{n}})$. 
\hfill \qed\medskip

\end{appendix}

\medskip\noindent
{\bf Acknowledgments.} 
This project started in June 2005, 
while the three authors were visiting ETH; they would like to thank 
Giovanni Felder for his kind invitations.

P.E. is deeply grateful to Victor Ginzburg 
for a lot of help with proofs of the main results 
in the part about equivariant $D$-modules, 
and for explanations on the Gan-Ginzburg functors. 
Without this help, this part could not have been written. 
P.E. is also very grateful to G. Lusztig for 
explaining to him the proof of Theorem \ref{semisi}, and 
to V. Ostrik and  D. Vogan for useful discussions.

The work of D.C. has been partially supported by the European Union 
through the FP6 Marie Curie RTN ENIGMA (contract number 
MRTN-CT-2004-5652).

\end{document}